\documentclass{amsbook}
\usepackage{amssymb}
\usepackage{amsfonts}
\usepackage{amsmath}
\usepackage{lastpage}
\usepackage[legalpaper,bookmarks=true,colorlinks=true,linkcolor=blue,citecolor=blue]{hyperref}
\usepackage{graphicx}%
\usepackage{fancyhdr}
\usepackage{color}
\usepackage[mathlines]{lineno}
\usepackage{lscape}
\usepackage{epsfig}
\usepackage[]{listings}

\newtheorem{theorem}{Theorem}
\theoremstyle{plain}

\newtheorem{corollary}{Corollary}

\newtheorem{definition}{Definition}

\newtheorem{lemma}{Lemma}

\newtheorem{proposition}{Proposition}
\newtheorem{remark}{Remark}

\numberwithin{equation}{section}

\begin{document}
\Large
\pagenumbering{roman}
\begin{center}

\huge \textbf{Gane Samb LO}\\
\bigskip \vskip 5cm
\Huge \textbf{Weak Convergence (IIA) - Functional and Random Aspects of the Univariate Extreme Value Theory} 
\vskip 5cm

 \huge \textbf{Calgary, Alberta. 2017}.\\

\bigskip \Large  \textbf{DOI} : http://dx.doi.org/10.16929/srm/2016.0009\\
\end{center}

\newpage
\newpage
\begin{center}
\huge \textbf{SPAS BOOKS SERIES}
\end{center}

\bigskip \bigskip

\Large

 \begin{center}
 \textbf{GENERAL EDITOR of SPAS EDITIONS}
 \end{center}

\bigskip
\noindent \textbf{Prof Gane Samb LO}\\
gane-samb.lo@ugb.edu.sn, gslo@ugb.edu.ng\\
Gaston Berger University (UGB), Saint-Louis, SENEGAL.\\
African University of Sciences and Technology, AUST, Abuja, Nigeria.\\

\bigskip

\begin{center}
\Large \textbf{ASSOCIATED EDITORS}
\end{center}

\bigskip
\noindent \textbf{Charles E CHIDUME}\\
cchidume@aust.edu.ng\\
African University of Sciences and Technology (AUST), Abuja, Nigeria\\

\bigskip
\noindent \textbf{KEHINDE DAHUD SHANGODOYIN}\\
shangodoyink@mopipi.ub.bw\\
UNIVERSITY Of BOTSWANA\\

\noindent \textbf{Blaise SOME}\\
some@univ-ouaga.bf\\
Chairman of LANIBIO, UFR/SEA\\
Ouaga I Pr Joseph Ki-Zerbo University.\\

\bigskip
\begin{center}
\Large \textbf{ADVISORS}
\end{center}

\bigskip

\noindent \textbf{Ahmadou Bamba SOW}\\
ahmadou-bamba.sow@ugb.edu.sn\\
Gaston Berger University, Senegal.\\

\noindent \textbf{Tchilabalo Abozou KPANZOU}\\
kpanzout@yahoo.fr\\
Kara University, Togo.\\

\newpage

\begin{center}
\LARGE \textbf{LIST OF SPAS EDITIONS BOOKS \\
PUBLISHED OR UNDER PRESS}.\\

\end{center}
\Large
\bigskip \noindent \textbf{COLLECTION TEXTBOOKS SPAS EDITIONS BOOKS PUBLISHED OR UNDER PRESS}.\\

\noindent  (1)\textit{Weak Convergence (IA) - Sequences of Random Variables}.\\
Authors : Gane Samb LO, Modou Ngom and Tchilabalo Atozou. KPANZOU.\\
ISSBN : 978-2-9559183-1-9. Doi : 10.16929/sbs/2016.0001.\\

\noindent (2) \textit{Convergence vague - Suite de variables al\'eatoires}.\\
Authors : Gane Samb LO, Modou Ngom and Tchilabalo Atozou. KPANZOU.\\
ISSBN : 978-2-9559183-1-9. Doi : 10.16929/sbs/2016.0002.\\

\noindent (3) \noindent \textit{A course on Elementary Probability Theory}.\\
Author : Gane Samb LO.\\
ISSBN : 978-2-9559183-3-3. Doi : 10.16929/sbs/2016.0003.\\

\newpage
\noindent \textbf{Library of Congress Cataloging-in-Publication Data}\\

\noindent Gane Samb LO, 1958-\\

\noindent Weak Convergence (IIA) - Functional and Random Aspects of the Univariate Extreme Value Theory.\\

\noindent SPAS Books Series, 2018.\\



\noindent \textit{DOI} : 10.16929/sbs/2016.0009\\

\noindent \textit{ISBN}  978-2-9559183-9-5

\newpage

\noindent \textbf{Author : Gane Samb LO}\\
\bigskip

\bigskip
\noindent \textbf{Emails}:\\
\noindent gane-samb.lo@ugb.edu.sn, ganesamblo@ganesamblo.net.\\

\bigskip
\noindent \textbf{Url's}:\\
\noindent www.ganesamblo@ganesamblo.net\\
\noindent www.statpas.net/cva.php?email.ganesamblo@yahoo.com.\\

\bigskip \noindent \textbf{Affiliations}.\\
Main affiliation : Gaston Berger University, UGB, SENEGAL.\\
African University of Sciences and Technology, AUST, ABuja, Nigeria.\\
Affiliated as a researcher to : LSTA, Pierre et Marie Curie University, Paris VI, France.\\

\noindent \textbf{Teaches or has taught} at the graduate level in the following universities:\\
Saint-Louis, Senegal (UGB)\\
Banjul, Gambia (TUG)\\
Bamako, Mali (USTTB)\\
Ouagadougou - Burkina Faso (UJK)\\
African Institute of Mathematical Sciences, Mbour, SENEGAL, AIMS.\\
Franceville, Gabon\\

\bigskip \noindent \textbf{Dedicatory}.\\

\noindent \textbf{To my wife Mbaye Ndaw Fall who is accompanying for decades with love and patience}

\bigskip \noindent \textbf{Acknowledgment of Funding}.\\

\noindent The author acknowledges continuous support of the World Bank Excellence Center in Mathematics, Computer Sciences and Intelligence Technology, CEA-MITIC. His research projects in 2014, 2015 and 2016 are funded by the University of Gaston Berger in different forms and by CEA-MITIC.

\newpage

\title{Weak Convergence (IIA) - Functional and Random Aspects of the Univariate Extreme Value Theory}

\noindent \textbf{Abstract}. The univariate extreme value theory deals with the convergence in type of powers of elements of sequences of cumulative distribution functions on the real line when the power index gets infinite. In terms of convergence of random variables, this amounts to the the weak convergence, in the sense of probability measures weak convergence, of the partial maximas of a sequence of independent and identically distributed random variables. In this monograph, this theory is comprehensively studied in the broad frame of weak convergence of random vectors as exposed in Lo et al.(2016). It has two main parts. The first is devoted to its nice mathematical foundation. Most of the materials of this part is taken from the most essential Lo\`eve(1936,177) and Haan (1970), based on
the stunning theory of regular, pi or gamma variation. To prepare the statistical applications, a number contributions I made in my PhD and my Doctorate of Sciences are added in the last chapter of the last chapter of that part. Our real concern is to put these materials together with others, among them those of the authors from his PhD dissertations and Science doctorate thesis, in a way to have an almost full coverage of the theory on the real line that may serve as a master course of one semester in our universities. As well, it will help the second part of the monograph. This second part will deal with statistical estimations problems related to extreme values. It addresses various estimation questions and should be considered as the beginning of a survey study to be updated progressively. Research questions are tackled therein. Many results of the author, either unpublished or not sufficiently known, are stated and/or updated therein.\\

\noindent \textbf{Keywords.} .\\

\noindent \textbf{AMS 2010 Classification Subjects :} 60XXX; 62G30; 62G32; 60E07.

\newpage
\noindent \textbf{R\'esum\'e}. La th\'eorie des valeurs extr\^emes univari\'ees traite de la convergence en type de puissances d'\'el\'ements de suites de fonctions de r\'epartition sur l'ensemble des nombres r\'eels lorsque la puissance devient infinie. En termes de convergence des variables al\'eatoires, cela se reformule par la convergence faible, dans le sens de celle de mesures de probabailit\'e, des maxima partiels d'une suite de variables al\'eatoires r\'eelles ind\'ependantes et identiquement r\'eparties. Dans ce monographe, cette th\'eorie est \'etudi\'ee de mani\`ere exhaustive, et plac\'ee dans le cadre g\'en\'eral de la faible convergence des vecteurs al\'eatoires expos\'ee dans Lo et al. (2016). Il comporte deux parties principales. La premi\`ere est consacr\'ee \`a sa tr\`es belle fondation math\'ematique autour des notions variation r\'eguli\`ere et de variation pi and gamma. La plupart des mat\'eriaux de cette partie proviennent de Lo\`eve(1963, 1977) and de Haan (1970). Afin de pr\'eparer les applications statistiques, nous avons introduit certaines contibution à la théorie dans le dernier chapitre de cette partie. Notre v\'eritable pr\'eoccupation, ici, est de mettre ces mat\'eriaux ensemble avec d'autres, parmi eux ceux de l'auteur issus des ses th\`eses de doctorat, de mani\`ere \`a aboutir \`a une couverture unifi\'ee et presque compl\`ete de la th\'eorie sur l'ensemble des nombres r\'eels dans un forme d'un texte pouvant servir pour un cours doctoral d'un semestre dans nos universit\'es. De plus, cela aidera pour la deuxi\`eme partie. Justement, cette deuxi\`eme partie traitera de probl\`emes d'estimation statistiques li\'es aux valeurs extr\^emes. Elle aborde diverses questions d'estimation et doit \^etre consid\'er\'ee comme le d\'ebut d'une revue de r\'esults du domaine qui devra \^etre mise \`a jour r\'eguli\`erement. Des questions de recherche y sont abord\'ees. De nombreux r\'esultats de l'auteur, non publi\'es ou pas suffisamment connus, y sont \'enonc\'es et / ou mis \`a jour.\\

\newpage
\begin{center}
\LARGE
jappa t\'e bagna bayi\\
Doh t\'e bagna taxaw\\
Yokka t\'e bagna wagni
\end{center}

\maketitle

\frontmatter
\tableofcontents
\mainmatter
\Large

\chapter*{General Preface}

\noindent \textbf{This textbook} is the first of series whose ambition is to cover broad part of Probability Theory and Statistics.  These textbooks are intended to help learners and readers, of all levels, to train themselves.\\

\noindent As well, they may constitute helpful documents for professors and teachers for both courses and exercises.  For more ambitious  people, they are only starting points towards more advanced and personalized books. So, these textbooks are kindly put at the disposal of professors and learners.

\bigskip \noindent \textbf{Our textbooks are classified into categories}.\\

\noindent \textbf{A series of introductory  books for beginners}. Books of this series are usually accessible to student of first year in 
universities. They do not require advanced mathematics.  Books on elementary probability theory and descriptive statistics are to be put in that category. Books of that kind are usually introductions to more advanced and mathematical versions of the same theory. The first prepare the applications of the second.\\

\noindent \textbf{A series of books oriented to applications}. Students or researchers in very related disciplines  such as Health studies, Hydrology, Finance, Economics, etc.  may be in need of Probability Theory or Statistics. They are not interested by these disciplines  by themselves.  Rather, the need to apply their findings as tools to solve their specific problems. So adapted books on Probability Theory and Statistics may be composed to on the applications of such fields. A perfect example concerns the need of mathematical statistics for economists who do not necessarily have a good background in Measure Theory.\\

\noindent \textbf{A series of specialized books on Probability theory and Statistics of high level}. This series begin with a book on Measure Theory, its counterpart of probability theory, and an introductory book on topology. On that basis, we will have, as much as possible,  a coherent presentation of branches of Probability theory and Statistics. We will try  to have a self-contained, as much as possible, so that anything we need will be in the series.\\

\noindent Finally, \textbf{research monographs} close this architecture. The architecture should be so large and deep that the readers of monographs booklets will find all needed theories and inputs in it.\\

\bigskip \noindent We conclude by saying that, with  only an undergraduate level, the reader will  open the door of anything in Probability theory and statistics with \textbf{Measure Theory and integration}. Once this course validated, eventually combined with two solid courses on topology and functional analysis, he will have all the means to get specialized in any branch in these disciplines.\\

\bigskip \noindent Our collaborators and former students are invited to make live this trend and to develop it so  that  the center of Saint-Louis becomes or continues to be a reknown mathematical school, especially in Probability Theory and Statistics.

\chapter*{General Preface of Our Series of Weak Convergence}

\noindent \textbf{The series Weak convergence} is an open project with three categories.\\

\noindent \textbf{The special series  Weak convergence I} consists of texts devoted to the core theory of weak convergence, each of them concentrated on the handling of one specific class of objects. The texts will have labels $A$, $B$, etc. Here are some examples.\\

\noindent (1) Weak convergence of Random Vectors (IA).\\

\noindent (2) Weak convergence of stochastic processes and empirical processes (IB).\\

\noindent (3) Weak convergence of random measures (IC).\\

\noindent (4) Etc.\\

\bigskip \noindent \textbf{The special series Weak convergence II} consists of textbooks related to the theory of weak convergence, each of them concentrated on one specialized field using weak convergence. Usually, these sub-fields are treated apart in the literature. Here, we want to put them in our general frame as continuations of the Weak Convergence Series I. Some examples are the following.\\

\noindent (1) Weak laws of sums of random variables.\\

\noindent (3) Univariate Extreme value Theory.\\

\noindent (4) Multivariate Extreme value Theory.\\

\noindent (5) Etc.\\

\bigskip \noindent \textbf{The special series Weak convergence III} consists of textbooks focusing on statistical applications of Parts of the Weak Convergence Series I and Weak Convergence Series II. Examples :\\

\noindent (1) A handbook of Gaussian Asymptotic Distribution Using the Functional Empirical Process.\\

\noindent (2) A handbook of Statistical Estimation of the Extreme Value index.\\

\noindent (1) etc.\\

\part{The portal}
\chapter{Univariate Extreme Value Theory : A portal} \label{portal}

\bigskip \noindent In this chapter, we give a global picture and basic notation in univariate extreme value theory in the independent and identical distribution setting. Theoretically, readers solely interested in the statistical estimation might just read this portal and go Part II reserved to Statistical aspects.\\

\noindent We will guide the reader how to use these component if he wishes to expertise himself in this theory. But he will
find all what he needs to go with us in the statistical and probability aspects we will be dealing with in the first part.\\

\noindent As this monograph is part of our Probability and Statistics series, it should be read after the basic introduction of weak convergence theory that was fully and broadly exposed in \cite{wc-srv-ang}. This monograph is cited all along this one. For example, the broad theory of weak convergence is rounded up the first section, with a special focus on specific tools \\

\noindent First of all, we will have to deal with a sequence $X_{1},X_{2},...$ of independent copies (s.i.c) of a real random variable ($rv$) $X$ with \textit{df} $F(x)=\mathbb{P}(X\leq x)$, all of them being defined on the same probability space ($\Omega ,A,\mathbb{P}).$ For each $n\geq 1,$ the order statistics associated with the sample $X_{1},X_{2},...,X_{n}$ is denoted as $%
X_{1,n}\leq X_{2,n}\leq ...\leq X_{n,n}.$\\

\noindent The support of the distribution function $F$ is $[lep(F),uep(F)]$ where%
\begin{equation*}
lep(F)=\inf \{x\in \mathbb{R},F(x)>0\}
\end{equation*}

is called the lower endpoint of $F$ (denoted \textit{lep}) and%
\begin{equation*}
upe(F)=\sup \{x\in \mathbb{R},F(x)<1\}
\end{equation*}

\noindent is its upper endpoint (denoted \textit{uep}).\\

\bigskip \noindent In the Univariate Extreme Value Theory (UEVT), particularly in the extreme value index estimation, we often
have to use the logarithm transformation $Y=\log ^{+}X,$ with $\log^{+}x=\left( \log x\right) 1_{(x>0)}$ where $1_{A}$ stands for the indicator
function of the set $A$.\\

\noindent It would be better and simpler to work with the logarithm function in place of $\log^{+}$. So we may add an additional assumption based on the following arguments. Indeed, we are mainly concerned by extreme values of the observations near $uep(F)$. If $uep(F)\leq 0$, we translate all the data by a real $t>0$ large enough to ensure that $uep(F)+t>0.$ The translated data are associated to the distribution function $F(\circ -t)$ with upper
endpoint $uep(F)+t>0$. In both cases $uep(F)\leq 0$ and $uep(F)>0$, we may assume that, eventually at the cost translating the data, that the
extreme values of the observations are positive. Based on these facts, we may and do assume from now that the random
variable is positive that is $X>0$ and we simply write $Y=\log X$.\\

\noindent Hence $Y_{1},Y_{2},...$ will be a sequence of independent copies of $Y$ with \textit{df} $G(y)=F(e^{y})$, $y\in \mathbb{R}$ and $Y_{1,n}\leq Y_{2,n}\leq ...\leq Y_{n,n}$ are the order statistics associated with $Y_{1},Y_{2},...,Y_{n}$.

\section{Weak convergence}

Extreme value theory in $\mathbb{R}$ begins with the knowledge of the asymptotic law of the sequence of the maxima $X_{n,n}$ when $n$ tends to
infinity under the hypothesis that the observations are independent and identically distributed. This means that the extreme value theory with this
respect is part of the weak convergence theory.\\

\bigskip \noindent The weak convergence theory in $\mathbb{R}^k$, $k\geq 1$ is exposed in details in Chapter 2 in \cite{wc-srv-ang}, and more specifically in Chapter 4 of the same textbook, for real random variables. All the needs of the \textit{UEVT} regarding weak convergence are to be found there. In this particular case, the link of this theory with generalized inverses in $\mathbb{R}$ is described details in Chapter 4, \cite{wc-srv-ang} ref{cv.R}, along with other useful and specific tools.\\

\bigskip \noindent The general of weak convergence in metric spaces is provided in Chapter 2 in \cite{wc-srv-ang}, this chapter being largely inspired by the book of \cite{billingsley} and that of \cite{vaart}.\\

\bigskip \noindent Here, we are going only to give the main tools we have need to ensure a linear reading of the book.\\

\bigskip \noindent To introduce the weak convergence in $\mathbb{R}$, consider a sequence of real random variables $Z_{1},Z_{2},...$
with distributions functions $H_{1},H_{2},...$

\begin{definition} \label{portal.def1} The sequence of real random variables $Z_{1},Z_{2},...$ with
distributions functions $H_{1},H_{2},...$ converges weakly to a real random
variable $Z$ with distribution function $H$ iff and only if one of these
equivalent assertions holds.\\

\noindent  (i) For any continuous and bounded functions $f:\mathbb{R}%
\rightarrow \mathbb{R},$%
\begin{equation*}
\lim_{n\rightarrow +\infty }\int fd\mathbb{P}_{Z_{n}}=\int fd\mathbb{P}_{Z}.
\end{equation*}

\noindent (ii) For any continuity point $x$ of $H,$%
\begin{equation*}
\lim_{n\rightarrow +\infty }H_{n}(x)=H(x).
\end{equation*}%

\noindent (iii) For any $u\in \mathbb{R},$%
\begin{equation*}
\lim_{n\rightarrow +\infty }\int e^{iux}d\mathbb{P}_{Z_{n}}(x)=\int e^{iux}fd%
\mathbb{P}_{Z}.
\end{equation*}
\end{definition}

\bigskip \noindent When $(Z_{n})_{n\geq 1}$ weakly converges to $Z$, we mainly use the notation%
\begin{equation*}
Z_{n}\rightsquigarrow Z\text{ or }H_{n}\rightsquigarrow H
\end{equation*}

\bigskip \noindent and we may also use $Z_{n}\rightarrow _{w}Z$ ($w$ standing for weakly) or $%
Z_{n}\rightarrow _{d}Z$ ($d$ standing for : in distribution). We also shift
to the distribution functions and say : $(H_{n})_{n\geq 1}$ weakly converges
to $H$.\\

\noindent Point (i) is the main definition. Points (ii) and (iii) are parts of what is called the Portmanteau Theorem 
(See  Theorem 2 for the general case and Theorem 3  for the particular case of $\mathbb{R}^k$, Chapter 2, \cite{wc-srv-ang}).\\

\noindent An interesting property is that the convergence of the distribution functions is uniform when $H$ is continuous, in the particular case of $\mathbb{R}$. This gives : \\

\begin{proposition} \label{portal.prop1} The sequence of real random variables $Z_{1},Z_{2},...$ with
distributions functions $H_{1},H_{2},...$ converges weakly to a real random
variable $Z$ with a continuous distribution function $H$ iff and only if
\begin{equation*}
\lim_{n\rightarrow +\infty} \sup_{x\in \mathbb{R}} |H_{n}(x)-H(x)|=0.
\end{equation*}
\end{proposition}

\bigskip \noindent The proof of this is given in Point (5) in Chapter 4, \cite{wc-srv-ang}.\\

\bigskip \noindent When all the $Z_{1},Z_{2},..$ and $Z$ are absolutely continuous with respect
to a $\sigma $-finite measure in $\mathbb{R}$, with Radon-Nikodym derivatives, denoted by $f_{Z_{n}},$ $n\geq 1$ and $%
f_{Z}$,  we have the following result.

\begin{theorem} \label{portal.theo1} (Sch\'{e}ff\'{e}). If for any $x\in \mathbb{R},$%
\begin{equation*}
\lim \lim_{n\rightarrow +\infty }f_{Z_{n}}(x)=f_{Z}(x),\text{ }\nu .a.e,
\end{equation*}

\bigskip \noindent then the sequence $(Z_{n})_{n\geq 1}$ weakly converges to $Z$ and 
\begin{equation*}
\sup_{B\in \mathcal{B}(\mathbb{R})}\left\vert \mathbb{P}_{Z_{n}}(B)-\mathbb{P%
}_{Z}(B)\right\vert =\frac{1}{2}\int \left\vert f_{Z_{n}}-f_{Z}\right\vert
d\nu \rightarrow 0\text{ as }n\rightarrow \infty .
\end{equation*}
\end{theorem}

\bigskip \noindent \textbf{Proof}. See proof of Theorem 4 in Chapter 2, in Lo \cite{wc-srv-ang}, or \cite{billingsley}.\\

\bigskip \noindent This theorem of Sch\'{e}ff\'{e} is very handy when one needs to find rate of
convergence with presence of probability densities.\\

\bigskip \noindent But weak convergence in $\mathbb{R}$ may be done entirely with inverses functions.
Define the inverse function of the distribution function H by
\begin{equation}
H^{-1}(u)=\inf \{x\in \mathbb{R},F(x)\geq u\},0\leq u\leq 1.
\label{inverseH}
\end{equation}

\bigskip \noindent If $H$ is right-continuous and nondecreasing, $H^{-1}$ is also nondecreasing but left-continuous, as proved in Point 3 in Chapter 4 in \cite{wc-srv-ang}. Let us point right now these two equivalence formulas :

\begin{equation}
H^{-1}(u)\leq t\Longleftrightarrow u\leq F(t)  \label{portal.invprop1}
\end{equation}%

\bigskip \noindent and
\begin{equation}
H^{-1}(u)>t\Longleftrightarrow u>F(t). \label{portal.invprop2}
\end{equation}

\bigskip \noindent In Chapter 4 in \cite{wc-srv-ang}, we expose and prove a long list of properties of
generalized function of nondecreasing funtions not necessarily with values in $[0,1]$. We will discover these properties when needed. Right now we need the Points 4 and 5 in just mentioned section. The first is

\begin{lemma} \label{portal.lem1} $\ $Let $(H_{n})_{n\geq 1}$ weakly converges to $H$, then $(H_{n}^{-1})_{n\geq 1}$ weakly converges to $H^{-1},$ that is for any continuity point of $H^{-1},$%
\begin{equation*}
H_{n}^{-1}(x)\rightarrow H(x)\text{ as }n\rightarrow +\infty .
\end{equation*}
\end{lemma}

\noindent The second is 

\begin{lemma} \label{portal.lem2} $(H_{n})_{n\geq 1}$ weakly converges to a continuous distribution
function $H$ of ans only if
\begin{equation*}
\sup_{x\in \mathbb{R}}\left\vert H_{n}(x)-H(x)\right\vert \rightarrow 0\text{
as }n\rightarrow +\infty .
\end{equation*}%
\end{lemma}

\bigskip \noindent The generalized inverse transform or the quantile transform is instrumental
in all parts of extreme value theory and in its statistical branch.\\

\bigskip \noindent In the frame we are constructing in this chapter, the logarithm transform $Y=\log X$ plays an important role, specially in the statistical part. It implies
\begin{equation}
G^{-1}(u)=\log F^{-1}(u),0\leq u\leq 1\text{.}  \label{logTransf}
\end{equation}

\bigskip \noindent Throughout the text, the $G$ stands for the distribution function of $Y$.

\section[The three nondegenerated extreme values distributions]{The three nondegenerated extreme values distributions and the generalized extreme value distribution}

\subsection{Convergence in type}

\begin{definition} \label{portal.def2} The sequence of real random variables $Z_{1},Z_{2},...$ with
distributions functions $H_{1},H_{2},...$ converges in type to a real random
variable $Z$ with distribution function $H$ iff and only if there exists a
sequence positive real numbers $(a_{n}>0)_{n\geq 0}$ and a sequence of real
numbers $(b_{n})_{n\geq 0}$ such that one of these assertions hold :\\

\noindent (i) We have

\begin{equation*}
\frac{Z_{n}-b_{n}}{a_{n}}\rightsquigarrow Z
\end{equation*}%

\noindent (ii) For any continuity point $x$ of $H,$%
\begin{equation*}
H_{n}(a_{n}x+b_{n})\longrightarrow H(x).
\end{equation*}
\end{definition}

\bigskip \noindent By Lemma \ref{portal.lem1} above, we have the following\\

\begin{lemma} \label{portal.lem3} The sequences of real random variables $Z_{1},Z_{2},...$ with distributions functions $%
H_{1},H_{2},...$ converges in type to $Z$ with distribution function $H,$
that is there exist a sequence positive real numbers $(a_{n}>0)_{n\geq 0}$
and a sequence of real numbers $(b_{n})_{n\geq 0}$(ii) such that for any continuity
point $x$ of $H,$%
\begin{equation*}
H_{n}(a_{n}x+b_{n})\longrightarrow H(x).
\end{equation*}%
Then for any continuity $x$ point $H^{-1},$ we also have%
\begin{equation*}
\frac{H_{n}^{-1}(x)-b_{n}}{a_{n}} \longrightarrow H^{-1}(x).
\end{equation*}
\end{lemma}

\bigskip \noindent The next lemma will allow to define the \textbf{convergence in type}, that is uded in UEVT.

\begin{lemma} \label{portal.lem4}   Let $(H_{n})_{n\geq 0}$ be a sequence of probability
distribution functions. Suppose there exist sequences $(a_{n}>0)_{n\geq 0},$ $%
(\alpha _{n}>0)_{n\geq 0},$ $(b_{n})_{n\geq 0}$ and $(\beta _{n})_{n\geq 0}$, and
probability distributions functions $H_{1}$ and $H_{2}$ such that%
\begin{equation}
\lim_{n\rightarrow \infty }F_{n}(a_{n}x+b_{n})=H_{1}(x),\text{ }x\in C(H_{1})
\label{conv1}
\end{equation}%

\bigskip \noindent and

\begin{equation}
\lim_{n\rightarrow \infty }F_{n}(\alpha _{n}x+\beta _{n})=H_{2}(x),\text{ }%
x\in C(H_{2}).  \label{conv2}
\end{equation}%

\bigskip \noindent Then there exist reals numbers $A>0$ and $B$ such that, as $n\rightarrow
\infty ,$

\begin{equation}
\alpha _{n}/a_{n}\rightarrow A\text{ and }(\beta
_{n}-b_{n})/a_{n}\rightarrow B,  \label{conditionType}
\end{equation}%
$as$ $n\rightarrow +\infty ,$and for any $x\in \mathbb{R},$%
\begin{equation}
H_{2}(x)=H_{1}(Ax+B).  \label{relType}
\end{equation}%

\bigskip \noindent Inversely, If (\ref{conv1}) and (\ref{conditionType}) hold both, then
(\ref{conv2}) is true, where $H_{2}$ defined in \ref{relType}.
\end{lemma}

\bigskip \noindent Formula (\ref{relType}) defines an equivalence class in the class of all real
probability distribution functions. Let us denote this equivalence relation by $\mathcal{R}_{type}$. And we say that $H_{1}$ and $H_{2}$ are of the same type if one is obtained from the other by affine transformation of the argument. The lemma says that the limit of convergence is unique in type,
meaning that all the possible limits in type are of the same type.\\

\bigskip \noindent The proof the lemma is given is Lemma \ref{evt.lem.1} in Chapter \ref{evt} following
the lines of \cite{resnick}.\\

\noindent Now, we are going to apply this to the special case of the sequences of the maxima $X_{n,n}$ of the samples.

\section{The three nondegenerated extreme values distributions}

\noindent Let us begin by the important Gnedenko theorem which actually covered more than the stability of the maximum but also concerned that of the sums.\\

\subsection{Gnedenko Theorem}


\bigskip \noindent Before we state the theorem, recall that random variable is degenerated if and only it is concentrated on a single point almost surely.

\begin{theorem} \label{portal.theo2} (Gnedenko) Let $X_{1},X_{2},...$ be independent copies (s.i.c) of a real random variable ($rv$) $X$ with \textit{df} $F(x)=\mathbb{P}(X\leq x)$, all of them
being defined on the same probability space ($\Omega ,A,\mathbb{P}).$
Define, for each $n\geq 1,$ $X_{n,n}=\max (X_{1},X_{2},..,X_{n})$.\\

\noindent Then $X_{n,n}$ in type to nondegenerated random variable $Z$ with distribution
function $H,$ that is : there exist sequences $(a_{n}>0)_{n\geq 0}$ and $%
(b_{n})_{n\geq 0}$ such that 
\begin{equation}
\frac{X_{n,n}-b_{n}}{a_{n}}\rightsquigarrow Z  \label{EXTDOM1}
\end{equation}

\bigskip \noindent or equivalently, for any continuity $x$ point of $H$,%
\begin{equation}
\lim_{n\rightarrow \infty }F^{n}(a_{n}x+b_{n})=H(x),  \label{EXTDOM2}
\end{equation}

\bigskip \noindent if and only $H$ is one of the three types of distributions :\\

\noindent \textbf{The type of Gumbel:} 
\begin{equation}
\Lambda (x)=\exp (-\exp (-x)),x\in R  \label{T1}
\end{equation}

\bigskip \noindent \textbf{The type of Fr\'{e}chet of parameter }$\alpha >0:$%
\begin{equation}
\varphi _{\alpha }(x)=\exp (-x^{-\alpha })1_{(x\geq 0)}.  \label{T2}
\end{equation}

\bigskip \noindent \bigskip \noindent \textbf{The type of Weibull of parameter }$\beta >0:$%
\begin{equation}
\psi _{\beta }(x)=\exp (-(-x)^{\beta })1_{(x<0)}+1_{(x\geq 0)}.  \label{T3}
\end{equation}
\end{theorem}

\bigskip \noindent Throughout this text, random variables respectively associated with the distribution functions given in (\ref{T1}), (\ref%
{T2}) and (\ref{T3}) will be denoted as $\Lambda$,  $F_{r}(\alpha)$, $\alpha>0$ and $W(\beta)$, $\alpha>0$, in the same order.\\

\bigskip \noindent Proof. The full proof is given in that of Theorem \ref{evt.theo1} in Chapter \ref{evt}. The proof is a direct rephrasing of the one given in \cite{resnick}, which itself is rather a classical one. But we complete the proofs by giving full details of the solutions the Hamel equations along with the principles of Littlewood in Chapter \ref{funct}.\\

\bigskip \noindent  
When  (\ref{EXTDOM1}) or (\ref{EXTDOM2}) hold, it is said that the maxima $%
X_{n,n}\ $attracted to $Z$ or $F$ is attracted to $H$ or $F$ is in the
attraction domain of $H$ denoted 
\begin{equation*}
F\in D(H).
\end{equation*}

\bigskip \noindent To avoid confusion with the attraction domains concerning infinitely divisible
laws using sums of independent random variables, we use the term of extreme
domain of attraction. We have this simple result that does not need to be
proved.\\

\begin{proposition}
\bigskip Two distributions of the same type have the same extreme domain of attraction.
\end{proposition}

\bigskip \noindent So the Theorem of Gnedenko says that the only three nondegenerated extreme
domains of attractions are those of (\ref{T1}), (\ref{T2}) and (\ref{T3}). We notice that these distributions functions are continuous. This implies, in view of Lemma \ref{portal.lem2} above, that the convergence 
(\ref{EXTDOM2}) holds uniformy in $x$ on $\mathbb{R}$. And this gives :\\

\begin{proposition}
$F\in D(H)$ where $H$ $\in \{\Lambda ,\varphi _{\alpha },\psi _{\beta }\}$
if and only if  there exist sequences $(a_{n}>0)_{n\geq 0}$ and $%
(b_{n})_{n\geq 0}$ such thator equivalently, for any continuity point of $H$,%
\begin{equation}
\lim_{n\rightarrow \infty }\sup_{x\in \mathbb{R}}\left\vert
F^{n}(a_{n}x+b_{n})-H(x)\right\vert =0.
\end{equation}
\end{proposition}

\subsection{Immediate and simple examples} \label{portal.subsec_examples_simples}

We are going to use the convergence of the probability distribution functions. Since the three possible limits the Gnedenko's Theorem are continuous, we have to check the convergence at any point of $\mathbb{R}$\newline

\noindent Let us consider three simple examples.\\

\bigskip \noindent \textbf{(a)} $X \sim \mathcal{E}$,\\

\noindent that is $X$ is standard exponential random variable with probability distribution function 
\begin{equation*}
F(x)=(1-\exp (-x))1_{(x\geq 0)}, \text{ } x\in \mathbb{R}.
\end{equation*}

\bigskip  By using the distribution functions, we want to prove that 

\begin{equation*}
M_{n}-\log n\overset{d}{\rightarrow }\Lambda \text{ as } x\rightarrow +\infty. \label{portal.example1}
\end{equation*}

\noindent Indeed, we have

\begin{equation*}
P(M_{n}-\log n\leq x)=P(M_{n}\leq x+\log n)=F(M_{n}\leq x+\log n)^{n}.
\end{equation*}

\bigskip \noindent But for any $x\in \mathbb{R},$ $x+\log n\geq 0$ for $n\geq
\exp (-x).$ Then for large values of $n$, $P(M_{n}\leq x+\log n)=(1-\exp
(-x-\log n))$ and next for any $x\in \mathbb{R}$ and for $n$ large enough,

\begin{equation*}
P(M_{n}-\log n\leq x)=(1-\frac{e^{-x}}{n})\rightarrow e^{-e^{-x}}=\Lambda(x).
\end{equation*}

\noindent We conclude that $X \in D(\Lambda)$.\\

\bigskip \noindent \textbf{(b)} $X \sim \mathcal{P}ar(\alpha)$, $\alpha>0$,\\

\noindent that is $X$ is a Pareto random variable with parameter $\alpha>0$ with probability distribution function
\begin{equation*}
F(x)=(1-x^{-\alpha }))1_{(x\geq 1)},  \text{ } x\in \mathbb{R}
\end{equation*}

\bigskip \noindent We have to prove that 

\begin{equation*}
n^{-1/\alpha }M_{n}\overset{d}{\rightarrow }C(\alpha) \text{ as } x\rightarrow +\infty. \label{portal.example2}
\end{equation*}

\bigskip \noindent The observation $X_i$ are nonnegative since the support of a $\mathcal{P}ar$ law is $\mathbb{R}_{+}$. So the maxima $M_{n}$ are nonnegative for any $n\geq 1$. We may discuss two cases.\newline

\bigskip \noindent Case $x\leq 0$. In this case, we have
\begin{equation*}
P(n^{-1/\alpha }M_{n}\leq 0)=0=\varphi _{\alpha }(x),
\end{equation*}

\bigskip \noindent and then \ref{portal.example2} holds.\newline

\bigskip \noindent Case $x>0$. In this case 
\begin{equation*}
P(n^{-1/\alpha }M_{n}\leq x)=P(M_{n}\leq n^{1/\alpha }x).
\end{equation*}

\bigskip \noindent For large values of $n$, we have $n^{1/\alpha }x>1$ (take $n\geq (1/x)^{-\alpha}$, to ensure that) and for these values, 
\begin{eqnarray*}
P(n^{-1/\alpha }M_{n} &\leq &x)=F(n^{1/\alpha }x)^{n}=(1-(n^{1/\alpha
}x)^{-\alpha })^{n} \\
&=&(1-\frac{x^{-\alpha }}{n})^{n}\rightarrow \exp (-x^{-\alpha }) \\
&=&\varphi _{\alpha }(x).
\end{eqnarray*}

\bigskip \noindent So \ref{portal.example2} also for this case. Then it holds for any $x\in \mathbb{R}$.\\

\noindent Conclusion : $X \in D(FR(\alpha))$.\\

\bigskip \noindent \textbf{(c)} $X \sim \mathcal{U}(0,1)$, 

\noindent that is $X$ is uniformly distributed on $(0,1)$ with probability distribution function : 
\begin{equation*}
F(x)=x1_{(0\leq x\leq 1)}+1_{(x\geq 1)}, \text{ } x\in \mathbb{R}.
\end{equation*}

\bigskip \noindent We want to prove that

\begin{equation*}
n(M_{n}-1)\overset{d}{\rightarrow }W(1) \text{ as } x\rightarrow +\infty. \label{portal.example3}
\end{equation*}

\bigskip \noindent We have
\begin{equation*}
P(n(M_{n}-1)\leq x)=F(1+\frac{x}{n})^{n}.
\end{equation*}

\noindent We have two cases.\newline

\bigskip \noindent Case $x\geq 0$. We see that $1+x/n$ is nonnegative $n\geq 1$ et 
\begin{equation*}
P(n(M_{n}-1)\leq x)=F(1+\frac{x}{n})^{n}=1=\psi _{1}(x)
\end{equation*}

\bigskip \noindent and we see that \label{portal.example3} holds.\\

\noindent Case $x<0$. For large values of $n$, we have $0\leq 1+x/n\leq 1$ (take $x\geq -n$ $i.e.$ $n\geq -(x)\geq 0$, to get it) and for these values of $n$, 
\begin{eqnarray*}
P(n(M_{n}-1) &\leq &x)=F(1+\frac{x}{n})^{n} \\
&=&(1+\frac{x}{n})^{n}\rightarrow e^{x}=\psi _{1}(x).
\end{eqnarray*}

\bigskip \noindent We get that for any $x\in \mathbb{R}$,%
\begin{equation*}
P(n(M_{n}-1)\leq x)\longrightarrow \psi _{1}(x).
\end{equation*}

\noindent Conclusion : $W \in D(W(1))$.\\

\bigskip \subsection{The Generalized Extreme Value (GEV) Distribution}

Now we know the three types of extreme value distributions from Theorem \ref{portal.theo1}, let
us try to gather them into one form with the help of Lemma \ref{portal.lem2}.

\bigskip \noindent Let us change our sequences  $(a_{n}>0)_{n\geq 1}$ and $(b_{n})_{n\geq 1}$
by $(\alpha _{n}>0)_{n\geq 1}$ and $(\beta _{n})_{n\geq 1}$.\\

\noindent For the case of $H_{1}=\varphi _{\alpha }$, take $A=\gamma =1/\alpha $ and $B=1$. We get $\gamma >0$ and
\begin{equation*}
H_{2}(x)=G_{\gamma }(x)=\exp (-(1+\gamma x)^{-1/\gamma })1_{(1+\gamma x\geq0)}.
\end{equation*}

\bigskip \noindent For $H_{1}=\psi _{\beta}$, choose $A=-\gamma =-1/\beta $ and $V=-1$. We get $\gamma <0$ and 
\begin{equation*}
H_{2}(x)=G_{\gamma }(x)=\exp (-(1+\gamma x)^{-1/\gamma })1_{(1+\gamma x\geq
0)}+1_{(1+\gamma x<0)}.
\end{equation*}

\bigskip \noindent Finally, for $H_1(x)=exp(-e^{0-x})$, we may use the following the limit  
\begin{equation*}
\Lambda(x)=\lim_{\gamma \rightarrow 0}\exp (-(1+\gamma x)^{-1/\gamma })
\end{equation*}

\noindent and define $\Lambda(x)$ as an extension of  $\exp (-(1+\gamma x)^{-1/\gamma })$ by continuity at $\gamma =0$, and write

\begin{equation*}
\Lambda (x)=G_{0}(x),x\in \mathbb{R}.
\end{equation*}

\bigskip \noindent We are now able to gather the whole extreme domain of attraction by one
parameterized extreme value distribution 
\begin{equation*}
G_{\gamma}=\exp (-(1+\gamma x)^{-1/\gamma }),\text{ } 1+\gamma x\geq 0,
\end{equation*}

\bigskip \noindent where $\gamma <0$ corresponds to the Frechet domain, $\gamma =0$ to the
Gumbel domain and $\gamma >0$ to Weibull domain. Throughout the text, we
will use these compact form.

\subsection{Malmquist and uniform representation} \label{portal.uniformrep}

\noindent We are introducing  important representation tools that can greatly help to handle the extremes.\\

\bigskip \noindent \textbf{A - Uniform representations}.\\

\noindent We already introduced the generalized inverse (see (\ref{inverseH}) above). We are
going to use it for studied sequence of real random variables $X_1$, $X_2$, ... By (\ref{portal.invprop1}) above, we
have this equivalence

\begin{equation}
\forall (t,s)\in \mathbb{R}\times \lbrack 0,1],\text{ }\left( F^{-1}(s)\leq
t\right) \Leftrightarrow \left( s\leq F(t\right) ).  \label{invEquiv}
\end{equation}

\bigskip \noindent Now, let $U\sim \mathcal{U}(0,1)$  following a uniform law on $(0,1)$
defined on ($\Omega ,\mathcal{A},\mathbb{P}).$ It iseasy to see that $V=1-U$
follows also a uniform law on $(0,1).$ By (\ref{invEquiv}), we have%
\begin{equation*}
\forall (x\in \mathbb{R}),\text{ }\mathbb{P}(F^{-1}(V)\leq x)=\mathbb{P}%
(V\leq F(x))=F(x).
\end{equation*}

\bigskip \noindent So, by the characterization of the law of a real random variable by its
distribution function, we see that $X$ and $F^{-1}(1-U)$ have the same distrution, denoted as 

\begin{equation*}
X=_{d}F^{-1}(1-U).
\end{equation*}

\bigskip \noindent This equality in laws enables to replace the whole sequence $X_{1},X_{2},...$
that is studied in this text by in this text by a sequence $F^{-1}(1-U_{1}),$
$F^{-1}(1-U_{2}),$..., where $U_{1},U_{2},...$ are independent and uniform
random variables on $(0,1)$. And we have the equalities in distribution of
stochastic processes : 
\begin{equation*}
\left\{ X_{j},j\geq 1\}=_{d}\{F^{-1}(1-U_{j}),j\geq 1\right\} ,
\end{equation*}

\begin{equation*}
\{\left\{ X_{1,n},X_{2,n},...X_{n,n}\right\} ,n\geq 1\}
\end{equation*}

\begin{equation*}
=_{d}\left\{
\{F^{-1}(1-U_{n,n}),F^{-1}(1-U_{n-1,n}),...,F^{-1}(1-U_{1,n})\},n\geq
1\right\} ,
\end{equation*}

\begin{equation*}
\left\{ Y_{j},j\geq 1\}=_{d}\{\log F^{-1}(1-U_{j}),j\geq 1\right\} ,
\end{equation*}

\bigskip \noindent and 
\begin{equation*}
\{\left\{ X_{1,n},X_{2,n},...X_{n,n}\right\} ,n\geq 1\}
\end{equation*}

\begin{equation*}
=_{d}\left\{ \{\log
F^{-1}(1-U_{n,n}),\log F^{-1}(1-U_{n-1,n}),...,\log
F^{-1}(1-U_{1,n})\},n\geq 1\right\} .
\end{equation*}

\bigskip \noindent With these representations that preserve the laws of the concerned
stochastic processes, all the exact laws and asympotic laws using solely
the probability laws are true for both processes $\left\{X_{j},j\geq 1\right\}$ and $\left\{F^{-1}(1-U_{j}),j\geq 1\right\} $.\\

\noindent To these elements, we add the Malmquist representation.\\

\noindent \textbf{B - Malmquist representation}.\\

\noindent As announced, we shall use the following Malmquist representation.

\begin{equation}
\{\log (\frac{U_{j+1,n}}{U_{j,n}})^{j},j=1,...,n\}=_{d}\{E_{1,n},...,E_{n,n}\}, \label{malquist}
\end{equation}

\bigskip \noindent where $E_{1,n},...,E_{n,n}$ is an array of independent standard exponential
random variables. We write $E_{i}$ instead of $E_{i,n}$ for simplicity sake.\\ 

\noindent Details and proof in Proposition 31, Chapter 4, \cite{wc-srv-ang}.\\

\subsection{Regular variation and $\pi$-variation}

\bigskip \noindent The notions of regular variation and $\pi$-variation are very useful tools in
extreme value theory and its applications. We may treat these notion in the
neighbourhood of zero or in that of $+\infty$. Both approaches are
equivalent and on may move from one to the other by the inverse transform. So
let
\begin{equation*}
\begin{tabular}{llll}
$S:$ & $(a_{0},+\infty \lbrack $ & $\longmapsto $ & $R$ \\ 
& $x$ & $\hookrightarrow $ & $S(x)$%
\end{tabular}%
\end{equation*}%

\noindent be a measurable function, that is integrable on compact sets, where $a_{0}>0,$
is some positive real number.\\
 
\subsubsection{Regular Variation}.

\begin{definition} \label{portal.def3} The function $S$ is said to be regularly varying with exponent $\rho \in R$
at infinity, and we denote $S\in RV(\rho,+\infty)$ if and only if for all $\mu \in
R_{+}$%
\begin{equation}
\lim_{x\rightarrow +\infty }\frac{S(\mu x)}{S(x)}=\mu ^{\rho }.
\label{RVDEF}
\end{equation}

\noindent If $\rho=0$, it is said that $S$ is slowly varying at infinty and we write $S \in SV(+\infty)$.
\end{definition}

\bigskip \noindent In this formula, only the final values of $x$ matter and these values are in 
$(a_{0},+\infty \lbrack $ for large values of $x$.\\

\bigskip \noindent In Chapter \ref{evt}, Section \ref{evt.sec.rvsv}, a complete theory of regularly varying functions is exposed there, at least all we want on them here.\\


\noindent Before we proceed further, suppose that we denote $s(u)=S(1/u)$, $0<u<\min (1,1/a_{0})$ We have the regularly variation in the neighborhood
at $0,$ from \ (\ref{RVDEF}), as follows
\begin{equation*}
\lim_{x\rightarrow +\infty }\frac{s(\mu x)}{s(x)}=\mu ^{-\rho },
\end{equation*}

\noindent and we write $s \ in RV(-\rho,0)$ and $s \ in SV(0)$ of $\rho=0$.\\

\bigskip \noindent Among interesting results, we will have these points

\bigskip \noindent \textbf{Karamata representation}. The function $S$ is a regularly varying function at +$\infty $
with exponent $\rho $ is and only if there exists a constant $c>0$ and
there\ exist functions $p(u)$ and $b(u)$ of $u\rightarrow $ $(a_{0},+\infty
\lbrack $ satisfying%
\begin{equation*}
(p(x),b(x))\rightarrow (0,0)\text{ as }x\rightarrow +\infty ,
\end{equation*}

\noindent such that $S$ admits the following representation of Karamata%
\begin{equation*}
S(x)=c(1+p(x))x^{\rho }\exp (\int_{x}^{+\infty }\frac{b(y)}{y}dy).
\end{equation*}

\bigskip \noindent \textbf{$\pi$-variation}.\\

\noindent This variation is more adapted for the neighboorhod of zero.\\

\bigskip \noindent Suppose that $s$ is slowly varying at infinity, which is equivalent to saying
that $s(\circ)$\ is slowly varying at zero.  Define
$$
U(t)=s(t) + \int_{0}^{u_0} \frac{s(u)}{du}, \ \ 0\leq t\leq u_0.
$$

\noindent By using the Karamata representation, one may readily prove, for any $\mu >0$ and for any $x >0,x \neq 1$, that 

\begin{equation}
\lim_{u\rightarrow 0} \frac{U(\mu u)-U(u)}{s(u)}=\log \mu   \label{PV1}
\end{equation}

\noindent and
\begin{equation}
\lim_{u\rightarrow 0} \frac{U(\lambda u)-U(u)}{s( xu)-s(u)}=\frac{\log \mu }{\log
x }.  \label{PV2}
\end{equation}

\bigskip \noindent This has been established in the proof of Theorem \ref{evt.dehaan.rep}, Section \ref{evt.sec.rvsv}, Chapter \ref{evt}, 
in the part $(c) \Rightarrow (d)$.\\

\bigskip \noindent Both of these formulas are particular cases of this definition.\\

\begin{definition} (Definition-Theorem) \label{portal.def4} A function $T(u)$ of $u\in (0,1)$ is of $\pi $-variation if and
only of one these two propositions holds.\\

\noindent (a) There exists a slowly varying function on $(0,1)$ at zero such that $%
\mu >0$ and for any $\kappa >0,$%
\begin{equation*}
\lim \frac{T(\mu u)-T(\kappa u)}{s(u)}=\log (\mu /\kappa ).
\end{equation*}
\end{definition}

\noindent (b) For any $\mu >0$ and for any $\kappa >0,\kappa \neq 1$,
\begin{equation}
\frac{s(\lambda u)-s(u)}{s(\kappa u)-s(u)}=\frac{\log \mu}{\log\kappa }.  \label{PV22}
\end{equation}

\bigskip \noindent Slowly varying functions have interesting uniform convergence properties, both
in deterministic and random frames. We have :\\

\begin{proposition} \label{portal.def5} Let  $S(u)$ be a function $u\in (0,1)$ that is slowly varying at zero. We have 
the following uniform convergence in deterministic and random versions.\\

\bigskip \noindent (a) Let $A(h)$ and $B(h)$ two functions of $h\in (0,+\infty \lbrack $ such that for each $h\in (0,+\infty \lbrack $, we have $%
0<A(h)\leq B(h)<+\infty $ and  $(A(h),B(h))\rightarrow (0,0)$ as $h\longrightarrow 0$. Suppose that there exist two real numbers $A$ and $B$
satisfying $0<C<D<+\infty$ such that 
\begin{equation}
C<\lim \inf_{h\rightarrow +\infty }A(h)/B(h)=A>0,\text{ }\lim
\sup_{h\rightarrow +\infty }B(h)/A(h)<B.  \label{C1}
\end{equation}

\noindent Then, we have

\begin{equation*}
\lim_{h\rightarrow +\rightarrow }\sup_{A(h)\leq u,v\leq B(h)}\left\vert 
\frac{S(u)}{S(v)}-1\right\vert =0.
\end{equation*}

\bigskip \noindent (b)  Let $A(h)$ and $B(h)$ two families, indexed by $h\in
(0,+\infty \lbrack ,$ of real-valued applications defined on a probability
space $(\Omega ,\mathcal{A},\mathbb{P})$  such that for each $h\in
(0,+\infty \lbrack $, we have $0<A(h)\leq B(h)<+\infty .$ Suppose that there
exist two families $A^{\ast }(h)$ and $B^{\ast }(h)$, indexed by $h\in
(0,+\infty \lbrack ,$ of $\mathbf{measurable}$ real-valued applications
defined on $(\Omega ,\mathcal{A},\mathbb{P})$ such that
for each  $h\in (0,+\infty \lbrack ,$ $A^{\ast }(h)\leq A(h)\leq B(h)\leq
B^{\ast }(h)$, and such that

\begin{equation}
\lim \sup_{h\rightarrow +\infty }\lim_{\lambda \rightarrow +\infty }\inf 
\mathbb{P}(B^{\ast }(h)/A(h)>\lambda )=0.  \label{C2R}
\end{equation}

\noindent and

\begin{equation}
\lim \sup_{h\rightarrow +\infty }\lim_{\lambda \rightarrow +\infty }\inf 
\mathbb{P}(A^{\ast }(h)/B^{\ast }(h)<1/\lambda )=0.  \label{C3R}
\end{equation}

\noindent We say that the family $\{B^{\ast }(h),h\in h\in (0,+\infty \lbrack \}$ is
asymptotically bounded in probability against $+\infty $ and the family $\{B^{\ast }(h),h\in h\in (0,+\infty \lbrack \}$ is asymptotically bounded in probability against $0$ and accordingly, we say that the family $\{B(h),h\in h\in (0,+\infty\lbrack \}$ is asymptotically bounded in outer probability against $+\infty$ and the family $\{A(h),h\in h\in (0,+\infty \lbrack \}$ is asymptotically bounded in outer probability against $0$.\\

\noindent Then any $\eta >0,$ for any $\delta >0$, there exists a measurable subset $\Delta (\delta )$ of such
that 
\begin{equation*}
\left( \sup_{A(h)\leq u,v\leq B(h)}\left\vert \frac{S(u)}{S(v)}-1\right\vert
>\eta \right) \subset \Delta (\delta ),
\end{equation*}

\noindent with
\begin{equation*}
\mathbb{P}(\Delta (\delta ))\leq \delta .
\end{equation*}

\noindent Consequently, if the quantities
\begin{equation*}
\sup_{A(h)\leq u,v\leq B(h)}\left\vert \frac{S(u)}{S(v)}-1\right\vert >\eta 
\end{equation*}

\noindent are measurable for $h\in h\in (0,+\infty \lbrack $, we have that%
\begin{equation*}
\sup_{A(h)\leq u,v\leq B(h)}\left\vert \frac{S(u)}{S(v)}-1\right\vert
\rightarrow _{\mathbb{P}}as\text{ }h\rightarrow +\infty .
\end{equation*}
\end{proposition}

\bigskip
\noindent \textbf{Proof}. See Lemma Lemma \ref{evt.rvsv.lemUnif} in Chapter \ref{evt}\\

\subsection{Theorem of Karamata and Theorem of de Haan}

\bigskip \noindent How to link extreme domains to regular or slowly variation? We have these
three characterizations.

\begin{proposition} \label{portal.rd}
We have the following characterizations for the three extremal domains.

\bigskip \noindent (a) $F\in D(H_{\gamma })$, $\gamma >0,$ if and only if there exist a
constant $c$ and functions $a(u)$ and $\ell (u)$ of $u\rightarrow $ $u\in
]0,1]$ satisfying
\begin{equation*}
(a(u),\ell (u))\rightarrow (0,0)\text{ as }u\rightarrow +\infty ,
\end{equation*}%
such that $F^{-1}$ admits the following representation of Karamata%
\begin{equation}
F^{-1}(1-u)=c(1+a(u))u^{-\gamma }\exp (\int_{u}^{1}\frac{\ell (t)}{t}dt). \label{portal.rdf}
\end{equation}

\bigskip \noindent (b) $F\in D(H_{\gamma })$, $\gamma <0,$ if and only if $uep(F)<+\infty $ and
there exist a constant $c$ and functions $a(u)$ and $\ell (u)$ of $u\in ]0,1]
$ satisfying
\begin{equation*}
(a(u),\ell (u))\rightarrow (0,0)\text{ as }u\rightarrow +\infty ,
\end{equation*}

\bigskip \noindent such that $F^{-1}$ admit the following representation of Karamata%
\begin{equation}
uep(F)-F^{-1}(1-u)=c(1+a(u))u^{-\gamma }\exp (\int_{u}^{1}\frac{\ell (t)}{t}
dt). \label{portal.rdw}
\end{equation}

\bigskip \noindent (c) $F\in D(H_{0})$ if and only if there exist a constant $d$ and a slowly
varying function $s(u)$ such that 
\begin{equation}
F^{-1}(1-u)=d+s(u)+\int_{u}^{1}\frac{s(t)}{t}dt,0<u<1, \label{portal.rdg}
\end{equation}

\bigskip \noindent and there exist a constant $c$ and functions $a(u)$ and $\ell (u)$ of $%
u\rightarrow $ $u\in ]0,1]$ satisfying
\begin{equation*}
(a(u),\ell (u))\rightarrow (0,0)\text{ as }u\rightarrow +\infty ,
\end{equation*}

\bigskip \noindent such that $s$ admits the representation%
\begin{equation}
s(u)=c(1+a(u)) \exp (\int_{u}^{1}\frac{\ell (t)}{t}dt). \label{portal.rdgs}
\end{equation}

\noindent Moreover, if $F^{-1}(1-u)$ is differentiable for small values of $s$ such
that $r(s)=-s(F^{-1}(1-s))^{\prime }=udF^{-1}(1-s)/ds$ is slowl varying at
zero, then \ref{portal.rdg} may be replaced by 
\begin{equation}
F^{-1}(1-u)=d+\int_{u}^{u_{0}}\frac{r(t)}{t}dt,0<u<u_{0}<1, \label{portal.rdgr}
\end{equation}

\noindent which will be called a reduced de Haan representation of $F^{-1}.$

\end{proposition}

\bigskip \noindent These representations are important. Proofs of them are to be found in the proofs of Proposition \ref{evt.sec.rvsv.repF}, Section \ref{evt.sec.rvsv} on Chapter \ref{evt} in the lines of proofs in \cite{loeve} or in \cite{dehaan}.\\

\noindent But we rather use the same representations but on $G$.\\

\subsubsection{Representations implied by the logarithm transformation}

\bigskip \noindent Since we deal with both $F$ and $G$ using the logarithm transform, we should
have the exact relation between them relatively to their belonging to extreme domains of attraction. From there,
we derive representations for $G^{-1}$.\\

\noindent Define the first asymptotic moment

\begin{equation*}
R(x,F)=\int_{x}^{uep(F)}\frac{1-F(t)}{1-F(x)}dt,x<uep(F).
\end{equation*}

\bigskip \noindent We have this result.

\begin{proposition}
We have the following equivalences.\\

\bigskip \noindent (a) Let $\gamma >0.$ Then $F\in D(H_{\gamma })$  $\Longleftrightarrow G\in
D(H_{0})$ and $R(x,G)\rightarrow \gamma $ as $x\rightarrow uep(G).$\\

\bigskip \noindent (b) $F\in D(H_{0})\Longleftrightarrow G\in D(H_{0})$ and $R(x,G)\rightarrow 0
$ as $x\rightarrow uep(G)$.\\

\noindent (c) Let $\gamma >0.$ Then $F\in D(H_{\gamma })\Longleftrightarrow F\in
D(H_{\gamma }).$
\end{proposition}

\bigskip \noindent \textbf{Proof}. See proof of Proposition ZZZ...\\

\bigskip \noindent Now based on this proposition, we expose the quantile representations for $G^{-1}$.

\begin{proposition} \label{portal.rdlog}
We have the following characterizations for the three extremal domains.\\

\bigskip \noindent (a) $F\in D(H_{\gamma })$, $\gamma >0,$ if and only if there exist a
constant $c$ and functions $p(u)$ and $b(u)$ of $u\rightarrow $ $u\in ]0,1]$
satisfying
\begin{equation*}
(p(u),b(u))\rightarrow (0,0)\text{ as }u\rightarrow +\infty ,
\end{equation*}

\noindent such that $xF^{-1}$ admit the following representation :
\begin{equation*}
G^{-1}(1-u)=c+\log (1+a(u)-\gamma \log u+\int_{u}^{1}\frac{b(t)}{t}dt.
\end{equation*}

\bigskip \noindent (b) $F\in D(H_{\gamma})$, $\gamma <0,$ if and only if $uep(G)<+\infty $ and
there exist a constant $c$ and functions $p(u)$ and $b(u)$ of $u\in ]0,1]$
satisfying
\begin{equation*}
(p(u),b(u))\rightarrow (0,0)\text{ as }u\rightarrow +\infty ,
\end{equation*}

\bigskip \noindent such that $G^{-1}$ admit the following representation :
\begin{equation*}
uep(G)-G^{-1}(1-u)=c(1+a(u))\text{ }u^{-\gamma }\exp (\int_{u}^{1}\frac{b(t)}{t}dt).
\end{equation*}

\bigskip \noindent (c) $F\in D(H_{0})$ if and only if there exist a constant $d$ and a slowly
varying function $s(u)$ such that 
\begin{equation*}
G^{-1}(1-u)=d+s(u)+\int_{u}^{1}\frac{s(t)}{t}dt,0<u<1,
\end{equation*}

\bigskip \noindent and there exist a constant $c$ and functions $p(u)$ and $b(u)$ of $%
u\rightarrow $ $u\in ]0,1]$ satisfying%
\begin{equation*}
(p(u),b(u))\rightarrow (0,0)\text{ as }u\rightarrow +\infty ,
\end{equation*}

\bigskip \noindent such that $s$ admits the representation%
\begin{equation*}
s(u)=c(1+p(u)) \ \exp (\int_{u}^{1}\frac{b(t)}{t}dt).
\end{equation*}
\end{proposition}

\subsection{General normalizing and centering sequences}

\bigskip \noindent From the representations of the quantile functions of distribution functions in the extreme domaon of attraction, we are are able to find general expressions of normalizing an centering coefficients as given below.\\

\begin{proposition} \label{portal.extdomain}
We have\\

\noindent (a) If $F\in D(G_{\gamma }),\gamma >0,$ then%
\begin{equation*}
\frac{X_{n,n}}{F^{-1}(1-1/n)}\leadsto Fr(1/\gamma ).
\end{equation*}

\noindent (b) If $F\in D(G_{\gamma }),\gamma <0,$ then $uep(F)<+\infty $ and%
\begin{equation*}
\frac{X_{n,n}-uep(F)}{uep(F)-F^{-1}(1-1/n)}\leadsto W(-1/\gamma ).
\end{equation*}

\noindent (c) If $F\in D(G_{0}),$ then
\begin{equation*}
\frac{X_{n,n}-F^{-1}(1-1/n)}{F^{-1}(1-1/(ne))-F^{-1}(1-1/n)}\leadsto \Lambda.
\end{equation*}
\end{proposition}

\bigskip \noindent \textbf{Proof}. Let us proceed case by case.\\

\noindent \textbf{Case \ $F\in D\left( G_{\gamma }\right) $ , $\gamma >0$}. By  (\ref{portal.rdf}), we have for $n\geq 1,$

\begin{equation*}
X_{n,n}=F^{-1}\left( 1-U_{1,n}\right) =c\left( 1+a(U_{1,n})\right) \left(
U_{1,n}\right) ^{-\gamma }\exp \left( \int_{U_{1,n}}^{1/n}\frac{(u)}{u}%
du\right) .
\end{equation*}

\noindent and 
\begin{equation*}
F^{-1}(1-1/n)=c\left( 1+a(1/n)\right) n^{\gamma }\exp \left( \int_{1/n}^{1}%
\frac{\ell (u)}{u}du\right) .
\end{equation*}

\bigskip \noindent We have

\bigskip 
\begin{equation*}
X_{n,n}/F^{-1}(1-1/n)=\frac{1+a(U_{1,n})}{1+a(1/n)}\left( nU_{1,n}\right) ^{-\gamma }\exp \left(
\int_{U_{1,n}}^{1/n}\frac{\ell (u)}{u}du\right) ,n\geq 1.
\end{equation*}

\bigskip \noindent Put 
\begin{equation*}
\ell _{n}=\sup \left\{ \left\vert b(t)\right\vert ;t\leq \max \left(
U_{1,n},1/n\right) \right\},
\end{equation*}

\bigskip \noindent 
\begin{equation*}
a_{n}=\sup \left\{ \left\vert a(t)\right\vert ;t\leq \max \left(U_{1,n},1/n\right) \right\}.
\end{equation*}

\bigskip \noindent Since  $U_{1,n}\rightarrow _{\mathbb{P}}0$ as $%
n\longrightarrow +\infty$, we have     
\begin{equation*}
(a_n \ell _{n})\longrightarrow _{\mathbb{P}} (0,0) \text{ as }n\longrightarrow +\infty. 
\end{equation*}

\noindent Now 

$$
\biggr| \frac{1+a(U_{1,n})}{1+a(1/n)}-1\biggr|\leq \frac{|a(1/n|+|a(U_{1,n}|}{1+a(1/n)}\equiv A_n,
$$

\bigskip \noindent and obviously, $A_n \longrightarrow_{\mathbb{P}} 0 \text{ as }n\longrightarrow +\infty$. We also have
\begin{eqnarray*}
\left\vert \exp \left( \int_{U_{1,n}}^{1/n}\frac{\ell (u)}{u}du\right)
-1\right\vert  &\leq &\left\vert \int_{U_{1,n}}^{1/n}\frac{\ell (u)}{u}%
du\right\vert  \\
&\leq &\ell _{n}\int_{\min \left( U_{1,n},1/n\right) }^{\max \left(
U_{1,n},1/n\right) }\frac{du}{u} \\
&\leq &\ell _{n}\log \left[ \frac{\max \left( U_{1,n},1/n\right) }{\min
\left( U_{1,n},1/n\right) }\right]  \\
&\leq &\ell _{n}\left\vert \log nU_{1,n}\right\vert.
\end{eqnarray*}

\noindent We know that $nU_{1,n}$ weakly converges to a standard exponential random
variable $E(1)$ and by the continuous mapping theorem (see Proposition 6, Chapter 2, \cite{wc-srv-ang}), $\left\vert \log \left( nU_{1,n}\right) \right\vert 
$ weakly converges to $\left\vert \Lambda \right\vert =\left\vert \log
E(1)\right\vert .$ This leads to $\left\vert \log \left( nU_{1,n}\right)
\right\vert =O_{\mathbb{P}}(1).$ Then 
\begin{equation*}
\exp \left( \int_{U_{1,n}}^{1/n}\frac{\ell (u)}{u}du\right) -1=O_{p}(\ell
_{n})=o_{p}(1)\longrightarrow 0
\end{equation*}%

\noindent and
\begin{equation*}
X_{n,n}/F^{-1}(1-1/n)=\left( nU_{1,n}\right) ^{-\gamma }\left(1+o_{p}(1)\right) (1 + O(A_n)) .
\end{equation*}

\noindent But, for $x\geq 0,$ 
\begin{eqnarray*}
\mathbb{P}(\left( nU_{1,n}\right) ^{-\gamma } &\leq &x)=\mathbb{P}(\left(
n(1-U_{n,n}\right) ^{-\gamma }\leq x) \\
&=&\mathbb{P}\left( U_{n,n}\leq \left( 1-\frac{x^{-1/\gamma }}{n}\right)
\right)  \\
&=&\left( 1-\frac{x^{-1/\gamma }}{n}\right) \text{ (for large values of }n)
\\
&\rightarrow &\exp (-x^{-1/\gamma })=\varphi _{1/\gamma }(x)
\end{eqnarray*}%

\noindent Since $\mathbb{P}(\left( nU_{1,n}\right) ^{-\gamma }\leq x)=0$ for $x\leq 0,$

\begin{equation*}
\left( nU_{1,n}\right) ^{-\gamma }\rightsquigarrow Fr\left( 1/\gamma \right).
\end{equation*}

\noindent and then%
\begin{equation*}
X_{n,n}/F^{-1}(1-1/n)\rightsquigarrow Fr\left( 1/\gamma \right) .
\end{equation*}

\bigskip \noindent  \textbf{Cas de Weibull $\protect\gamma <0$}. By (\ref{portal.rdw}), we have for $n\geq 1,$

\begin{equation*}
x_{0}\left( F\right) -F^{-1}\left( 1-U_{1,n}\right) =c\left(
1+a(U_{1,n})\right) \left( U_{1,n}\right) ^{-\gamma }\exp \left(
\int_{U_{1,n}}^{1/n}\frac{\ell (u)}{u}du\right) .
\end{equation*}

\noindent and 
\begin{equation*}
x_{0}\left( F\right) -F^{-1}\left( 1-U_{1,n}\right) =c\left( 1+a(1/n)\right)
n^{\gamma }\exp \left( \int_{U_{1,n}}^{1/n}\frac{\ell (u)}{u}du\right) .
\end{equation*}

\noindent This leads, for $n\geq 1,$ to%
\begin{eqnarray*}
\frac{X_{n,n}-x_{0}\left( F\right) }{x_{0}\left( F\right) -F^{-1}\left(
1-U_{1,n}\right) } &=&-(1+O(A_n))\left( nU_{1,n}\right) ^{-\gamma }\exp \left(
\int_{U_{1,n}}^{1/n}\frac{\ell (u)}{u}du\right)(1+O(A_n))  \\
&=&-\left( nU_{1,n}\right) ^{-\gamma } \left( 1+O(A_n)\right)\left( 1+O(\ell_n)\right)\\
&=&-\left( nU_{1,n}\right) ^{-\gamma }\left( 1+o_{\mathbb{P}}(1)\right)
\end{eqnarray*}

\noindent by the computations done before. And for $x\in \mathbb{R},$ 
\begin{eqnarray*}
\mathbb{P}(-\left( nU_{1,n}\right) ^{\gamma } &\leq &x)=\mathbb{P}(-\left(
n(1-U_{n,n}\right) ^{-\gamma }\leq x) \\
&=&\mathbb{P}\left( U_{n,n}\leq \left( 1-\frac{(-x)^{-1/\gamma }}{n}\right)
\right) .
\end{eqnarray*}

\noindent If $x\geq 0,$ $\mathbb{P}(U_{n,n}\leq (1-\frac{(-x)^{1/\gamma }}{n}))=1$. If 
$x\leq 0,$ we have for large values of $n,$
\begin{equation*}
\mathbb{P}(U_{n,n}\leq (1-\frac{(-x)^{1/\gamma }}{n}))=\left( 1-\frac{%
(-x)^{-1/\gamma }}{n}\right) ^{n}\rightarrow \exp (-(-x)^{-1/\gamma }).
\end{equation*}%

\noindent Then 
\begin{eqnarray*}
\frac{X_{n,n}-x_{0}\left( F\right) }{x_{0}\left( F\right) -F^{-1}\left(
1-U_{1,n}\right) } &=&-\left( nU_{1,n}\right) ^{-\gamma }\left(
1+o_{p}(1)\right)  \\
&\rightsquigarrow &W(-1/\gamma ).
\end{eqnarray*}

\noindent \textbf{Case of Gumbel $\gamma=0$}. By (\ref{portal.rdg}), we have for $n\geq 1$%
\begin{equation*}
X_{n,n}=c-s\left( U_{1,n}\right) +\int_{U_{1,n}}^{1}\frac{s(t)}{t}dt.
\end{equation*}

\noindent and%
\begin{equation*}
F^{-1}(1-1/n)=c-s\left( 1/n\right) +\int_{1/n}^{1}\frac{s(t)}{t}dt.
\end{equation*}

\noindent Then for $n\geq 1,$

\begin{equation*}
\frac{X_{n,n}-F^{-1}(1-1/n)}{s(1/n)}=\frac{s(1/n)-s\left( U_{1,n}\right) }{%
s(1/n)}+\int_{U_{1,n}}^{1/n}\frac{s(t)}{s(1/n)}\times \frac{1}{t}dt.
\end{equation*}

\noindent By Proposition \ref{portal.def5}, we have 
\begin{equation*}
\sup \{\left\vert \frac{s(t)}{s(1/n)}-1\right\vert ,U_{1,n}\wedge 1/n\leq
t\leq \text{ }U_{1,n}\vee 1/n\}=d_{n}\longrightarrow _{\mathbb{P}}0.
\end{equation*}%

\noindent and 
\begin{equation*}
\sup \left\{\left\vert \frac{s(t)}{s(1/n)}-1\right\vert ,1/(ne)\leq t\leq \text{ }%
1/n\right\}=c_{n}\longrightarrow 0.
\end{equation*}%

\noindent Then,
\begin{equation*}
\left\vert \frac{s(1/n)-s\left( U_{1,n}\right) }{s(1/n)}\right\vert \leq
d_{n}\longrightarrow _{\mathbb{P}}0
\end{equation*}

\noindent and 

\begin{equation*}
\left\vert \int_{U_{1,n}}^{1/n}\left( \frac{s(t)}{s(1/n)}-1\right) \frac{dt}{%
t}\right\vert \leq d_{n}\left\vert \log \left( nU_{1,n}\right) \right\vert
\rightarrow 0.
\end{equation*}

\bigskip \noindent We conclude that
\begin{eqnarray*}
\frac{X_{n,n}-F^{-1}(1-1/n)}{s(1/n)} &=&o_{P}(1)+\int_{U_{1,n}}^{1/n}\frac{1%
}{t}dt \\
&=&\log (nU_{1,n})+o_{\mathbb{P}}(1) \\
&\rightsquigarrow &\Lambda .
\end{eqnarray*}

\noindent We also have 
\begin{equation*}
\frac{X_{n,n}-b_{n}}{a_{n}}=o_{p}(1)+\log \left( nU_{1,n}\right) \leadsto
\Lambda .
\end{equation*}

\noindent further 
\begin{eqnarray*}
\frac{F^{-1}\left( 1-1/(ne)\right) -F^{-1}\left( 1-1/n\right) }{s(1/n)} &=&%
\frac{s(1/n)-s(1/(ne))}{s(1/n)}+\int_{1/n_{e}}^{1/n}\frac{s(t)}{s(1/n)}\frac{%
dt}{t} \\
&=&\frac{s(1/n)-s(1/(ne))}{s(1/n)}+\int_{1/(ne)}^{1/n}\left( \frac{s(t)}{%
s(1/n)}-1\right) \frac{dt}{t} \\
&&+\int_{1/(ne)}^{1/n}\frac{dt}{t} \\
&=&O(c_{n})+\int_{1/(ne)}^{1/n}\frac{dt}{t} \\
&=&o(1)+1.
\end{eqnarray*}

\noindent By Lemma \ref{portal.lem4},

\begin{equation*}
\frac{X_{n,n}-F^{-1}(1-1/n)}{s(1/n)}\rightsquigarrow \Lambda. \label{portal.gumbel1}
\end{equation*}

\noindent We easily see that if (\ref{portal.rdgr}) holds, similar but less heavy computations also lead to 
\begin{equation*}
\frac{X_{n,n}-F^{-1}(1-1/n)}{r(1/n)}\rightsquigarrow \Lambda. \label{portal.gumbel2}
\end{equation*}

\bigskip \noindent Let us give a few number of applications.

\subsection{Examples of more complex normalizing and centering coefficients} \label{portal.ss.examples} $ $\\

\bigskip \subsubsection{Upper extreme extreme of a standard Gaussian random variable.}

\noindent Suppose that $X\sim N(0,1)$, with probability density function%
\begin{equation*}
\phi _{d}(t)=\frac{1}{\sqrt{2\pi }}e^{-t^{2}/2},-\infty <t<+\infty 
\end{equation*}

\noindent and distribution function  
\begin{equation*}
\phi (x)=\frac{1}{\sqrt{2\pi }}\int_{-\infty
}^{x}e^{-t^{2}/2}dt=\int_{-\infty }^{x}\phi _{d}(t)dt.
\end{equation*}%

\noindent Let us give some useful expansion of functions related to this law, all of
them demonstrated in Section \ref{evt.techniq} of Chapter \ref{evt}. From
these formalae, we will be able to clearly expose the law of of $X_{n,n}.$
First, $1-\phi (x)$ admits the expansion  for $x>0,$%
\begin{equation}
C\left\{ \frac{1}{x}-\frac{1}{x^{2}}\right\} e^{-x^{2}/2}\leq 1-\phi (x)\leq 
\frac{Ce^{-x^{2}/2}}{x},  \label{portal.Gauss.boundPhi}
\end{equation}

\noindent where $C=1/\sqrt{2\pi}$. From this, the quantile $\phi ^{-1}(1-s)$ is expanded for $s\downarrow 0$, as%
\begin{equation}
\phi ^{-1}(1-s)=\left\{ (2\log (1/s))^{1/2}-\frac{\log 4\pi +\log \log (1/s)%
}{2(2\log (1/s))^{1/2}}+O((\log \log (1/s)^{2}(\log 1/s)^{-1/2}))\right\} .
\label{portal.Gauss.quantile}
\end{equation}

\noindent The derivative of $Q(1-s)=\phi ^{-1}(1-s)$ is, as $s\downarrow 0,$%

\begin{eqnarray}
\phi ^{-1}(1-s) &=&(2\log (1/s))^{1/2}-\frac{\log 4\pi +\log \log (1/s)}{%
2(2\log (1/s))^{1/2}}  \label{portal.Gauss.quantile} \\
&&+O((\log \log (1/s)^{2}(\log 1/s)^{-1/2})).
\end{eqnarray}

\bigskip \noindent  These three formulae respectively correspond to (\ref{evt.techniq.boundPhi}), \ (\ref{evt.techniq.quantile}) and to 
(\ref{evt.techniq.derivative}) in Section \ref{evt.techniq} of Chapter \ref{evt}. Now, we see that

\begin{equation*}
r(s)=-s(Q(1-s))^{\prime }=(2\log 1/s)^{1/2}(1+o(1))
\end{equation*}

\noindent \noindent is a slowly varying function at zero. By Proposition \ref{portal.rd}, we
have the representation (\ref{portal.rdgr}). So by (\ref{portal.gumbel2}), we have
\begin{equation*}
\frac{X_{n,n}-Q(1-1/n)}{r(1/n)}\rightsquigarrow \Lambda,
\end{equation*}

\noindent that is,

\begin{equation*}
\frac{X_{n,n}-Q(1-1/n)}{(2\log n)^{1/2}}\rightsquigarrow \Lambda .
\end{equation*}

\noindent Since, we know that $\phi \in D(\Lambda )$, we have that

\begin{equation*}
\frac{X_{n,n}-Q(1-1/n)}{Q(1-1/(ne))-Q(1-1/n)}\rightsquigarrow \Lambda .
\end{equation*}

\noindent Set, for $n\geq 2$,

\begin{equation*}
b_{n}=(2\log n)^{1/2}-\frac{\log 4\pi +\log \log n}{2(2\log n)^{1/2}}.
\end{equation*}%

\noindent It is easily got from  (\ref{portal.Gauss.quantile}) that%
\begin{equation*}
\frac{Q(1-1/n)-b_{n}}{(2\log n)^{1/2}}\rightarrow 0,\text{ as }n\rightarrow
\infty .
\end{equation*}%

\noindent we conclude that
\begin{equation*}
\frac{X_{n,n} - (2\log n)^{1/2} - \left( \frac{\log 4\pi +\log \log n}{2(2\log n)^{1/2}} \right)}{(2\log n)^{1/2}} \rightsquigarrow \Lambda .
\end{equation*}%


\noindent We get the following by-result, as a consequence of Lemma \ref{portal.lem4}
\begin{equation*}
\frac{Q(1-1/(ne))-Q(1-1/n)}{(2\log n)^{1/2}}\rightarrow 1\text{, as }%
n\rightarrow \infty .
\end{equation*}%

\noindent However, a direct proof of this result is given in (\ref{evt.techniq.A}) in  Section \ref{evt.techniq} of Chapter \ref{evt}.

\section{The Cs\"org\H{o} et al. space}

We already mentioned in Subsection \ref{portal.uniformrep}, that we will
replace the sequences of random variables $X_{1},X_{2},...$ by their uniform
representations $X_{1}=F^{-1}(U_{1}),X_{2}=F^{-1}(U_{2}),...$ In many
situations, we will suppose that the sequence $U_{1},U_{2},...$ are defined
on the Cs\"org\H{o}-Cs\"org\H{o}-Horv\`ath-Mason space. This space is described in  lemma
\ref{portal.cchm} below. Before we state it, let us introduce some notation.\\

\bigskip \noindent For a sequence of  a sequence $U_{1},U_{2},..$ of independent random
variables uniformly distributed on (0,1), we define, for each $n\geq 1,$ the
empirical quantile function%
\begin{equation*}
U_{n}(s)=\frac{j}{n}\text{ }U_{1,n}\leq s<U_{j+1,n},j=1,...,n
\end{equation*}

\noindent and $zero$ elsewhere, with the convention that $U_{0,n}=0$ and $U_{n+1,n}=1.$%
The uniform quantile process is defined for $n\geq 1$ by 
\begin{equation*}
V_{n}(s)=\left\{ 
\begin{tabular}{l}
$U_{j,n}\text{ \ for }\frac{j=1}{n}<s\leq \frac{j}{n},j=1,...,n$ \\ 
U$_{1,n}$ for $s=0.$%
\end{tabular}%
\right. 
\end{equation*}

\noindent The uniform empirical process is given for  for $n\geq 1$ 
\begin{equation*}
\mathbb{G}_{n}(s)=\{\sqrt{n}(U_{n}(s)-s),0\leq s\leq 1\}
\end{equation*}

\noindent and the quantile process is%
\begin{equation*}
\mathbb{V}_{n}(s)=\{\sqrt{n}(U_{n}(s)-s),0\leq s\leq 1\}.
\end{equation*}

\noindent Finally a standard brownian bridge $\{B(s),0\leq s\leq 1\}$ is a centered
continuous Gaussian process with variance and co-variance function%
\begin{equation*}
\mathbb{E}(B(s)B(t))=\min (s,t)-st,\text{ }(s,t)\in \lbrack 0,1]^{2}.
\end{equation*}

\bigskip \noindent Now, here is the so-important theorem of Cs\"org\H{o}-Cs\"org\H{o}-Horv\`ath-Mason (1986).

\begin{theorem} \label{portal.cchm} (See \cite{cchm}) There exists a
probability space holding a sequence of independent uniform random variables 
$U_{1},\ U_{2},$\ ... and a sequence of Brownian bridges $B_{1},B_{2},...$\
such that for each $0<\nu <1/4$, $as$\ $n\rightarrow \infty ,$%
\begin{equation}
\underset{1/n\leq s\leq 1-1/n}{\sup }\frac{\left\vert \sqrt{n}%
(U_{n}(s)-s)-B_{n}(s)\right\vert }{(s(1-s))^{1/2-\nu }}=O_{p}(n^{-\nu })
\label{2.1}
\end{equation}

\noindent and
 
\begin{equation}
\underset{1/n\leq s\leq 1-1/n}{\sup }\frac{\left\vert B_{n}(s)-\sqrt{n}%
(s-V_{n}(s))\right\vert }{(s(1-s))^{1/2-\nu }}=O_{p}(n^{-\nu }),  \label{2.2}
\end{equation}

\noindent where for each $n\geq 1$, $U_{n}(s)=j/n$ for $U_{j,n}\leq s<U_{j+1,n}$ is
the uniform empirical $df$ and $V_{n}(s)=U_{j,n}$ for $(j-1)/n<s\leq j/n,$
and $V_{n}(0)=U_{1,n},$ is the uniform quantile function and, finalyy, $%
U_{1,n}\leq ...\leq U_{n,n\text{ }}$are the order statistics of $%
U_{1},...,U_{n}$ with by convention $U_{0,n}=0=1-U_{n+1,n}.$ 
\end{theorem}


\part{Functional Aspects of Univariate Extreme Value Theory}

\chapter{Mathematical Background of Extreme Value Theory} \label{evt}

In this paper we will expose the very beginning of Extreme value
Theory, that is the extreme values domains of attraction for independent and
identically distributed randoms variables (\textit{iid}).\\

\noindent We begin this chapter by defining the notion of convergence in type, which
is a sub-domain of weak convergence.

\bigskip

\section{Convergence in type}

Basic and main notations concerning the general theory of weak convergence
are provided in the monograph \cite{wc-srv-ang} of this series. In particular, in this text, the most
used definition of weak convergence concerns the distribution functions in a
broad sense : a real distribution function is a function $F$ from $\mathbb{R}$ to $\mathbb{R}$,
right-continuous and non-decreasing. If $F(-\infty )=0$ and $F(+\infty )=1$, it becomes a probability distribution function. Next by the Probability Foundation Theorem of Kolmogorov, there exists a probability space holding a real random variable $X$ such that on that probability space, we have $F(x)=\mathbb{P}(X\leq x)$, $x \in \mathbb{R}$.\\

\noindent  We remind the definition of weak convergence just to make the reader ready to go through. In case he wants more details, he will be free to go back to \cite{wc-srv-ang} which is entirely devoted to that theory.\\

\begin{definition} \label{cv.def.001} The sequence of random vectors $X_{n} :(\Omega_{n},\mathcal{A}_n,\mathbb{P}_{n})\mapsto (\mathbb{R},\mathbb{B}(\mathbb{R}))$ weakly converges to $X : (\Omega_{\infty},\mathcal{A}_{\infty},\mathbb{P}_{\infty})\mapsto (\mathbb{R},\mathbb{B}(\mathbb{R}^k))$ 
if and only if for any point $x \in \mathbb{R}$ such that $\mathbb{P}(X=x)=0$, that is $x$ is 
a continuity point of $F_{X}(x)=\mathbb{P}(X_n \leq x)$,
\begin{equation}
F_{X_n}(x)=\mathbb{P}(X_n \leq x) \rightarrow F_{X}(x)=\mathbb{P}(X \leq x) \text{ as } n\rightarrow +\infty. \label{evt.000b}
\end{equation}
\end{definition}

\begin{remark} \label{evt.rem.001} By this definition, we see that the convergence of distribution concerns only the distribution. If you replace in (\ref{evt.000b}), $(X_{n})_{n\geq 0}$ by another sequence $(\widetilde{X}_{n})_{n\geq 0}$ such that $X_{n}=_{d}\widetilde{X}_{n}$ for each $n\geq 0$, it remains true. Also, the weak limit is unique in distribution because any random variable is fully determined by its its probability distribution function.
\end{remark}

\bigskip

\noindent The theory of extreme value deals with a particular case of weak convergence named as \textit{convergence in type} and defined below.

\bigskip

\begin{definition} \label{cv.def.002} A sequence of probability distribution functions $%
(F_{n})_{n\geq 0}$ converges in type to the probability distribution
functions $H$ if and only if there exist a sequence of positive real numbers 
$(a_{n}>0)_{n\geq 0}$ and a sequence of real numbers $(b_{n})_{n\geq 0}$
such that the sequence of distributions functions $F_{n}(a_{n}x+b_{n})$
weakly converges to a probability distribution function $H(x)$, that is, for any
continuity point of $H$, 
\begin{equation}
\lim_{n\rightarrow \infty }F_{n}(a_{n}x+b_{n})=H(x).  \label{evt.001}
\end{equation}
\end{definition}

\bigskip \noindent We may rephrase this definition into a different version. We may use the
Kolmogorov Theorem to place ourselves in a probability space $(\Omega,\mathcal{A},\mathbb{P})$ holding a sequence a random variable $Z$\ and a
sequence of random variables $(X_{n})_{n\geq 0}$ such that $X$ has the probability distribution function $H$ and each $X_{n}$ has the probability
distribution function $F_{n}$. With this representation, (\ref{evt.001}) is equivalent to 
\begin{equation*}
\frac{X_{n}-b_{n}}{a_{n}}\rightsquigarrow Z\text{ }\Longleftrightarrow \text{
}\frac{X_{n}-b_{n}}{a_{n}}\rightarrow _{d}Z,\text{ }as\text{ }n\rightarrow
+\infty ,
\end{equation*}

\noindent where the symbols respectively denote the weak convergence and the
convergence in distribution. These two convergences being equivalent, except
when we deal with non-measurability.\\

\bigskip

\noindent In the definition, the sequences $(a_{n}>0)_{n\geq 0}$ and $(b_{n})_{n\geq 0}
$ are not unique, nor is $H$. But one can only change them by satisfying the
relations below.

\bigskip In the sequel, we will use this notation : $C(H)$ denotes the set of continuity points of distribution function $H$.\\

\begin{lemma} \label{evt.lem.1} Let $(F_{n})_{n\geq 0}$ be a sequence of probability distribution functions Suppose there exist sequences 
$(a_{n}>0)_{n\geq 0}$, $(\alpha_{n}>0)_{n\geq 0}$, $(b_{n})_{n\geq 0}$ and $(\beta _{n})_{n\geq 0}$, probability distribution functions $H_{1}$ and $H_{2}$ such that

\begin{equation}
\lim_{n\rightarrow \infty }F_{n}(a_{n}x+b_{n})=H_{1}(x),\text{ }x\in C(H_{1})
 \label{evt.001a}
\end{equation}

\noindent and
\begin{equation}
\lim_{n\rightarrow \infty }F_{n}(\alpha _{n}x+\beta _{n})=H_{2}(x),\text{ }x\in C(H_{2}).  \label{evt.001b}
\end{equation}

\noindent Then there exist reals numbers $A>0$ and $B$ such that, as $n\rightarrow +\infty$,
\begin{equation}
\alpha _{n}/a_{n}\rightarrow A\text{ and }(\beta
_{n}-b_{n})/a_{n}\rightarrow B, \text{ as } n\rightarrow +\infty. \label{evt.001c}
\end{equation}

\noindent and for any $x\in \mathbb{R}$
\begin{equation}
H_{2}(x)=H_{1}(Ax+B).  \label{evt.002}
\end{equation}

\bigskip \noindent Reversely, If (\ref{evt.001a}) and (\ref{evt.001c}) hold both, then
(\ref{evt.001b}) is true, where $H_{2}$ defined in (\ref{evt.002}).
\end{lemma}

\bigskip

\noindent Formula (\ref{evt.002}) defines an equivalence class in the class of all real
probability distribution functions. Let us denote this equivalence relation
by $\mathcal{R}_{type}$. And we say that $H_{1}$ and $H_{2}$ are of the same
type if one is obtained from the other by a non-constant affine transformation of the
argument. The lemma says that the limit of convergence is unique in type,
meaning that all the possible limits in type are of the same type.\\

\bigskip \noindent \textbf{Proof}. We will use the Skorohod representation theorem for weak limits on on $\mathbb{R}$ (see Theorem 11, Section 4, \cite{wc-srv-ang}, in this series).\\

\noindent First assume that (\ref{evt.001a}), (\ref{evt.001b}) and (\ref{evt.001c}) hold. We are going to
prove (\ref{evt.002}). Placing ourselves in the right space with the help of Kolmogorov Theorem, we say that $F_{n}$
is the probability distribution function of $X_{n}$, $H_{1}$ is the probability distribution function
\ of $Z_{1}$ and $H_{2}$ is the probability distribution function of $Z_{2}$, all of these random variables being defined on the same probability space and as $n\rightarrow +\infty$, 

\begin{equation*}
U_{n}=\frac{X_{n}-b_{n}}{a_{n}}\rightsquigarrow Z_{1}\text{ and }V_{n}=\frac{%
X_{n}-\beta _{n}}{\alpha _{n}}\rightsquigarrow Z_{2}.
\end{equation*}

\noindent By the Skorohod Theorem we mentioned at the opening of the proof, there is also a probability space holding random
variables $\widetilde{Z}_{1}$, $\widetilde{Z}_{2},$ $\widetilde{U}_{n}$, $%
\widetilde{V}_{n},$ $n\geq 0,$ such that we have the following equalities in
distribution%
\begin{equation*}
\widetilde{Z}_{1}=_{d}Z_{1},\text{ }\widetilde{Z}_{2}=_{d}Z_{2},\text{ }%
\widetilde{U}_{n}=_{d}U_{n},\text{ }\widetilde{V}_{n}=_{d}V_{n},\text{ for }%
n\geq 0
\end{equation*}

\noindent and, $as$ $n\rightarrow +\infty ,$ 
\begin{equation*}
\widetilde{U}_{n}\rightarrow \widetilde{Z}_{1}\text{ }a.s.\text{ and }%
\widetilde{V}_{n}\text{=}\frac{X_{n}-\beta _{n}}{\alpha _{n}}\rightarrow
Z_{2}\text{ }a.s.
\end{equation*}

\noindent Here we only need the convergence of probability, that is, as $n\rightarrow +\infty$,
\begin{equation*}
\widetilde{U}_{n}\rightarrow _{P}\widetilde{Z}_{1}\text{ and }\widetilde{V}%
_{n}\text{=}\frac{X_{n}-\beta _{n}}{\alpha _{n}}\rightarrow _{P}Z_{2}.
\end{equation*}

\noindent It is evident that for any $n\geq 1$,
\begin{equation*}
a_{n}\widetilde{U}_{n}+b_{n}=_{d}\alpha _{n}\widetilde{V}_{n}+\beta _{n}.
\end{equation*}

\noindent Denote, for each $n\geq 1$
\begin{equation*}
\widetilde{X}_{n}(1)=a_{n}\widetilde{U}_{n}+b_{n}=_{d}X_{n}(2)=\alpha _{n}%
\widetilde{V}_{n}+\beta _{n}=_{d}X_{n}.
\end{equation*}

\noindent Then for each $\ngeq 1$,

\begin{equation*}
\widetilde{X}_{n}(1)=_{d}\widetilde{X}_{n}(2),
\end{equation*}

\noindent and, as $n\rightarrow +\infty$,

\begin{equation*}
\frac{\widetilde{X}_{n}(1)-b_{n}}{a_{n}}\rightarrow _{P}Z_{1}\text{ and }%
V_{n}=\frac{\widetilde{X}_{n}(2)-\beta _{n}}{\alpha _{n}}\rightarrow
_{P}Z_{2}.
\end{equation*}

\noindent This leads to 
\begin{equation*}
(\widetilde{X}_{n}(1)-\beta _{n})/\alpha _{n}=_{d}\left\{ (\widetilde{X}%
_{n}(2)-b_{n})/a_{n}\right\} (\alpha _{n}/a_{n})+\left\{
(b_{n}-b_{n})/a_{n}\right\} (\alpha _{n}/a_{n}).
\end{equation*}

\noindent The right member converges to $Z_{2}$ in probability, then in distribution by Proposition 12, Section 6 in \cite{wc-srv-ang}. The second member converges in probability to $A(Z_{1}-B),$ then also in
distribution. By Remark \ref{evt.rem.001} above, we have
\begin{equation*}
Z_{2}=A(Z_{1}-B)
\end{equation*}

\noindent and this implies (\ref{evt.002}) and we have the first part of the proof.\\

\noindent To complete the proof, suppose that (\ref{evt.001a}) and (\ref{evt.001b}) hold. We are going to
prove that \ref{evt.001c} and (\ref{evt.002}) also holds when both $H_{1}$ and $H_{2}$ are non-degenerated. If so, their generalized inverses $H_{1}^{-1}$ and $H_{2}^{-1}$ are also non
degenerated. Remind the definition of $H_{1}^{-1}$%
\begin{equation*}
H_{1}^{-1}(u)=\inf {x\in \mathbb{R},\text{ }H_{1}(x)\geq u},\text{ 
}u\in \lbrack 0,1],
\end{equation*}

\noindent which is a left-continuous and non-decreasing function. Properties of
generalized inverses are fully exposed in Chapter 4 in \cite{wc-srv-ang}. By using Point 8 in the mentioned  section, we may find to continuity points $u_{1}$ and $u_{2}$ of both $H_{1}^{-1}$ and $H_{2}^{-1}$ such that
\begin{equation*}
u_{1}<u_{2},\text{ \ \ \ \ }H_{1}^{-1}(u_{1})<H_{1}^{-1}(u_{2}),\text{ \ \ }%
H_{1}^{-1}(u_{1})<H_{1}^{-1}(u_{2}).
\end{equation*}

\noindent Next, also by Point 4 in in Chapter 4 in \cite{wc-srv-ang}, the weak convergence of $L_{n}$,
defined by $L_{n}(x)=F_{n}(a_{n}x$ $+b_{n})$, $x\in \mathbb{R},$ to $H_{1}$ in
\ref{evt.001a} implies convergence of their generalized inverses : $%
L_{n}^{-1}\rightsquigarrow H_{1}^{-1}$, that for any $u\in C(H_{1}),$%
\begin{equation*}
\frac{F_{n}^{-1}(u)-b_{n}}{a_{n}}\rightarrow H_{1}^{-1}(u).
\end{equation*}

\noindent Applying this to $u_{i},$ $i=1,2$. We get, as $n\rightarrow +\infty$,
\begin{equation}
\frac{F_{n}^{-1}(u_{i})-b_{n}}{a_{n}}\rightarrow H_{1}^{-1}(u_{i}),\text{ }%
i=1,2.  \label{evt.002b}
\end{equation}%

\noindent By taking the difference in each of the two formulas for $u_1$ and $u_2$, we have as $n\rightarrow +\infty$,
\begin{equation}
\frac{F_{n}^{-1}(u_{2})-F_{n}^{-1}(u_{1})}{a_{n}}\rightarrow
H_{1}^{-1}(u_{2})-H_{1}^{-1}(u_{1})  \label{evt.002a}
\end{equation}%

\noindent and
\begin{equation*}
\frac{F_{n}^{-1}(u_{2})-F_{n}^{-1}(u_{1})}{\alpha _{n}}\rightarrow
H_{2}^{-1}(u_{2})-H_{2}^{-1}(u_{1}).
\end{equation*}

\noindent By taking the ratio of the two last formula, we arrive at, as $n\rightarrow +\infty$,
\begin{equation*}
\alpha _{n}/a_{n}\rightarrow \frac{H_{2}^{-1}(u_{2})-H_{2}^{-1}(u_{1})}{%
H_{1}^{-1}(u_{2})-H_{1}^{-1}(u_{1})}:A>0.
\end{equation*}

\noindent From
\begin{equation*}
\frac{F_{n}^{-1}(u_{1})-b_{n}}{a_{n}}\rightarrow H_{1}^{-1}(u_{1})
\end{equation*}

\noindent and 
\begin{equation*}
\frac{F_{n}^{-1}(u_{1})-\beta _{n}}{a_{n}}=\frac{F_{n}^{-1}(u_{1})-\beta _{n}%
}{\alpha _{n}}(\alpha _{n}/a_{n})\rightarrow AH_{2}^{-1}(u_{1})
\end{equation*}

\noindent and by taking their difference, we get as $n\rightarrow +\infty$,
\begin{equation*}
\frac{\beta _{n}-b_{n}}{a_{n}}\rightarrow
H_{1}^{-1}(u_{1})-AH_{2}^{-1}(u_{1})=B.
\end{equation*}

\noindent The proof is complete.\\

\bigskip \noindent At this stage, we want to make two remarks.\\

\begin{remark} \label{evt.rem.1} $H$ is degenerated if and only $H^{-1}$ is degenerated, that is they have only one increase points at the same time. Next if (\ref{evt.002}) holds, $H_1$ is degenerated if and only if $H_2$ is. To see that, suppose that $H_{1}$ is
degenerated to the constant $a$, that is
\begin{equation*}
H_{1}(x)=\left\{ 
\begin{tabular}{lll}
$1$ & $if$ & $x\geq a$ \\ 
$0$ & $if$ & $x<a$%
\end{tabular}
\right. .
\end{equation*}
Then%
\begin{equation*}
H_{2}(x)=H_{1}(\frac{x-B}{A})=\left\{ 
\begin{tabular}{lll}
$1$ & $if$ & $x\geq (a-B)/A$ \\ 
$0$ & $if$ & $x<(a-B)/A.$%
\end{tabular}
\right. 
\end{equation*}

\noindent Thus $H_{2}$ is also degenerated. In this reasoning, we may exchange the roles of $H_1$ and $H_2$ and consider the inverse relation 
$H_{2}(x)=H_{1}(\frac{x-B}{A})$, $x \in \mathbb{R}$. We conclude that $H_1$ is degenerated if and only if $H_2$ is.\\
\end{remark}

\begin{remark} \label{evt.rem.2} Once we have the convergence in type, say (\ref{evt.001a}),
we may change the sequences by setting first $\gamma_{n}=F_{n}^{-1}(u_{2})-F_{n}^{-1}(u_{1})$ so that by (\ref{evt.002a}), 
\begin{equation*}
\gamma _{n}/\alpha _{n}\rightarrow C=H_{1}^{-1}(u_{2})-H_{1}^{-1}(u_{1})>0.
\end{equation*}

\bigskip \noindent Next, set $\delta _{n}=F_{n}^{-1}(u_{1})$. By (\ref{evt.002b}), we get
\begin{equation*}
\frac{\delta _{n}-b_{n}}{a_{n}}\rightarrow H_{1}^{-1}(u_{1})=D.
\end{equation*}
\end{remark}

\noindent By applying the theorem, we have for $x\in C(H_{1})$%
\begin{equation*}
F_{n}(\gamma _{n}x+\delta _{n})\rightarrow H(Cx+D).
\end{equation*}

\noindent We will keep in mind this choice for the extreme value theory :
\begin{equation*}
\left\{ 
\begin{tabular}{l}
$\gamma _{n}=F_{n}^{-1}(u_{2})-F_{n}^{-1}(u_{1})$ \\ 
$\delta _{n}=F_{n}^{-1}(u_{1})$%
\end{tabular}%
\right. 
\end{equation*}

\bigskip \noindent We have also 

\begin{lemma} \label{evt.lem.2} Let $F$ be a nondegenerated probability distribution function. If for any $x\in \mathbb{R}$, we have $F(ax+b)=F(cx+d)$, for real numbers $a>0$, $c>0$, $b$ and $d$, then $a=c$ and $b=d$.
\end{lemma}

\noindent \textbf{Proof of Lemma \ref{evt.lem.2}}. Consider the probability space $([0,1],\mathcal{B}([0,1]),\lambda ),$ where $\lambda $
the standard Lebesgue measure which is a probability. The canonical injection for $[0,1]$ to $\mathbb{R}$ follows the uniform standard law. Then $X=F^{-1}(U)$ has the probability distribution function $F$. Also $Z_{1}=(X-b)/a$ has probability distribution function $F(ax+b)$ and $Z_{2}=(X-d)/c$ the probability distribution function $F(cx+d)$. The equality $F(ax+b)=F(cx+d)$, $x\in \mathbb{R}$, implies for any $u\in \lbrack 0,1]$,

$$
(F^{-1}{u}-b)/a=(F^{-1}{u}-b)/c.
$$

\noindent Then
\begin{equation*}
\lambda (Z_{1}=Z_{2})=\lambda ({u\in \lbrack
0,1],(F^{-1}{u}-b)/a=(F^{-1}{u}-b)/c})=1.
\end{equation*}

\noindent Hence $Z_{1}=Z_{2}$ $a.s$ which leads do $aZ_{1}-cZ_{1}=c-b$ $a.s$ If $a\neq c,$ 
$Z_{1}=(c-d)/(a-c)=A$ $a.s$ and next, $X=aA+b$ $a.s$. So, unless $X$ is degenerated, we
have $a=c$ and next the equality $c=d$ easily follows.\\

\section{The different non-degenerated type limits in extreme value theory} \label{evt.types}

Let $X_{1}$, $X_{2}$, ... a sequence of independent and identically distributed (\textit{iid}) random variables (\textit{rv}'s) defined on a probability space $(\Omega,\mathcal{A},\mathbb{P})$, of common probability distribution function $F$. The univariate and classical extreme Value Theory is related on the characterization of the possible limits in type of the sequence of maxima
\begin{equation*}
M_{n}=\max (X_{1},...,X_{n}),\text{ }n\geq 1.
\end{equation*}

\noindent If the sequence $(M_{n})_{n\geq 0}$ converges in type to a random variable $Z$ of probability distribution function $H$, we say that the sequence is attracted to $Z$. This notion is more precise if we use the probability distribution functions, to say that $F$ is attracted to $H$, in other terms : $F$ lies in the extreme value domain of attraction of $H$, denoted $F \in D(H),$ if and only there exist sequences $(a_{n}>0)_{n\geq 0}$ and $%
(b_{n})_{n\geq 0}$ such that
\begin{equation*}
\frac{M_{n}-b_{n}}{a_{n}}\rightsquigarrow Z \text{ as } n \rightarrow +\infty,
\end{equation*}

\noindent that is, for any $x\in C(H)$%
\begin{equation}
F^{n}(a_{n}x+b_{n})\rightarrow H(x) \text{ as } n \rightarrow +\infty (\label{evt.00G})
\end{equation}

\bigskip
\noindent  Right now, we say that by the pioneering work of \cite{gnedenko43} and others who
will be cited that there are only three types non-degenerated possible
limits. This is one the most beautiful result of that theory.

\begin{theorem} \label{evt.theo1} (Fisher-Tippet (1928), Gnedenko (1943)). Suppose that $F \in D(H)$, where $H$ is nondegenerated. Then only three nondegenerated possible types of $H$ are the following.\\

\noindent \textbf{The Gumbel type of distribution function} :

\begin{equation}
\Lambda (x)=\exp (-\exp (-x)),\text{ }x\in \mathbb{R},  \label{dl05a}
\end{equation}%

\bigskip \noindent \textbf{The Fr\'{e}chet type \textit{df} of parameter $\gamma >0$} :

\begin{equation}
\phi _{\gamma }(x)=\exp (-x^{-\gamma })\mathbb{I}_{\left[ 0,+\infty \right[
}(x),\text{ }x\in \mathbb{R}\   \label{dl05b}
\end{equation}

\bigskip \noindent \textbf{The Weibull type \textit{df} of parameter $\gamma <0$} :

\begin{equation}
\psi _{\gamma }(x)=\exp (-(x)^{-\gamma })\mathbb{I}_{\left] -\infty ,0\right]
}(x)+(1-1_{\left] -\infty ,0\right] }(x)), \ x\in \mathbb{R},\ 
\label{dl05c}
\end{equation}

\bigskip\noindent where $I_{A}$ denotes the indicator function of the set $A$.
\end{theorem}

\bigskip \noindent We may the following notation in the sequel. $D(\phi)=\cup _{\gamma >0}D(\phi _{\gamma }),$ $D(\psi )=\cup _{\gamma >0}D(\psi_{\gamma }),$ and $\Gamma =D(\phi )\cup D(\psi )\cup D(\Lambda )$.\newline

\noindent \textbf{Proof of Theorem \ref{evt.theo1}}. Suppose that (\ref{evt.00G}) holds. For a fixed $t>0$,  we apply that formula for the indices $[nt]$ where $[.]$ stands for the integer part of a real number, that is the greatest integer less or equal to that number. We get sequences of real numbers $(a_{n}(t)>0)_{n\geq 0}$ and $(b_{n}(t))_{n\geq 0}$ such that

\begin{equation}
F^{[nt]}(a_{[nt]}x+b_{[nt]}(t))\rightarrow H(x)  \label{evt.003a}
\end{equation}

\noindent for any $x$ $\in C(H).$ Right here, we remark that the function $a$ which
mapps $t$ to $a_{[nt]}$ and the function $b$ which associates $t$ to $b_{[nt]}$ are measurable, since for exemple for the first case, for any Borel set $B$ of $\mathbb{R}$,
\begin{equation*}
(a\in B)=\sum_{a_{k}\in B}\left\{ t,[nt]=k\right\} =\sum_{a_{k}\in B}[\frac{k%
}{n},\frac{k+1}{n}[,
\end{equation*}

\noindent which is a Borel set. Further, (\ref{evt.003a}) implies for $x$ $\in C(H)$%
\begin{equation}
F^{[nt]}(a_{n}x+b_{n}(t))=(F^{n}(a_{n}x+b_{n}(t)))^{([nt]/n)}\rightarrow
H^{t}(x).  \label{evt.003b}
\end{equation}

\bigskip\noindent By Lemma \ref{evt.lem.1} and by comparing (\ref{evt.003a}) and (\ref{evt.003b})
there exists real numbers $\alpha (t)>0$ and $\beta (t),$ such that
\begin{equation}
a_{[nt]}/a_{n}\rightarrow \alpha (t)\text{ and (}b_{n}-b_{[nt]})/a_{[nt]}\rightarrow \beta (t)  \label{evt.003d}
\end{equation}

\noindent and for any $x\in \mathbb{R}$
\begin{equation}
H(\alpha (t)x+\beta (t))=H^{t}(x).  \label{evt.003c}
\end{equation}

\bigskip \noindent Let us apply this latter to the product $st$ to have%
\begin{equation*}
H^{st}(x)=H(\alpha (st)x+\beta (st))
\end{equation*}

\noindent But also
\begin{equation}
H^{st}(x)=(H^{s}(x))^{t}=H(\alpha (s)x+\beta (s))^{t}.  \label{evt.003e}
\end{equation}

\noindent We apply (\ref{evt.003c}) to $x=\alpha (s)x+\beta (s)$ to get

\begin{equation}
H^{st}(x)=H(\alpha (t)\left\{ \alpha (s)x+\beta (s)\right\} +\beta (t)).
\label{evt.003i}
\end{equation}

\noindent We arrive at the equality for any $x\in \mathbb{R}$, 
\begin{equation*}
H(\alpha (st)x+\beta (st))=H(\alpha (s)\alpha (t)x+\alpha (t)\beta (s)+\beta
(t)).
\end{equation*}

\noindent From this, and since $H$ is nondegenerated, we get by Lemma \ref{evt.lem.2}, 

\begin{equation}
\forall (s,t)\in (\mathbb{R}_{+}\setminus {0})^{2},\alpha (st)=\alpha
(s)\alpha (t)  \label{evt.003f}
\end{equation}

\noindent and
\begin{equation}
\forall (s,t)\in (\mathbb{R}_{+}\setminus {0})^{2},\beta (st)=\alpha (t)\beta
(s)+\beta (t).  \label{evt.003g}
\end{equation}

\noindent Formula (\ref{evt.003f}) is a Hamel-Cauchy equation. The solutions of such equations are given in Chapter \ref{funct}. By Corrolary \ref{funct.cor.1} of that chapter, we get that there exists $\rho
\in \mathbb{R}$ such that%
\begin{equation*}
\alpha (t)=t^{\rho },t>0.
\end{equation*}

\noindent Let us consider the three cases corresponding to three signs of $\rho$.\\

\noindent \textbf{Case $\rho =0$}.\\

\noindent Then $\alpha (t)=1$ and (\ref{evt.003f}) implies

\begin{equation*}
\forall (s,t)\in (\mathbb{R}_{+}\setminus {0})^{2},\beta (st)=\beta (s)+\beta (t).
\end{equation*}

\noindent By \ref{funct.cor.1} of Chapter \ref{funct}, for $c=-\beta (e),$ we have
\begin{equation*}
\beta (st)=-c\log t,t>0.
\end{equation*}

\noindent Formula (\ref{evt.003c}) gives

\begin{equation}
H^{t}(x)=H(x-c\log t),t>0.  \label{evt.003h}
\end{equation}

\noindent We are going to make a number of remarks and implications of this fact.
First, for a fixed $x\in \mathbb{R}$, $H^{t}(x)$ in non-increasing in $t>0,$ and since 
$H$ is non-decreasing, we necessarily have that $c>0$. Next, we are going to
see that $H$ has an unbounded above support, that is $H(x)<1,$ for all $x\in \mathbb{R}$. Otherwise, there exists $x_{0}$ such that $H(x_{0})=1$. This would imply, from (\ref{evt.003h}),

\begin{equation*}
1=H(x_{0}-c\log t),t>0.
\end{equation*}

\noindent But $\{x_{0}-c\log t),t>0\}=\mathbb{R}$ so that we would get $H(x)=1$ for all $x\in \mathbb{R}$
and then by right continuity $H(-\infty )=1.$ This is a contradiction. So $%
H(x)<1,$ for all $x\in \mathbb{R}$. By the same argument $H(0)=0$ implies that $H(x)=0
$ for all $x\in \mathbb{R}$. Now let
\begin{equation*}
p=-\log (-\log H(0)),
\end{equation*}

\noindent that is
\begin{equation*}
H(0)=\exp (-e^{-p}).
\end{equation*}

\noindent Next from  (\ref{evt.003h}), we get

\begin{equation*}
H^{t}(0)=H(-c\log t),t>0.
\end{equation*}

\noindent From there, make the change of variance $x=-c\log t$, that is $t=\exp (-x/c)$
and next

\begin{eqnarray*}
H(x)=(H(0))^{\exp (-x/c)}&=&\exp (-\exp (-x/c)\exp (-p))\\
&=&\exp (-\exp (-\frac{x+cp}{c})),x\in \mathbb{R}.
\end{eqnarray*}

\bigskip
\noindent By letting $a=c>0$ and $b=-cp$, we have

\begin{equation*}
H(x)=\exp (-\exp (-\frac{x-b}{a})),x\in \mathbb{R}
\end{equation*}

\noindent which is of type of
\begin{equation*}
\Lambda (x)=\exp (-e^{-x}),x\in \mathbb{R}.
\end{equation*}

\noindent \textbf{Case $\rho <0$}.\\

\noindent By (\ref{evt.003g}) and by the symmetry of the roles of $s$ and $t,$%
\begin{equation*}
\alpha (t)\beta (s)+\beta (t)=\alpha (s)\beta (t)+\beta (s),
\end{equation*}

\noindent which implies
\begin{equation*}
\frac{\beta (s)}{1-\alpha (s)}=\frac{\beta (s)}{1-\alpha (s)},
\end{equation*}

\noindent for any $s>0$ and $t>0$ with $s\neq 1$ and $t\neq 1.$ Then $\beta
(t)/(1-\alpha (t)=c$ is a constant $c$ for $t\neq 1$ and next 

\begin{equation*}
\beta (t)=c(1-t^{\rho }),t>0,t\neq 1,
\end{equation*}

\noindent and \ref{evt.003c} becomes for $t>0,t\neq 1$%
\begin{equation*}
H^{t}(x)=H(t^{\rho }(x-c)+c),x\in \mathbb{R},
\end{equation*}

\noindent that is for $t>0,t\neq 1$%
\begin{equation*}
H^{t}(x+c)=H(xt^{\rho }+c),x\in \mathbb{R}.
\end{equation*}

\noindent Let us change $H(x)$ into the type $H_{1}(x)=H(x+c)$ which is of the same
type of $H$ to get

\begin{equation}
H_{1}^{t}(x)=H_{1}(xt^{\rho }),x\in \mathbb{R}.  \label{evt.003j}
\end{equation}

\noindent As for the precedent case, let us see that the fact $H_{1}(0)=1$ is impossible. Let us assume that $H_{1}(0)=1$ and set
$$
a=\inf \{x<0, H_{1}(x)=1\}
$$

\noindent If $a=-\infty$, we would get that $H_1=1$ on the whole real line, which is absurd. Then $a$ is finite. If we have $H_1(x)=0$ for all $x<a$, $H_1$ would be degenerated. Finally, it would exists
$x_0<a$ such that $H_1(x_0)>0$. By the definition of $a$, we would have $H_1(x_0)<1$ and hence, $0<H_1(x_0)<1$. By (\ref{evt.003j})), we would have
\begin{equation*}
H_{1}^{t}(x_{0})=H(x_{0}t^{\rho }),t>0,t\neq 1.
\end{equation*}

\noindent is impossible since, as $t$ increases (while avoiding the
value $1$), the left-side member of this equality is decreasing to zero while
the right-side member is nondecreasing. Then we exclude the fact that $H_{1}(x)=0.$ Now, application of \ref{evt.003j} to $x=0$ gives%
\begin{equation*}
H_{1}^{t}(0)=H_{1}(0),t>0,t\neq 1
\end{equation*}

\noindent Since $0\leq H_{1}(0)\leq 1,$ we have only two possibilities $H_{1}(0)=1$ or 
$H_{1}(0)=0.$ The case $H_{1}(0)=1$ has already been excluded. We keep $H_{1}(0)=0.$ Now
we are going to see that $H_{1}(1)$ is different of $0$ and of $1.$ Indeed, $%
H_{1}(1)=0$ implies through (\ref{evt.003j}) that%
\begin{equation*}
0=H_{1}(t^{\rho }),t>0,t\neq 1.
\end{equation*}

\noindent Then for $t$ small enough, $H_{1}(t^{\rho })=0$ and $H_{1}(x)=0$ for $x\leq
t^{\rho }.$ By letting $t\downarrow 0,$ we have $H_{1}(x)=0$ for $x\leq
t^{\rho }\nearrow +\infty $ so that $H_{1}(x)=0,x\in \mathbb{R}$. This is impossible.
Next   $H_{1}(1)=1$ implies%

\begin{equation*}
1=H_{1}(t^{\rho }),t>0,t\neq 1.
\end{equation*}

\noindent Then for $t$ large enough,%
\begin{equation*}
H_{1}(x)\text{ for }x\geq t^{\rho }\searrow 0\text{ as }t\nearrow +\infty .
\end{equation*}
Then%
\begin{equation*}
H_{1}(x)=1,x\geq 0.
\end{equation*}

\noindent By combining this with the fact that $H_{1}(0)=0,$ we arrive at the conclusion
that $H_{1}$ concentrated at zero, which is in contraction with our
assumption. Then we have
\begin{equation*}
H_{1}(0)=0\text{ and }0<H_{1}(1)<1.
\end{equation*}

\noindent Finally, (\ref{evt.003j}) yields for $t>0,t\neq 1,$ 
\begin{equation*}
H_{1}^{t}(1)=H_{1}(t^{\rho })
\end{equation*}

\noindent and by change of variables $x=t^{\rho }>0,x\neq 1$%
\begin{equation*}
H_{1}(x)=H_{1}(1)^{x^{1/\rho }}=\exp (-x^{1/\rho }(-\log H_{1}(1)))=\exp
(-ax^{-\gamma }),
\end{equation*}

\noindent where $a=-\log H_{1}(1)>0,\gamma =-1/\rho .$ By putting together all what
preceeds and by using right continuity to handle the point $x=1$, we have%

\begin{equation*}
H_{1}(x)=\left\{ 
\begin{tabular}{lll}
$\exp (-ax^{-\gamma })$ & if  & $x\geq 0$ \\ 
0 & if  & $x<0$,
\end{tabular}
\right. 
\end{equation*}

\noindent which is of the type of%
\begin{equation*}
\Phi _{\gamma }(x)=\left\{ 
\begin{tabular}{lll}
$\exp (-x^{-\gamma })$ & if  & $x\geq 0$ \\ 
0 & if  & $x<0$.
\end{tabular}
\right. 
\end{equation*}

\noindent \textbf{Case $\rho>0$}.\\

\noindent In the proof of Case $\rho <0$, the part from the beginning and to Formula (\ref{evt.003g}) does not depend on the value of $\rho \neq 0.$ So we may start this proof from \ref{evt.003j}. We are going to prove the upper endpoint of $H_{1}(x)$ is non-positive, that is $uep(H_1)\leq 0$. Suppose that $uep(H_1)>0$.  If we have $H_1(x)=0$ for all $x<uep(H_1)$, $H_1$ would be degenerated. Hence, it would exists $x_1<uep(H_1)$ such that $H_1(x_1)>0$. Of we choose $x_0$ such that $\max(0,x_1)<x_0<uep(H_1)$, we would get $0<H_1(x_0)<1$. We would finally have 

\begin{equation}
H_{1}^{t}(x_{0})=H_{1}(x_{0}t^{\rho }),t>0,t\neq 1.
\end{equation}

\noindent The left-side member decreases in $t$ while the right-side is nondecreasing in $t$. This impossibility combined with the first one and the right-continuity at zero, allows to conclude that 
$uep(H_1)\leq 0$ and then

\begin{equation*}
H_{1}(x)=1,x\geq 0.
\end{equation*}

\noindent As well, since $H_{1}$ is nondegenerated, there exists $x_{1}<0$ such that $0<H_{1}(x_{1})<1$ and (\ref{evt.003j}) implies

\begin{equation}
H_{1}^{t}(x_{1})=H_{1}(x_{1}t^{\rho }),t>0,t\neq 1.
\end{equation}

\noindent Set $x=x_{1}t^{\rho }<0$ to get

\begin{eqnarray*}
H_{1}(x)&=&H_{1}(x_{1})^{(x/x_{1})^{1/\rho }}\\
&=&\exp ((x/x_{1})^{1/\rho }\log H_{1}(x_{1}))=\exp (-b(-x)^{\gamma }),x\neq -1,
\end{eqnarray*}

\noindent where $b=-(-x)^{1/\rho }\log H_{1}(x_{1})>0$. Put together all what precedes
and use the right-continuity at $-1$ to get that

\begin{equation*}
H_{1}(x)=\left\{ 
\begin{tabular}{lll}
1 & if  & $x\geq 0$ \\ 
$\exp (-b(-x)^{\gamma })$ & if  & $x<0$,
\end{tabular}
\right. 
\end{equation*}

\noindent which is of the type of
\begin{equation*}
\Psi _{\gamma }(x)=\left\{ 
\begin{tabular}{lll}
1 & if  & $x\geq 0$ \\ 
$\exp (-b(-x)^{\gamma })$ & if  & $x<0$%
\end{tabular}%
\right. ,\gamma >0.
\end{equation*}

\begin{remark}
\noindent Threorem \ref{evt.theo1} does not say that these three limits in type effectively occur. It only says that if $H$ is non-degenerated and is limit of type of $F^n$, then it is necessarily of the these three types.\\

\noindent In the portal chapter \ref{portal}, in Subsection \ref{portal.subsec_examples_simples}, we already gave simple examples leading to these three types.
\end{remark}

\noindent We want to finish with a generalization of the result in the lines just after Remark \ref{evt.rem.2}. This result is a by-product of the proof of Lemma \label{evt.lem.1}. It allows to find arbitrary normalizing and centering coefficients once we know that a sequence of distribution functions $F_n$ converges in type to a non-degenerated distribution $H$. Is is stated as follows.\\

\begin{lemma} \label{evt.lem.1plus} Let $F_{n}$ be a sequence of probability distribution functions weakly converging to a non-degenerated distribution function $H$. Consider any other distribution of type of $H$ of the form $G(x)=H(Ax+B)$, $x\in \mathbb{R}$ with $A>0$. Choose 
$0<u_{1} <u_{2} <1$ such that $u_{1}$ and $u_{2}$ are continuity points of $H^{-1}$ and $G^{-1}$. Set  for $n \geq 1$

$$
\gamma _{n}=F_{n}^{-1}(u_{2})-F_{n}^{-1}(u_{1}),
$$

$$
\delta _{n}=F_{n}^{-1}(u_{1}),
$$

\noindent and 

\begin{equation*}
a=\frac{G^{-1}(u_{2})-G^{-1}(u_{1})}{H^{-1}(u_{2})-H^{-1}(u_{1})}
\end{equation*}

\noindent and 

\begin{equation*}
b=H^{-1}(u_{1})-a^{-1}(u_{1}).
\end{equation*}

\noindent Then $F_n(\gamma_{n}x+\delta_n) \rightarrow H(ax+b)$ for any continuity point of $H$.\\
\end{lemma} 

\noindent We are going to state a useful continuous version of that result, which is very important for the characterization of the Gumbel extreme domain of attraction.\\

\begin{lemma} \label{evt.lem.pregumbel} \bigskip Let $H$ be a strictly increasing and continous distribution
function on $S(H)=[lep(H),uep(H)]$. Let $\left( F_{n}\right) _{n\geq 1}$ be
a sequence of probability distribution which weakly converges to $H.$ Choose
two real numbers  $u_{1}$ and $u_{2}$ such that $0<u_{1}<u_{2}<1$. Let $%
\left( u_{1}(n)\right) _{n\geq 1}$ and $\left( u_{2}(n)\right) _{n\geq 1}$
be two sequences of numbers in $]0,1[$ such that $u_{1}(n)\longrightarrow
u_{1}$ and $u_{1}(n)\longrightarrow u_{2}$ as $n\rightarrow +\infty .$ Denote%
\begin{equation*}
\left\{ 
\begin{tabular}{l}
$b=H^{-1}(u_{1})$ \\ 
$a=H^{-1}(u_{2})-H^{-1}(u_{2})$%
\end{tabular}
\right. 
\end{equation*}

\noindent and define the sequences
\begin{equation*}
\left\{ 
\begin{tabular}{l}
$b_{n}=F_{n}^{-1}(u_{1}(n))$ \\ 
$a_{n}=F_{n}^{-1}(u_{2}(n))-F_{n}^{-1}(u_{2}(n))$%
\end{tabular}%
\right. ,n\geq 1.
\end{equation*}

\bigskip \noindent  Put $A=a-b,B=b$ and $A_{n}=a_{n}-b_{n},$ $B_{n}=b_{n}$, $n\geq 1$.\\

\noindent We have\\

\noindent (I) $A>0$, $A_{n}>0$ for large values of $n$, and for any $x\in S(H)$,
\begin{equation*}
F_{n}(A_{n}x+B_{n})\rightarrow H(Ax+B)\text{ as }n\rightarrow \infty .
\end{equation*}

\noindent (II) If for all $x\in D$ such that 
$$
F_{n}(c_{n}x+d_{n})\rightarrow H(x),
$$ 

\noindent for sequences $\left(c_{n}>0\right) _{n\geq 1}$ and $\left( d_{n}\right) _{n\geq 1}$, then we have

$$
\frac{F_{n}^{-1}(u_{2}(n))-F_{n}^{-1}(u_{1}(n))}{c_{n}} \rightarrow H^{-1}(a)-H^{-1}(b) \text{ as } n\rightarrow +\infty.
$$

\noindent and

$$
\frac{F_{n}^{-1}(u_{1}(n))-d_{n}}{c_{n}} \rightarrow 0 \text{ as } n \rightarrow +\infty.
$$

\end{lemma}

\bigskip \noindent \textbf{Proof}. First, we fix $\delta >0$ such that $\delta <(u_{2}-u_{1})/2.$ So for large
values of $n$ we have
\begin{equation}
\normalsize
0<u_{1}-\delta /2<u_{1}(n)<u_{1}+\delta /2<u_{2}-\delta
/2<u_{2}(n)<u_{2}+\delta /2<1.  \label{evt.lem.12}
\end{equation}

\noindent By continuity of $H$, there exist $\varepsilon _{0}$ such that for $0<\varepsilon <$ $\varepsilon _{0}$ , 
\begin{equation*}
\min \{H(y+\varepsilon )-H(y),H(y)-H(y-\varepsilon ),y\in \{a,b\}\}<\delta
/2.
\end{equation*}

\noindent This implies, for $0<\varepsilon <\varepsilon _{0}$, that 
\begin{equation*}
H(y+\varepsilon )-\delta /2<H(y)<H(y+\varepsilon )+\delta /2
\end{equation*}

\noindent Since $F_{n}\rightsquigarrow H$, there exist sequences of real numbers $\left(
c_{n}>0\right) _{n\geq 1}$ and $\left( d_{n}\right) _{n\geq 1}$ such that
for each $x\in D,$ $F_{n}(c_{n}x+d_{n})\rightarrow H(x).$ Then for $y\in
\{a,b\},$ we have 
\begin{equation*}
F_{n}(c_{n}(y-\varepsilon )+d_{n})\rightarrow H(y-\varepsilon)\text{ and } F_{n}(c_{n}(y+\varepsilon )+d_{n})\rightarrow H(y+\varepsilon).
\end{equation*}

\noindent Then for $n$ large enough
\begin{equation}
F_{n}(c_{n}(y-\varepsilon )+d_{n})<H(y-\varepsilon )+\delta 
\label{evt.lem.13A}
\end{equation}

\noindent and

\begin{equation}
F_{n}(c_{n}(y+\varepsilon )+d_{n})>H(y+\varepsilon )-\delta. \label{evt.lem.13B}
\end{equation}

\noindent From this point, we handle the cases $y=a$ and $y=b$ one after the other. We
have for large values of $n$,

\begin{equation*}
b_{n}=F_{n}^{-1}(u_{1}(n))\Longrightarrow F_{n}(b_{n}-)<u_{1}(n)\leq
F_{n}(b_{n})
\end{equation*}

\bigskip \noindent Fix $n_{0}$ such that (\ref{evt.lem.12}), (\ref{evt.lem.13A}) and (\ref{evt.lem.13B}) hold for $n\geq n_{0}.$ For $n\geq n_{0},$ we have

\begin{equation*}
F_{n}(b_{n})\geq u_{2}(n)>u_{1}-\delta /2=H(b)-\delta >H(b+\varepsilon
)-\delta >F_{n}(c_{n}(b+\varepsilon )+d_{n}).
\end{equation*}

\bigskip \noindent And \ $F_{n}(b_{n})>F_{n}(c_{n}(y+\varepsilon )+d_{n})$ for $n\geq n_{0}$
implies

\begin{equation*}
b_{n}>c_{n}(y+\varepsilon )+d_{n},\text{ for }n\geq n_{0}.
\end{equation*}

\noindent As well, for $n\geq n_{0}$,

\begin{equation*}
F_{n}(b_{n}-)<u_{1}(n)<u_{1}+\delta =H(b)+\delta /2<H(y+\varepsilon )+\delta
<F_{n}(c_{n}(b+\varepsilon )+d_{n})
\end{equation*}

\noindent which implies that
\begin{equation*}
F_{n}(b_{n}-)<F_{n}(c_{n}(b+\varepsilon )+d_{n}).
\end{equation*}

\noindent By the definition of the left limit, for each $n\geq n_{0}$, there exists $h>0$ such that

\begin{equation*}
F_{n}(b_{n}-)\leq F_{n}(b_{n}-h)<F_{n}(c_{n}(b+\varepsilon )+d_{n}),
\end{equation*}

\noindent which ensures that
\begin{equation*}
b_{n}-h<c_{n}(b+\varepsilon )+d_{n}.
\end{equation*}

\noindent By combining all these results, wa have for $n\geq n_{0}$ and for $%
0<\varepsilon <\varepsilon _{0}$,
\begin{equation}
c_{n}(b-\varepsilon )+d_{n}<b_{n}<c_{n}(b+\varepsilon )+d_{n},
\label{evt.lem.14A}
\end{equation}

\noindent that is
\begin{equation*}
\left\vert \frac{b_{n}-d_{n}}{c_{n}}\right\vert <\varepsilon .
\end{equation*}

\noindent Taking the limit superior and next letting $\varepsilon \downarrow 0$ together give
\begin{equation*}
\frac{b_{n}-d_{n}}{c_{n}}\rightarrow 0\text{ as }n\rightarrow \infty .
\end{equation*}

\noindent As well, simular considerations that led to (\ref{evt.lem.14A}) yield : for $%
n\geq n_{0}$ and for $0<\varepsilon <\varepsilon _{0},$%
\begin{equation}
c_{n}(a-\varepsilon )+d_{n}<a_{n}<c_{n}(a+\varepsilon )+d_{n}.
\label{evt.lem.14B}
\end{equation}

\noindent Formulas (\ref{evt.lem.14A}) and (\ref{evt.lem.14B}) together ensure that
for $n\geq n_{0}$ and for $0<\varepsilon <\varepsilon _{0},$%
\begin{equation*}
c_{n}(a-b-2\varepsilon )<a_{n}-b_{n}<c_{n}(a-b+2\varepsilon ),
\end{equation*}

\noindent that is for $n\geq n_{0}$ \ (recall that $a_{n}\geq b_{n}$ for all $n\geq 1$
and $a>b).$ 
\begin{equation*}
\left\vert \frac{a_{n}-b_{n}}{c_{n}}-(a-b)\right\vert <2\varepsilon .
\end{equation*}

\noindent Similarly, taking the limit superior and next letting $\varepsilon
\downarrow 0$ gives%
\begin{equation*}
\frac{a_{n}-b_{n}}{c_{n}}\rightarrow (a-b)\text{ as }n\rightarrow \infty .
\end{equation*}
 
\noindent Part (II) of the lemma is the summary of the formulas proved by the computations above. Part (I) is the result of the application of Lemma \ref{evt.lem.1} to these formula.\\

\noindent Before we begin the characterizations of distribution functions in the extreme value domain, we have to take a serious and deep tour of the important classes of functions that are in the heart of univariate extreme value theory.\\

\section{Regularly and $\pi$-variation} \label{evt.sec.rvsv}

\noindent This section is devoted to representations of functions involved in Extreme
Value Theory, in particular to representations of Regularly varying functions and $\pi$-variation functions. The Karamata Representation Theorem and that of \text{de} Haan will be our special guests.\\

\bigskip \noindent Let us begin by the concept of regular variation.\\

\noindent \textbf{(A) Regular and Slow Variation}.\\

\bigskip \noindent Throughout this section, we will deal with functions  $U:\mathbb{R}%
_{+}\longrightarrow \mathbb{R}_{+}$ that are measurable and Lebesgue
locally integrable and not vanishing in the neighborhood of $+\infty ,$
that is, $(\forall 0\leq a\leq b<+\infty)$, $U\in L([a,b],\lambda )$. This will
ensure in particular that $U$ is continuous $a.e.$ and the formula
\begin{equation*}
\left( \int_{0}^{x}U(t)d\lambda (t)\right) ^{\prime }=U(x),\text{ }\lambda
-a.e.\text{ on }[0,+\infty \lbrack .
\end{equation*}

\noindent holds.\\

\subsection{Definitions}

\bigskip \noindent We begin to define the regular and the slow variation.

\begin{definition}
The function $U$ is regularly varying at $+\infty $ with exponent $\rho \in 
\mathbb{R}$ if and only for any $\gamma >0,$

\begin{equation}
\lim_{x\rightarrow +\infty }U(\gamma x)/U(x)=\gamma ^{\rho }.
\label{evt.rvsv.DEFRV}
\end{equation}

\noindent When (\ref{evt.rvsv.DEFRV}) holds, we use the notation

\begin{equation*}
U\in RV(\rho ,+\infty ).
\end{equation*}

\noindent If $\rho =0,$ we say that $U$ is slowly varying at $+\infty ,$ denoted as%
\begin{equation*}
U\in U\in RV(0,+\infty )\text{ or }U\in SV(+\infty ).
\end{equation*}
\end{definition}

\bigskip \noindent We are going to work with limits at $+\infty$. We also frequently work in
neighborhoods of $zero.$ So, we have to adapt the definition above for
functions $u:[0,u_{0}]\longrightarrow \mathbb{R}_{+}$, where $u_{0}>0$, defined in a right-neighborhood of zero. 
Using the transform

\begin{equation*}
u(s)=U(1/s),s\in \lbrack 0,u_{0}],
\end{equation*}

\noindent allows to transfer the definition in the following way.\\

\begin{definition}
The function $u:[0,u_{0}]\longrightarrow \mathbb{R}_{+}$, where $u_{0}>0,$
is regularly varying at $zero$ with exponent $\rho \in \mathbb{R}$ if and
only for any $\gamma >0$,
\begin{equation}
\lim_{s\rightarrow 0}u(\gamma s)/u(s)=\gamma ^{\rho }.
\label{evt.rvsv.DEFRV0}
\end{equation}

\noindent When (\ref{evt.rvsv.DEFRV}), we use the notation

\begin{equation*}
u\in RV(\rho ,0).
\end{equation*}

\noindent If $\rho =0,$ we say that $u$ is slowly varying at $zero$, denoted as

\begin{equation*}
u\in RV(0,0)\text{ or }u\in SV(0).
\end{equation*}
\end{definition}

\bigskip \noindent We will be able to move from one of these two versions to the other, simply by
remarking that : for $u(s)=U(1/s),$ $u\in RV(\rho ,0)$ if and only if $U\in
RV(-\rho ,+\infty ).$ So, the theory will be made for one the version and
translated to the other if needed.\\

\bigskip \noindent At the very beginning, let us notice these immediate properties.

\begin{lemma} \label{evt.rvsv.charac} Let $U:\mathbb{R}_{+}\longrightarrow \mathbb{R}_{+}$ be a
measurable function Lebesgue locally integrable and not vanishing in the
neighborhood of $+\infty ,$ that is, $\forall\text{ }(0\leq a\leq b<+\infty)$, $U\in L([a,b],\lambda )$ where $\lambda $ is the Lebesgue measure. The function $U
$ is reguarly varying if and only if for any $x>0,$%
\begin{equation}
\lim_{x\rightarrow +\infty }U(\gamma x)/U(x)\text{ exists in }\mathbb{R}_{+}.
\label{evt.rvsv.DEFRVG}
\end{equation}
\end{lemma}

\bigskip \noindent \textbf{Proof}. It is clear that (\ref{evt.rvsv.DEFRVG}) holds, if $U$ is regularly varying. Now, suppose now that (\ref{evt.rvsv.DEFRVG}). Denote,
for $x>0,$%
\begin{equation*}
\lim_{t\rightarrow +\infty }U(tx)/U(t)=h(x).
\end{equation*}

\noindent Now for any $x>0,y>0,$%
\begin{eqnarray*}
h(xy) &=&\lim_{t\rightarrow +\infty }U(txy)/U(t) \\
&=&\lim_{t\rightarrow +\infty }\left\{ \frac{U(txy)}{U(tx)}\right\} \left\{ 
\frac{U(tx)}{U(t)}\right\} =h(x)h(y).
\end{eqnarray*}

\noindent The function $h:]0,+\infty \lbrack $ is measurable and satisfies the Hamel
Equation : $h(xy)=h(x)h(y),$ $x>0,y>0$. By Corollary \label{funct.cor.1} in Chapter \ref{funct}, the unique solution is%
\begin{equation*}
h(x)=x^{h(e)},x>0.
\end{equation*}

\noindent So $U\in RV(h(e),+\infty ).$

\bigskip \noindent  The next lemma gives some algebras on regular varying functions.\\

\begin{lemma} \label{evt.rv.alg}
If $U\in RV(\rho _{1},+\infty )$ and $V\in RV(\rho _{2},+\infty )$ then 

\noindent (1) $UV\in RV(\rho _{1}+\rho _{2},+\infty)$.\\

\noindent (2) $UV^{-1}\in RV(\rho _{1}-\rho_{2},+\infty)$.
\end{lemma}

\bigskip \noindent \textbf{Proof}. The proofs are immediate.\\

\bigskip Regular variation is used in Extreme value Theory by the so
important Karamata representation we are going to introduce in the next subsection.

\bigskip 

\subsection{Karamata Representation Theorem} \label{evt.subsec.karamata}

\bigskip \noindent Let us begin by these two lemmas.

\begin{lemma} \label{evt.karamata.lem1} (\text{de} Haan, 1970). Let $U:\mathbb{R}_{+}\longrightarrow 
\mathbb{R}_{+}$ be a measurable function Lebesgue locally integrable and
not vanishing in the neighborhood of $+\infty ,$ that is, $\forall\text{ }(0\leq a\leq b<+\infty)$, $U\in L([a,b],\lambda)$ where $\lambda $ is the Lebesgue measure. The following assertions holds.\\

\noindent (A) If $U\in RV(\rho ,+\infty ),\rho >-1$, $U^{\ast
}(x)=\int_{0}^{x}U(t)d\lambda (t)$ satisfies $U^{\ast }(+\infty )=+\infty $
and $U^{\ast }\in RV(\rho +1,+\infty)$.\\

\noindent (B) If $U\in RV(\rho ,+\infty ),\rho <-1$, $U_{\ast }(x)=\int_{x}^{+\infty
}U(t)d\lambda (t)$ satisfies $U_{\ast }(+\infty )<+\infty $ and $U_{\ast
}\in RV(\rho +1,+\infty)$.\\

\noindent (C) If $U\in RV(-1,+\infty ),$then  $U^{\ast }\in RV(0,+\infty )$. If
further $U_{\ast }(+\infty )<+\infty $, then $U_{\ast }\in RV(1,+\infty ).$
\end{lemma}

\noindent \textbf{Proof of Lemma \ref{evt.karamata.lem1}}.\\

\noindent \textbf{Part (A)}. Let $U\in RV(\rho ,+\infty ),\rho >-1$. We have to prove that $U^{\ast
}(+\infty )=+\infty .$ Let $c>1.$ We have $U(ct)/U(t)\rightarrow c^{\rho
}>c^{-1}.$ Let $\varepsilon >0$ such that $\rho -\varepsilon >-1$. Then, there
exists $x_{0}>1,$ such that for $t\geq x_{0},$%
\begin{equation*}
U(ct)/U(t)>c^{-1+\varepsilon }.
\end{equation*}
For $n$ such that $c^{n}\geq x_{0},$%
\begin{equation*}
\int_{c^{n}}^{c^{n+1}}U(t)d\lambda (t)=\int_{c^{n-1}}^{c^{n}}cU(ct)d\lambda
(t)>c^{\varepsilon }\int_{c^{n-1}}^{c^{n}}U(t)d\lambda (t).
\end{equation*}

\noindent Now let $n_{0}$ the first integer such that $c^{n_{0}}\geq x_{0}$ and $%
\int_{c^{n_{0}}}^{c^{n_{0}+1}}U(t)d\lambda (t)\neq 0.$ For any $n>n_{0},$ we
have
\begin{eqnarray*}
\int_{c^{n}}^{c^{n+1}}U(t)d\lambda (t) &>&c^{\varepsilon
}\int_{c^{n-1}}^{c^{n}}U(t)d\lambda (t)>c^{2\varepsilon
}\int_{c^{n-2}}^{c^{n-1}}U(t)d\lambda (t) \\
... &>&c^{(n-n_{0})\varepsilon }\int_{c^{n_{0}}}^{c^{n_{0}+1}}U(t)d\lambda
(t).
\end{eqnarray*}

\noindent Next
\begin{eqnarray*}
\int_{0}^{+\infty }U(t)d\lambda (t) &\geq &\int_{x_{0}}^{+\infty
}U(t)d\lambda (t)\geq \sum_{n\geq n_{0}}\int_{c^{n}}^{c^{n+1}}U(t)d\lambda
(t) \\
&\geq &\left\{ \int_{c^{n_{0}}}^{c^{n_{0}+1}}U(t)d\lambda (t)\right\}
\left\{ \sum_{n\geq n_{0}}c^{(n-n_{0})\varepsilon }\right\} .=+\infty .
\end{eqnarray*}

\noindent To complete the proof, we have to check that we can find $n_{0}$ such that $%
c^{n_{0}}\geq x_{0}$ and $\int_{c^{n_{0}}}^{c^{n_{0}+1}}U(t)d\lambda (t)\neq
0.$ If we cannot, then $U=0$ $a.e$. on $[N_{0},+\infty \lbrack $ where $N_{0}
$ is the first integer such that $c^{N_{0}}\geq 1.$ This is impossible
because of the assumption. Thus 
\begin{equation*}
U_{\ast }(+\infty )=+\infty .
\end{equation*}
For $\gamma >0,$ we have\bigskip 
\begin{equation*}
\frac{U^{\ast }(\gamma x)}{U^{\ast }(x)}=\frac{\int_{0}^{\gamma
x}U(t)d\lambda (t)}{\int_{0}^{x}U(t)d\lambda (t)}=\frac{\gamma
\int_{0}^{x}U(\gamma t)d\lambda (t)}{\int_{0}^{x}U(t)d\lambda (t)}.
\end{equation*}

\noindent Now we may apply Part (A) of \ Lemma \ref{evt.technical.hospital}, for $%
f(t)=U(\gamma t)$ and $g(t)=U(t)$ to get that%
\begin{equation}
\lim_{x\rightarrow +\infty }\frac{U^{\ast }(\gamma x)}{U^{\ast }(x)}=\gamma
\lim_{x\rightarrow +\infty }\frac{U(\gamma x)}{U(x)}=\gamma ^{\rho +1}.
\label{evt.technical.hosp1}
\end{equation}

\bigskip \noindent \textbf{Part (B)}. Let $U\in RV(\rho ,+\infty ),\rho <-1.$ Let us prove that $U_{\ast
}(+\infty )<\infty .$ Let $c>1.$ We have $U(ct)/U(t)\rightarrow c^{\rho
}<c^{-1}.$ Le $\varepsilon >0$ such that $\rho +\varepsilon <-1.$Then, there
exists $x_{0}>1,$ such that for $t\geq x_{0},$%
\begin{equation*}
U(ct)/U(t)<c^{-1-\varepsilon }=\delta <1.
\end{equation*}

\noindent For $n$ such that $c^{n}\geq x_{0},$%
\begin{equation*}
\int_{c^{n}}^{c^{n+1}}U(t)d\lambda (t)=\int_{c^{n-1}}^{c^{n}}cU(ct)d\lambda
(t)<\delta \int_{c^{n-1}}^{c^{n}}U(t)d\lambda (t).
\end{equation*}

\noindent Now let $n_{0}$ the first integer such that $c^{n_{0}}\geq x_{0}.$ For any $%
n>n_{0},$ we have
\begin{eqnarray*}
\int_{c^{n}}^{c^{n+1}}U(t)d\lambda (t) &<&\delta
\int_{c^{n-1}}^{c^{n}}U(t)d\lambda (t)>\delta ^{\varepsilon
}\int_{c^{n-2}}^{c^{n-1}}U(t)d\lambda (t) \\
... &<&\delta ^{n-n_{0}}\int_{c^{n_{0}}}^{c^{n_{0}+1}}U(t)d\lambda (t).
\end{eqnarray*}

\noindent Next
\begin{eqnarray*}
\int_{c^{n_{0}}}^{+\infty }U(t)d\lambda (t) &=&\sum_{n\geq
n_{0}}\int_{c^{n}}^{c^{n+1}}U(t)d\lambda (t) \\
&\leq &\left\{ \int_{c^{n_{0}}}^{c^{n_{0}+1}}U(t)d\lambda (t)\right\}
\left\{ \sum_{n\geq n_{0}}\delta ^{n-n_{0}}\right\}  \\
&=&(1-\delta )^{-1}\left\{ \int_{c^{n_{0}}}^{c^{n_{0}+1}}U(t)d\lambda
(t)\right\} <+\infty .
\end{eqnarray*}

\noindent Next

\begin{equation*}
\int_{0}^{+\infty }U(t)d\lambda (t)=\int_{0}^{c^{n_{0}}}U(t)d\lambda
(t)+\int_{c^{n_{0}}}^{+\infty }U(t)d\lambda (t)<+\infty .
\end{equation*}

\noindent For $\gamma >0,$ we have\bigskip 
\begin{equation*}
\frac{U_{\ast }(\gamma x)}{U_{\ast }(\gamma x)}=\frac{\int_{\gamma
x}^{+\infty }U(t)d\lambda (t)}{\int_{x}^{+\infty }U(t)d\lambda (t)}=\frac{%
\gamma \int_{x}^{+\infty }U(\gamma t)d\lambda (t)}{\int_{x}^{+\infty
}U(t)d\lambda (t)}.
\end{equation*}

\noindent Now we may apply Part (B) of  Lemma \ref{evt.technical.hospital}, for $%
f(t)=U(\gamma t)$ and $g(t)=U(t)$ to get that
\begin{equation}
\lim_{x\rightarrow +\infty }\frac{U^{\ast }(\gamma x)}{U^{\ast }(x)}=\gamma
\lim_{x\rightarrow +\infty }\frac{U(\gamma x)}{U(x)}=\gamma ^{\rho +1}.
\label{evt.technical.hosp2}
\end{equation}

\bigskip  \noindent \textbf{Part (C)}. Let $U\in RV(0,+\infty )$. As for $U^{\ast}$, either $U^{\ast}(+\infty )=+\infty$ and we use again the proof of Part (A) above from (\ref{evt.technical.hosp1}) that needed only that $U_{\ast }$ is infinite at $%
+\infty .$ Either $U_{\ast }(+\infty )=\ell >0$ finite and then for $\gamma
>0$,

\begin{equation*}
\frac{U^{\ast }(\gamma x)}{U^{\ast }(x)}\rightarrow \ell /\ell =1=\gamma
^{0}.
\end{equation*}

\noindent As for $U_{\ast }$, the condition $U_{\ast }(+\infty )<+\infty$ leads to (\ref{evt.technical.hosp2}) [ in the proof of Part 2] to get the same conclusion. In both case,
\begin{equation*}
U_{\ast }\in RV(0,+\infty ).
\end{equation*}

\bigskip The next lemma is the main Karamata Theorem. We expose it as a lemma and next give the induced representation as the theorem in this book.\\

\begin{lemma} \label{evt.karamata.lem2} We have the following assertions.\\

\noindent (A) If $U\in RV(\rho ,+\infty ),\rho \geq -1$, then
\begin{equation*}
b(x)=\frac{xU(x)}{\int_{0}^{x}U(t)d\lambda (t)}\rightarrow \rho +1\text{ as }%
x\rightarrow +\infty 
\end{equation*}

\noindent and there is a constant $c>0$,  for $x\geq 0$,
\begin{equation}
U(x)=cx^{-1}b(x)\exp \left( -\int_{1}^{x}t^{-1}b(t)dt\right) 
\end{equation}

\noindent (B) $U\in RV(\rho ,+\infty ),\rho <-1$. We have
\begin{equation*}
B(x)=\frac{xU(x)}{\int_{x}^{\infty }U(t)d\lambda (t)}\rightarrow -(\rho +1)%
\text{ as }x\rightarrow +\infty ,
\end{equation*}

\noindent and there is a constant $c>0$, for $x\geq 0$,
\begin{equation}
U(x)=cx^{-1}B(x)\exp \left( -\int_{1}^{x}t^{-1}B(t)dt\right) .
\end{equation}

\noindent The functions $b$ and $B$ are bounded.
\end{lemma}

\bigskip \noindent \textbf{Proof of Lemma \ref{evt.karamata.lem2}}.\\

\noindent \textbf{Part (A)}. We begin to remark that if $U\in RV(\rho ,+\infty )$ with $\rho
\geq -1$, then $xU(x)\in RV(\rho +1,+\infty )$ and $\int_{0}^{x}U(t)d\lambda
(t)\in RV(\rho +1,+\infty)$ by Lemma \ref{evt.karamata.lem1}. Applying
Lemma \ref{evt.rv.alg},
\begin{equation*}
b(\circ )\in VR(0,+\infty ).\end{equation*}

\noindent Now, by Lemma \ref{evt.karamata.lem1}, $b(\circ )\in RV(\rho+1)$. And
\begin{eqnarray*}
x^{-1}b(x) &=&\frac{U(x)}{\int_{0}^{x}U(t)d\lambda (t)}=\left(
\int_{0}^{x}U(t)d\lambda (t)\right) ^{\prime }/\left(
\int_{0}^{x}U(t)d\lambda (t)\right)  \\
&=&\left( \log \int_{0}^{x}U(t)d\lambda (t)\right) ^{\prime }.
\end{eqnarray*}

\noindent This gives, for $x>0$, 
\begin{equation*}
\log \int_{0}^{x}U(t)d\lambda (t)=\int_{1}^{x}t^{-1}b(t)dt+c_{1}, \text{ } a.e.
\end{equation*}

\noindent where $c_{1}$ is some constant. Then, for some constant $c>0,$ for for $x>0,$%
\begin{equation}
\int_{0}^{x}U(t)d\lambda (t)=c\exp \left( \int_{1}^{x}t^{-1}b(t)dt\right), \text{ } a.e.,
\label{evt.technical.rep0}
\end{equation}%

\noindent and next, for $x>$,
\begin{equation}
U(x)=cx^{-1}b(x)\exp \left( \int_{1}^{x}t^{-1}b(t)dt\right), \text{ } a.e., 
\label{evt.technical.rep1}
\end{equation}

\bigskip \noindent By making the change of variable $s=tx$, we have 
\begin{equation*}
\int_{0}^{x}U(s)d\lambda (s)=\int_{0}^{1}xU(tx)d\lambda (t).
\end{equation*}

\noindent Using this and the Fatou-Lebesgue Theorem, we have
\begin{eqnarray*}
\liminf_{x\rightarrow \infty }b(x)^{-1} &=&\liminf_{x\rightarrow \infty }
\frac{\int_{0}^{x}U(t)d\lambda (t)}{xU(x)} \\
&=&\liminf_{x\rightarrow \infty }\int_{0}^{1}\frac{U(tx)d\lambda (t)}{U(x)%
} \\
&\geq &\int_{0}^{1}\liminf_{x\rightarrow \infty }\left\{ \frac{U(tx)}{U(x)%
}\right\} d\lambda (t) \\
&=&\int_{0}^{1}t^{\rho}d\lambda (t)=(1+\rho )^{-1}.
\end{eqnarray*}

\noindent If $\rho =-1,$ we have $\liminf_{x\rightarrow \infty }b(x)^{-1}=+\infty $.
Hence $\limsup_{x\rightarrow +\infty }b(x)=0$ and then%
\begin{equation*}
b(x)\rightarrow 1+\rho \text{ for }\rho =-1.
\end{equation*}

\noindent What happens for $\rho >-1?$ We recall that%
\begin{equation*}
\liminf_{x\rightarrow \infty }b(x)^{-1}=\lim_{x\uparrow +\infty }\inf
\{b(t)^{-1},t\geq x\}.
\end{equation*}

\noindent Since $\inf \{b(t)^{-1},t\geq x\}\nearrow 1+\rho >0,$ then for some $%
x_{0}>0,\inf \{b(t)^{-1},t\geq x\}>(1+\rho )/2.$ This implies that $%
b(t)^{-1}>(1+\rho )/2=M>0$ for all $t\geq x_{0}.$ It comes that $b(\circ)$ is bounded on $[x_{0},+\infty \lbrack .$ Thus, for $x\geq x_{0},$ for $%
t\geq 1,$ for $s\geq 1,$ 
\begin{equation*}
\left\vert \frac{b(tx)-b(x)}{t}\right\vert \leq Mt^{-1}\in L([1,s],\lambda ).
\end{equation*}

\noindent Further, since $b(\circ )\in RV(0,+\infty ),$ we have for any fixed $t\in
\lbrack 1,s],$%
\begin{equation*}
\frac{b(tx)-b(x)}{t}\rightarrow 0\text{ as }x\rightarrow \infty .
\end{equation*}

\noindent We may apply the Dominated Lebesgue Theorem to get
\begin{equation*}
\int_{1}^{s}\frac{b(tx)-b(x)}{t}d\lambda (t)\rightarrow 0.
\end{equation*}

$\forall\text{ }(0\leq a\leq b<+\infty)$, $U\in L([a,b],\lambda)$

\noindent But, as $x\in +\infty$,
\begin{equation}
\int_{1}^{s}\frac{b(tx)-b(x)}{t}d\lambda (t)=\int_{1}^{s}t^{-1}b(tx)d\lambda
(t)-b(x)\log s\rightarrow 0. 
\label{evt.technical.limA}
\end{equation}

\noindent Now from (\ref{evt.technical.rep0}) and from Part (A) of Lemma \ref%
{evt.karamata.lem1}, $h(x)=\exp \left( \int_{1}^{x}t^{-1}b(t)dt\right) \in
RV(\rho +1,+\infty )$. Then
\begin{eqnarray*}
\frac{h(sx)}{h(x)} &=&\exp (\int_{x}^{xs}t^{-1}b(t)dt) \\
&=&\exp (\int_{1}^{s}u^{-1}b(xu)du).
\end{eqnarray*}

Then
\begin{equation}
\int_{1}^{s}t^{-1}b(tx)dt=\log \left\{ \frac{h(sx)}{h(x)}\right\}
\rightarrow \log s^{(1+\rho )}\text{ as }x\rightarrow \infty .
\label{evt.technical.limB}
\end{equation}

\bigskip \noindent By comparing (\ref{evt.technical.limA}) and (\ref{evt.technical.limB}), we get

\begin{equation*}
\int_{1}^{s}t^{-1}b(tx)dt-b(x)\log s\rightarrow 0\text{ and }%
\int_{1}^{s}t^{-1}b(tx)dt-(1+\rho )\log s\rightarrow 0\text{ .}
\end{equation*}

\noindent This leads to 
\begin{equation*}
b(x)\rightarrow 1+\rho \text{ as }x\rightarrow +\infty. 
\end{equation*}

\noindent This finishes the proof of Part A.\\

\bigskip \noindent \textbf{Part (B)}. As in Part (A), if $U\in RV(\rho ,+\infty )$ with $\rho
<-1$, then $xU(x)\in RV(\rho +1,+\infty )$ and $\int_{x}^{+\infty
}U(t)d\lambda (t)\in RV(\rho +1,+\infty )$ by Lemma \ref{evt.karamata.lem1}.
By applying Lemma \ref{evt.rv.alg}, we get
\begin{equation*}
B(\circ )\in VR(0,+\infty ).
\end{equation*}

\bigskip \noindent 
Since $A=\int_{0}^{\infty }U(t)d\lambda (t)$ is finite by Lemma \ref{evt.karamata.lem1}, we have
\begin{eqnarray*}
x^{-1}B(x) &=&\frac{U(x)}{A-\int_{0}^{x}U(t)d\lambda (t)}=-\left(
A-\int_{0}^{x}U(t)d\lambda (t)\right) ^{\prime }/\left(
A-\int_{0}^{x}U(t)d\lambda (t)\right)  \\
&=&-\left( \log \left( A-\int_{0}^{x}U(t)d\lambda (t)\right) \right)
^{\prime },\text{ }a.e.,
\end{eqnarray*}

\noindent This gives, for $x>0$,
\begin{equation*}
\log \left( A-\int_{0}^{x}U(t)d\lambda (t)\right)
=-\int_{1}^{x}t^{-1}B(t)dt+c_{1}, \text{ } a.e.
\end{equation*}

\noindent where $c_{1}$ is some constant. Then, for some constant $c>0$,  for $x>0$,
\begin{equation} \normalsize
\int_{x}^{+\infty }U(t)d\lambda (t)=A-\int_{0}^{x}U(t)d\lambda (t)=c\exp
\left( -\int_{1}^{x}t^{-1}B(t)dt\right), \text{ } a.e.  \label{evt.technical.rep3}
\end{equation}

\noindent and next, for $x>0$,
\begin{equation}
U(x)=cx^{-1}B(x)\exp \left( -\int_{1}^{x}t^{-1}B(t)dt\right). 
\end{equation}

\noindent From here, we will be able to conclude by re-conducting the same methods as in Part (A), if we prove the analogues of 
(\ref{evt.technical.limA}) and (\ref{evt.technical.limB}) for $B(\circ)$, that is, as $x \rightarrow +\infty$,

\begin{equation}
\int_{1}^{s}t^{-1}B(tx)dt - B(x) \log s \rightarrow 0, \label{evt.technical.limAA}
\end{equation}

\noindent and 

\begin{equation}
\int_{1}^{s}t^{-1}B(tx)dt + (1+\rho ) \log s \rightarrow 0. \label{evt.technical.limAB}
\end{equation}

\noindent To get (\ref{evt.technical.limAA}), we may get it as as we did for \ref{evt.technical.limA} only by showing that $B(\circ)$ is ultimately bounded. To get this, we also use the Fatou-Lebesgue Theorem to have
\begin{eqnarray*}
\liminf_{x\rightarrow \infty }b(x)^{-1} &=&\liminf_{x\rightarrow +\infty}
\frac{\int_{x}^{+\infty}U(t)d\lambda (t)}{xU(x)} \\
&=&\liminf_{x\rightarrow \infty }\int_{1}^{+\infty}\frac{U(tx)d\lambda (t)}{U(x)%
} \\
&\geq &\int_{1}^{+\infty}\liminf_{x\rightarrow \infty }\left\{ \frac{U(tx)}{U(x)%
}\right\} d\lambda (t) \\
&=&\int_{0}^{1}t^{\rho}d\lambda (t)=-(1+\rho)^{-1},
\end{eqnarray*}

\noindent where we took into account that $1+\rho < 0$. As in Part (B), this ensures that $B(\circ)$ and this ensures (\ref{evt.technical.limAA}).\\
\\

\noindent To establish (\ref{evt.technical.limAB}), we may see from (\ref{evt.technical.rep3}) and from Part (B) of Lemma \ref%
{evt.karamata.lem1} that  
$$
h(x)=\exp \left( -\int_{1}^{x}t^{-1}B(t)dt\right) \in RV(\rho +1,+\infty).
$$ 

\noindent Then

\begin{equation*}
\frac{h(sx)}{h(x)}=\exp (-\int_{x}^{xs}t^{-1}B(t)dt)=\exp
(-\int_{1}^{s}t^{-1}B(xt)dt)
\end{equation*}

\noindent and
\begin{equation}
\int_{1}^{s}t^{-1}B(tx)dt=\log \left\{ \frac{h(sx)}{h(x)}\right\}
\rightarrow -(1+\rho )\log s\text{ as }x\rightarrow \infty ,
\label{evt.technical.limC}
\end{equation}

\noindent and this implies  
\begin{equation*}
\int_{1}^{s}t^{-1}(B(tx)+(1+\rho ))dt\rightarrow 0\text{ as }x\rightarrow
+\infty,
\end{equation*}

\noindent which is (\ref{evt.technical.limAB}). Besides, it is clear from the details of the proof that the functions $b$ and $c$ and $B$ are bounded.\\

\noindent The proof of this important theorem is now complete.\\

\bigskip \noindent The following representation result is a key tool in Extreme value Theory.\\ 

\begin{theorem} \label{evt.karamata.theo} (Karamata's Theorem) Let $U:\mathbb{R}_{+}\longrightarrow \mathbb{R%
}_{+}$ be a measurable function Lebesgue locally integrable and not
vanishing in the neighborhood of $+\infty ,$ that is, $\forall \text{}(0\leq a\leq
b<+\infty)$, $U\in L([a,b],\lambda )$ where $\lambda $ is the Lebesgue
measure.\\

\noindent Then $U\in RV(\rho ,+\infty )$ if and only there exist two measurables functions 
$a(x)$ and $\ell(x)$ of $x\in \mathbb{R}$ and a constant $c>0$ such that $(p(x),\ell (x))\rightarrow (0,0)$ as $x\rightarrow \infty$, for any $x\geq
0$
\begin{equation*}
U(x)=c(1+a(x))\exp (\int_{1}^{x}t^{-1}\ell (t)d\lambda (t)).
\end{equation*}

\noindent Besides, the functions $a(x)$ and $\ell(x)$ of $x\in \mathbb{R}$ are bounded in neighborhood of $+\infty$.
\end{theorem}

\bigskip \noindent P\textbf{roof}. Let $U\in RV(\rho ,+\infty ).$ Let first $\rho \neq -1.$ Use
Lemma \ref{evt.karamata.lem2}. Take $m=b$ in Part (A) and $m=-B$ in Part (B)
so that $\ell (x)=m(x)-(1-\rho )\rightarrow 0$ as $x\rightarrow \infty $ and
next%
\begin{equation*}
U(x)=c\left\vert 1+\rho \right\vert (1+p(x))x^{-1}\exp (\int_{1}^{x}\frac{%
(1-\rho )}{t}d\lambda (t))\exp (\int_{1}^{x}t^{-1}\ell (t)d\lambda (t)),
\end{equation*}
where $p(x)\rightarrow 0$ as $x\rightarrow \infty .$ For $C=c\left\vert
1+\rho \right\vert ,$ we get for $x>0,$%
\begin{equation}
U(x)=C(1+p(x))x^{\rho }\exp (\int_{1}^{x}t^{-1}\ell (t)d\lambda (t)).
\label{evt.karamataG}
\end{equation}

\noindent If $\rho =-1,$ then $xU(x)\in RV(0,+\infty ).$ Then $xU(x)$ admits a
representation
\begin{equation*}
xU(x)=C(1+p(x))\exp (\int_{1}^{x}t^{-1}\ell (t)d\lambda (t))
\end{equation*}

\noindent which implies

\begin{equation*}
U(x)=C(1+p(x))x^{-1}\exp (\int_{1}^{x}t^{-1}\ell (t)d\lambda (t)).
\end{equation*}

\noindent The reverse direction is quite straightforward. Suppose that  (\ref%
{evt.karamataG}) holds. Then for any fixed $\gamma >0,$ we have%
\begin{equation*}
U(\gamma x)/U(x)=\frac{1+p(\gamma x)}{1+p(x)}\gamma ^{\rho }\exp
(\int_{x}^{\gamma x}t^{-1}\ell (t)d\lambda (t)).
\end{equation*}

\noindent For any $\varepsilon >0,$ we have for $x>0,$ large enough, $0<1-\varepsilon
\leq 1+p(\gamma x),1+p(x)\leq 1+\varepsilon $ and $\max \{\left\vert
b(t)\right\vert ,\left\vert t\right\vert \leq \max (x,\gamma x)\}\leq
\varepsilon $ and then

\begin{equation*}
\frac{1-\varepsilon }{1+\varepsilon }\gamma ^{\rho -\varepsilon }\leq \frac{%
U(\gamma x)}{U(x)}\leq \frac{1+\varepsilon }{1-\varepsilon }\gamma ^{\rho
+\varepsilon }.
\end{equation*}

\noindent Hence
\begin{equation*}
\frac{1-\varepsilon }{1+\varepsilon }\gamma ^{\rho -\varepsilon }\leq \liminf_{x\rightarrow +\infty }\frac{U(\gamma x)}{U(x)}\leq \limsup_{x\rightarrow +\infty }\frac{U(\gamma x)}{U(x)}\leq \frac{1+\varepsilon 
}{1-\varepsilon }\gamma ^{\rho +\varepsilon }.
\end{equation*}

\noindent Letting $\varepsilon \downarrow 0$ leads to

\begin{equation*}
\lim_{x\rightarrow +\infty }\frac{U(\gamma x)}{U(x)}=\gamma ^{\rho }.
\end{equation*}

\noindent The measurability and the boundedness of $p$ and $\ell$ come from that of $m$, that is that of $b$ and $B$.\\

\noindent The proof is now complete.\\

\bigskip \noindent \textbf{Remark}. We may see in this proof that the representation (\ref%
{evt.karamataG}) is true, whenever $b(x)\rightarrow \lambda =\rho +1\in
]0,+\infty \lbrack $ or $B(x)\rightarrow \lambda =-(\rho +1)\in ]0,+\infty
\lbrack $ and (\ref{evt.karamataG}) ensures that $U\in RV(\rho ,+\infty ).$
We then get this inverse lemma to Lemma \ref{evt.karamata.lem2}

\begin{lemma} \label{evt.karamata.lem3} Let $U:\mathbb{R}_{+}\longrightarrow \mathbb{R}_{+}$ be
a measurable function Lebesgue locally integrable and not vanishing in the
neighborhood of $+\infty$ that is, $\forall 0\leq a\leq b<+\infty ,U\in
L([a,b],\lambda )$ where $\lambda $ is the Lebesgue measure. If, for some $\rho \in \mathbb{R}$,

\begin{equation*}
b(x)=\frac{xU(x)}{\int_{0}^{x}U(t)d\lambda (t)}\rightarrow \lambda \in
]0,+\infty \lbrack \text{ as }x\rightarrow +\infty ,
\end{equation*}

\noindent then $\ U\in RV(\lambda -1,+\infty )$.\\

\noindent If $U$ in integrable over subsets of the form $[a,+\infty[$, $0<a\in \mathbb{R}$ and
\begin{equation*}
b(x)=\frac{xU(x)}{\int_{x}^{+\infty }U(t)d\lambda (t)}\rightarrow \lambda
\in ]0,+\infty \lbrack \text{ as }x\rightarrow +\infty ,
\end{equation*}

\noindent then $U\in RV(-\lambda -1,+\infty ).$
\end{lemma}

\subsection{Weak conditions for regular variation}

It may be handy to have weaker conditions for establishing that some function is regularly varying. This subsection provides such weaker conditions.  
 
\begin{lemma} \label{evt.rvsv.VRWEAK} Let $U:\mathbb{R}_{+}\longrightarrow \mathbb{R}_{+}$ be a monotone function
not vanishing in the neighborhood of $+\infty .$ Then $U$ is regularly
varying if and only if there exist two sequences $(\lambda _{n})_{n\geq 1}$
and $(a_{n})_{n\geq 1}$ such that%
\begin{equation*}
\lim_{n\rightarrow +\infty }\lambda _{n+1}/\lambda _{n}=1\text{ and }%
\lim_{n+\infty }a_{n}=+\infty 
\end{equation*}

\noindent and such that for all $x>0$, $\lambda _{n}U(a_{n}x)$ admits a limit in $\mathbb{R}_{+}\setminus \{0\},$ denoted by%
\begin{equation*}
\lim_{n\rightarrow +\infty }\lambda _{n}U(a_{n}x)=\ell (x).
\end{equation*}

\noindent Moreover, the function $\ell (\circ )$ is regularly varying with the same
exponent as $U(\circ ).$
\end{lemma}

\bigskip \noindent \textbf{Proof}. Let us prove it for $U$ non-increasing first. Then%
\begin{equation*}
\inf U=U(+\infty )=\lim_{x\uparrow +\infty }U(x)\leq U(0).
\end{equation*}

\noindent \textbf{Case $U(+\infty )>0$}. Then $U(0+)\in ]0,U(0)[.$ We have for any $x>0$,%
\begin{equation*}
\lim_{t\rightarrow +\infty }\frac{U(tx)}{U(t)}=\frac{U(+\infty )}{U(+\infty)}=1,
\end{equation*}

\noindent that is : $U\in RV(0,+\infty ).$ At the same time, we have, by taking $\lambda
_{n}=1$ and $a_{n}=n,$ that%
\begin{equation*}
\lim_{n\rightarrow +\infty }\lambda _{n}U(a_{n}x)=U(+\infty )=\ell (x)\text{
positive and finite,}
\end{equation*}

\noindent which is also in \ $RV(0,+\infty ).$ So the equivalence holds for $U(0+)>0.$\\

\bigskip \noindent \textbf{Case $U(+\infty)=0$}.\\

\noindent \textbf{Part (A)}. Suppose that there exist two sequences $(\lambda _{n})_{n\geq 1}$
and $(a_{n})_{n\geq 1}$ such that%
\begin{equation*}
\lim_{n\rightarrow +\infty }\lambda _{n+1}/\lambda _{n}=1\text{ and }%
\lim_{n+\infty }a_{n}=+\infty .
\end{equation*}

\noindent Define for $t>0,$%
\begin{equation*}
n(t)=\inf \{n\geq 0,a_{n+1}>t\},
\end{equation*}

\noindent which implies for any $t>0$,
\begin{equation*}
a_{n(t)}\leq t<a_{n(t)+1}.
\end{equation*}

\noindent It is clear that $n(t)\rightarrow \infty $ as $t\rightarrow +\infty$. Using the monotonicity of $U$, we have%
\begin{eqnarray*}
\frac{U(tx)}{U(t)} &\leq &\frac{U(xa_{n(t)+1})}{U(a_{n(t)})} \\
&=&\frac{\lambda _{n(t)+1}U(xa_{n(t)+1})}{\lambda _{n(t)}U(a_{n(t)})}\left\{ 
\frac{\lambda _{n(t)}}{\lambda _{n(t)+1}}\right\} \rightarrow \frac{\ell (x)%
}{\ell (1)}.
\end{eqnarray*}

\noindent Also

\begin{eqnarray*}
\frac{U(tx)}{U(t)} &\geq &\frac{U(xa_{n(t)})}{U(a_{n(t)+1})} \\
&=&\frac{\lambda _{n(t)}U(xa_{n(t)})}{\lambda _{n(t)+1}U(a_{n(t)+1})}\left\{ 
\frac{\lambda _{n(t)+1}}{\lambda _{n(t)}}\right\} \rightarrow \frac{\ell (x)%
}{\ell (1)}.
\end{eqnarray*}%

\noindent Then 
\begin{equation*}
\frac{U(tx)}{U(t)}\rightarrow h(x)=\ell (x)/\ell (1)\text{ finite.}
\end{equation*}

\noindent By Lemma \ref{evt.rvsv.charac}, we conclude that $U\in RV(\ell (e)/\ell
(1),+\infty )$.\\

\noindent \textbf{Part (B)}. Now, let us suppose that $U$ is regularly varying. Recall that $U$
is non-increasing. We may define the generalized inverse%
\begin{equation*}
V(y)=\inf \{x>0,U(x)\leq y\}.
\end{equation*}

\noindent Let us show that for any $y>0,$%
\begin{equation}
U(V(y)+)\leq y\leq U(V(y)-) \text{ }y>0.  \label{evt.rvsv.Uinv}
\end{equation}

\bigskip \noindent By definition of the infimum, there exists a sequence $(x_{n})_{n\geq 1}$ such
that

\begin{equation*}
\forall (n\geq 1),U(x_{n})\leq y\text{ and }x_{n}\downarrow V(y)\text{ as }%
n\uparrow +\infty .
\end{equation*}

\noindent Then by letting $n\uparrow +\infty $ in $U(x_{n})\leq y$ we get $U(V(y)+)\leq y$. Next, for any $\eta >0,V(y)-\eta $ cannot satisfy $U(V(y)-\eta )\leq y$, otherwise $V(y)$ would not be the infimum of $\{x>0,U(x)\leq y\}$. So  $U(V(y)-\eta )\geq y$ for all $\eta>0$. By letting 
$\eta \downarrow 0$, we get $U(V(y)-)\geq y$. Since $U(\circ )$ does not take
the null value in the neighborhood of zero and since $U(+\infty)=0$, it becomes
clear that  
\begin{equation}
\lim_{y\rightarrow +\infty }V(y)=+\infty .  \label{evt.rvsv.inf}
\end{equation}

\noindent Suppose now that $U \in RV(\rho ,+\infty)$. Let us also show that%
\begin{equation}
\lim_{t\rightarrow +\infty }\frac{U(t-)}{U(t)}=\lim_{t\rightarrow +\infty }%
\frac{U(t+)}{U(t)}=1.  \label{evt.rvsv.ratio3}
\end{equation}

\noindent To see this, let $0<x<1.$ For any $t>0$ fixed, for any $s$ such that $tx<s<t,
$ $U(t)\leq U(s)\leq U(tx)$. Then, we have

\begin{equation*}
1\leq \frac{U(s)}{U(t)}\leq \frac{U(tx)}{U(t)}.
\end{equation*}

\noindent By letting $s\uparrow 0$ and next $t\uparrow +\infty $, we get 
\begin{equation}
1\leq \limsup_{t\rightarrow +\infty }\frac{U(t-)}{U(t)}\leq
\lim_{t\rightarrow +\infty }\frac{U(tx)}{U(t)}=x^{\rho }.
\label{evt.rvsv.ratio1}
\end{equation}

\noindent If $\rho =0$ we have 
\begin{equation}
\lim_{t\rightarrow +\infty }\frac{U(t-)}{U(t)}=1.  \label{evt.rvsv.ratio2}
\end{equation}

\noindent If $\rho \neq 0,$ we let $x\uparrow 1$ in (\ref{evt.rvsv.ratio1}) to get (\ref%
{evt.rvsv.ratio2}). This proves the first limit in (\ref{evt.rvsv.ratio3}).
To get the second limit, consider $x>1,t>0,$ $t<s<tx,$ to get

\begin{equation*}
\frac{U(tx)}{U(t)}\leq \frac{U(s)}{U(t)}\leq 1.
\end{equation*}

\noindent We will able to conclude by letting successively $s\downarrow t$, $t\uparrow +\infty$ and $x\downarrow 1$.\\

\noindent Combining (\ref{evt.rvsv.Uinv}), (\ref{evt.rvsv.inf}) and (\ref{evt.rvsv.ratio3}) leads to
\begin{equation*}
\lim_{y\rightarrow +\infty }U(V(y))/y=1.
\end{equation*}

\noindent Finally take $\lambda _{n}=n$ and $a_{n}=V(1/n)$ for $n\geq 1.$ We get%
\begin{eqnarray*}
\lim_{n\rightarrow +\infty }\lambda _{n}U(a_{n}x) &=&\lim_{n\rightarrow
+\infty }\left\{ \lambda _{n}U(a_{n}x)\right\} \left\{ \frac{U(a_{n}x)}{%
U(a_{n})}\right\}  \\
&=&\lim_{n\rightarrow +\infty }\left\{ \frac{U(V(1/n))}{(1/n)}\right\}
\left\{ \frac{U(a_{n}x)}{U(a_{n})}\right\} =x^{\rho }.
\end{eqnarray*}

\noindent The proof of the lemma is complete for $U$ non-increasing. To extend this
a non-increasing fonction $U$, we use the transform $V=1/U$ and take into account the three equivalences :\\

\noindent (1) :  ($U$ non-increasing, finite and nonzero in the neighberhood of $+\infty
)\Longleftrightarrow$ ($1/U$ non-decreasing finite and nonzero in the
neighborhood of $+\infty$).\\

\noindent (2) : $\left( U\in RV(\rho ,+\infty ),\rho \neq 0\right) \Longleftrightarrow
\left( 1/U\in RV(-\rho ,+\infty ),\rho \neq 0\right)$.\\

\noindent (3) There exist two sequences $(\lambda_{n})_{n\geq 1}$ and $(a_{n})_{n\geq 1}$
such that 
$$
\lim_{n\rightarrow +\infty }\lambda _{n+1}/\lambda _{n}=1 \text{ and } \lim_{n+\infty }a_{n}=+\infty,
$$ 

\noindent and such that for all $x>0$ 

$$
\lim_{n\rightarrow +\infty}\lambda _{n}U(a_{n}x)=\ell (x)
$$ 

\noindent is equivalent to : \\

\noindent there exist two sequences $(\lambda _{n}^{\prime })_{n\geq 1}$ and $(a_{n}^{\prime })_{n\geq
1}$ such that 
$$
\lim_{n\rightarrow +\infty }\lambda _{n+1}^{\prime }/\lambda
_{n}^{\prime }=1 \text{ and} \lim_{n+\infty }a_{n}^{\prime }=+\infty,
$$

\noindent and such that for all $x>0$,

$$
\lim_{n\rightarrow +\infty }\lambda _{n}U(a_{n}x)=\ell^{\prime }(x)
$$ 

\noindent with

\begin{equation*}
\lambda _{n}^{\prime }=1/\lambda _{n}^{^{\prime }},\text{ }a_{n}^{\prime
}=a_{n},\text{ }\ell ^{\prime }(x)=1/\ell (x),x>0.
\end{equation*}

\noindent The proof is now complete with these remarks.\\

\bigskip  \bigskip \noindent  Inspired by the proof of Theorem 1.1.2 in \cite{dehaan}, we get this
result, which is more simpler.

\begin{lemma} \label{evt.rvsv.gslo} Let $U:\mathbb{R}_{+}\longrightarrow \mathbb{R}_{+}$ be a monotone function not
vanishing in the neighborhood of $+\infty$. Then $U\in RV(\rho ,+\infty )$
if and only if for any integer $p>0$,
\begin{equation}
\lim_{t \rightarrow +\infty }U(pt)/U(t)=p^{\rho}. \label{evt.rvsv.DEFRVN}
\end{equation}
\end{lemma}

\noindent \textbf{Proof}. We have to prove that (\ref{evt.rvsv.DEFRVN}) implies that $U\in
RV(\rho ,+\infty)$. So assume that (\ref{evt.rvsv.DEFRVN}) holds and at the
first place that $U$ is non-decreasing. Then we necessarily have $\rho
\geq 0$. Let us begin to prove that 
\begin{equation}
\lim_{n\rightarrow +1}U(n+1)/U(n)=1.  \label{evt.rvsv.ratio4}
\end{equation}

\noindent Let $c>1$. For large enough value of $n$, we have $(n+1)/n<c$, and for such values of $n$,
\begin{equation*}
1\leq \frac{U(n+1)}{U(n)}\leq \frac{U(cn)}{U(n)}
\end{equation*}

\noindent and

\begin{equation*}
1\leq \lim_{n\rightarrow +\infty }\sup \frac{U(n+1)}{U(n)}\leq
\lim_{n\rightarrow +\infty }\frac{U(cn)}{U(n)}=c^{\rho }.
\end{equation*}

\bigskip

\noindent Letting $c\downarrow 1$ proves (\ref{evt.rvsv.ratio4}). Now let $x>0$ fixed.
We want to extend (\ref{evt.rvsv.ratio4}) to $p=x \mathbb{R}_{+}$. Let $\varepsilon >0.$ By
density of the set $\mathbb{Q}_{+}$\ of rational positive numbers in $%
\mathbb{R}_{+},$ we may find $r^{\prime }$ and $r^{\prime \prime }$ such
that $x-\varepsilon <r^{\prime }<x<r"<x+\varepsilon .$ \ By writing $%
r^{\prime }$ and $r^{\prime \prime }$ with a common denominator, wet get
that there exists non-negative integers $p$, $q$ and $r\neq 0$ such that    
\begin{equation*}
x-\varepsilon <(p/r)<x<(q/r)<x+\varepsilon .
\end{equation*}

\noindent Take for any $t>0,n(t)=[t/r]$, that is $n(t)$ is integer and 
\begin{equation*}
n(t)r\leq t<r(n(t)+1.
\end{equation*}

\noindent We easily check that $n(t)\rightarrow +\infty $ as $t\rightarrow +\infty .$ By
combining the last two double inequalities, we have

\begin{equation*}
pn(t)\leq tx\leq q(n(t)+1)
\end{equation*}

\noindent By non-decreasingness of $U$, we have 
\begin{eqnarray*}
\frac{U(tx)}{U(t)} &\leq &\frac{U(q(n(t)+1))}{U(rn(t))} \\
&=&\left\{ \frac{U(q(n(t)+1))}{U(n(t)+1)}\right\} \left\{ \frac{U(n(t))}{%
U(rn(t))}\right\} \left\{ \frac{U(n(t)+1)}{U(n(t))}\right\}  \\
&\rightarrow &q^{\rho }r^{\rho }=(q/r)^{\rho }\leq (x+\varepsilon )^{\rho },
\end{eqnarray*}

\noindent since $\rho \geq 0$. Likely, we have 
\begin{eqnarray*}
\frac{U(tx)}{U(t)} &\geq &\frac{U(p(n(t)+1))}{U(rn(t))} \\
&=&\left\{ \frac{U(p(n(t)+1))}{U(n(t)+1)}\right\} \left\{ \frac{U(n(t))}{%
U(rn(t))}\right\} \left\{ \frac{U(n(t)+1)}{U(n(t))}\right\}  \\
&\rightarrow &p^{\rho }r^{\rho }=(p/r)^{\rho }\leq (x-\varepsilon )^{\rho }.
\end{eqnarray*}%
\begin{equation*}
(x-\varepsilon )^{\rho }\leq \lim_{x+\infty }\inf \frac{U(tx)}{U(t)}\leq
\lim_{x+\infty }\sup \frac{U(tx)}{U(t)}\leq (x+\varepsilon )^{\rho },
\end{equation*}

\noindent for any $\varepsilon >0.$ We get the searched result by letting $\varepsilon
\downarrow 0.$\\

\bigskip \noindent To finish, we have to give the proof for $U$ non-increasing. We easily get
the proof by using the transform $1/U$. 

\bigskip 

\subsubsection{Regular variation and generalized inverses}

For the needs of Extremen value Theory for example, we frequently use the
quantile function which is a generalized inverse of the distribution
function. This lemma will greatly help.

\bigskip 

\begin{proposition} \label{evt.rvsv.Inverse} Let $U:\mathbb{R}_{+}\longrightarrow \mathbb{R}_{+}$ be a
a non-constant and monotone function not vanishing in the
neighborhood of $+\infty $ such that $U(+\infty )=+\infty .$ If $U\in
RV(\rho ,+\infty ),$ $\rho\neq0$, then the generalized inverse of $U$ defined when $U$ is non-decreasing by
\begin{equation}
V(y)=\inf \{x>0,U(x)\geq x\}
\end{equation}

\noindent and when $U$ is non-increasing by

\begin{equation}
V(y)=\inf \{x>0,U(x)\leq x\}
\end{equation}

\noindent is in $RV(1/\rho ,+\infty).$
\end{proposition}

\noindent To prove this, we need this lemma.

\begin{lemma} \label{evt.rvsv.lem7} Let $U:\mathbb{R}_{+}\longrightarrow \mathbb{R}_{+}$ $\rho -$%
regularly varying. Let $(a_{n})_{n\geq 0}$ and $(b_{n})_{n\geq 0}$ be two
sequences of positive real numbers such that  $a_{n}/b_{n}\rightarrow c$  as 
$n\rightarrow +\infty .$ If $c\in \mathbb{R}_{+}\setminus \{0\},$ then%
\begin{equation}
\lim_{x\rightarrow +\infty }\frac{U(a_{n})}{U(b_{n})}=c^{\rho }.
\label{evt.rvsv.ratioT}
\end{equation}

\noindent If $\rho\neq0$, then for any $c\in \mathbb{R}_{+}\cup \{+\infty \},$ (\ref{evt.rvsv.ratioT}) holds.
\end{lemma}

\bigskip 

\noindent \textbf{Proof of \ref{evt.rvsv.lem7}}. First, let $c\in \mathbb{R}_{+}\setminus
\{0\}$. We are going to use of Karamata representation in Theorem \ref{evt.karamata.theo} to write

\begin{equation*}
U(x)=c(1+p(x))x^{\rho }\exp (\int_{1}^{x}t^{-1}\ell (t)d\lambda (t)),
\end{equation*}
with ($p(x),\ell (x))\rightarrow (0,0)$ as $x\rightarrow +\infty .$ We have%
\begin{eqnarray}
\frac{U(a_{n})}{U(b_{n})} &=&\left\{ \frac{1+p(a_{n})}{1+p(b_{n})}\right\}
\left( \frac{a_{n}}{b_{n}}\right) ^{\rho }\exp
(\int_{b_{n}}^{a_{n}}t^{-1}\ell (t)d\lambda (t))  \label{evt.rvsv.ratioT1} \\
&=&(1+o(1))\left( \frac{a_{n}}{b_{n}}\right) ^{\rho }\exp
(\int_{b_{n}}^{a_{n}}t^{-1}\ell (t)d\lambda (t)).  \notag
\end{eqnarray}

\noindent Set 
\begin{equation*}
\varepsilon _{n}=\sup \{\left\vert b(t)\right\vert ,0\leq t\leq \max
(a_{n},b_{n})\}
\end{equation*}

\noindent and remark that, $a_{n}/b_{n}\rightarrow c$ finite and positive, as
$n\rightarrow +\infty$. Hence,  for any fixed $\eta >0$ with $\eta <c,$
there exists $n_{0}$ such that for any $n\geq n_{0},$ $b_{n}(c-\eta )\leq
a_{n}\leq b_{n}(c+\eta )$. We get from all this that

\begin{equation*}
\left\vert \int_{b_{n}}^{a_{n}}t^{-1}\ell (t)d\lambda (t)\right\vert \leq
\varepsilon _{n}\log \left\{ \frac{\max (a_{n},b_{n})}{\max (a_{n},b_{n})}%
\right\} \leq \varepsilon _{n}\log \frac{c+\eta }{c-\eta }=A_{n}\rightarrow
0.
\end{equation*}

\noindent Then
\begin{equation*}
\frac{U(a_{n})}{U(b_{n})}=\left\{ \frac{1+p(a_{n})}{1+p(b_{n})}\right\}
\left( \frac{a_{n}}{b_{n}}\right) ^{\rho }(1+O(A_{n}))\rightarrow c^{\rho }.
\end{equation*}

\noindent So (\ref{evt.rvsv.ratioT}) is true for $c\in \mathbb{R}_{+}\setminus \{0\}.$\\

\noindent From there, we prove the second statement of the lemma for $\rho >0$. If the statement holds for $\rho>0$ and if
$U \in RV(-\rho,+\infty)$, we apply it to $1/U$ which is in $U \in RV(rho,+\infty)$ and get if for $U$.\\ 

\bigskip \noindent \textbf{So, in our second step}, we let $\rho >0$ and suppose that $c=0.$ We may write $a_{n}=b_{n}r_{n}$
where $r_{n}=a_{n}/b_{n}\rightarrow 0.$ We surely have $a_{n}<b_{n}$ for
large values. For those values,%
\begin{equation*}
\left\vert \int_{b_{n}}^{a_{n}}t^{-1}\ell (t)d\lambda (t)\right\vert \leq
\varepsilon _{n}\log \frac{b_{n}}{a_{n}}=-\varepsilon _{n}\log
r_{n}=r_{n}^{-\varepsilon _{n}}.
\end{equation*}

\noindent Plugging this in (\ref{evt.rvsv.ratioT1}) leads to  
\begin{equation*}
(1+o(1))r_{n}{}^{\rho +\varepsilon _{n}}\leq \frac{U(a_{n})}{U(b_{n})}\leq
(1+o(1))r_{n}{}^{\rho -\varepsilon _{n}},
\end{equation*}

\noindent which implies for $n$ large so that $\rho /2\leq \rho -r_{n}\leq 3\rho /2$
and $0<r_{n}<1,$\bigskip\ 
\begin{equation*}
(1+o(1))r_{n}{}^{\rho /2}\leq \frac{U(a_{n})}{U(b_{n})}\leq
(1+o(1))r_{n}{}^{3\rho /2},
\end{equation*}

\noindent and this ensures that

\begin{equation*}
\frac{U(a_{n})}{U(b_{n})}\rightarrow 0=0^{\rho }\text{ as }n\rightarrow
+\infty .
\end{equation*}

\noindent \textbf{In a third and last step}, let $\rho >0$ and $c=+\infty$ Use the notation of
the second step where $r_{n}\rightarrow +\infty .$ We have 
\begin{equation*}
\frac{U(b_{n})}{U(a_{n})}=(1+o(1))r_{n}^{-\rho }\exp
(\int_{a_{n}}^{b_{n}}t^{-1}\ell (t)d\lambda (t)),
\end{equation*}

\noindent where, for large values of $n$,
\begin{equation*}
\left\vert \int_{a_{n}}^{b_{n}}t^{-1}\ell (t)d\lambda (t)\right\vert \leq
\varepsilon _{n}\log (a_{n}/b_{n})=r_{n}^{\varepsilon _{n}}.
\end{equation*}

\noindent This leads to, for large values of $n$, 
\begin{equation*}
(1+o(1))r_{n}^{-\rho +\varepsilon _{n}}\leq \frac{U(b_{n})}{U(a_{n})}\leq
(1+o(1))r_{n}^{-\rho +\varepsilon _{n}}.
\end{equation*}

\noindent Since $r_{n}\rightarrow +\infty$, $\varepsilon _{n}\rightarrow 0$
and $\rho>0$, it comes that
 
\begin{equation*}
\frac{U(b_{n})}{U(a_{n})}\rightarrow 0\Longrightarrow \frac{U(a_{n})}{U(b_{n})}\rightarrow +\infty =(+\infty )^{\rho }.
\end{equation*}

\noindent The proof is now complete.\\

\bigskip \noindent Now, we are able to prove Proposition \ref{evt.rvsv.Inverse}, following the lines of \cite{dehaan}, who quoted a personal communication of
W. Vervaart.\\

\bigskip \noindent \textbf{Proof of Proposition \ref{evt.rvsv.Inverse}}. We are going to give the proof only for
$U$ non-decreasing and $\rho>0$. The other case of $U$ non-increasing and $\rho<0$ is obtained by using the transform $1/U$. Let us assume the conditions
and the notations of the proposition. We are going to prove
\begin{equation}
\forall (x>0),\lim_{t\rightarrow \infty }\frac{V(tx)}{V(t)}=x^{\rho }.
\label{evt.rvsv.RVINS}
\end{equation}

\noindent Suppose (\ref{evt.rvsv.RVINS}) false. It is easy to see that $V(+\infty
)=+\infty $ and $V$ is also non-decreasing. Besides, an analogue to Formula (\ref{evt.rvsv.Uinv}) which we established for the generalized function of a non-increasing can be easily derived, in a very similar way, for the
generalized inverse of our non-drecreasing function, in the form%
\begin{equation}
U(V(y)-)\leq y\leq U(V(y)+), \text{ }y>0. \label{evt.rvsv.InvInv}
\end{equation}

\noindent Formula  (\ref{evt.rvsv.RVINS}) is false if and only for some $x_{0}>0,$
either 
\begin{equation*}
\lim_{t\rightarrow \infty }\inf \frac{V(tx_{0})}{V(t)}=c\neq x_{0}^{\rho }
\end{equation*}

\noindent or
\begin{equation*}
\lim_{t\rightarrow \infty }\sup \frac{V(tx_{0})}{V(t)}=c\neq x_{0}^{\rho }.
\end{equation*}

\noindent In both cases, there exists a sequence $(t_{n})_{n\geq 1}$ such that $%
t_{n}\rightarrow \infty $ and 
\begin{equation}
\lim_{t\rightarrow \infty }\frac{V(t_{n}x_{0})}{V(t_{n})}=c\neq x_{0}^{\rho
}.  \label{evt.rvsv.ratio7}
\end{equation}

\noindent Let us consider the three following cases.\\

\bigskip \noindent \textbf{Case 1}. $c$ is finite and positive.  We get, by using Inequality (\ref%
{evt.rvsv.InvInv}) , the non-decreasingness and the definitions of right
and left limits :
\begin{equation}
x_{0}=\frac{t_{n}x_{0}}{t_{n}}\leq \frac{U(V(t_{n}x_{0})+)}{U(V(t_{n})-)}%
\leq \frac{U(V(t_{n}x_{0})+1)}{U(V(t_{n})-1)}.  \label{evt.rvsv.ratio8}
\end{equation}

\noindent Let us take $a_{n}=V(t_{n}x_{0})+1$ and $b_{n}=V(t_{n})-1.$ Since $%
a_{n}/b_{n}\rightarrow c$ finite and positive by (\ref{evt.rvsv.ratio7}),
and ($a_{n},b_{n})\rightarrow (+\infty ,+\infty )$, we are in a position to apply Lemma \ref{evt.rvsv.Inverse} and let $n\rightarrow +\infty $ in (\ref{evt.rvsv.ratio8} )  to obtain
\begin{equation*}
x_{0}\leq c^{\rho }.
\end{equation*}

\noindent In a similar way, we also have
\begin{equation}
x_{0}=\frac{t_{n}x_{0}}{t_{n}}\geq \frac{U(V(t_{n}x_{0})-)}{U(V(t_{n})+)}%
\leq \frac{U(V(t_{n}x_{0})-1)}{U(V(t_{n})+1)}.  \label{evt.rvsv.ratio9}
\end{equation}

\noindent But the same method, we also get $x_{0}\geq c^{\rho }.$ The final conclusion
is $x_{0}=c^{\rho },$ which contradicts our supposition. Hence (\ref%
{evt.rvsv.RVINS}) is true for $c$ finite and positive.\\

\noindent \textbf{Case 2} : $c=0$. Use $a_{n}=V(t_{n}x_{0})+1$ and $b_{n}=V(t_{n})-1.$ Since $%
a_{n}/b_{n}\rightarrow 0$ and since $\rho >0$ in the present case, we let $%
n\rightarrow +\infty $ in (\ref{evt.rvsv.ratio8}) to get $x_{0}=0.$\\ 

\noindent \textbf{Case 3} : $c=+\infty .$Use $a_{n}=V(t_{n}x_{0})-1$ and $b_{n}=V(t_{n})+1.$
Since $a_{n}/b_{n}\rightarrow +\infty $ and since $\rho >0$ in the present
case, we let $n\rightarrow +\infty $ in (\ref{evt.rvsv.ratio9}) to get $%
x_{0}=+\infty .$\\

\noindent When we put all this together, we say that if (\ref{evt.rvsv.RVINS}) is
false, then (\ref{evt.rvsv.ratio7}) holds for some number $c$ and some
finite and positive $x_{0}$. The cases 2 and 3 above showed that $c$ is
necessarily positive and finite. And the case 1 showed the last $c$ cannot
be positive and finite. The conclusion is that (\ref{evt.rvsv.RVINS}) is
true.\\

\bigskip \noindent Before, we close the door, let us show a nice uniform convergence property of slowly varying functions.

\begin{lemma} \label{evt.rvsv.lemUnif} Let  $S(u)$ be a function $u\in (0,1)$ that is slowly varying at zero. We have 
the following uniform convergence in deterministic and random versions.\\

\bigskip \noindent (a) Let $A(h)$ and $B(h)$ two functions of $h\in (0,+\infty \lbrack $ such that for each $h\in (0,+\infty \lbrack $, we have $%
0<A(h)\leq B(h)<+\infty $ and  $(A(h),B(h))\rightarrow (0,0)$ as $h\longrightarrow 0$. Suppose that there exist two real numbers $0<C<D<+\infty$ such that 
\begin{equation}
C<\lim \inf_{h\rightarrow +\infty }A(h)/B(h),\text{ }\lim
\sup_{h\rightarrow +\infty }B(h)/A(h)<B.  \label{C1}
\end{equation}

\noindent Then, we have

\begin{equation*}
\lim_{h\rightarrow +\rightarrow }\sup_{A(h)\leq u,v\leq B(h)}\left\vert 
\frac{S(u)}{S(v)}-1\right\vert =0.
\end{equation*}

\bigskip \noindent (b)  Let $A(h)$ and $B(h)$ two families, indexed by $h\in
(0,+\infty \lbrack ,$ of real-valued applications defined on a probability
space $(\Omega ,\mathcal{A},\mathbb{P})$  such that for each $h\in
(0,+\infty \lbrack $, we have $0<A(h)\leq B(h)<+\infty .$ Suppose that there
exist two families $A^{\ast }(h)$ and $B^{\ast }(h)$, indexed by $h\in
(0,+\infty \lbrack ,$ of $\mathbf{measurable}$ real-valued applications
defined on $(\Omega ,\mathcal{A},\mathbb{P})$ such that
for each  $h\in (0,+\infty \lbrack ,$ $A^{\ast }(h)\leq A(h)\leq B(h)\leq
B^{\ast }(h)$, and such that

\begin{equation}
\lim \sup_{h\rightarrow +\infty }\lim_{\lambda \rightarrow +\infty }\inf 
\mathbb{P}(B^{\ast }(h)/A(h)>\lambda )=0.  \label{C2R}
\end{equation}

\noindent and

\begin{equation}
\lim \sup_{h\rightarrow +\infty }\lim_{\lambda \rightarrow +\infty }\inf 
\mathbb{P}(A^{\ast }(h)/B^{\ast }(h)<1/\lambda )=0.  \label{C3R}
\end{equation}

\noindent We say that the family $\{B^{\ast }(h),h\in h\in (0,+\infty \lbrack \}$ is
asymptotically bounded in probability against $+\infty $ and the family $\{B^{\ast }(h),h\in h\in (0,+\infty \lbrack \}$ is asymptotically bounded in probability against $0$ and accordingly, we say that the family $\{B(h),h\in h\in (0,+\infty\lbrack \}$ is asymptotically bounded in outer probability against $+\infty$ and the family $\{A(h),h\in h\in (0,+\infty \lbrack \}$ is asymptotically bounded in outer probability against $0$.\\

\noindent Then any $\eta >0,$ for any $\delta >0$, there exists a measurable subset $\Delta (\delta )$ of such
that 
\begin{equation*}
\left( \sup_{A(h)\leq u,v\leq B(h)}\left\vert \frac{S(u)}{S(v)}-1\right\vert
>\eta \right) \subset \Delta (\delta ),
\end{equation*}

\noindent with
\begin{equation*}
\mathbb{P}(\Delta (\delta ))\leq \delta .
\end{equation*}

\noindent Consequently, if the quantities
\begin{equation*}
\sup_{A(h)\leq u,v\leq B(h)}\left\vert \frac{S(u)}{S(v)}-1\right\vert >\eta 
\end{equation*}

\noindent are measurable for $h\in h\in (0,+\infty \lbrack $, we have that%
\begin{equation*}
\sup_{A(h)\leq u,v\leq B(h)}\left\vert \frac{S(u)}{S(v)}-1\right\vert
\rightarrow _{\mathbb{P}}as\text{ }h\rightarrow +\infty .
\end{equation*}\end{lemma}

\bigskip

\noindent \textbf{Proof}. \bigskip \noindent Let us use the Kamarata Representation Theorem \ref{evt.karamata.theo} of $S$ :
there exist a constant $c$ and functions $a(u)$ and $b(u)$ of $u\in ]0,1]$ satisfying
\begin{equation*}
(a(u),b(u))\rightarrow (0,0)\text{ as }u\rightarrow +\infty ,
\end{equation*}

\bigskip \noindent such that $S$ is written as%
\begin{equation}
S(u)=c(1+a(u))\exp (\int_{u}^{1}\frac{b(t)}{t}dt).
\label{evt.katamata.unif.lem}
\end{equation}

\bigskip \noindent \textbf{Proof of Point (a)}. Suppose that Condition (\ref{C1}) holds. Let $\varepsilon >0$
such that $\varepsilon <1.$ Then two functions $h_{1}(\varepsilon )$ and $h_{2}(\varepsilon )$ of $\varepsilon \in ]0,1[,$%
\begin{equation*}
h_{1}(\varepsilon )=C^{\varepsilon }\frac{1-\varepsilon }{1-\varepsilon }
\text{ and }h_{2}(\varepsilon )=D^{\varepsilon }\frac{1+\varepsilon }{%
1-\varepsilon }
\end{equation*}

\noindent both converge to $0$ as $\varepsilon \downarrow 0$. So for any $\eta >0,$
there exist $\varepsilon _{0},$ $0<\varepsilon <\varepsilon _{0}<1,$  such that 
\begin{equation}
1-\eta \leq h_{1}(\varepsilon ),h_{2}(\varepsilon )\leq 1+\eta .
\label{evt.borneinf}
\end{equation}

\noindent So, let $\eta >0$ and let $\varepsilon _{0}<1$ such that (\ref{evt.borneinf}) holds. Fix $\varepsilon ,$ $0<\varepsilon <\varepsilon _{0}.$ Now, by the assumptions on the functions $b$ and $p,$ there exists $t_{0}$ such that for 
$0\leq t\leq t_{0}$
\begin{equation*}
\max (\left\vert p(t)\right\vert ,\left\vert b(t)\right\vert )\leq
\varepsilon .
\end{equation*}

\noindent Since $B(h)\rightarrow 0$ as $h\rightarrow +\infty ,$ and since (\ref{C1})
holds, there is a value $h_{0}>0$ such that $h\geq h_{0}$ implies that $%
0\leq B(h)\leq t_{0}$ and 
\begin{equation*}
C\leq A(h)/B(h)\text{ and }B(h)/A(h)\leq D.
\end{equation*}

\noindent Then for  $h\geq h_{0}$ and $(u,v)\in \lbrack A(h),B(h)]^{2},$ we have 
\begin{eqnarray}
\frac{S(u)}{S(v)} &=&\frac{1+p(u)}{1+p(v)}\exp \left( \int_{u}^{v}\frac{b(t)%
}{t}dt\right)  \\
&\leq &\frac{1+\varepsilon }{1-\varepsilon }\exp \left( \sup_{0\leq t\leq
t_{0}}\left\vert b(t)\right\vert \int_{u}^{v}\frac{dt}{t}\right) 
\label{evt.upperbound} \\
&\leq &\frac{1+\varepsilon }{1-\varepsilon }\exp \left( \varepsilon \log
\left\{ \frac{\max (u,v)}{\min (u,v)}\right\} \right) .  \notag
\end{eqnarray}

\noindent Thus  for $h\geq h_{0}$ and $(u,v)\in \lbrack A(h),B(h)]^{2},$ we have 
\begin{equation*}
\frac{S(u)}{S(v)}\leq D^{\varepsilon }\frac{1+\varepsilon }{1-\varepsilon }.
\end{equation*}

\noindent By using lower bounds on place of upper bounds in (\ref{evt.upperbound}), we
also have $h\geq h_{0}$ and $(u,v)\in \lbrack A(h),B(h)]^{2},$%
\begin{equation*}
\frac{S(u)}{S(v)}\geq C^{\varepsilon }\frac{1-\varepsilon }{1-\varepsilon }.
\end{equation*}

\noindent By putting together the previous facts, we have, for $h\geq h_{0}$ and $%
(u,v)\in \lbrack A(h),B(h)]^{2}$ 
\begin{equation}
1-\eta \leq \frac{S(u)}{S(v)}\leq 1+\eta .  \label{evt.conPart}
\end{equation}

\noindent This implies that for for any $\eta >0,$ we have found $h_{0}$ such that for 
$h\geq h_{0}$ and $(u,v)\in \lbrack A(h),B(h)]^{2},$ we have
\begin{equation*}
\sup_{A(h)\leq u,v\leq B(h)}\left\vert \frac{S(u)}{S(v)}-1\right\vert \leq \eta .
\end{equation*}

\noindent Thus
\begin{equation*}
\lim_{h\rightarrow +\infty }\sup_{A\leq u,v\leq B}\left\vert \frac{S(u)}{S(v)%
}-1\right\vert =\lim_{h\rightarrow +\infty }\sup_{A\leq u,v\leq B}\left\vert 
\frac{S(u)}{S(v)}-1\right\vert =0.
\end{equation*}

\bigskip \noindent \textbf{Proof of Point (b)}. Suppose that the conditions of this point hold. For
any $\delta >0,$ there exist a real number $h_{1}>0$, and a number $\lambda >0
$ $(\lambda =\lambda (\delta )$ and a real number $h_{1}=h_{1}(\delta )$ both depend $on$ $%
\delta )$ such that

\begin{equation*}
\mathbb{P}(\lambda ^{-1}\leq A^{\ast }(h)/B^{\ast }(h),B^{\ast }(h)/A^{\ast
}(h)\leq \lambda )\geq 1-\delta /2.
\end{equation*}

\noindent If the latter property holds for a number $\lambda >0,$ it also holds for
any greater number. So, we may and do choose $\lambda >1.$
Put  $C=\lambda ^{-1}$ and $D=\lambda$.  From here, we follow partially use the proof of Point (a). Let $\eta >0$ and consider $\varepsilon _{0},$ $%
0<\varepsilon _{0}<1$ such that for any  $0<\varepsilon <\varepsilon _{0}<1,$
\ whe have 
\begin{equation}
1-\eta \leq h_{1}(\varepsilon ),h_{2}(\varepsilon )\leq 1+\eta .
\end{equation}

\noindent And let $t_{0}>0$ such that for any $0\leq t\leq t_{0}$%
\begin{equation*}
\max (\left\vert p(t)\right\vert ,\left\vert b(t)\right\vert )\leq
\varepsilon .
\end{equation*}

\noindent Since $B^{\ast }(h)\rightarrow _{P}0$ as $h\rightarrow 0,$ there exists a
value $h_{2}>0$ such that for any $h\geq h_{2},$ we have 
\begin{equation*}
\mathbb{P}(B^{\ast }(h)>t_{0})<\delta /2.
\end{equation*}

\noindent Denote $h_{0}=\max (h_{1,}h_{2}).$ The conditions under which (\ref{evt.conPart}) was proved are satisfied on the event  $(\lambda ^{-1}\leq A^{\ast
}(h)/B^{\ast }(h),B^{\ast }(h)/A^{\ast }(h)\leq \lambda )\cap (B^{\ast
}(h)<t_{0}),$ $h\geq h_{0}.$ Hence, we have on $(\lambda ^{-1}\leq A^{\ast
}(h)/B^{\ast }(h),B^{\ast }(h)/A^{\ast }(h)\leq \lambda )\cap (B^{\ast
}(h)<t_{0}),$ for $h\geq h_{0}$  
\begin{equation*}
\sup_{A^{\ast }(h)\leq u,v\leq B^{\ast }(h)}\left\vert \frac{S(u)}{S(v)}%
-1\right\vert \leq \eta .
\end{equation*}

\noindent Let us denote
\begin{equation*}
\Delta (\delta ,h)^{c}=(\lambda ^{-1}\leq A^{\ast }(h)/B^{\ast }(h),B^{\ast
}(h)/A^{\ast }(h)\leq \lambda )\cap (B^{\ast }(h)\leq t_{0}).
\end{equation*}

\noindent We have for $h\geq h_{0},$%
\begin{eqnarray*}
\mathbb{P}(\Delta (\delta ,h)) &\leq &\mathbb{P}((\lambda ^{-1}\leq A^{\ast
}(h)/B^{\ast }(h),B^{\ast }(h)/A^{\ast }(h)\leq \lambda )^{c})+\mathbb{P}%
(B^{\ast }(h)<t_{0}) \\
&\leq &\delta /2+\delta /2=\delta .
\end{eqnarray*}

\noindent We also have
\begin{equation*}
\Delta (\delta )^{c}\subset \left( \sup_{A^{\ast }(h)\leq u,v\leq B^{\ast
}(h)}\left\vert \frac{S(u)}{S(v)}-1\right\vert \leq \eta \right) 
\end{equation*}

\noindent This gives for $h\geq h_{0},$%
\begin{equation*}
\left( \sup_{A(h)\leq u,v\leq B(h)}\left\vert \frac{S(u)}{S(v)}-1\right\vert
>\eta \right) \subset \left( \sup_{A(h)\leq u,v\leq B^{\ast }(h)}\left\vert 
\frac{S(u)}{S(v)}-1\right\vert >\eta \right) \subset \Delta (\delta )
\end{equation*}

\bigskip \noindent Thus for any $\eta >0,$ for any $\delta >0,$ we have found $h_{0}>0$
such that for $h\geq h_{0},$  
\begin{equation*}
\left( \sup_{A(h)\leq u,v\leq B(h)}\left\vert \frac{S(u)}{S(v)}-1\right\vert
>\eta \right) \subset \Delta (\delta ),
\end{equation*}

\noindent with $\mathbb{P}(\Delta (\delta ))\leq \delta$.\\

\noindent The proof is complete $\blacksquare$

\newpage
\bigskip \noindent\textbf{$\pi$-variation}.\\

\noindent When dealing with $\pi $-varying functions, we adopt the former assumptions,
that is we work on functions $U:\mathbb{R}_{+}\longrightarrow \mathbb{R}_{+}$
that are measurable and Lebesgue locally integrable and not vanishing in
the neighborhood of $+\infty ,$ that is, $\forall (0\leq a\leq b<+\infty)$, $U\in L([a,b],\lambda ).$ But we add the conditions :\\ 

\noindent (PV1) $U$ is non-decreasing.\\

\noindent (PV2) For all $t>0,$ for any $x>0,x\neq 1,$ there exist $y_{0}(t,x)$ such
that 
\begin{equation*}
(t \geq y_{0}(t,x)) \Longrightarrow (\left\vert U(tx)-U(x)\right\vert >0).
\end{equation*}

\bigskip \noindent Condition (PV2) is automatically implied by the increasingness of $U$ in a
neighborhood of $+\infty$.\\ 

\bigskip \noindent We will see $\pi $-variation is strongly linked to slowly variation. that We
are able to define the $\pi $-variation.\\

\bigskip 

\begin{definition} \label{evt.defPivar}
A function $U:\mathbb{R}_{+}\longrightarrow \mathbb{R}_{+}$ satisfying (PV1)
and (PV2) is $\pi $-varying at infinity, denoted $U\in \Pi (+\infty )$ if
and only if 
\begin{equation*}
\forall (y>0),\forall (x>0,x\neq 1),\lim_{t\rightarrow +\infty} \frac{U(ty)-U(t)}{U(tx)-U(t)}=\frac{%
\log y}{\log x}.
\end{equation*}
\end{definition}

\bigskip 

\noindent \textbf{Remark}. We understand now why the condition (PV2) is imposed as a necessary
condition to develop the theory of $\pi $-variation.\\

\bigskip \noindent Before giving interesting properties of $\pi $-variation in connection of
slow variation, let us state the important representation of \text{de} Haan(1970).

\bigskip 

\subsection{Theorem of \text{de} Haan}

\noindent Before we proceed any further, let give this simple and useful identity.

\bigskip 

\begin{proposition} \label{evt.karamata.lem4} Let $U:\mathbb{R}_{+}\longrightarrow \mathbb{R}_{+}$ be
a measurable function Lebesgue locally integrable on [0,+$\infty \lbrack .$
Then for%
\begin{equation*}
g(x)=U(x)-\frac{1}{x}\int_{1}^{x}U(t)dt,
\end{equation*}

\noindent we have
\begin{equation*}
U(x)=g(x)+\int_{1}^{x}t^{-1}g(t)dt,a.e.
\end{equation*}
\end{proposition}

\bigskip \noindent \textbf{Proof}. Define $g(\circ )$ as in the statement. Then%
\begin{equation}
\int_{1}^{x}t^{-1}g(t)dt=\int_{1}^{x}t^{-1}U(t)dt-\int_{1}^{x}\left(
\int_{1}^{t}U(s)ds\right) t^{-2}dt.  \label{evt.karamata.ident1}
\end{equation}

\noindent By integration by parts, we have for all $x>0,$%
\begin{eqnarray*}
\int_{1}^{x}\left( \int_{1}^{t}U(s)ds\right) t^{-2}dt &=&\int_{1}^{x}\left(
\int_{1}^{t}U(s)ds\right) d(-t^{-1}) \\
&=&\left[ -t^{-1}\left( \int_{1}^{t}U(s)ds\right) \right] _{1}^{x}+%
\int_{1}^{x}t^{-1}U(t)dt \\
&=&-\frac{1}{x}\int_{1}^{x}U(t)dt+\int_{1}^{x}t^{-1}U(t)dt.
\end{eqnarray*}

\noindent By combining this (\ref{evt.karamata.ident1}), we get%
\begin{equation*}
\int_{1}^{x}t^{-1}g(t)dt=\frac{1}{x}\int_{1}^{x}U(t)dt,
\end{equation*}

\noindent which, by definition of $g(\circ )$, leads to%
\begin{equation*}
\int_{1}^{x}t^{-1}g(t)dt=U(x)-g(x).
\end{equation*}

\noindent This puts an end to the proof.\\

\bigskip 
 
\noindent Let us \textit{}give the Theorem of \text{de} Haan(1970) which includes the representation we
need.\\

\begin{theorem} \label{evt.dehaan.rep}
Let $U:\mathbb{R}_{+}\longrightarrow \mathbb{R}_{+}$ be a function that is
measurable and Lebesgue locally integrable and not vanishing in the
neighborhood of $+\infty ,$ that is, $\forall 0\leq a\leq b<+\infty ,U\in
L([a,b],\lambda )$ and such that (PV1) and (PV2) hold. Then the following
assertions are equivalent.\\

\noindent \textbf{(a)} U$\in \Pi (+\infty)$.\\

\noindent \textbf{(b)} The function 
\begin{equation*}
g(x)=U(x)-\frac{1}{x}\int_{1}^{x}U(t)dt,
\end{equation*}

\noindent is slowly varying at $+\infty$.\\

\noindent \textbf{(c)} 
\begin{equation*}
\lim_{x\rightarrow +\infty }\frac{x\int_{1}^{x}U(t)dt-2\int_{1}^{x}%
\int_{1}^{y}U(y)dydt}{x^{2}U(x)-x\int_{1}^{x}U(t)dt}=1/2.
\end{equation*}

\noindent There exist a slowly varying function $g(\circ)$ and a real constant $c$
such that, for $x>0$,
\begin{equation*}
U(x)=c+g(x)+\int_{1}^{x}t^{-1}g(t)dt.
\end{equation*}

\noindent \textbf{(e)} There exists a positive slowly varying $g(x)$ of $x>0$ such that for each 
$x>0,$%
\begin{equation*}
\lim_{t\rightarrow +\infty }\frac{U(tx)-U(t)}{g(t)}=\log x.
\end{equation*}
\end{theorem}

\bigskip 

\noindent \noindent \textbf{Proof}. \textbf{We begin to prove that (a)$\Longrightarrow (b)$}. Assume (a) holds. Fix $x>0$ and $x\neq 1$. Let us
show that the function in $t>0,$ 
\begin{equation*}
h(t)=U(t)-U(tx)
\end{equation*}

\noindent is slowly varying. Remark that by (PV1), $h$ is non-negative. Indeed, we
have
\begin{eqnarray*}
\frac{h(ty)}{h(t)} &=&\frac{U(ty)-U(txy)}{U(t)-U(tx)} \\
&=&\frac{U(t)-U(txy)}{U(t)-U(tx)}-\frac{U(t)-U(ty)}{U(t)-U(tx)} \\
&\rightarrow &\frac{\log xy}{\log x}-\frac{\log y}{\log x}=1.
\end{eqnarray*}

\noindent Next, by Lemma \ref{evt.karamata.lem2}, we have%
\begin{equation*}
th(t)/\int_{0}^{t}h(s)ds\rightarrow 1,
\end{equation*}

\noindent which implies by taking the inverses $1/\circ$, as t $\rightarrow +\infty$,
\begin{equation*}
\frac{\int_{0}^{t}U(s)dt-\int_{0}^{t}U(xs)ds}{t(U(tx)-U(t))}\rightarrow 1 \text{ as } t \rightarrow +\infty.
\end{equation*}

\noindent Check that we way use integrals on $[1,+\infty \lbrack $ instead of
integrals on $[0,+\infty \lbrack $ in this formula. Next use right change of
variables to have

\begin{equation*}
C(t)=\frac{\frac{1}{t}\int_{1}^{t}U(xs)dt-\frac{1}{tx}\int_{1}^{tx}U(s)ds}{%
U(t)-U(tx)}\rightarrow 1\text{ as }t\rightarrow +\infty .
\end{equation*}

\noindent Let us explain the coming computations. We want to prove that $g(\circ )$ is
slowly varying, that is for $y>0$ fixed, $g(ty)/g(t)\rightarrow 1$
which equivalent to
\begin{equation*}
\frac{g(ty)-g(t)}{g(t)}\rightarrow 0\text{ as }t\rightarrow +\infty .
\end{equation*}

\noindent All the coming details are intended to express the quantity 
\begin{equation*}
\frac{g(ty)-g(t)}{g(t)}
\end{equation*}

\noindent and to have it bounded by a quantity that tends to zero. To this purpose, remark that
\begin{equation*}
C(t)-1=C(t)-\frac{U(t)-U(tx)}{U(t)-U(tx)}\rightarrow 0.
\end{equation*}

\noindent But the left member becomes, by putting the terms together,%
\begin{eqnarray*}
&=&-\frac{\left( U(t)-\frac{1}{t}\int_{1}^{t}U(xs)dt\right) -\left( U(tx)-%
\frac{1}{tx}\int_{1}^{tx}U(s)ds\right) }{U(t)-U(tx)} \\
&=&\frac{g(tx)-g(t)}{U(t)-U(tx)}\rightarrow 0.
\end{eqnarray*}

\noindent We get for any $0<x<1,$%
\begin{equation}
\frac{g(tx)-g(t)}{U(t)-U(tx)}\rightarrow 0\text{ as }t\rightarrow \infty .
\label{evt.karamata.Pi1}
\end{equation}

\bigskip \noindent We are almost done if we are able to replace $U(t)-U(tx)$ by $g(t).$ So we
are going to compare $U(t)-U(tx)$ with $g(t)$. First, remark that, by a change of variable,
\begin{equation*}
\frac{1}{t}\int_{1}^{t}U(s)ds=\int_{1/t}^{1}U(ts)ds
\end{equation*}

\noindent and next,
\begin{eqnarray*}
\frac{g(t)}{U(t)-U(tx)} &=&\frac{U(t)-\int_{1/t}^{1}U(ts)ds}{U(t)-U(tx)} \\
&=&\frac{U(t)-\left( \int_{1/t}^{1/2}U(ts)ds+\int_{1/2}^{1}U(ts)ds\right) }{%
U(t)-U(tx)} \\
&=&\frac{\frac{1}{2}U(t)-\int_{1/t}^{1/2}U(ts)ds+\frac{1}{2}%
U(t)-\int_{1/2}^{1}U(ts)ds}{U(t)-U(tx)}.
\end{eqnarray*}

\noindent Since
\begin{equation*}
\frac{1}{2}U(t)=\int_{1/2}^{1}U(t)ds,
\end{equation*}

\noindent we arrive at 
\begin{equation*}
\frac{g(t)}{U(t)-U(tx)}=\frac{\frac{1}{2}U(t)-\int_{1/t}^{1/2}U(ts)ds}{%
U(t)-U(tx)}+\int_{1/2}^{1}\frac{U(t)-U(ts)}{U(t)-U(tx)}ds.
\end{equation*}

\noindent Next, using the non-decreasingness of $U$ leads to, for $t\geq 1,$ 
\begin{eqnarray*}
\frac{\frac{1}{2}U(t)-\int_{1/t}^{1/2}U(ts)ds}{U(t)-U(tx)} &\geq &\frac{%
\frac{1}{2}U(t)-\int_{1/t}^{1/2}U(t/2)ds}{U(t)-U(tx)} \\
&=&\frac{\frac{1}{2}U(t)-(\frac{1}{2}-1/t)U(t/2)}{U(t)-U(tx)} \\
&\geq &\frac{1}{2}\frac{U(t)-U(t/2)}{U(t)-U(tx)}.
\end{eqnarray*}

\noindent We arrive at
\begin{equation*}
\frac{g(t)}{U(t)-U(tx)}\geq \frac{1}{2}\frac{U(t)-U(t/2)}{U(t)-U(tx)}%
+\int_{1/2}^{1}\frac{U(t)-U(ts)}{U(t)-U(tx)}ds.
\end{equation*}

\noindent Remark that the function in the integral of the second term in the left
member is non-negative on $[0,1]$. Apply Point (a) and Fatou-Lebesgue
theorem (in the last integral) to get
\begin{eqnarray*}
\liminf_{t\rightarrow \infty }\frac{g(t)}{U(t)-U(tx)} &\geq &\frac{1}{_{2}}%
\left( \frac{-\log 2}{\log x}\right) +\int_{1/2}^{1}\left( \frac{\log s}{%
\log x}\right) ds \\
&=&\frac{1}{_{2}}\left( \frac{-\log 2}{\log x}\right) +\frac{1}{\log x}\left[
s\log s-s\right] _{1/2}^{1} \\
&=&\frac{1}{_{2}}\left( \frac{-\log 2}{\log x}\right) +\frac{-1-\left( -\log
2\right) /2+1/2}{\log x} \\
&=&1/\left\{ 2(-\log x)\right\} >0.
\end{eqnarray*}

\noindent This implies that

\begin{equation*}
\limsup_{t\rightarrow +\infty }\frac{U(t)-U(tx)}{g(t)}\leq 2(-\log x).
\end{equation*}

\noindent By combining this with (\ref{evt.karamata.Pi1}), we get get 
\begin{equation*}
\frac{\left\vert g(tx)-g(t)\right\vert }{g(t)}\leq \frac{\left\vert
g(tx)-g(t)\right\vert }{U(t)-U(tx)}\times \frac{U(t)-U(tx)}{g(t)}
\end{equation*}

\noindent and next, since $\log (1/x)$ is finite and nonzero, 
\begin{equation}
\lim_{t\rightarrow +\infty }\frac{\left\vert g(tx)-g(t)\right\vert }{g(t)}=0,
\label{evt.karamata.Pi2}
\end{equation}

\noindent for any $x>0,x\neq 1$. To finish the proof, remark that (\ref%
{evt.karamata.Pi2}) is obvious for $x=1$ and for $x>1,$ we have for $%
y=1/x\in ]0,1[$,

\begin{eqnarray*}
\lim_{t\rightarrow +\infty }\frac{g(tx)-g(t)}{g(t)} &=&\lim_{s\rightarrow
+\infty }\frac{g(s)-g(sy)}{g(sy)} \\
&=&\lim_{s\rightarrow +\infty }\frac{g(s)-g(sy)}{g(y)}\times \frac{g(y)}{%
g(sy)}=0,
\end{eqnarray*}

\noindent since by (\ref{evt.karamata.Pi2}), $(g(s)-g(sy))/g(y)\rightarrow 0$ and $%
g(y)/g(sy)\rightarrow 1$ as $s\rightarrow \infty .$    

\bigskip \noindent \textbf{Let us prove that $(b)\Longleftrightarrow (c)$}. It is easy to see that $g(\circ )$ $\in
RV(0,+\infty )$ $\Longleftrightarrow q(x)=xg(x)\in RV(1,+\infty )$ which, by
Lemma \ref{evt.karamata.lem3} is equivalent to
\begin{equation}
\frac{xq(x)}{\int_{1}^{x}q(t)dt}\rightarrow 2\text{ as }x\rightarrow \infty .
\label{evt.karamara.Pi3}
\end{equation}

\noindent But, by partial integration,
\begin{eqnarray*}
\int_{1}^{x}tU(t)dt &=&\int_{1}^{x}t\left( \int_{1}^{t}U(s)ds\right)
^{\prime }dt \\
&=&\left[ t\int_{1}^{t}U(s)ds\right] _{1}^{x}-\int_{1}^{x}%
\int_{1}^{t}U(s)dsdt. \\
&=&x\int_{1}^{x}U(s)ds-\int_{1}^{x}\int_{1}^{t}U(s)dsdt.
\end{eqnarray*}

\noindent Then

\begin{eqnarray*}
\int_{1}^{x}q(t)dt &=&\int_{1}^{x}tU(t)dt-\int_{1}^{x}\left(
\int_{1}^{t}U(s)ds\right) dt \\
&=&x\int_{1}^{x}U(s)ds-2\int_{1}^{x}\int_{1}^{t}U(s)dsdt
\end{eqnarray*}

\noindent and, by (\ref{evt.karamara.Pi3}),
\begin{equation*}
\frac{x\int_{1}^{x}U(s)ds-2\int_{1}^{x}\left( \int_{1}^{t}U(s)ds\right) dt}{%
x^{2}U(x)-x\int_{1}^{x}U(t)dt}\rightarrow 1/2.
\end{equation*}

\noindent Then (b) is equivalent to 
\begin{equation*}
\lim_{x\rightarrow +\infty }\frac{x\int_{1}^{x}U(s)ds-2\int_{1}^{x}\left(
\int_{1}^{t}U(s)ds\right) dt}{x^{2}U(x)-x\int_{1}^{x}U(t)dt}=1/2.
\end{equation*}

\noindent \textbf{Let us prove that $(c)\Longrightarrow (d)$}. Assume (c) holds. Then $(b)$ is true and $g$
is slowly varying. By Lemma \ref{evt.karamata.lem4}, the relation 
\begin{equation*}
g(x)=U(x)-\frac{1}{x}\int_{1}^{x}U(t)dt,
\end{equation*}

\noindent leads to 
\begin{equation*}
U(x)=g(x)+\int_{1}^{x}t^{-1}g(t)dt,a.e.
\end{equation*}

\noindent and Point $(d)$ is true.

\bigskip \noindent \textbf{Let us prove that $(d)\Longrightarrow (e)$}. Suppose $(d)$ holds. Then for any $t>0$ and $x>0
$
\begin{eqnarray*}
\frac{U(tx)-U(t)}{g(t)} &=&\left\{ \frac{g(tx)}{g(t)}-1\right\}
+\int_{x}^{tx}s^{-1}\left\{ \frac{g(s}{g(t)}\right\}  \\
&=&\int_{x}^{tx}s^{-1}ds+\left\{ \frac{g(tx)}{g(t)}-1\right\}
+\int_{x}^{tx}\left\{ \frac{g(s)}{g(t)}-1\right\} ds
\end{eqnarray*}

\noindent Since $g(\circ)$ is slowly varying, by Lemma \ref{evt.rvsv.lemUnif} above, implies that , 
\begin{equation*}
\sup_{s\in \lbrack x\min (1,t),x\max (1,t)]}\left\vert \frac{g(s)}{g(t)}%
-1\right\vert =\varepsilon (x,t)\rightarrow 0\text{ as }x\rightarrow \infty .
\end{equation*}

\noindent Then we have 
\begin{equation*}
\left\vert \frac{U(tx)-U(t)}{g(t)}-\log x\right\vert \leq \varepsilon
(x,t)(1+\left\vert \log x\right\vert )\rightarrow 0\text{ as }x\rightarrow
+\infty .
\end{equation*}

\bigskip 

\noindent \textbf{Let us prove that $(e)\Longrightarrow (a)$}. Suppose that $(e)$ holds. Then for $%
y>0,x>0,x\neq 1$,
\begin{equation*}
\frac{U(ty)-U(t)}{U(tx)-U(t)}=\frac{U(ty)-U(t)}{g(t)}\times \frac{g(t)}{%
U(tx)-U(t)}\rightarrow \frac{\log y}{\log x}\text{ as }x\rightarrow \infty .
\end{equation*}

\noindent And $(a)$ holds.\\

\newpage
\section[Characterizations of the extreme domain]{Characterizations of distribution functions in the extreme domain of attraction} \label{evt.charac}

\subsection{Characterization of the Frechet domain}

The Frechet domain is characterized by this theorem.\\

\begin{theorem} \label{evt.extdom.charac.frechet} Let $F$ be a distribution function.\\

\noindent (a) We have the following equivalence :\\

\noindent (a1) $F$ is in the domain of attraction of the Frechet
type distribution function of parameter $\alpha >0,$%
\begin{equation*}
\varphi _{\alpha }(x)=\exp (-x^{-\alpha })1_{(x\geq 0)},
\end{equation*}

\noindent that is

\begin{equation}
F\in D(\varphi _{\alpha })  \label{evt.extdomC1}
\end{equation}

\noindent if and only if \\

\noindent (a2) the survival function  $\overline{F}(\circ )=1-F(\circ )$ is $%
(-\alpha)$-regularly varying at $+\infty $, that is 
\begin{equation}
\forall (\gamma >0),\lim_{x\rightarrow +\infty }\frac{1-F(\gamma x)}{1-F(x)}%
=\gamma ^{-\alpha }.  \label{evt.extdomC2}
\end{equation}

\noindent (b) Besides, if $F$ is the domain of attraction of a Frechet distribution,
then the upper endpoint is infinite, that is 
\begin{equation}
uep(F)=+\infty.   \label{evt.extdomNormFR}
\end{equation}

\noindent and for 
\begin{equation*}
a_{n}=F^{-1}(1-1/n),n\geq 1.
\end{equation*}

\noindent we have for any $x\in \mathbb{R},$%
\begin{equation}
F^{n}(a_{n}x)\rightarrow \varphi _{\alpha }(x).  \label{evt.extdomConv1}
\end{equation}
\end{theorem}

\bigskip \noindent \textbf{Proof}. Suppose that $1-F(\circ )\in RV(-\alpha ,+\infty )$\ . In Proposition %
\ref{portal.extdomain}, we already showed that $F\in D(\varphi _{\alpha })$
if we have the representation

\begin{equation}
F^{-1}(1-s)=c(1+p(x))x^{-1/\alpha }\exp (\int_{s}^{1}s^{-1}\ell (s)d\lambda
(s)),  \label{evt.rvsv.repQF}
\end{equation}

\noindent with ($p(s),\ell (s))\rightarrow (0,0)$ as $x\rightarrow +\infty $ and $%
F^{n}(a_{n}x)\rightarrow \varphi _{\alpha }(x)$ for $a_{n}=F^{-1}(1-1/n),n%
\geq 1.$ But if $1-F(\circ )\in RV(-\alpha ,+\infty ),$ we have by
Proposition \ref{evt.rvsv.Inverse} that ($1-F)^{-1}(1/x)=F^{-1}(1-1/x)\in
RV(1/\alpha ,+\infty ).$ We may write the Karamata representation for $%
F^{-1}(1-1/x)$ which gives (\ref{evt.rvsv.repQF}) for $s=1/x$.\\

\noindent The new thing to prove is that $F\in D(\varphi _{\alpha })$ implies that $%
1-F(\circ )\in RV(-\alpha ,+\infty ).$ Now suppose that $F\in D(\varphi
_{\alpha }).$ We are going to put ourselves in a position to use Lemma \ref%
{evt.rvsv.VRWEAK}. Since $F\in D(\varphi _{\alpha }),$ there exist two sequences 
$(a_{n}>0)_{n\geq 0}$ and $(b_{n})_{n\geq 0}$ such that for any $x\in \mathbb{R}$,

\begin{equation*}
F^{n}(a_{n}x+b_{n})\rightarrow \varphi _{\alpha }(x)\text{ as }n\rightarrow
+\infty .
\end{equation*}

\noindent Let $s>1$. We have, as $n\rightarrow +\infty$,
 
\begin{equation}
F^{[ns]}(a_{[ns]}x+b_{[ns]})=\left\{ F^{n}(a_{[ns]}x+b_{[ns]})\right\}
^{[ns]/n}\rightarrow \varphi _{\alpha }(x), 
\label{evt.domFR00}
\end{equation}

\noindent and next

\begin{equation}
F^{n}(a_{[ns]}x+b_{[ns]})\rightarrow \varphi _{\alpha }(x)^{1/s}=\varphi
_{\alpha }((1/s)^{-1/\alpha }x)\text{ as }n\rightarrow +\infty .
\label{evt.domFR0}
\end{equation}

\noindent By \ Lemma \ref{evt.lem.1}, we have 

\begin{equation}
\frac{b_{[ns]}-b_{n}}{a_{n}}\rightarrow 0\text{ and }\frac{a_{[ns]}}{a_{n}}%
\rightarrow (1/s)^{-1/\alpha }=\rho >1\text{ as }n\rightarrow +\infty .
\label{evt.domFR1}
\end{equation}

\noindent Set $n(1)=[(s/(s-1)]$ and $n(i+1)=[n(i)s],i\geq 1$. We remark that%

\begin{equation*}
n(i)\geq n(i-1)s\geq ...\geq s^{i-1}n(1)\rightarrow +\infty \text{ as }%
i\rightarrow \infty 
\end{equation*}

\noindent and

\begin{equation*}
\frac{n(i+1)}{n(i)}=\frac{[n(i)s]}{n(i)}\rightarrow s\text{ as }i\rightarrow
+\infty .
\end{equation*}

\bigskip \noindent Replacing $n$ by $n(i)$ is Formula (\ref{evt.domFR1}) gives%
\begin{equation*}
\frac{b_{n(i+1)}-b_{n(i)}}{a_{n(i)}}\rightarrow 0\text{ and }\frac{a_{n(i+1)}%
}{a_{n(i)}}\rightarrow \rho >1\text{ as }i\rightarrow +\infty 
\end{equation*}

\noindent \noindent which implies that $a_{n(i)}\rightarrow +\infty $ as $i\rightarrow +\infty .$
We have to prove also that

\begin{equation}
\frac{b_{n(i)}}{a_{n(i)}}\rightarrow 0\text{ as }i\rightarrow +\infty .
\label{evt.domFR2}
\end{equation}

\bigskip \noindent To see that, put $c_{n}(i)=b_{n(i+1)}-b_{n(i)},i\geq 0$ with the convention
that $b_{n(0)}=0.$ We have 
\begin{equation*}
\frac{c_{n(i)}}{a_{n(i)}}\rightarrow 0\text{ as }i\rightarrow +\infty 
\end{equation*}

\bigskip \noindent  Let $\varepsilon >0$ such $r=(1/\rho )+\varepsilon <1.$ There exists $i_{0}$
such that for $i\geq i_{0},$

\begin{equation*}
\left\vert \frac{c_{n(i)}}{a_{n(i)}}\right\vert \leq \varepsilon \text{ and }%
\frac{a_{n(i)}}{a_{n(i+1)}}\leq r.
\end{equation*}

\bigskip \noindent We have for $i>i_{0}$

\bigskip 
\begin{eqnarray*}
\frac{b_{n(i+1)}}{a_{n(i)}} &=&\frac{1}{a_{n(i)}}\sum\limits_{k=1}^{i}c_{n}(k) \\
. &=&\frac{1}{a_{n(i)}}\sum\limits_{k=1}^{i_{0}}c_{n}(k)+\sum%
\limits_{i_{0}}^{i}\frac{c_{n}(k)}{a_{n(i)}}.
\end{eqnarray*}

\bigskip Since $\sum\limits_{k=1}^{i_{0}}c_{n}(k)$ is fixed, we have

\begin{equation*}
\frac{1}{a_{n(i)}}\sum\limits_{k=1}^{i_{0}}c_{n}(i)\rightarrow 0\text{ as }%
i\rightarrow +\infty .
\end{equation*}

\bigskip Next, we obtain

\begin{eqnarray*}
\sum\limits_{i_{0}}^{i}\frac{c_{n}(k)}{a_{n(i)}} &=&\sum\limits_{i_{0}}^{i}%
\frac{c_{n}(k)}{a_{n(k)}}\left\{ \frac{a_{n(k)}}{a_{n(i)}}\right\}  \\
&\leq &\varepsilon \sum\limits_{i_{0}}^{i}\left\{ \frac{a_{n(k)}}{a_{n(k+1)}%
}\right\} \left\{ \frac{a_{n(k+1)}}{a_{n(k+2)}}\right\} ...\left\{ \frac{%
a_{n(i-1)}}{a_{n(i)}}\right\}  \\
&\leq &\varepsilon \sum\limits_{i_{0}}^{i}r^{i-k}=\varepsilon
r^{i_{0}-k}\sum\limits_{p=0}^{i-i_{0}}r^{p}\leq \frac{\varepsilon
r^{i_{0}-k}}{1-r}.
\end{eqnarray*}

\bigskip \noindent So, for any $\varepsilon >0$, 
\begin{equation*}
\limsup_{i\rightarrow +\infty }\frac{b_{n(i+1)}}{a_{n(i)}}\leq \frac{\varepsilon r^{i_{0}-k}}{1-r},
\end{equation*}

\bigskip which , by letting $\varepsilon \downarrow 0$, leads to 
\begin{equation*}
\limsup_{i\rightarrow +\infty }\frac{b_{n(i+1)}}{a_{n(i+1)}}=\limsup_{i\rightarrow +\infty }\left\{ \frac{b_{n(i+1)}}{a_{n(i)}}\right\} \left\{ \frac{a_{n(i)}}{a_{n(i+1)}}\right\} =0.
\end{equation*}

\bigskip \noindent This gives (\ref{evt.domFR2}). This combined with the application of (\ref{evt.domFR00}) for $n=n(i-1)$, that is
\begin{equation*}
F^{n(i)}(a_{n(i)}x+b_{n(i)})\rightarrow \varphi _{\alpha }(x)\text{ as }%
i\rightarrow +\infty .
\end{equation*}

\bigskip \noindent ensures, via Lemma \ref{evt.lem.1}, that%
\begin{equation*}
F^{n}(a_{n(i)}x)\rightarrow \varphi _{\alpha }(x)\text{ as }i\rightarrow
+\infty .
\end{equation*}

\bigskip \noindent This implies that for any $x>0$%
\begin{equation*}
n(i)\log F(a_{n(i)}x)\rightarrow \log \varphi _{\alpha }((1/s)^{-1/\alpha
}x)=-x^{-\alpha }/s\text{ as }i\rightarrow +\infty .
\end{equation*}

\bigskip \noindent This is possible if and only if $F(a_{n(i)}x)\rightarrow 1$ as $i\rightarrow
+\infty $ and 
\begin{equation}
n(i)(1-F(a_{n(i)}x)\rightarrow x^{-\alpha }\text{ as }i\rightarrow +\infty .
\label{evt.extdomFR5}
\end{equation}

\bigskip \noindent By Lemma \ref{evt.rvsv.VRWEAK}, we may conclude that $1-F\in D(\varphi
_{\alpha })$.\\

\noindent It remains to prove (\ref{evt.extdomNormFR}) in Point (b). Suppose that that 
$uep(F)$ is finite. Consider $x>\inf (0,uep(F)).$ Formula (\ref%
{evt.extdomFR5}) would implies for $i$ large enough 
\begin{equation*}
n(i)(1-F(a_{n(i)}x)=0\rightarrow x^{-\alpha },
\end{equation*}

\bigskip \noindent which is absurd.\\

\subsection{Characterization of the Weibull domain}

\bigskip 

The Weibull domain is characterized by this theorem.

\begin{theorem} \label{evt.extdom.charac.weibull} Let $F$ be a probability distribution function.\\

\noindent (a) We have the following equivalence.

\noindent (a.1) $F$ is in the domain of attraction of a Weibull type of distribution function of parameter $\beta >0$,
\begin{equation*}
\psi _{\beta }(x)=\exp (-(-x)^{\beta })1_{(x\geq 0)}+1_{(x\geq 0)},
\end{equation*}

\noindent that is
\begin{equation}
F\in D(\psi _{\beta })  \label{evt.extdomWEC1}
\end{equation}

\noindent if and only if

\begin{equation}
uep(F)<+\infty 
\end{equation}

\noindent and 
\begin{equation}
x\hookrightarrow F^{\ast }(x)=F(uep(F)-\frac{1}{x})  \label{evt.extdomWEC2A}
\end{equation}

\noindent is $(-\beta )$-regularly varying, that is F$^{\ast }\in D(\varphi _{\alpha })$ :  
\begin{equation}
\forall (\gamma >0),\lim_{x\rightarrow +\infty }\frac{1-F^{\ast }(\gamma x)}{1-F^{\ast }(x)}=\gamma ^{-\beta }.  \label{evt.extdomWEC2B}
\end{equation}

\bigskip \noindent (b) Besides, if $F$ is the domain of attraction of a Weibull distribution of
parameter $\beta >0,$ for 
\begin{equation*}
a_{n}=uep(F)-F^{-1}(1-1/n),n\geq 1.
\end{equation*}

\noindent we have for any $x\in \mathbb{R}$,

\begin{equation}
F^{n}(a_{n}x+uep(F))\rightarrow \psi _{\beta}(x).
\end{equation}
\end{theorem}

\bigskip \noindent  \textbf{Proof}. As in the Frechet's case, we are going to focus the implication of (\ref{evt.extdomWEC2A}) and (\ref{evt.extdomWEC2B}) by (\ref{evt.extdomWEC1}). The other reverses implications are sufficiently handled in Chapter \ref{portal}. We begin similarly as in the case of Frechet's domain. Suppose that $F\in D(\psi_{\beta})$. We are going to put ourselves in a position to use again
Lemma \ref{evt.rvsv.VRWEAK}. There exist sequences $(a_{n}>0)_{n\geq 0}$ and $(b_{n})_{n\geq 0}$ such that for any $x\in \mathbb{R}$,
\begin{equation*}
F^{n}(a_{n}x+b_{n})\rightarrow \psi _{\beta }(x).
\end{equation*}

\noindent Let $s>1$. We obtain

\begin{equation}
F^{n}(a_{[ns]}x+b_{[ns]})\rightarrow \psi _{\beta }(x)^{1/s}=\psi _{\beta
}((1/s)^{1/\beta }x).
\end{equation}

\noindent By  Lemma \ref{evt.lem.1}, we have 
\begin{equation}
\frac{b_{[ns]}-b_{n}}{a_{n}}\rightarrow 0\text{ and }\frac{a_{[ns]}}{a_{n}}%
\rightarrow (1/s)^{1/\beta }=\rho <1.
\end{equation}

\noindent We are going to re-use the sequence $n(i)_{i\geq 1}$ introduced above to establish that 
\begin{equation*}
b_{n(i)}\rightarrow uep(F)\text{ and }\frac{a_{n(i+1)}}{a_{n(i)}}\rightarrow
\rho <1,
\end{equation*}

\noindent which implies that $a_{n(i)}\rightarrow 0$ as $i\rightarrow +\infty $. We
have to prove also that

\begin{equation}
\frac{b_{n(i+1)}-b_{n(i)}}{a_{n(i)}}\rightarrow 0\rightarrow uep(F)\text{
and }\frac{uep(F)-b_{n(i)}}{a_{n(i)}}\rightarrow 0.
\end{equation}

\noindent Let us prove that  $(b_{n(i)})_{i\geq 1}$ is a Cauchy sequence. Put again $%
c_{n}(i)=b_{n(i+1)}-b_{n(i)},i\geq 0$ with the convention that $b_{n(0)}=0.$
We have 

\begin{equation*}
\frac{c_{n(i)}}{a_{n(i)}}\rightarrow 0\text{ as }i\rightarrow +\infty .
\end{equation*}

\noindent From these limits, we may find, for any $\varepsilon >0$ such $r=\rho
+\varepsilon <1$, an integer $i_{0}$ such that for $i\geq i_{0}$,

\begin{equation*}
\left\vert \frac{c_{n(i)}}{a_{n(i)}}\right\vert \leq \varepsilon \text{ and }%
\frac{a_{n(i+1)}}{a_{n(i)}}\leq r.
\end{equation*}

\noindent This ensures that for any $i>i_{0}$%
\begin{equation*}
a_{n(i)}\leq ra_{n(i-1)}\leq r^{2}a_{n(i-2)}\leq ...\leq
r^{i-i_{0}}a_{n(i_{0})}
\end{equation*}

\noindent We have for $i_{0}<i<j$,

\begin{eqnarray*}
b_{n(j+1)}-b_{n(i)} &=&\sum\limits_{k=i}^{j}c_{n}(k) \\
&=&\sum\limits_{k=i}^{j}\left\{ \frac{c_{n}(k)}{a_{n(k)}}\right\} a_{n(k)}
\\
&=&\varepsilon \sum\limits_{k=i}^{j}a_{n(k)} \\
&\leq &\varepsilon r^{-i_{0}}a_{n(i_{0})}\sum\limits_{k=i}^{j}r^{k} \\
&\leq &\frac{r^{-i_{0}}a_{n(i_{0})}}{1-r}\varepsilon .
\end{eqnarray*}

\noindent Hence, we obtain

\begin{equation}
\limsup_{i\rightarrow +\infty ,j\rightarrow +\infty }b_{n(j+1)}-b_{n(i)}\leq \frac{r^{-i_{0}}a_{n(i_{0})}}{1-r}\varepsilon .
\label{evt.extdomSuiteBn1}
\end{equation}

\noindent \noindent By letting $\varepsilon \downarrow 0$, we get that $(b_{n(i)})_{i\geq 1}$ is
a Cauchy sequence and hence, there exist a real number $b$ such that 
\begin{equation*}
b_{n(i)}\rightarrow b.
\end{equation*}

\noindent Let us use again these formulas to see that

\begin{eqnarray*}
\frac{b_{n(j+1)}-b_{n(i)}}{a_{n(i)}} &=&\sum\limits_{k=i}^{j}\left\{ \frac{c_{n}(k)}{%
a_{n(k)}}\right\} \left\{ \frac{a_{n(k)}}{a_{n(i)}}\right\}  \\
&\leq &\varepsilon \sum\limits_{k=i}^{j}r^{k-i}=\varepsilon
r^{-i}\sum\limits_{k=i}^{j}r^{k}.
\end{eqnarray*}

\noindent When $j\rightarrow \infty $, we get

\begin{equation*}
\frac{b-b_{n(i)}}{a_{n(i)}}\leq \varepsilon
r^{-i}\sum\limits_{k=i}^{+\infty }r^{k}=\frac{\varepsilon }{1-r}
\end{equation*}

\noindent for any $\varepsilon \downarrow 0.$ This is enough to prove that 
\begin{equation*}
\frac{b-b_{n(i)}}{a_{n(i)}}\rightarrow 0.
\end{equation*}

\bigskip \noindent By Lemma \ref{evt.lem.1}, we have for any $x\in \mathbb{R}$  

\begin{equation}
F^{n(i)}(a_{n(i)}x+b)\rightarrow \psi _{\beta }(x),\text{as }i\rightarrow
+\infty .  \label{evt.extdomWE3}
\end{equation}

\noindent \bigskip For $x=0$, (\ref{evt.extdomWE3}) implies that%
\begin{equation*}
F(b)=1.
\end{equation*}

\bigskip \noindent Let $h>0$. Since $a_{n(i)}\rightarrow 0$, we have for large values of $i$, $b-h\leq b-a_{n(i)}$ and 
\begin{equation*}
F^{n(i)}(b-h)\leq F^{n(i)}(b-a_{n(i)})\rightarrow \psi _{\beta }(-1)=1/e%
\text{ as }i\rightarrow +\infty .
\end{equation*}

\noindent This implies for large values of $i$,

\begin{equation*}
F(b-h)\leq \left\{ \frac{1,01}{e}\right\} ^{1/n(i)}<1.
\end{equation*}

\noindent Hence 
\begin{equation*}
b=uep(F).
\end{equation*}

\noindent \bigskip Now we have for $x>0$,

\begin{equation*}
F^{n(i)}(b-a_{n(i)}x)\rightarrow \psi _{\beta }(-x)=\exp (x^{\beta}),
\end{equation*}

\noindent which implies that
\begin{equation*}
n(i)(1-F(b-a_{n(i)}x))\rightarrow x^{\beta },
\end{equation*}

\noindent and by change of variable $x\longrightarrow 1/x$,

\begin{equation*}
n(i)(1-F(b-\frac{a_{n(i)}}{x}))\rightarrow x^{-\beta }
\end{equation*}

\noindent which is exactly 
\begin{equation*}
n(i)(1-F^{\ast }(b-a_{n(i)}^{\ast }x))\rightarrow x^{-\beta},
\end{equation*}

\noindent where $a_{n(i)}^{\ast }=1/a_{n}(i)\rightarrow +\infty .$ Then by Lemma \ref%
{evt.rvsv.VRWEAK}, we conclude that $1-F^{\ast }\in RV(-\beta ).$

\bigskip

\subsection{Characterization of the Gumbel domain}

We give the main characterization of $D(\Lambda)$ here. Before we give the
characterization, let us highlight the following implication. Suppose that
some distribution function $F$ satisfies :

\begin{equation*}
\forall (0<s,t\neq 1<1),\lim_{u\downarrow 0}\frac{F^{-1}(1-su)-F^{-1}(1-u)}{%
F^{-1}(1-tu)-F^{-1}(1-u)}=\frac{\log s}{\log t},
\end{equation*}

\noindent that is $F^{-1}(1-\circ )$ is $\pi $-varying at $zero$. Put $%
U(x)=F^{-1}(1-1/x),x>1.$ Then $U$ is $\pi $-varying at $+\infty .$ We may
use the \text{de} Haan representation given in Theorem \ref{evt.dehaan.rep} : there
exists a slowly varying function $g$ and a constant $c$ such that%
\begin{equation*}
U(x)=d+g(x)+\int_{1}^{x}t^{-1}g(t)dt,x>1.
\end{equation*}

\noindent Put $s(u)=g(1/u),0<u<1.$ We get the representation 
\begin{equation}
F^{-1}(1-u)=d+s(u)+\int_{u}^{1}t^{-1}s(t)dt,0<u<1, \label{evt.rdg}
\end{equation}

\noindent where $s(\circ)$ is still slowly varying at $zero.$ We proved in Point (c)
of Proposition \ref{portal.extdomain} in Chapter \ref{portal}\ \ that $F\in
D(\Lambda )$ whenever  the latter representation (\ref{evt.rdg}) holds. The
next theorem says that the reverse implication is also true.

\bigskip 

\begin{proposition} \label{evt.extdom.charac.gumbel} Let $F$ be a  probability distribution function. The following propositions
are equivalent.\\

\noindent (a) $F\in D(\Lambda )$.\\

\noindent (b) For any $0<s,t\neq 1<1,$\\
\begin{equation}
\lim_{u\downarrow 0}\frac{F^{-1}(1-su)-F^{-1}(1-u)}{F^{-1}(1-tu)-F^{-1}(1-u)}=\frac{\log s}{\log t}. \label{evt.gumbel.PIVAR})
\end{equation}

\bigskip \noindent (c) There exist a slowly varying function $s$ and a constant $d$ such that
for $0<u<1,$ 
\begin{equation}
F^{-1}(1-u)=d+s(u)+\int_{u}^{1}t^{-1}s(t)dt.  \label{evt.dehaan}
\end{equation}
\end{proposition}

\bigskip \noindent \textbf{Proof}. Based on the arguments given in the introduction and Theorem \ref%
{evt.dehaan.rep}, we see that we only have to prove $(a)\Longrightarrow (b)$. So, suppose $(a)$ holds that is : there exist sequences $(c_{n}>0)_{n\geq 1}$ and $(d_{n})_{n\geq 1}$ such that for any $x\in \mathbb{R}$,

\begin{equation*}
F^{n}(c_{n}x+d_{n})\rightarrow \exp (-e^{-x}).
\end{equation*}

\bigskip \noindent We are going to apply Lemma \ref{evt.lem.pregumbel} since $H(x)=exp(-e^{-x}),x\in \mathbb{R},$ is strictly increasing on its support which is $\mathbb{R}$. Consider for $x>1$, 
\begin{equation*}
\left\{ 
\begin{tabular}{lll}
$u_{1}(n)=(1-1/n)^{n}$ & $\rightarrow $ & $u_{1}=e^{-1}.$ \\ 
$u_{2}(n)=(1-/(nz))^{n}$ & $\rightarrow $ & $u_{2}=e^{-(1/z)}$%
\end{tabular}
\right..
\end{equation*}

\bigskip \noindent Remark that for any $0<s<1,$ $(F^{n})^{-1}(s)=F^{-1}(s^{1/n})$ and $%
H^{-1}(e^{-(1/x)})=\log x$ and $H^{-1}(e^{-1})=0$ $.$ Applying Part (II) of
Lemma \ref{evt.lem.pregumbel} gives

\begin{equation*}
\frac{F^{-1}(1-1/(nx))-F^{-1}(1-1/n)}{c_{n}}\rightarrow
H^{-1}(e^{-(1/x)})-H^{-1}(e^{-1})\log z.
\end{equation*}

\noindent By doing the same with $y>1,$ we get a similar formula in $y$ and by dividing the to
two formulas leads to : For any $x>1,y>1,$%
\begin{equation}
\frac{F^{-1}(1-1/(nx))-F^{-1}(1-1/n)}{F^{-1}(1-1/(ny))-F^{-1}(1-1/n)}%
\rightarrow \frac{\log x}{\log y}.  \label{evt.gumbel.01}
\end{equation}

\noindent Put $U(x)=F^{-1}(1-1/x),x>1.$ We are going to prove that for any $x>1,$ 
\begin{equation}
\frac{U(n+1)-U(n)}{U(nx)-U(n)}\rightarrow 0,  \label{evt.gumbel.02}
\end{equation}

\noindent and use it to conclude. Let  $y>1$. We have $n+1\leq yn$ and for $n$ large
enough, and since $U$ is non-decreasing, 
\begin{equation*}
0\leq \frac{U(n+1)-U(n)}{U(nx)-U(n)}\leq \frac{U(ny)-U(n)}{U(nx)-U(n)},
\end{equation*}

\bigskip \noindent and then, by (\ref{evt.gumbel.01}),
 
\begin{equation*}
0\leq \limsup_{n\rightarrow +\infty} \frac{U(n+1)-U(n)}{U(nx)-U(n)}\leq \frac{\log y}{\log x}.
\end{equation*}

\noindent We get the result by letting $y\downarrow 1.$ Now fix $x>1.$ We are going to
use that for any $\varepsilon >0,$ we have $[z+1]/[z]<1+\varepsilon $ for $t$
large enough. Fix $\varepsilon >0$ such that $(1+\varepsilon )x>1.$ We have
for $s$ large enough,%
\begin{eqnarray*}
1-\frac{U([z]+1)-U([z])}{U([t]x)-U([z])} &=&\frac{U([z]x)-U([z]+1)}{%
U([z]x)-U([z])} \\
&\leq &\frac{U(zx)-U(z)}{U([z]x)-U([z])} \\
&\leq &\frac{U(([z]+1)x)-U([z])}{U([z]x)-U([z])} \\
&\leq &\frac{U([z](1+\varepsilon )x)-U([z])}{U([z]x)-U([z]).}.
\end{eqnarray*}

\bigskip \noindent Taking the limits of the extreme members as $t\rightarrow +\infty$ and applying (\ref{evt.gumbel.01})
and (\ref{evt.gumbel.02}), we get%
\begin{equation*}
1\leq \limsup \frac{U(zx)-U(z)}{U([z]x)-U([z])}\leq \frac{\log x+\log
(1+\varepsilon )}{\log x}.
\end{equation*}

\noindent By letting $\varepsilon \downarrow 1,$ we get for any $x>0,$%
\begin{equation*}
\frac{U(zx)-U(z)}{U([z]x)-U([z])}\rightarrow 1.
\end{equation*}

\noindent By doing the same with $y>1,$ we get a similar formula and by dividing
the to two formulas and by applying (\ref{evt.gumbel.01}), which uses limits
over the integers, we get : $x>1,y>1,$

\bigskip 
\begin{eqnarray*}
\frac{U(zx)-U(z)}{U([zy)-U([z])} &=&\frac{U(zx)-U(z)}{U([z]x)-U([z])}\times 
\frac{U([z]y)-U([z])}{U([z]y)-U([z])}\times \frac{U([z]x)-U([z])}{
U([z]y)-U([z])} \\
&\rightarrow &1\times 1\times \frac{\log x}{\log y}.
\end{eqnarray*}

\noindent Finally setting $s=1/x$ and $t=1/y$ and $u=1/z$, we arrive at (\ref{evt.gumbel.PIVAR}), which was our target.

\newpage

\subsection{Quantile representations} \label{evt.quantile.rep}

We are going to particularize the representations we obtained above to
probability distributions functions $F$. Next for the needs of statistical
estimation on the extreme domain of attraction, we will also give
representations for the functions $G(x)=F(e^{x}),x\in \mathbb{R}$. 

\subsubsection{Direct representation of the quantile of a distribution in the extreme domain of attraction}.

\bigskip \noindent Let us work case by case.\\

\noindent \textbf{Case $F\in D(\varphi_{\alpha}),\alpha >0$}.\\

\noindent By Theorem \ref{evt.extdom.charac.frechet}, we have that $1-F\in RV(\alpha
,+\infty ).$ Then $U=1/(1-F)\in RV(\alpha ,+\infty )$ with $U(+\infty
)=U(+\infty ).$ By Proposition\ \ref{evt.rvsv.Inverse}, $U^{-1}\in
RV(1/\alpha ,+\infty ).$ But $U^{-1}(x)=F^{-1}(1-1/x),x>1$. By Theorem \ref{evt.karamata.theo}, we have the representation%
\begin{equation*}
F^{-1}(1-1/x)=c(1+a_{1}(x))x^{1/\alpha }\exp (\int_{1}^{x}t^{-1}\ell
_{1}(t)dt),x>1
\end{equation*}

\noindent where $(a_{1}(x),\ell _{1}(x))\rightarrow (0,0)$ as $x\rightarrow +\infty .$
Put $s=1/x,0<s<1,$ $a(s)=a_{1}(1/s),\ell (s)=-\ell _{1}(1/s)$ to get%
\begin{equation*}
F^{-1}(1-s)=c(1+a(s))s^{-1/\alpha }\exp (\int_{s}^{1}t^{-1}\ell
(t)dt),x>1,0<s<1,
\end{equation*}

\noindent where $(a(s),\ell (s))$ $\rightarrow $ $(0,0)$ as $s\rightarrow 0.$

\bigskip \noindent \textbf{Case $F\in D(\psi _{\beta }),\beta>0$}.

\noindent  By Theorem \ref{evt.extdom.charac.frechet}, $uep(F)<+\infty $ and $F^{\ast
}\in RV(-\beta )$ where $F^{\ast }(x)=F(uep(F)-1/x),x\in \mathbb{R}.$ By
using the above case, We have the representation 
\begin{equation*}
\left( F^{\ast }\right) ^{-1}(1-s)=c_{1}(1+a_{1}(s))s^{-1/\beta }\exp
(\int_{s}^{1}t^{-1}\ell _{1}(t)dt),x>1,0<s<1,
\end{equation*}

\noindent where $(a_{1}(s),\ell _{1}(s))$ $\rightarrow $ $(0,0)$ as $s\rightarrow 0.$
But it is easy to see that $\left( F^{\ast }\right)
^{-1}(1-s)=1/(uep(F)-F^{-1}(1-s)),0<s<1,$ with gives%
\begin{equation*}
uep(F)-F^{-1}(1-s))=c(1+a(s))s^{1/\beta }\exp (\int_{s}^{1}t^{-1}\ell
(t)dt),x>1,0<s<1,
\end{equation*}

\noindent where $c=1/c_{1},(a(s),\ell (s))=$ $(a_{1}(s)/(1+a_{1}(s)),-\ell
(s))\rightarrow $ $(0,0)$ as $s\rightarrow 0$.\\

\bigskip \noindent \textbf{Case $F\in D(\Lambda)$}.\\

\noindent we already saw in Lemma \ref{evt.extdom.charac.gumbel} that if $F \in D(\Lambda)$, there exist a slowly varying function $s$ and a constant $d$ such that
for $0<u<1,$ 
\begin{equation}
F^{-1}(1-u)=d-s(u)+\int_{u}^{1}t^{-1}s(t)dt.  \label{evt.dehaan}
\end{equation}

\bigskip \noindent We conclude by this :

\begin{proposition} \label{evt.sec.rvsv.repF}
We have the following characterizations for the three extremal domains.

\bigskip \noindent (a) $F\in D(H_{\gamma })$, $\gamma >0,$ if and only if there exist a
constant $c$ and functions $a(u)$ and $\ell (u)$ of $u\in ]0,1]$ satisfying
\begin{equation*}
(a(u),\ell (u))\rightarrow (0,0)\text{ as }u\rightarrow +\infty ,
\end{equation*}

\bigskip \noindent such that $F^{-1}$ admits the following representation of Karamata

\begin{equation*}
F^{-1}(1-u)=c(1+a(u))u^{-\gamma }\exp (\int_{u}^{1}\frac{\ell (u)}{u}du). \label{evt.rdf}
\end{equation*}

\bigskip \noindent (b) $F\in D(H_{\gamma })$, $\gamma <0,$ if and only if $uep(F)<+\infty $ and
there exist a constant $c$ and functions $a(u)$ and $\ell (u)$ of $u\in ]0,1]
$ satisfying
\begin{equation*}
(a(u),\ell (u))\rightarrow (0,0)\text{ as }u\rightarrow +\infty ,
\end{equation*}

\bigskip \noindent such that $F^{-1}$ admit the following representation of Karamata%
\begin{equation*}
uep(F)-F^{-1}(1-u)=c(1+a(u))u^{-\gamma }\exp (\int_{u}^{1}\frac{\ell (u)}{u}
du). \label{evt.rdw}
\end{equation*}

\bigskip \noindent (c) $F\in D(H_{0})$ if and only if there exist a constant $d$ and a slowly
varying function $s(u)$ such that 
\begin{equation*}
F^{-1}(1-u)=d+s(u)+\int_{u}^{1}\frac{s(u)}{u}du,0<u<1, \label{evt.rdg}
\end{equation*}

\bigskip \noindent and there exist a constant $c$ and functions $a(u)$ and $\ell (u)$ of $%
u\rightarrow $ $u\in ]0,1]$ satisfying%
\begin{equation*}
(a(u),\ell (u))\rightarrow (0,0)\text{ as }u\rightarrow +\infty ,
\end{equation*}

\bigskip \noindent such that $s$ admits the representation%
\begin{equation*}
s(u)=c(1+a(u))u^{-\gamma }\exp (\int_{u}^{1}\frac{\ell (u)}{u}du).
\end{equation*}
\end{proposition}

\bigskip \noindent As to the distribution fuction $G$, defined earlier by $G(y)=F(e^y)$, $y\in \mathbb{R}$ with $G^{-1}(1-s)=\log F^{-1}(1-s)$, $0<s<1$, its representations come from the combination of the just given ones and the points of Lemma ZZZ which ensure that 
$$
F \in D(G_0) \Rightarrow G\in D(D(G_0))
$$

\noindent and
$$
(F \in D(G_{\gamma}), \text{}{\gamma}<0) \Rightarrow (F \in D(G_{\gamma}), \text{}{\gamma}<0).
$$

\noindent This gives :

\begin{proposition} \label{evt.sec.rvsv.repG}
We have the following characterizations for the three extremal domains.\\

\bigskip \noindent (a) $F\in D(H_{\gamma})$, $\gamma >0,$ if and only if there exist a
constant $c$ and functions $p(u)$ and $b(u)$ of $u\rightarrow $ $u\in ]0,1]$
satisfying
\begin{equation*}
(p(u),b(u))\rightarrow (0,0)\text{ as }u\rightarrow +\infty ,
\end{equation*}

\bigskip \noindent such that $G^{-1}$ admit the following representation of Karamata

\begin{equation*}
G^{-1}(1-u)=c+\log (1+a(u)-\gamma \log u+\int_{u}^{1}\frac{b(u)}{u}du.
\end{equation*}

\bigskip \noindent (b) $F\in D(H_{\gamma })$, $\gamma <0,$ if and only if $uep(F)<+\infty $ and
there exist a constant $c$ and functions $p(u)$ and $b(u)$ of $u\in ]0,1]$
satisfying
\begin{equation*}
(p(u),b(u))\rightarrow (0,0)\text{ as }u\rightarrow +\infty ,
\end{equation*}

\bigskip \noindent such that $G^{-1}$ admit the following representation of Karamata

\begin{equation*}
uep(G)-G^{-1}(1-u)=c(1+a(u))\text{ }u^{-\gamma }\exp (\int_{u}^{1}\frac{b(u)}{u}du).
\end{equation*}

\bigskip \noindent (c) $F\in D(H_{0})$ if and only if there exist a constant $d$ and a slowly
varying function $s(u)$ such that 
\begin{equation*}
G^{-1}(1-u)=d+s(u)+\int_{u}^{1}\frac{s(u)}{u}du,0<u<1,
\end{equation*}

\bigskip \noindent and there exist a constant $c$ and functions $p(u)$ and $b(u)$ of $%
u\rightarrow $ $u\in ]0,1]$ satisfying

\begin{equation*}
(p(u),b(u))\rightarrow (0,0)\text{ as }u\rightarrow +\infty ,
\end{equation*}

\bigskip \noindent such that $s$ admits the representation%
\begin{equation*}
s(u)=c(1+p(u))\text{ }u^{-\gamma }\exp (\int_{u}^{1}\frac{b(u)}{u}du).
\end{equation*}
\end{proposition}

\section{EVT Criteria Based of the Derivatives of the Distribution Function} \label{evt.vonmises}

\section{Other Technical computations} \label{evt.techniq}

\subsection{Gaussian Extremes}

\noindent This subsection is devoted to a number of expansions results for the distribution function $\Phi(x)$, $x>0$ and the quantile function 
$\Phi^{-1}(s)$, $0<s<1$ of a standard normal random variable. These expressions are useful to describe the extremes of samples from Gaussian variables. The results exposed here are used for example in Subsection \ref{portal.ss.examples} of Chapter \ref{portal}.

\subsubsection{Gaussian tails} Let 
\begin{equation*}
\phi (x)=\frac{1}{\sqrt{2\pi }}\int_{-\infty
}^{x}e^{-t^{2}/2}dt=\int_{-\infty }^{x}\phi _{d}(t)dt
\end{equation*}

\noindent where
\begin{equation*}
\phi _{d}(t)=\frac{1}{\sqrt{2\pi }}e^{-t^{2}/2},-\infty <t<+\infty .
\end{equation*}

\noindent Put $C=1/\sqrt{2\pi }.$ We have for $x>0,$%
\begin{eqnarray*}
C^{-1}(1-\phi (x)) &=&\int_{x}^{+\infty }t^{-1}\text{ }d(-e^{-t^{2}/2}) \\
&=&\left[ -t^{-1}e^{-t^{2}/2}\right] _{t=x}^{t=+\infty }-\int_{x}^{+\infty
}t^{-2}\text{ }e^{-t^{2}/2}dt \\
&=&x^{-1}e^{-x^{2}/2}-\int_{x}^{+\infty }t^{-2}\text{ }te^{-t^{2}/2}dt\leq
x^{-1}.
\end{eqnarray*}%

\noindent But, also
\begin{eqnarray*}
\int_{x}^{+\infty }t^{-2}\text{ }te^{-t^{2}/2}dt &=&\int_{x}^{+\infty }t^{-2}%
\text{ }d(-e^{-t^{2}/2}) \\
&=&\left[ -t^{-2}e^{-t^{2}/2}\right] _{t=x}^{t=+\infty }-\int_{x}^{+\infty
}t^{-3}e^{-t^{2}/2}dt \\
&=&x^{-2}e^{-x^{2}/2}-\int_{x}^{+\infty }t^{-3}e^{-t^{2}/2}dt.
\end{eqnarray*}

\noindent This, combined with the first formulae, gives%
\begin{eqnarray*}
C^{-1}(1-\phi (x)) &=&x^{-1}-\int_{x}^{+\infty }t^{-2}\text{ }te^{-t^{2}/2}dt
\\
&=&x^{-1}e^{-x^{2}/2}-x^{-2}e^{-x^{2}/2}+\int_{x}^{+\infty
}t^{-3}e^{-t^{2}/2}dt \\
&\geq &e^{-x^{2}/2}(x^{-1}-x^{-2}).
\end{eqnarray*}

\noindent We arrive at, for $x>0,$%
\begin{equation}
C\left\{ \frac{1}{x}-\frac{1}{x^{2}}\right\} e^{-x^{2}/2}\leq 1-\phi (x)\leq 
\frac{Ce^{-x^{2}/2}}{x},  \label{evt.techniq.boundPhi}
\end{equation}%

\noindent and next,
\begin{equation*}
1-\phi (x)=\frac{Ce^{-x^{2}/2}}{x}(1+O(1/x)),
\end{equation*}

\noindent as $x\rightarrow \infty .$

\subsubsection{Tail quantile} Let 
\begin{equation*}
s=1-\phi (x).
\end{equation*}

\noindent We get as $s\rightarrow 0,$ which implies $x\rightarrow \infty ,$%
\begin{equation*}
s=\frac{Ce^{-x^{2}/2}}{x}(1+O(1/x))
\end{equation*}

\noindent and
\begin{equation*}
\log s=\log C-x^{2}/2-\log x+\log (1+O(1/x))
\end{equation*}

\noindent that is
\begin{equation*}
\log s=\log C-x^{2}/2-\log x+O(1/x).
\end{equation*}

\noindent We get
\begin{equation}
2\log (1/s)/x^{2}=1+(2\log C)/x^{2}+(2\log x)/x^{2}+O(1/x^{3}).
\label{form1}
\end{equation}%

\noindent Hence
\begin{equation*}
2\log (1/s)/x^{2}=1+O(x^{-2}\log x)=1+o(1),
\end{equation*}

\noindent meaning, as $s\rightarrow 0,$ 
\begin{equation*}
x=\phi ^{-1}(1-s)\sim (2\log 1/s)^{1/2}.
\end{equation*}

\noindent By going back to \ref{form1}, we get
\begin{eqnarray*}
x &=&(2\log (1/s))^{1/2}\{1+(2\log C)/x^{2}+(2\log
x)/x^{2}+O(1/x^{3})\}^{-1/2} \\
&=&(2\log (1/s))^{1/2}\{1-\log C)/x^{2}-\log x/x^{2}+O(x^{-2}\log x)\} \\
&=&(2\log (1/s))^{1/2}\left\{ 1-\frac{\log C+\log x+O(x^{-2}\log x)}{2\log
(1/s)}\left( 1+O(x^{-2}\log x)\right) \right\}  \\
&=&(2\log (1/s))^{1/2}\left\{ 1-\frac{2\log C+2\log x+O(x^{-2}\log x)}{4\log
(1/s)}\left( 1+O(x^{-2}\log x)\right) \right\} .
\end{eqnarray*}

\noindent By denoting $\varepsilon (s)=O((\log 1/s)^{-1}\log \log 1/s)$, we have%
\begin{eqnarray*}
&=&(2\log (1/s))^{1/2}\left\{ 1-\frac{\log 2\pi +\log 2+(\log \log
1/s)+\varepsilon (s)}{4\log (1/s)}(1+\varepsilon (s))\right\}  \\
&=&(2\log (1/s))^{1/2}\left\{ 1-\frac{\log 4\pi +\log \log (1/s)+\varepsilon
(s)}{4\log (1/s)}(1+\varepsilon (s))\right\} ,
\end{eqnarray*}

\noindent where $\varepsilon (s)=O((\log 1/s)^{-1/2}).$ We also have%
\begin{eqnarray*}
x &=&(2\log (1/s))^{1/2}\left\{ 1-\frac{\log 4\pi +\log \log
(1/s)+\varepsilon (s)}{4\log (1/s)}\right.  \\
&&-\left. \frac{\log 4\pi +\log \log (1/s)+\varepsilon (s)}{4\log (1/s)}%
\varepsilon (s)\right\} 
\end{eqnarray*}

\noindent This leads to
\begin{eqnarray*}
x &=&(2\log (1/s))^{1/2}-\frac{\log 4\pi +\log \log (1/s)}{2(2\log
(1/s))^{1/2}} \\
&&+O((\log 1/s)^{-3/2}(\log \log (1/s))+O((\log \log (1/s)^{2}(\log
1/s)^{-1})) \\
&=&\left\{ (2\log (1/s))^{1/2}-\frac{\log 2\pi +\log \log (1/s)}{2(2\log
(1/s))^{1/2}}+O((\log \log (1/s)^{2}(\log 1/s)^{-1}))\right\} .
\end{eqnarray*}

\noindent We conclude by%
\begin{equation}
\phi ^{-1}(1-s)=\left\{ (2\log (1/s))^{1/2}-\frac{\log 4\pi +\log \log (1/s)%
}{2(2\log (1/s))^{1/2}}+O((\log \log (1/s)^{2}(\log 1/s)^{-1/2}))\right\} .
\label{evt.techniq.quantile}
\end{equation}

\bigskip 

\subsubsection{Derivative at +$\infty$} Set%
\begin{equation*}
r(u)=Q(1-s)=\phi ^{-1}(1-s),0<s<1.
\end{equation*}

\noindent We have
\begin{eqnarray*}
\left( Q(1-s)\right) ^{\prime } &=&-\phi _{d}(Q(1-s))^{-1} \\
&=&-\frac{1}{\sqrt{2\pi }}\exp (\frac{1}{2}Q(1-s)^{2}). \\
&=&-\frac{1}{\sqrt{2\pi }}\exp \left( -\frac{1}{2}\left( 2\log (1/s))\left\{
1-\frac{\log 4\pi +\log \log (1/s)+\varepsilon (s)}{4\log (1/s)}%
(1+\varepsilon (s))\right\} ^{2}\right) \right) .
\end{eqnarray*}

\noindent Now, use the properties of $\varepsilon (s)$, when $s$ is near zero, to see
that
\begin{eqnarray*}
\left( Q(1-s)\right) ^{\prime } &=&-\frac{1}{\sqrt{2\pi }}\exp \left( -\frac{%
1}{2}\left( 2\log (1/s))-\log 4\pi -\log \log (1/s)+o(1)\right) \right)  \\
&=&-\sqrt{2}s^{-1}(\log 1/s)^{1/2}(1+o(1)) \\
&=&-s^{-1}(2\log 1/s)^{1/2}(1+o(1)).
\end{eqnarray*}

\noindent We obtained, as $s\downarrow 0,$%
\begin{equation}
\left( Q(1-s)\right) ^{\prime }=-s^{-1}(2\log 1/s)^{1/2}(1+o(1)).
\label{evt.techniq.derivative}
\end{equation}

\subsubsection{Normalizing and centering coefficients} Denote $\varepsilon _{n}=(\log n)^{-1}(\log \log n)^{2}$. We get

\begin{equation}
\phi ^{-1}(1-1/n)=\left\{ (2\log n)^{1/2}-\frac{\log 4\pi +\log \log n}{%
2(2\log n)^{1/2}}+O(\varepsilon _{n})\right\}   \label{inverseN}
\end{equation}

\noindent and

\begin{equation*}
\phi ^{-1}(1-1/ne)=\left\{ (2+2\log n)^{1/2}-\frac{\log 4\pi +\log (2+\log n)%
}{2(2+\log n)^{1/2}}+O(\varepsilon _{n})\right\} ,
\end{equation*}

\noindent that is
\begin{eqnarray*}
a_{n} &=&\phi ^{-1}(1-1/ne)-\phi ^{-1}(1-1/n) \\
&=&(2+2\log n)^{1/2}-(2\log n)^{1/2} \\
&&+\frac{\log 4\pi }{2}\left\{ \frac{1}{(2\log n)^{1/2}}-\frac{1}{(2+2\log
n)^{1/2}}+O(\varepsilon _{n})\right\}  \\
&&+\frac{1}{2}\left\{ \frac{\log \log n}{(2\log n)^{1/2}}-\frac{\log (1+\log
n)}{(2+2\log n)^{1/2}})\right\}  \\
&=&A_{n}+B_{n}+C_{n}.
\end{eqnarray*}

\noindent We surely have%
\begin{eqnarray*}
B_{n}/\alpha _{n} &=&\frac{\log 4\pi }{2}\left\{ \frac{(2\log n)^{1/2}}{%
(2\log n)^{1/2}}-\frac{(2\log n)^{1/2}}{(2+2\log n)^{1/2}}+O((\log
n)^{-1/2}(\log \log n)^{2})\right\}  \\
&\rightarrow &0.
\end{eqnarray*}%

\noindent Next
\begin{eqnarray*}
C_{n} &=&\frac{1}{2}\left\{ \frac{\log \log n}{(2\log n)^{1/2}}-\frac{\log
(2+\log n)}{(2+2\log n)^{1/2}})\right\}  \\
&=&\frac{1}{2}\left\{ \frac{\log \log n}{(2\log n)^{1/2}}-\frac{\log (\log
n(1+1/\log n))}{(2\log n)^{1/2}(1+1/\log n)^{1/2}})\right\}  \\
&=&\frac{1}{2}\left\{ \frac{\log \log n}{(2\log n)^{1/2}}-\frac{\log \log
n+\log (1+1/\log n)}{(2\log n)^{1/2}}(1+1/\log n)^{-1/2}\right\}  \\
&=&\frac{1}{2}\frac{\log \log n}{(2\log n)^{1/2}}-\frac{\log \log n+\log
(1+1/\log n)}{(2\log n)^{1/2}}(1-\frac{1}{2\log n}+O(\log n)^{-2})) \\
&=&\frac{1}{2}\left( \frac{\log \log n}{(2\log n)^{1/2}}-\frac{\log \log
n+\log (1+1/\log n)}{(2\log n)^{1/2}}\right)  \\
&&-\frac{1}{2}\left( \frac{\log \log n+\log (1+1/\log n)}{(2\log n)^{1/2}}%
\right) \left( -\frac{1}{2\log n}+O(\log n)^{-2})\right) .
\end{eqnarray*}%

\noindent Then
\begin{eqnarray*}
C_{n}/\alpha _{n} &=&-\frac{1}{2}\log (1+1/\log n) \\
&&-\frac{1}{2}\left( \log \log n+\log (1+1/\log n)\right) \left( -\frac{1}{%
2\log n}+O(\log n)^{-2})\right)  \\
&\rightarrow &0.
\end{eqnarray*}%

\noindent Finally
\begin{eqnarray*}
A_{n}/\alpha _{n} &=&(2\log n)(1+1/\log n)^{1/2}-(2\log n) \\
&=&(2\log n)(1+\frac{1}{2\log n}+O(\left( \log n\right) ^{-2})-(2\log n) \\
&=&1+O((\log n)^{-1})\rightarrow 1.
\end{eqnarray*}

\noindent Formally, we got%
\begin{equation}
A_{n}/\alpha _{n}\rightarrow 1\text{ as }n\rightarrow +\infty .
\label{evt.techniq.A}
\end{equation}

\noindent Finally, from (\ref{inverseN}), 
\begin{equation*}
\frac{\beta _{n}-b_{n}}{\alpha _{n}}=O(\varepsilon _{n}/\alpha _{n})=O((\log
n)^{-1/2}(\log \log n)^{2}).
\end{equation*}

\noindent Hence
\begin{equation}
\frac{\beta _{n}-b_{n}}{\alpha _{n}}\rightarrow 0\text{ as }n\rightarrow
+\infty.  \label{evt.techniq.B}
\end{equation}

\subsection{L'Hospital Type Rules}

We expose a result of \text{de} Haan (1970), that is useful to prove the Karamata representation theorem \ref{evt.karamata.theo}.\\

\begin{lemma} \label{evt.technical.hospital} Let $f,g:\mathbb{R}_{+}\longrightarrow \mathbb{R}%
_{+}$ be two measurable function Lebesgue locally integrable, that is, $%
\forall 0\leq a\leq b<+\infty ,f,g\in L([a,b],\lambda )$ where $\lambda $ is
the Lebesgue measure.\\

\noindent (A) Suppose that%
\begin{equation*}
\int_{0}^{+\infty }f(x)d\lambda (x)=\lim_{x\rightarrow +\infty
}\int_{0}^{x}f(x)d\lambda (x)=+\infty \text{ and }\int_{0}^{+\infty
}g(x)d\lambda (x)=+\infty 
\end{equation*}

\noindent and
\begin{equation*}
\lim_{x\rightarrow +\infty }f(x)/g(x)=c\in \lbrack 0,+\infty ].
\end{equation*}

\noindent Then
\begin{equation*}
\lim_{x\rightarrow +\infty }\left( \int_{0}^{x}g(t)d\lambda (t)\right)
\left( \int_{0}^{x}g(t)d\lambda (t)\right) ^{-1}=c.
\end{equation*}

\noindent (B) Suppose that 
\begin{equation*}
\int_{0}^{+\infty }f(x)d\lambda (x)<+\infty \text{ and }\int_{0}^{+\infty
}g(x)d\lambda (x)<+\infty 
\end{equation*}

\noindent and
\begin{equation*}
\lim_{x\rightarrow +\infty }f(x)/g(x)=c\in \lbrack 0,+\infty ].
\end{equation*}

\noindent Then
\begin{equation*}
\lim_{x\rightarrow +\infty }\left( \int_{x}^{+\infty }g(t)d\lambda
(t)\right) \left( \int_{x}^{+\infty }g(t)d\lambda (t)\right) ^{-1}=c.
\end{equation*}
\end{lemma}

\bigskip \noindent \textbf{Proof of (A)}. (\text{de} Haan, 1970). Suppose that%
\begin{equation*}
\lim_{x\rightarrow +\infty }f(x)/g(x)=+\infty .
\end{equation*}

\noindent Then for any $A>0,$ there exists $x_{0}>0$ such that for any $t\geq x_{0},$%
\begin{equation*}
f(t)/g(t)\geq A\Longrightarrow f(t)\geq 2Ag(t).
\end{equation*}

\noindent For $x_{0}$ fixed, $\int_{x_{0}}^{x}g(t)d\lambda (t)\rightarrow \infty $ as $%
x\rightarrow \infty $ and $\int_{0}^{x_{0}}g(t)d\lambda (t)$ is finite by
locally integrability of $g.$ So for some $x_{1}>x_{0},$ we have for $x\geq
x_{1},$
\begin{equation*}
\int_{0}^{x_{0}}g(t)d\lambda (t)/\int_{x_{0}}^{x}g(t)d\lambda (t)\leq 1
\end{equation*}
 \noindent and next 
\begin{eqnarray*}
\frac{\int_{0}^{x}f(t)d\lambda (t)}{\int_{0}^{x}g(t)d\lambda (t)} &=&\frac{%
\int_{0}^{x_{0}}f(t)d\lambda (t)+\int_{x_{0}}^{x}f(t)d\lambda (t)}{%
\int_{0}^{x_{0}}g(t)d\lambda (t)+\int_{x_{0}}^{x}g(t)d\lambda (t)} \\
&\geq &\frac{\int_{x_{0}}^{x}f(t)d\lambda (t)}{\int_{0}^{x_{0}}g(t)d\lambda
(t)+\int_{x_{0}}^{x}g(t)d\lambda (t)} \\
&\geq &\frac{2A\int_{x_{0}}^{x}g(t)d\lambda (t)}{\int_{0}^{x_{0}}g(t)d%
\lambda (t)+\int_{x_{0}}^{x}g(t)d\lambda (t)} \\
&=&2A(1+\int_{0}^{x_{0}}g(t)d\lambda (t)/\int_{x_{0}}^{x}g(t)d\lambda
(t))^{-1} \\
&\geq &A.
\end{eqnarray*}

\noindent This leads to 
\begin{equation*}
\frac{\int_{0}^{x}f(t)d\lambda (t)}{\int_{0}^{x}g(t)d\lambda (t)}\rightarrow
+\infty \text{ as }x\rightarrow +\infty .
\end{equation*}

\noindent Suppose now that%
\begin{equation*}
\lim_{x\rightarrow +\infty }f(x)/g(x)=c\in \lbrack 0,+\infty \lbrack .
\end{equation*}

\noindent Then for any $\varepsilon >0,$ there exists $x_{0}>0$ such that for any $%
t\geq x_{0},$%
\begin{equation*}
c-\varepsilon \leq f(t)/g(t)\leq c+\varepsilon \Longrightarrow
(c-\varepsilon )g(t)\leq f(t)\leq (c+\varepsilon )g(t).
\end{equation*}

\noindent Then
\begin{eqnarray*}
(c-\varepsilon )&=&\frac{(c-\varepsilon )\int_{x_{0}}^{x}g(t)d\lambda (t)}{%
\int_{x_{0}}^{x}g(t)d\lambda (t)}\leq \frac{(c-\varepsilon
)\int_{x_{0}}^{x}g(t)d\lambda (t)}{\int_{x_{0}}^{x}g(t)d\lambda (t)}\\
&\leq&\frac{\int_{x_{0}}^{x}f(t)d\lambda (t)}{\int_{x_{0}}^{x}g(t)d\lambda (t)}%
\leq \frac{(c+\varepsilon )\int_{x_{0}}^{x}g(t)d\lambda (t)}{%
\int_{x_{0}}^{x}g(t)d\lambda (t)}=(c+\varepsilon )
\end{eqnarray*}

\noindent For $x_{0}$ fixed, $\int_{x_{0}}^{x}g(t)d\lambda (t)\rightarrow \infty $ and 
$\int_{x_{0}}^{x}g(t)d\lambda (t)\rightarrow \infty $ as $x\rightarrow
\infty .$ So for some $x_{1}>x_{0},$ we have for $x\geq x_{1},$%
\begin{equation*}
0\leq \int_{0}^{x_{0}}f(t)d\lambda (t)/\int_{x_{0}}^{x}f(t)d\lambda (t)\leq
\varepsilon 
\end{equation*}

\noindent and 
\begin{equation*}
\int_{0}^{x_{0}}g(t)d\lambda (t)/\int_{x_{0}}^{x}g(t)d\lambda (t)\leq
\varepsilon .
\end{equation*}

\noindent Then for $x\geq x_{1},$

\begin{eqnarray*}
\frac{\int_{0}^{x}g(t)d\lambda (t)}{\int_{0}^{x}g(t)d\lambda (t)} &=&\frac{%
\int_{0}^{x_{0}}f(t)d\lambda (t)+\int_{x_{0}}^{x}f(t)d\lambda (t)}{%
\int_{0}^{x_{0}}g(t)d\lambda (t)+\int_{x_{0}}^{x}g(t)d\lambda (t)} \\
&=&\frac{\int_{x_{0}}^{x}f(t)d\lambda (t)}{\int_{x_{0}}^{x}g(t)d\lambda (t)}%
\times \frac{(1+\int_{0}^{x_{0}}f(t)d\lambda
(t)/\int_{x_{0}}^{x}f(t)d\lambda (t))}{(1+\int_{0}^{x_{0}}g(t)d\lambda
(t)/\int_{x_{0}}^{x}g(t)d\lambda (t))}
\end{eqnarray*}

\noindent This leads to $x\geq x_{1}$%
\begin{equation*}
(c-\varepsilon )(1+\varepsilon )^{-1}\leq \frac{\int_{0}^{x}f(t)d\lambda (t)%
}{\int_{0}^{x}g(t)d\lambda (t)}\leq (c+\varepsilon )(1+\varepsilon ).
\end{equation*}

\noindent Hence for any $\varepsilon >0,$%
\begin{equation*}
(c-\varepsilon )(1+\varepsilon )^{-1}\leq \liminf_{x\rightarrow \infty }%
\frac{\int_{0}^{x}f(t)d\lambda (t)}{\int_{0}^{x}g(t)d\lambda (t)}\leq \lim
\sup_{x\rightarrow +\infty }\frac{\int_{0}^{x}f(t)d\lambda (t)}{%
\int_{0}^{x}g(t)d\lambda (t)}\leq (c+\varepsilon )(1+\varepsilon ).
\end{equation*}

\noindent By letting $\varepsilon \downarrow 0$, we conclude that%
\begin{equation*}
\lim_{x\rightarrow +\infty }\frac{\int_{0}^{x}f(t)d\lambda (t)}{%
\int_{0}^{x}g(t)d\lambda (t)}=c.
\end{equation*}

\noindent \textbf{Proof of (B)}. We use very similar but easier ways. Suppose that%
\begin{equation*}
\lim_{x\rightarrow +\infty }f(x)/g(x)=+\infty .
\end{equation*}

\noindent Then for any $A>0,$ there exists $x_{0}>0$ such that for any $t\geq x_{0},$%
\begin{equation*}
f(t)/g(t)\geq A\Longrightarrow f(t)\geq Ag(t).
\end{equation*}

\noindent For  $x\geq x_{0},$ $\int_{x}^{+\infty }f(t)d\lambda (t)$ and $%
\int_{x}^{+\infty }g(t)d\lambda (t)$ are finite and we may writeand next 
\begin{equation*}
\frac{\int_{x}^{+\infty }f(t)d\lambda (t)}{\int_{x}^{+\infty }g(t)d\lambda
(t)}\geq \frac{\int_{x}^{+\infty }Ag(t)d\lambda (t)}{\int_{x}^{+\infty
}g(t)d\lambda (t)}\geq A.
\end{equation*}

\noindent This leads to 
\begin{equation*}
\frac{\int_{x}^{+\infty }f(t)d\lambda (t)}{\int_{x}^{+\infty }g(t)d\lambda
(t)}\rightarrow +\infty \text{ as }x\rightarrow +\infty .
\end{equation*}

\bigskip  \noindent Suppose now that
\begin{equation*}
\lim_{x\rightarrow +\infty }f(x)/g(x)=c\in \lbrack 0,+\infty \lbrack .
\end{equation*}

\bigskip
\noindent Then for any $\varepsilon >0,$ there exists $x_{0}>0$ such that for any $%
t\geq x_{0},$

\begin{equation*}
c-\varepsilon \leq f(t)/g(t)\leq c+\varepsilon \Longrightarrow
(c-\varepsilon )g(t)\leq f(t)\leq (c+\varepsilon )g(t).
\end{equation*}

\bigskip \noindent Then $x\geq x_{0}$, we have
\begin{eqnarray*}
(c-\varepsilon )&=&(c-\varepsilon )\frac{\int_{x}^{+\infty }g(t)d\lambda (t)}{%
\int_{x}^{+\infty }g(t)d\lambda (t)}\leq \frac{\int_{x}^{+\infty
}f(t)d\lambda (t)}{\int_{x}^{+\infty }g(t)d\lambda (t)}\\
&\leq& \frac{(c+\varepsilon )\int_{x}^{+\infty }g(t)d\lambda (t)}{\int_{x}^{+\infty
}g(t)d\lambda (t)}=(c+\varepsilon )
\end{eqnarray*}

\noindent This leads to%
\begin{equation*}
\frac{\int_{x}^{+\infty }f(t)d\lambda (t)}{\int_{x}^{+\infty }g(t)d\lambda
(t)}\rightarrow c\text{ as }x\rightarrow \infty .
\end{equation*}

\chapter[Advanced Characterizations]{Advanced Characterizations of the Univariate Extreme Value Domain} \label{evtp}

\section{Introduction}  
\noindent In this part, we suppose that the reader already knows the results of chapter \ref{evt} particularly, of Section \ref{evt.charac} where the first and main characterization are given, alongside with the representation of Karamata for a distribution $F\in D(H_\gamma),\gamma \neq 0$ and the de Haan's one
for $F\in D(H_0)$ through the quantile transformation $F^{-1}(1-u)$, $u\in(0,1)$.\\

\noindent This part, as it is a free restitution  of the reading of Sections $2.5-2.8$ in \cite{dehaan}, will be based on a definition of the extreme domain using limit on $t\rightarrow\infty$ on $\mathbb{R}$, rather than on $n\rightarrow +\infty$ on $\mathbb{N}$.\\

\noindent We say that the probability distribution function $F$ is in the domain of attraction of $H$, in the sense of extreme value theory, if and only if :\\

\noindent \textbf{$(D_1)$} There exist two functions $a(t)$ and $b(t)$ of $t\in \mathbb{R}_{+}$ such that $a(t)\geq 0$ for all  $t\geq 0$ and 
\begin{eqnarray}
\forall x\in C(H), F^t(a(t)x+b(t))\rightarrow H(x), as  t\rightarrow\infty. \label{evtp.01}
\end{eqnarray}

\bigskip \noindent Throughout this text, we will apply, without mentioning it, the discretization of limits for $t\rightarrow t_0$ in the following way :\\

\noindent Let $t_0$ be an adherent point of set $D\in \mathbb{R}$. A numerical function $f(t)$ of $t \in D$  converges to $a \in \overline{\mathbb{R}}$ as $t\rightarrow t_0$ if and only if for any sequence $(\ t_n)_{n\ge 0} \subset D$ such that $t_n\rightarrow t_0$, as $n\rightarrow +\infty$, we have $f(t_n)\rightarrow a$ as $n\rightarrow +\infty$.\\

\noindent From this, let us quickly show that (\ref{evtp.01}) is equivalent to the classical definition.\\

\noindent \textbf{$(D_2)$} There exist two sequence of real numbers $(a_n)_{n\ge 0}$ and $(b_n)_{n\ge 0}$ such that

\begin{eqnarray}
\forall~ x\in C(H),~ F^n (a_n x + b_n x )\rightarrow H(x),~ as~ t\rightarrow +\infty {\label{evtp.02}}.
\end{eqnarray}

\bigskip \noindent First, that (\ref{evtp.01}) implies (\ref{evtp.02}) is a simple consequence of the discretization principle of continuous limits.\\

\noindent Secondly, if (\ref{evtp.02}) holds, we have for any $t\geq 1$, $a_t = a_{[t]}>0$ and $b_t=b_{[t]}$ and for any $x\in C(H)$
\begin{equation*}
F^t(a_t x + b_t)= \big(F^{[t]}(a_{[t]} x + b_{[t]})\big)^{\frac{t}{[t]}}\rightarrow H(x)~as~t\rightarrow +\infty \ (as~ [t]\rightarrow +\infty).
\end{equation*}

\bigskip \noindent Based on this, we will work with Definition $(D_1)$. By using the discretization principle, we also have that $\alpha(t)>0$, $a(t)>0$, $b(t)$ and $\beta(t)$ are functions of $t>0$ such that, for any $x\in C(H)$, as $t\rightarrow +\infty$,
\begin{equation}\label{evtp.03} 
F^t(a_t x + b_t)\rightarrow H(x)~ as~t\rightarrow +\infty
\end{equation}

\noindent and, as $t\rightarrow +\infty$,

\begin{equation}
\frac{\alpha(t)}{a(t)}\rightarrow A>0 \text{ and } \frac{\beta(t)-b(t)}{a(t)}\rightarrow B\in \mathbb{R}\label{evtp.04}.
\end{equation}

\noindent Then, for any $ x\in C(H_{A,B})$, where $H_{A,B}(x)=H(Ax+B)$ for all $x \in \mathbb{R}$,
$$
F^t(a_{t} x + b_{t})\rightarrow H(Ax + B)~ as~ t\rightarrow+\infty.
$$

\bigskip \noindent As well if for any $x\in C(H_1)$, 
$$
F^t(a_t x + b_t)\rightarrow H_1 (x) ~as~ t \rightarrow +\infty
$$
\noindent and for any $x\in C(H_2),$

\begin{eqnarray*}
F^{t}(\alpha (t)x+ \beta (t)) \rightarrow H_2(x)\quad   as\quad t\rightarrow +\infty.
\end{eqnarray*}

\bigskip \noindent Then, (\ref{evtp.04}) holds and
\begin{eqnarray*}
\forall x\in \mathbb{R}, H_2(x)=H_1(Ax+B).
\end{eqnarray*}

\bigskip \noindent \textbf{Remark R1}. This merely means that Lemma (\ref{evt.lem.1}) remains valid in the frame of definition (D2) above.\\

\noindent Now, it remains to introduce another version of expressing (\ref{evtp.01}. By taking the logarithm of the quantities involved in that limit, we have :\\

\noindent (a) $\forall x\in \big[lep(H), uep(H)\big]\cap C(H)=:C^*(H)$,

$$
t\log F(a(t)x+b(t))\rightarrow \log H(x).
$$

\bigskip \noindent (b) for any $x < lep(H)$ and $x> uep(F)$, $x$ is a continuity point of $H$ ($H$ is constant on $]-\infty, lep(H)[$ and equal to zero, $H=1$ on  $]uep(H), +\infty[$, whenever these intervals are not empty). We also have on these intervals 
\begin{eqnarray*}
t\log F(a(t)x+b(t))\rightarrow \log H(x).
\end{eqnarray*}

\bigskip \noindent In summary, for any $x\in C(H),$
\begin{eqnarray*}
t\log F(a(t)x+b(t)) \rightarrow \log H(x).
\end{eqnarray*}

\bigskip \noindent For $x\in C^*(H)$, since $\log H(x)\in (-\infty, 0)$ and $t \rightarrow +\infty$, it follows that $\log F(a(t)x+b(t)) \rightarrow 0,$ that is 
$$
F(a(t)x+b(t)) \rightarrow 1,
$$ 

\noindent and then 

\begin{eqnarray*}
\log F(a(t)x+b(t))=\log \big(1+(F(a(t)x+b(t))+1)\big)\sim - \big(1-F(a(t)x+b(t))\big)
\end{eqnarray*}

\bigskip \noindent Finally,
\begin{eqnarray*}
t\log F(a(t)x+b(t))\sim - \big(1-F(a(t)x+b(t))\big)\rightarrow t\log H(x).
\end{eqnarray*}

\bigskip \noindent Then, for $x\in C^*(H)$,
\begin{eqnarray}
t(1-F(a(t)x+b(t)) \rightarrow -\log H(x). \label{evtp.05}
\end{eqnarray}

\bigskip \noindent Suppose (\ref{evtp.05}) hold in turn. Then for $x\in C^*(H)$,
$$
t\big(\log F(a(t)x+b(t))\big) (1+o(11))= (+\log H(x))+o(1),
$$

\noindent which implies

$$
t\log F(a(t)x+b(t))=\frac{1}{1+o(1)}(+\log H(x))+o(1),
$$

\noindent which implies
$$
\log F^t(a(t)x+b(t))=+\log H(x)+o(1),
$$

\bigskip \noindent since $\log H(x)\in (-\infty, 0)$. Hence, for $x\in C^*(H),$
\begin{eqnarray*}
 F^t\big(a(t)+b(t)\big)=H(x)\exp (o(1))\rightarrow H(x).
\end{eqnarray*}

\bigskip \noindent Finally, we have the equivalence, as $t\rightarrow +\infty,~\forall x\in C^*(H)$ : 
\begin{eqnarray*}
F^t\big(a(t)+b(t)\big)\rightarrow H(x)
\iff
\forall x\in C^*(H),~~ t\big(1-F(a(t)x+b(t)\big)\rightarrow H(x).
\end{eqnarray*}

\bigskip  \noindent In the special case where  $H(x)=H_0(x)=\exp(-e^{-x})$, we get\\

\begin{lemma} \label{evtp.lem.00}
\noindent $F\in D(H_0)$ and only if and only if : there exist two functions $a(t)>0$ and $b(t)$ of $t>0$ such that :
 
\begin{equation}
\forall x\in \mathbb{R}, ~t(1-F(a(t)x+b(t)) \rightarrow e^{-x},~ as~ t\rightarrow+\infty. \label{evtp.05a}
\end{equation}
\end{lemma}

\bigskip \noindent Now we may begin our round up of the chapter.\\
\newpage

\section{$\Gamma$-variation}

\bigskip 

\begin{definition} We introduce the following definitions.\\

\noindent (1) Let $F$ be a distribution function, $\alpha (t)>0$ and $\beta (t)$ two
function of $t<uep(F).$ We denote%
\begin{equation}
\Gamma (F,x,\alpha (t),\beta (t))=\frac{1-F(\alpha (t)x+\beta (t))}{1-F(t)}, \ t<uep(F).
\label{evtp.06}
\end{equation}

\bigskip \noindent (2) A distribution function $F$ is said to be of $\Gamma $-variation if and
only if there exist two $\alpha (t)>0$ and $\beta (t)$ of $t<uep(F)$ such
that
\begin{equation}
(\forall x\in \mathbb{R)},\text{ }\Gamma (F,x,\alpha (t),\beta
(t))\rightarrow \exp (-x)\text{ as }t\rightarrow uep(F)^{-}.  \label{evtp.06a}
\end{equation}
\end{definition}

\bigskip \noindent We have the following facts

\bigskip 

\begin{lemma} \label{evtp.lem.01} Let $F$ be distribution function such that (\ref{evtp.06}) holds. We have \\

\noindent (A) There exist $t_{0}(x)<uep(F)$  such that \\   
\begin{equation}
t_{0}(x)<t<uep(F)\Longrightarrow \alpha (t)x+\beta (t)<uep(F). \label{evtp.06c}
\end{equation}

$$
\alpha (t)x+\beta (t) \rightarrow uep(F),  \ as \ t\rightarrow uep(F)^{-}.
$$

\noindent Precisely, we have :\\

\noindent (a) If $uep(F)=+\infty$ the for any $x\in \mathbb{R},$ 
\begin{equation*}
\alpha (t)x+\beta (t)\rightarrow +\infty \text{ as }t\rightarrow uep(F)^{-}.
\end{equation*}

\noindent (b) If $uep(F)<\infty$,  then for any $x\in \mathbb{R}$, there exists a real
number $t_{0}(x)<uep(F)$  such that \\   
\begin{equation}
t_{0}(x)<t<uep(F)\Longrightarrow \alpha (t)x+\beta (t)<uep(F).
\label{evtp.06c}
\end{equation}

\noindent (a) If $uep(F)<+\infty$ the for any $x\in \mathbb{R},$ 
\begin{equation*}
\alpha (t)x+\beta (t)\rightarrow +\infty \text{ as }t\rightarrow uep(F)^{-}.
\end{equation*}

\bigskip \noindent (B) If Formula \ref{evtp.06}) holds such that for $\beta(t)>0$ for $t<uep(F)$ in the case where $uep(F)=+\infty$, and  
$\beta(t)<uep(F)$ for $t<uep(F)$ in the case where $uep(F)<+\infty$, then we have, as $t\rightarrow uep(F)^{-}$

\begin{equation} \label{evtp.X01bis}
\alpha(t)/\beta(t) \rightarrow 0,
\end{equation}

\noindent when $uep(F)^{-}=+\infty$ and
 
\begin{equation} \label{evtp.X02bis}
\alpha(t)/(uep(F)-\beta(t)) \rightarrow 0.
\end{equation}

\noindent when $uep(F) <+\infty$.\\

\bigskip \noindent (C) If the functions $\alpha(t)>0$ and $\beta(t)<uep(F)$ of $t<uep(F)$ are finite and if Formulas \ref{evtp.X02bis} and \ref{evtp.X02bis} hold, then for all $x\in \mathbb{R}$, there exists $t(x)<uep(F)$ such that for $t(x)\leq t <uep(F)$, we have $\alpha(t)x+\beta(t)<uep(F)$.
\end{lemma}

\bigskip \bigskip \noindent \textbf{Proof of Lemma \ref{evtp.lem.01}}. Let  $F$ be distribution function such that (\ref{evtp.06}) holds.\\

\noindent \textbf{Proof of Part A}.\\

\noindent (a) Suppose that $uep(F)=+\infty $ and fix $x\in \mathbb{R}$. Let $\ell$ be an arbitrary adherent point of  $\alpha (t)x+\beta (t),$ as $t\rightarrow
uep(F)$. So there exists a sequence $(t_{n})_{n\geq 1}$ such that 
\begin{equation*}
\text{ }t_{n}\rightarrow uep(F)^{-}\text{ \ and }\alpha (t_{n})x+\beta
(t_{n})\rightarrow \ell \text{ as }n\rightarrow +\infty .
\end{equation*}

\noindent Put $A=\ell +1$ if $\ell $ is finite and consider an arbitrary real number $A>0$ if $\ell =-\infty$ . So, by definition of the limit, we may find a value $%
n_{0}\geq 0$ such that for any $n\geq n_{0}$, we have 

\begin{equation*}
\alpha (t_{n})x+\beta (t_{n})\leq A.
\end{equation*}

\noindent It follows that for any $n_{0}\geq 0,$ we have 
\begin{equation*}
\Gamma (F,x,\alpha (t_{n}),\beta (t_{n}))\geq \frac{A}{1-F(t_{n})}.
\end{equation*}

\noindent Hence,  
\begin{equation*}
\Gamma (F,x,\alpha (t),\beta (t))+\infty \text{  as }n\rightarrow +\infty ,
\end{equation*}

\noindent which contradicts (\ref{evtp.06a}). Thus, $\ell =+\infty$ is the unique
adherent point of  $\alpha (t)x+\beta (t),$ as $t\rightarrow uep(F)^{-}$. Hence,
Point (a) us proved.\\

\bigskip \noindent (b)  Let us suppose that $uep(F)<+\infty$ and let $x\in \mathbb{R}$ be fixed. If (\ref{evtp.06c}), we
would be able to construct a sequence $(t_{n})_{n\geq 1}$ such that

\begin{equation*}
\text{ }t_{n}\uparrow uep(F)^{-}\text{ \ and for all }n\geq 0,\text{ }\alpha
(t_{n})x+\beta (t_{n})\geq uep(F).
\end{equation*}

\noindent We would have
\begin{equation*}
\Gamma (F,x,\alpha (t_{n}),\beta (t_{n}))=\frac{0}{1-F(t_{n})}\rightarrow 0%
\text{ as }n\rightarrow +\infty .
\end{equation*}

\noindent This would  contradict (\ref{evtp.06a}). Thus (b) is proved.\\

\bigskip \noindent (c). Suppose  that $uep(F)<+\infty$. By Point (b), $\alpha (t)x+\beta (t)<uep(F)$ for $t$ near enough $uep(F)$. Suppose that 
$\alpha (t)x+\beta (t)$ does not converge to $uep(F)$ as $t\rightarrow uep(F)^{-}$. Thus, there exist $\varepsilon>0$ and $(t_n)_{n\geq 0}$ such that 
$t_n\rightarrow uep(F)^{-}$ as $n\rightarrow +\infty$ and for all $n\geq $, we have $\alpha(t_n)x+\beta (t_n)<uep(F)-\varepsilon$. It follows that

$$
\liminf_{n\rightarrow +\infty} \Gamma(F,x,\alpha(t_n),\beta(t_n)) \geq \liminf_{n\rightarrow +\infty} \frac{1-F(uep(F)-\varepsilon)}{1-F(t_n)}=+\infty, 
$$

\noindent which is impossible since the right-hand is finite and equal to $e^{x}$.\\

\noindent The three points (a), (b) and (c) are summarized in the statement of Point (A).\\

\noindent \textbf{Proof of Part B}.\\

\noindent Let us suppose that $uep(F)=+\infty$. Let $x<0$. Since $\alpha (t)x+\beta (t)\rightarrow 0$ as $t\rightarrow +\infty$, there exists $t_0>0$ such that : $t\geq t_0$ implies that  $x\alpha (t)x+\beta (t)$, that is (since $\beta (t)\geq 0$),

$$
\alpha (t)/\beta (t) \leq -1/x,
$$

\noindent for all $x<0$, for all $t\geq t_0$. We get the result by letting successively $x\downarrow -\infty$ and $t\rightarrow +\infty$.\\

\noindent Let us suppose that $uep(F)<\infty$. Let $x>0$. Since $\alpha (t)x+\beta (t) <uep(F)$ for $t$ near $uep(F)$, that is $t_0\leq <uep(F)$ for some $t_0<uep(F)$, and since $\beta(t)<uep(F)$ for $t<uep(F)$, we have
$$
\alpha (t)/(uep(F)-\beta (t)) \leq 1/x,
$$

\noindent for all $t>0$, for all $t_0\leq t <1$. We get the result by letting successively $x\uparrow +\infty$ and $t\rightarrow uep(F)^{-}$.\\

\noindent \textbf{Proof of Part C}. There is nothing to do for $uep(F)=+\infty$. Now, let us suppose now that $uep(F)$ is finite. Let us begin with a positive $x$. We have

$$
\alpha(t)x+\beta(t) =(uep(F)-\beta(t)) \left( x\frac{\alpha(t)}{uep(F)-\beta(t)} - 1\right) +uep(F), \ t<uep(F).
$$

\noindent Hence for any $\varepsilon >0$ such that $\varepsilon<1$, there exists $t(x)<uep(F)$ such that for $t_1(x)\leq t <uep(F)$, we have

$$
\alpha(t)x+\beta(t) <(uep(F)-\beta(t)) (- 1+\varepsilon) +uep(F)=\beta(t)+\varepsilon (uep(F)-\beta(t))<uep(F).
$$

\noindent $\blacksquare$\\

\bigskip \noindent As well, the following fact will be useful for the sequel.\\

\begin{lemma} \label{evtp.charac.Extra} Let $a$ be finite or an infinite non-negative number. Let $r(x)$ and $b(x)$ be two functions of $x<a$ and $b(x)\rightarrow 0$ as $x\rightarrow a^{-}$. Suppose that we have for $a=+\infty$, 

\begin{equation} \label{evtp.Exp_r_general}
r(x)=c +\int_{x_1}^{x} b(t) dt
\end{equation}

\noindent where $c$ is a constant and for $a<+\infty$,
\begin{equation} \label{evtp.Exp_r_finite}
r(x)=-\int_{x}^{a} b(t) dt.
\end{equation}

\noindent Then, \\

\noindent (a) we have, as $t\rightarrow a^{-}$, $r(t)/t \rightarrow 0$ if $a=+\infty$ and $r(t)/(a-t) \rightarrow 0$ if $a<+\infty$\\

\noindent and\\

\noindent (b) for all $x\in \mathbb{R}$, $xr(t)+t \rightarrow a$ as $t\rightarrow a^{-}$ and there exists $t(x)<a$ such that $xr(t)+t<a$ for 
$t(x)\leq t<a$.\\

\noindent These results still hold if we replace Conditions (\ref{evtp.Exp_r_general}) and (\ref{evtp.Exp_r_finite}) by the following one :  there exists $x_0<a$ such that the derivative function $r^{\prime}(x)$ exists for $x_0\leq x <a$, $r^{\prime}(x) \rightarrow 0$ as $x\rightarrow a^{-}$, and 
$r(x) \rightarrow 0$ as $x\rightarrow a^{-}$ 
\end{lemma}

\bigskip \noindent \textbf{Proof of Lemma \ref{evtp.charac.Extra}}.\\

\noindent \textbf{Proof of Part (a)}. Let us prove from Formula \ref{evtp.Exp_r_general}, that 
\begin{equation}
r(x)/x \rightarrow 0 \ as \ x\rightarrow a, \label{evtp.X1bis}
\end{equation}

\noindent when $a=+\infty$. Indeed, for any $\varepsilon>0$, there exists $x_2>x_1$ such that $|b(t)|\leq \varepsilon$ for $x\leq x_2$. Thus

\begin{eqnarray*}
0\leq \limsup_{x\rightarrow +\infty, \ x>x_2} r(x)/x&=&\limsup_{x\rightarrow +\infty} \frac{1}{x}\left\vert r(x_1)+\int_{x_1}^{x_2} b(t) dt +\int_{x_2}^{x} b(t) dt\right\vert\\
&\leq& \limsup_{x\rightarrow +\infty, x>x_2} \frac{1}{x}\left\vert r(x_1)+\int_{x_1}^{x_2} b(t) dt +\varepsilon (x-x_2) \right\vert\\
&\leq& \varepsilon.
\end{eqnarray*}

\noindent We get the desired result by letting $\varepsilon\downarrow 0$.\\

\noindent We have to prove, from Formula (\ref{evtp.Exp_r_finite}), that  
\begin{equation}
r(x)/(uep(F)-x) \rightarrow 0 \ as \ x\rightarrow a, \label{evtp.X2bis},
\end{equation}

\noindent when $a<+\infty$. But the proof is similar and still based on the fact that $b(t)\rightarrow 0$ as $x\rightarrow uep(F)$.\\

\noindent \textbf{Proof of Part (a)}. Now, it is clear that for $a=+\infty$, we have for all real $x$, $r(t)x+t=t(1+xr(t)/t) \rightarrow +\infty$ as 
$t \rightarrow a=+\infty$. If $a$ is finite, we have, as $a\rightarrow a^{-}$,
$$
xr(t)+t=xr(t)+t-a+a=(a-t)\left(x\frac{r(t)}{a-1}-1\right) + a \rightarrow a.
$$

\noindent The last part of Part (b) follows exactly as in the proof of Part (C) in the above Lemma \ref{evtp.lem.01}.\\

\noindent To finish, let us suppose that  there exists $x_0<a$ such that the derivative function $r^{\prime}(x)$ exists for $x_0\leq x <a$, $r^{\prime}(x) \rightarrow 0$ as $x\rightarrow a^{-}$, and $r(x) \rightarrow 0$ as $x\rightarrow a^{-}$ for $a$ finite. Here, we put $b(x)=r^{\prime}(x)$ and Formula (\ref{evtp.Exp_r_general}) holds in general. If $a$ is finite, we may find the value of the constant by letting $x\rightarrow a^{-}$ in Formula (\ref{evtp.Exp_r_general}), and thus get Formula (\ref{evtp.Exp_r_finite}). From there, the proof above may be reproduced word by word. $\blacksquare$\\

\bigskip \noindent Before we close this section, let us mention that the convergence in the $\Gamma$-variation is in fact continuous and locally compact. Precisely, we have

\begin{proposition} \label{evtp.prop.a01} Let $F\in D(H_0)$. Let $\alpha (t)>0$ and $\beta(t)$ be two functions of $t<uep(F)$ such for all  $x\in \mathbb{R}$, as $t\rightarrow uep(F)^{-}$, we have
\begin{equation}
\Gamma (F,x,\alpha (t),\beta(t)) \rightarrow \exp (-x). \label{evtp.a01}
\end{equation}

\noindent Then for any $x\in \mathbb{R}$, for any function $x(t)>0$ of $t<uep(F)$ such for all  $x(t) \rightarrow x$ as $t\rightarrow uep(F)^{-}$, we have
$$
\Gamma (F,x(t),\alpha (t),\beta(t)) \rightarrow \exp (-x),
$$

\noindent as $t\rightarrow uep(F)^{-}$.\\

\bigskip \noindent Also, for any $-\infty < a < b <+\infty$, we have, $t\rightarrow uep(F)^{-}$,

\begin{equation}
\sup_{x \in [a,b]} \left\vert \Gamma (F,x,\alpha (t),\beta(t)) -\exp (-x)\right\vert \rightarrow 0. \label{evtp.a02}
\end{equation}
\end{proposition}

\bigskip \noindent \textbf{Proof of Proposition \ref{evtp.prop.a01}}. Assume the hypotheses of the proposition hold. Fix $x \in \mathbb{R}$. Let $x(t)>0$ be a function of $t<uep(F)$ such for all  $x(t) \rightarrow x$ as $t\rightarrow uep(F)^{-}$. For any $\varepsilon>0$ fixed, there exists $t_0<upe(F)$ such that
$$
(t_0 \leq t <uep(F)) \Rightarrow x-\varepsilon \leq x x+\varepsilon.
$$

\noindent Thus, we have, for $t_0 \leq t <uep(F)$, that

$$
\Gamma(F,x+\varepsilon,\alpha (t),\beta(t)) \leq \Gamma(F,x+\varepsilon,\alpha (t),\beta(t))\leq \Gamma(F,x-\varepsilon,\alpha (t),\beta(t)).
$$

\noindent Then, as $t\rightarrow uep(F)^{-}$, we have 

\begin{eqnarray}
e^{-(x+\varepsilon)}&=&\Gamma(F,x+\varepsilon,\alpha (t),\beta(t))\\
&\leq&\liminf_{t\rightarrow uep(F)^{-}} \Gamma(F,x(t),\alpha (t),\beta(t))\\
&\leq&\limsup_{t\rightarrow uep(F)^{-}} \Gamma(F,x(x),\alpha (t),\beta(t))\\
&\leq&\limsup_{t\rightarrow uep(F)^{-}} \Gamma(F,x-\varepsilon,\alpha (t),\beta(t))\\
&=&e^{-(x-\varepsilon)}.
\end{eqnarray}

\noindent We may conclude by letting $\varepsilon \downarrow 0$ to get (\ref{evtp.a01}). The proves the first part.\\

\noindent To prove the second part, that is (\ref{evtp.a02}), we proceed by contradiction. Suppose that (\ref{evtp.a02}) is false. This suppose that

\begin{equation*}
\limsup_{t\rightarrow uep(F)^{-}} \sup_{x \in [a,b]} \left\vert \Gamma (F,x,\alpha (t),\beta(t)) -\exp (-x)\right\vert >0.
\end{equation*}

\noindent This implies that there exists $\varepsilon>0$ such that there exists a sequence $(t_n)_{n\geq 0}$ such that ${t_n\rightarrow uep(F)^{-}}$ as $n\rightarrow +\infty$ for which, for all $\geq 1$,

$$
\sup_{x \in [a,b]} \left\vert \Gamma (F,x,\alpha (t_n),\beta(t_n)) -\exp (-x)\right\vert >\varepsilon.
$$

\noindent Then for each $n\geq 1$, there exists $x_n$ such that

\begin{equation}
\left\vert \Gamma (F,x_n,\alpha (t_n),\beta(t_n)) -\exp (-x_n)\right\vert >\varepsilon. \label{evtp.a03}
\end{equation}

\bigskip \noindent Now, the the sequence $(x_n)_{n\geq 1}$ is $[a,b]$. By the Bolzano-Weierstrass, it admits a sub-sequence $(x_{n_k})_{k\geq 1}$ converging to a point $x\in [a,b]$. By applying (\ref{evtp.a03}) to the sequence $(x_{n_k})_{k\geq 1}$, as $n\rightarrow 0$, together with (\ref{evtp.a01}) and the continuity of the exponential function at $x$ lead the the contradiction that zero is greater than $\varepsilon>0$. This proves that (\ref{evtp.a02}) holds.$\blacksquare$\\

\newpage
\section{Gumbel Extreme Value Domain and $\Gamma$-variation}

In this section, we unveil the equivalence between for a distribution $F$ between lying in the Gumbel Extreme Value Domain and satisfying the $\Gamma$-variation. This equivalence will imply a number of interesting formulas and properties that will help in getting new characterizations in the Extreme Value domain of attraction. We begin with

\begin{proposition} \label{evtp.prop.01}
Let $F$ be a \textit{cdf} such that $F\in D(H_0)$. Then $F$ is of $\Gamma$-variation.
\end{proposition}

\bigskip \noindent \textbf{Proof of \ref{evtp.prop.01}}. Suppose that $F\in D(H_0)$. Then by definition, there exist two $a(t)>0$ and $b(s)$ of $s\mathbb{R}$ such for all  $x\in \mathbb{R}$, as $s\rightarrow +\infty$, we have $s\rightarrow +\infty$,
$$
F^{s}(a(s),x+b(s)) \rightarrow H_0(x)
$$

\noindent Now, by Formula (\ref{evtp.05a}), we have for all  $x\in \mathbb{R}$, as $s\rightarrow +\infty$,

\begin{equation}
s(1-F(a(s)x+b(s)) \rightarrow e^{-x} \label{evtp.07}
\end{equation}

\noindent Let us use this formula with
$$
s(t)=1/(1-F(t))
$$

\noindent so that $s(t)\rightarrow +\infty$ as $t\rightarrow uep(F)^{-}$. Put

$$
a(s(t))=\alpha(t) \text{ and } b(s(t))=\beta(t), \ \ t<uep(F).
$$

\noindent We get, for all  $x\in \mathbb{R}$, as $t\rightarrow uep(F)^{-}$, we have
$$
\Gamma (F,x,\alpha (t),\beta(t)) \rightarrow \exp (-x),
$$

\noindent which establishes the $\Gamma$-variation of $F$.$\blacksquare$\\

\noindent The next propositions show how we could change the functions in the definition of $\Gamma$-variation.\\

\begin{proposition} \label{evtp.prop.02}
Let $F$ be a \textit{cdf}. Suppose that $F$ is $\Gamma$-varying such that for some function $\alpha (t)>0$ of $t<uep(F)$, we have, $x\in \mathbb{R}$, as $t\rightarrow uep(F)^{-}$,

\begin{equation}
\Gamma (F,x,\alpha (t), \beta(t)) \rightarrow \exp (-x), \label{evtp.08}
\end{equation}

\noindent Then $F \in D(H_0)$.\\

\noindent Precisely, if there exist two functions $\alpha (t)>0$ and $\beta (t)$ of $t<uep(F)$ such that for all  $x\in \mathbb{R}$, as $t\rightarrow uep(F)^{-}$, we have
$$
\Gamma (F,x,\alpha (t),\beta(t)) \rightarrow \exp (-x),
$$

\noindent then, by denoting $t(s)=(F^{-1}(1-1/s)$, $a(s)=\alpha(t-(s))$ and  $b(s)=\beta(t(s))$ for for $s>1$, we have for all $x\in \mathbb{R}$,
$$
s(1-F(a(s)x+b(s)))\rightarrow e^{-x}, \ as \ s\rightarrow +\infty.
$$

\end{proposition}

\bigskip \noindent \textbf{Proof of Proposition \ref{evtp.prop.02}}. Assume that Formula (\ref{evtp.08}) hold for a distribution $F$. We may use it for
$$
t(s)=F^{-1}(a-1/s), \ \ s>1,
$$ 

\noindent so that $t(s)\rightarrow uep(F)^{-}$ as $s\rightarrow +\infty$. We get, for all $x\in \mathbb{R}$, as $s\rightarrow +\infty$,
\begin{equation}
\frac{1-F(\alpha (t(s))x+\beta(t(s)))}{1-F(t(s))} \rightarrow e^{-x}, \label{evtp.08a}
\end{equation}

\noindent with $a(s)=\alpha(t(s))$ and $b(s)=\beta(t(s))$, $s>1$. Now, let us prove that

\begin{equation}
s(1-F(t(s)) \rightarrow 1 \text{ as } \rightarrow +\infty. \label{evtp.08b}
\end{equation}

\bigskip \noindent To this end, let us use the properties of generalized inverses functions (see for example, Chapter 4 in \cite{wc-srv-ang} of the current series), in particular Point 9 in the above cited chapter, that is

$$
\forall u \in (0,1), \ F(F^{-1}(u)-) \leq u \leq F(F^{-1}(u)), 
$$

\bigskip \noindent since $F$ is right-continuous. Applying this to $u=1-1/s$, $s>1$, leads to

$$
\forall s>1, \ F(t(s)-)\leq 1-1/s \leq F(t(s)). 
$$

\noindent Since $\alpha(t(s))>0$ for $t<uep(F)$, then for any $\varepsilon>0$, for any $s>0$, we have by the finition of the left-hand limit

$$
\forall s>1, \ F(t(s)-\varepsilon \alpha(t(s)))\ 1-1/s \leq F(t(s)),
$$

\noindent which my be rewritten as

$$
\forall s>1, \ 1\leq \frac{1}{s(1-F(t(s))} \leq \frac{1-F(\alpha (t(s))\varepsilon x+t(s) )}{1-F(t(s))}.
$$

\noindent By applying Formula (\ref{evtp.08}), we have for any $\varepsilon>0$,

\begin{eqnarray}
1 &\leq& \liminf_{s\rightarrow +\infty} s(1-F(t(s))\\
	&\leq& \limsup_{s\rightarrow +\infty} s(1-F(t(s))\\
  &\leq& \limsup_{s\rightarrow +\infty} \frac{1-F(\alpha (t(s))\varepsilon x+t(s) )}{1-F(t(s))}\\
	&=& \lim_{s\rightarrow +\infty} \frac{1-F(\alpha (t(s))\varepsilon x+t(s) )}{1-F(t(s))}\\
	&=& e^{\varepsilon}.
\end{eqnarray}

\noindent By letting $\varepsilon \downarrow 0$, we get (\ref{evtp.08b}). By plugging this in Formula (\ref{evtp.08}) yields : for all  $x\in \mathbb{R}$, as $s\rightarrow +\infty$,

$$
s(1-F(a(s)x+b(s)))\rightarrow e^{-x},
$$

\noindent which, by the preliminary remarks of this chapter, implies that $F\in D(H_0))$.$\blacksquare$.\\

\bigskip \noindent The following proposition shows how the limit is affected when the auxiliary functions $\alpha(t)$ and $\beta(t)$ of $t<uep(F)$ are changed in a precise way.

\begin{proposition} \label{evtp.prop.03A}  Let $F$ be a \textit{cdf}. Suppose that there exist two functions $\alpha (t)>0$ and $\beta (t)$ of $t<uep(F)$ such that for all  $x\in \mathbb{R}$, as $t\rightarrow uep(F)^{-}$, we have
\begin{equation}
\Gamma (F,x,a(t),b(t)) \rightarrow \exp (-x). \label{evtp.11a}
\end{equation}

\noindent The following assertions hold.\\

\noindent (A) Suppose that there exist two functions $\alpha (t)>0$ and $\beta (t)$ of $t<uep(F)$ such that 
$s\rightarrow +\infty$,

\begin{equation}
\alpha(t)/a(t) \rightarrow A>0, \ A \in \mathbb{R} \text{ and } \frac{\beta(t)-b(t)}{a(t)} \rightarrow B \in \mathbb{R}, \label{evtp.11ab}
\end{equation}

\noindent then, for all  $x\in \mathbb{R}$, as $t\rightarrow uep(F)^{-}$, we have
\begin{equation}
\Gamma (F,x,\alpha(t),\beta(t)) \rightarrow \exp(-(Ax+B)). \label{evtp.11b}
\end{equation}

\noindent (B) If Formula (\ref{evtp.11b}) holds on top of Formula \ref{evtp.11a} with $A>0$ and $B\in \mathbb{R}$, then for for $t(s)=F^{-1}(1-1/s)$ with $s>1$, we have, as $s\rightarrow +\infty$ (which is equivalent to $t(s)\rightarrow uep(F)^{-}$)

\begin{equation}
\alpha(t(s))/a(t(s)) \rightarrow A>0, \ A \in \mathbb{R} \text{ and } \frac{\beta(t(s))-b(t(s))}{a(t(s))} \rightarrow B \in \mathbb{R}, \label{evtp.11b}
\end{equation}
\end{proposition}

\bigskip \noindent \textbf{Proof of Proposition \ref{evtp.prop.03A}}.\\

\noindent \textbf{Proof of Part (A)}. Suppose that Formulas (\ref{evtp.11a}) and (\ref{evtp.11ab}) hold. Then, there exists $t_0<uprp(F)$ such that for $t_0<t<uep(F)$,

$$
a(t)(A-\varepsilon)\leq \alpha(t) \leq a(t) (A+\varepsilon)
$$

\noindent and

$$
\alpha(t)(B-\varepsilon)+b(t) \leq \beta(t) \leq \alpha(t) (B+\varepsilon)+b(t).
$$

\noindent Now, for $x\geq 0$, by denoting $x_1(\varepsilon)=(A-\varepsilon)x+(B-\varepsilon)$ and $x_2(\varepsilon)=(A+\varepsilon)x+(B+\varepsilon)$, we have for  
$t_0<t<uep(F)$,

\begin{eqnarray}
a(t) x_1(\varepsilon)+b(t) \leq& \alpha(t)+\beta(t) \leq a(t) x_2(\varepsilon)+b(t). \label{evtp.11abc}
\end{eqnarray}

\noindent and then $t_0<t<uep(F)$,

$$
\Gamma (F,x_2(\varepsilon),a(t),b(t)) \leq \Gamma (F,x,\alpha(t),\beta(t)) \leq \Gamma (F,x_1(\varepsilon),a(t),b(t))
$$

\noindent Now, by applying (\ref{evtp.11a}) and (\ref{evtp.11b}) in the latter formula, we get as $t\rightarrow uep(F)^{-}$,

\begin{eqnarray*}
e^{-((A+\varepsilon)x+B+\varepsilon)}&\leq& \liminf_{t\rightarrow uep(F)^{-}} \Gamma (F,x,\alpha(t),\beta(t))\\
&\leq& \limsup_{t\rightarrow uep(F)^{-}} \Gamma (F,x,\alpha(t),\beta(t))\\
&\leq&  e^{-((A-\varepsilon)x+B-\varepsilon)}.
\end{eqnarray*}

\noindent Now by letting $\varepsilon \uparrow 0$, we get (\ref{evtp.11b}). For $x<0$, we proceed by the same manner and by reversing the inequalities in 
Formula (\ref{evtp.11abc}).\\

\noindent \textbf{Proof of Part B}. Let us suppose that Formula ({evtp.11a}). By denoting $t(s)=(F^{-1}(1-1/s)$, $a(s)=\alpha(t-(s))$ and  $b(s)=\beta(t(s))$ for for $s>1$, we have by Proposition \ref{evtp.prop.02}, that for all $x\in \mathbb{R}$

$$
s(1-F(a(s)x+b(s)) \rightarrow e^{-x}, \ as \ s\rightarrow +\infty.
$$

\noindent Now, by Lemma \ref{evtp.prop.02}, we get for all $x\in \mathbb{R}$

$$
F^{s}(a(s)x+b(s)) \rightarrow H_0(x), \ as \ s\rightarrow +\infty.
$$

\noindent Then, by Remark (R1) in the preliminaries of this section, we have that Formula (\ref{evtp.11b}) holds.$\blacksquare$\\

\bigskip \noindent But we want to go further and to show that Formula (\ref{evtp.11b}) in Lemma \ref{evtp.prop.03A}, actually holds for all $t$ converging to $uep(F)^{-}$ and not for the special case of $t(s)$. Indeed, we have the final form of the lemma as in :

\begin{proposition} \label{evtp.prop.03B}  Let $F$ be a \textit{cdf}. Suppose that there exist two functions $\alpha (t)>0$ and $\beta (t)$ of $t<uep(F)$ such that for all  $x\in \mathbb{R}$, as $t\rightarrow uep(F)^{-}$, we have
\begin{equation}
\Gamma (F,x,a(t),b(t)) \rightarrow \exp (-x). \label{evtp.11na}
\end{equation}

\noindent The following assertions hold.\\

\noindent (A) Suppose that there exist two functions $\alpha (t)>0$ and $\beta (t)$ of $t<uep(F)$ such that 
$s\rightarrow +\infty$,

\begin{equation}
\alpha(t)/a(t) \rightarrow A>0, \ A \in \mathbb{R} \text{ and } \frac{\beta(t)-b(t)}{a(t)} \rightarrow B \in \mathbb{R}, \label{evtp.11nab}
\end{equation}

\noindent then, for all  $x\in \mathbb{R}$, as $t\rightarrow uep(F)^{-}$, we have
\begin{equation}
\Gamma (F,x,\alpha(t),\beta(t)) \rightarrow \exp(-(Ax+B)). \label{evtp.11nb}
\end{equation}

\noindent (B) If Formula (\ref{evtp.11nb}) holds on top of Formula \ref{evtp.11na} with $A>0$ and $B\in \mathbb{R}$, then Formula (\ref{evtp.11nab}) holds.
\end{proposition}

\bigskip \noindent \textbf{Proof of Proposition \ref{evtp.prop.03B}}. Here, we only need to prove the part B. Theoof of Part B is based on a reasoning by contradiction, which is lengthy. We propose it in details in the Appendix Section \ref{evtp.sec.append}.\\

\bigskip \noindent The main door of the characterization is the following Theorem. Before we present it, let us introduce the asymptotic moments of $F$. At the first order, we define :

\begin{equation}
R(F,x)=(1-F(x))^{-1} \int_{x}^{uep(F)} 1-F(t)\ dt, \ x< uep(F). \label{evtp.masymp01}
\end{equation}

\noindent The $p$-th moment of $F$ is defined for $x< uep(F)$ by

\begin{equation}
R_p(F,x)=\frac{1}{1-F(x)} \int_{x}^{uep(F)} \int_{u_1}^{uep(F)}... \int_{u_{p-1}}^{uep(F)} 1-F(t)\ du_{1} \ ... \ 1-F(t)\ du_{p-1}dt. \label{evtp.masymp0p}
\end{equation}

\noindent By convention, we usually write $R_1(F,.)=R(F,.)$ and $R_2(F,.)=W(F,.)$. When there is no reason of confusion, we drop the symbol $F$.\\

\noindent We have :\\

\begin{theorem} \label{evtp.theo3} Let $F \in D(H_0)$. Then, $R(x)$ is finite for $x<uep(F)$, and  such for all  $x\in \mathbb{R}$, we have, as $t\rightarrow uep(F)^{-}$, 
$$
\Gamma(F,x,R(F,t), t) \rightarrow \exp (-x).
$$
\end{theorem}

\bigskip \noindent \textbf{Proof of \ref{evtp.theo3}}. To come soon.\\

\bigskip \noindent This theorem is used with the following result as a main tool.\\

\begin{proposition} \label{evtp.prop.04A}
Let $F\in D(H_0)$. Then there exists $x_0$ such that 
\begin{equation}
c=1 - \int_{x_0}^{uep(F)} 1-F(t) dt>0. \label{evtp.prop040}
\end{equation}

\noindent Then function defined by
$$
F_1(x)=\left(1 - \int_{x}^{uep(F)} 1-F(t) dt \right) 1_{(x\geq x_0)}
$$ 

\noindent is a probability distribution function lying in $D(H_0)$ such that
$$
R(F,t)/R(F_1,t) \rightarrow 1 \ as \ t\rightarrow uep(F)^{-}. 
$$
\end{proposition} 

\bigskip \noindent \textbf{Proof of \ref{evtp.prop.04A}}. By Part (B) of Lemma \ref{evtp.lem.01}, 
$$
A(x)=\int_{x}^{uep(F)} 1-F(t) dt
$$

\noindent is finite when $x$ is near $uep(F)$. Let $x_1$ a point $x_1<uep(F)$ such that $A(x_1)$ finite. Then $A(x)$ converges to zero, as $x\rightarrow uep(F)^{-}$ by the Dominated Convergence Theorem. So there exists $x_0<uep(F)$ such that Formula (\ref{evtp.prop040}) holds. Hence $F_1$ is well-defined as a probability distribution function. We are going to prove the points of the proposition by establishing that for all $x\in \mathbb{R}$, as $t\rightarrow uep(F)^{-}$,

\begin{equation}
\Gamma(F_1, x, R(F,t), t) \rightarrow e^{-x}.
\end{equation}

\noindent Let $y$ be an arbitrary real number. By Lemma \ref{evtp.lem.01}, $R(F,u)y+y \rightarrow e^{-y}$ for $u\rightarrow uep(F)^{-}$. Then by Theorem \ref{evtp.theo3}, we have for all $y\mathbb{R}$, for all $x\mathbb{R}$, as $u\rightarrow uep(F)^{-}$, 

\begin{equation}
\Gamma(F,x R(F, R(F, u)y+u),R(F, u)y+u) \rightarrow e^{-x}, \label{evtp.015a}
\end{equation}

\noindent meaning

\begin{equation}
\frac{1-F(R(F, R(F, u)y+u)x+R(F, u)y+u)}{1-F(R(F, u)y+u)} \rightarrow e^{-x}. \label{evtp.015b}
\end{equation}

\noindent But, still by \ref{evtp.theo3}, we have, as $u\rightarrow uep(F)^{-}$,
 
\begin{equation}
\frac{1-F(R(F, u)y+u)}{1-F(u)} \rightarrow e^{-y}. \label{evtp.015c}
\end{equation}

\noindent By combining Formulas (\ref{evtp.015a})-(\ref{evtp.015c}), we get for $t(u)=R(F, u)y+u$, for all $x\mathbb{R}$,

$$
\Gamma(F, x, R(F,t(u), t(u)) \rightarrow e^{-(x-y)},
$$

\noindent as $t(u)\rightarrow uep(F)^{-}$. Hence by proposition \ref{evtp.prop.03B},

$$
R(F,t(u))/R(F,u)\rightarrow 1 \ as \ u\rightarrow uep(F)^{-}.
$$

\noindent Now, remark that for $x<uep(F)$ near enough $uep(F)$, that is for some $x_0<uep(F)$, we have 

$$
R(F,x)=\frac{1-F_1(x)}{1-F(x)}, \ x_0\leq x<uep(F).
$$

\noindent Thus

$$
\left(\frac{1-F_1(t(u))}{1-F_1(u)}\right) \left( \frac{1-F(t(u))}{1-F(u)}\right)^{-1} \rightarrow 1 \ as \ u\rightarrow uep(F)^{-}.
$$

\noindent Hence, since $t(u)=R(F, u)y+u$, we have

\begin{eqnarray*}
&&\lim_{u\rightarrow uep(F)^{-}}\left(\frac{1-F_1(t(u))}{1-F_1(u)}\right)\\
&=& \lim_{u\rightarrow uep(F)^{-}} \left( \frac{1-F(t(u))}{1-F(u)}\right)=e^{-y}.
\end{eqnarray*}

\noindent It follows that $F_1$ is of $\Gamma$-variation with $\alpha(t)=R(F,t)$, $t<uep(F)$. It follows that $F_1 \in D(H_0)$ and by the way, $uep(F)=uep(F_1)$.\\

\bigskip \noindent The just discovered property is very important. Let us rephrase it again : $F_1$, which is a transformation of $F\in D(H_0)$, is also in $D(H_0)$ and is of $\Gamma$-variation. Then, there two functions $\alpha(t)>0$ and $\beta(t)$ of $t<uep(F)=uep(F_1)$, and for all  $x\in \mathbb{R}$, as $t\rightarrow uep(F)^{-}$, we have
$$
\Gamma (F,x,\alpha (t),\beta(t)) \rightarrow \exp (-x), \label{evtp.03F}
$$

 \begin{equation}
1-F(x)= c(x) \exp \left( \int_{x_1}^{uep(F)} -\frac{a(t)}{r(t)} dt\right), \ x_1\geq x <uep(F),
\end{equation}

\bigskip \noindent where $c(x)$, $a(x)$ and $r(x)$ are functions $x<uep(F)$ such that :

\noindent $(c(x),a(x))\rightarrow (c,1)$ as $x\rightarrow uep(F)^{-}$, and $r(x)$ is positive and differentiable and $r^{\prime}(x) \rightarrow 0$ as $x\rightarrow uep(F)^{-}$. Moreover, the function $r(x)$ of $x<uep(F)$ may be taken as : for $uep(F)=+\infty$,

$$
r(x)=c_1+\int_{x_1}^{x} b(t) dt, \ x_1\leq x <uep(F)
$$ 

\noindent and for $uep(F)<+\infty$, and we may replace $\alpha (t)$ by $R(F_1,t)$ or by $R(F,t)$, $t<uep(F)$. This means that $R(F_1,x)/R(F,x)\rightarrow 1$ as $x\rightarrow uep(F)^{-}$. But, by the definition of the $R(F,.)$ in Formula (\ref{evtp.masymp01}) and that of $F_1$, we have, for $x<uep(F)$,

\begin{equation}
R(F_1,x)=\left( \int_{x}^{uep(F)} \int_{u}^{uep(F)} 1-F(t)dt du \right)/\left( \int_{x}^{uep(F)} 1-F(t)dt \right). \label{evtp.quot41}
\end{equation}

\noindent and next, for $x<uep(F)$,

\begin{equation}
H(x,F)=:\frac{R(F_1,t)}{R(F,t)}=\frac{(1-F(x))\left( \int_{x}^{uep(F)} \int_{u}^{uep(F)} 1-F(t)dt du \right)}{\left( \int_{x}^{uep(F)} 1-F(t)dt \right)^2}. \label{evtp.quot21a}
\end{equation}

\noindent By by Formula (\ref{evtp.masymp0p}) and by the notation $W(F,.)=R_2(F,.)$, we have

\begin{equation}
H(F,x)=:H(x)=W(F,x)/R(F,x)^2, \ x<uep(F). \label{evtp.quot21b}
\end{equation}

\bigskip \noindent In summary, the lines above yield the following result.\\

\begin{lemma} \label{evtp.lem.quot21} Let $F$ be a \textit{cdf} on $\mathbb{R}$, such that $F\in D(H_0)$. Then the function 
$W(F,x)R(F,x)^{-2}$ of $x<uep(F)$ is defined and finite and we have

$$
W(F,x)R(F,x)^{-2} \rightarrow 1 \ as \ x\rightarrow uep(F)^{-}.
$$
\end{lemma}

\bigskip \noindent Follows a characterization of $D(H_{0})$ which reverses the result in Lemma \ref{evtp.lem.quot21} and contains a useful representation.\\

\begin{theorem} \label{evtp.theo.08} Let $F$ be a \textit{cdf} on $\mathbb{R}$. The following assertions are equivalent.\\

\noindent (a) $F \in D(H_0)$.\\

\noindent (b) The functions $R(F,x)$ and $W(F,x)R(F,x)^{-2}$ of $x<uep(F)$ is defined and finite and we have

$$
W(F,x)R(F,x)^{-2} \rightarrow 1 \ as \ x\rightarrow uep(F)^{-}.
$$

\noindent is defined for $x<uep(F)$ and $H(x)\rightarrow 1$ as $x\rightarrow uep(F)^{-}$.

\noindent (c) $F$ admits the following representation : there exists $x_1<uep(F)$ and a constant $c>0$ such that  :

\begin{equation}
1-F(x)=c(x) \exp \left( - \int_{x_1}^{x}\frac{a(t)}{r(t)} dt\right), \ x_1\leq x <uep(F), \label{evtp.dehaanD}
\end{equation}

\noindent $c(x)$, $a(x)$ and $r(x)$ are measurable functions such that $(c(x),a(x)) \rightarrow (c,1)$ as $x\rightarrow uep(F)^{-}$ and $r(.)$ is positive and differentiable satisfying :  $r^{\prime}(x) \rightarrow 0$ as $x\rightarrow uep(F)^{-}$ and $r(x) \rightarrow 0$ as $x\rightarrow uep(F)^{-}$ if $uep(F)$ is finite.\\

\noindent Moreover, we may take for $uep(F)=+\infty$

$$
r(x)=c_1 + \int_{x_1}^{x} b(t) dt, \ x_1\leq x <uep(F),
$$ 

\noindent and $uep(F)<+\infty$

$$
r(x)=-\int_{x}^{uep(F)} b(t) dt, \ x_1\leq x <uep(F),
$$ 

\noindent where $c_1$ is a constant and $r^{\prime}(x)=b(x) \rightarrow 0$ as $x\rightarrow uep(F)^{-}$. 
\end{theorem}

\bigskip \noindent \textbf{Proof of Theorem \ref{evtp.theo.08}}. Let use a circular method by proving : $(a)\Rightarrow (b)$, $(b)\Rightarrow (c)$ and $(c)\Rightarrow (a)$.\\

\noindent \textit{Proof of $(a)\Rightarrow (b)$}. This is exactly Lemma \ref{evtp.lem.quot21}.\\

\bigskip \noindent \textit{Proof of $(b)\Rightarrow (c)$}.\\

\noindent We begin to set $b(x)=-1+H(x)$, with $H(x)=W(F,x)R(F,x)^{-2}$, $x<uep(F)$. By the assumptions, $b(x) \rightarrow 0$ as $x\rightarrow uep(F)^{-}$. Also, the function $r(.)=R(F_1,.)$ (see Formula \ref{evtp.quot41})
$$
r(x)=\frac{\int_{x}^{uep(F)} \int_{u}^{uep(F)} 1-F(t)dt du}{\int_{x}^{uep(F)} 1-F(t)dt}, \ x<uep(F),
$$

\noindent is finite and positive. A simple computation shows that $L^{\prime}(x)=b(x)$ $a.e.$ Thus, for all $x_1<x<uep(F)$,

\begin{equation} \label{evtp.Exp_r_general}
r(x)=\int_{x_1}^{x} b(t) dt +r(x_1).
\end{equation}

\noindent If $uep(F)<+\infty$, we have for $x<uep(F)$,
\begin{eqnarray*}
\int_{x}^{uep(F)} \int_{u}^{uep(F)} 1-F(t)dt du &\leq& \int_{x}^{uep(F)} \int_{x}^{uep(F)} 1-F(t)dt du\\
&\leq& ({uep(F)}-x) \int_{x}^{uep(F)} 1-F(t)dt du\\
\end{eqnarray*}

\noindent and next, $0\leq r(x)\leq ({uep(F)}-x)\rightarrow 0$ as $x\rightarrow uep(F)^{-}$. By plugging this in the Formula \ref{evtp.Exp_r_general}, we get
$$
\int_{x_1}^{uep(F)} b(t) dt =-r(x_1)
$$

\noindent and next

\begin{equation} \label{evtp.Exp_r_finite}
r(x)=-\int_{x}^{uep(F)} b(t) dt, \ x<uep(F).
\end{equation}

\noindent By putting 
$$
w(x)=\int_{x}^{uep(F)} \int_{u}^{uep(F)} 1-F(t)dt du, \ x<uep(F),
$$

\noindent we see that $r^{-1}(x)dx=-dw(x)/w(x)$ and thus, for all $x_1<x<uep(F)$,

$$
\int_{x_1}^{x} r^{-1}(t) dt= - \log w(x) + \log w(x_1).
$$

\noindent which leads to

$$
w(x)=w(x_1)^{-1} \exp\left(- \int_{x_1}^{x} r^{-}(t) dt \right).
$$

\noindent From there, we notice first that

$$
L(x)=r(x)/R(F,x) \rightarrow 1 \ as \ 
$$

\noindent and next

$$
(1-F(x))L(x)=R(F,x)^{-2} w(x), \ x<upe(F).
$$

\noindent It follows from the three last equations that

$$
1-F(x)=d(x) w(x_1)^{-1} r(x)^{-2} \exp\left(- \int_{x_1}^{x} r^{-1}(t) dt \right),\ x<uep(F).
$$

\noindent where $d(x)\rightarrow 1$ as $x\rightarrow uep(F)$. Now, we have from both Formulas \ref{evtp.Exp_r_finite} and \ref{evtp.Exp_r_general} that $dr(x)=b(x)dx$, and thus for $x_1<x<uep(F)$,

\begin{eqnarray*}
\exp\left( \int_{x_1}^{x} \frac{b(t)}{r(t)} dt \right)&=&\exp\left( \int_{x_1}^{x} \frac{dr(t)}{r(t)} \right)\\
&=&\exp\left( \log r(x) - \log r(x_1) \right)\\
&=&r(x)/r(x_1).
\end{eqnarray*}

\noindent By combining the two last equations, we get

$$
1-F(x)=d(x) w(x_1)^{-1} r(x_1)^{2} \exp\left(- \int_{x_1}^{x} \frac{1+2b(t)}{r(t)} dt \right),\ x<uep(F).
$$

\noindent By putting $c(x)=d(x) w(x_1)^{-1} r(x_1)^{2}$, $c=w(x_1)^{-1} r(x_1)^{2}$ and $a(x)=1+2b(x)$, we arrive at Formula \ref{evtp.dehaanD}. We still have a few number of points to check. First, we obviously have, from Formulas \ref{evtp.Exp_r_finite} and \ref{evtp.Exp_r_general}, that $r^{\prime}(x)=b(x)\rightarrow 0$ as  $x\rightarrow uep(F)$.\\

\bigskip \noindent \textit{Proof of $(c)\Rightarrow (a)$}. Suppose that $(c)$ holds. Let us use $(a)$ by showing that $F$ is of $\Gamma$-variation. From the assumptions of (b), and by Lemma \ref{evtp.charac.Extra}, we have, as $t\rightarrow uep(F)^{-}$, that  $r(t)/t \rightarrow 0$ if $uep(F)=+\infty$ and $r(t)/(uep(F)-t) \rightarrow 0$ if $uep(F)<+\infty$, and for all $x\in \mathbb{R}$, $xr(t)+t \rightarrow uep(F)^{-}$ as $t\rightarrow uep(F)^{-}$ and there exists $t(x)<uep(F)$ such that $xr(t)+t<a$ for $t(x)\leq t<a$. Then for $x$ fixed, for  $t(x)\leq t<a$, we may apply Formula (\ref{evtp.dehaanD}) to $t$ and to $xr(t)+t$. Thus, we have 

$$
c(xr(t)+t)/c(x) \rightarrow 1, \ as \ t\rightarrow uep(F)^{-}.
$$

\noindent and then,

\begin{equation} \label{evtp.quot40}
\Gamma(F,x, r(t),t)=(1+o(1)) \left( -\int_{t}^{xr(t)+t}\frac{a(s)}{r(s)} ds\right).
\end{equation}

\noindent \noindent If $x=0$, we get $\Gamma(F,0, r(t),t)=1=e^{0}$ and there is nothing to prove.\\

\noindent Let us proceed with $|x|>0$.  For any $|y|\leq |x|$, we have $-|x|r(t)+t \leq yr(t)+t\leq |x|r(t)+t$, which implies that $yr(t)+t \rightarrow uep(F)^{-}$ uniformly in $|y|\leq |x|$. Since $r^{\prime}(t)\rightarrow 0$ as $t \rightarrow uep(F)^{-}$, we may find for an arbitrary $\varepsilon>0$, a value $t_0$ such that $t(-|x|)\vee t(|x|) \leq t_0 <uep(F)$ and $|r^{\prime}(s)|\leq \varepsilon$ for $t_0\leq s <uep(F)$. hence, for $t_0\leq t <uep(F)$

$$
|r(yr(t)+t)-r(t)|= \int_{t}^{yr(t)+t} |r'(s)|ds \leq \varepsilon |x| r(t).
$$

\noindent that is

$$
\sup_{|y|\leq |x|} \left\vert \frac{r(yr(t)+t)}{r(t)} -1\right\vert \leq \varepsilon |x|.
$$

\noindent Hence
$$
A(t)=\sup_{|y|\leq |x|} \left\vert \frac{r(yr(t)+t)}{r(t)} -1\right\vert \rightarrow 0 \ as \ t \rightarrow uep(F)^{-}.
$$

\noindent As well, we may see that 
$$
\sup_{s\in [(xr(t)+t)\wedge t, (xr(t)+t)\vee t]} |a(s)-1| \rightarrow 0 \ as \ t \rightarrow uep(F)^{-}.
$$

\noindent By combining all this, we are able to find, for an arbitrary $\eta>0$, a value $t_1<uep(F)$ such that for $t_1\leq t <uep(F)$, we have

$$
B(t)=\sup_{s\in [(xr(t)+t)\wedge t, (xr(t)+t)\vee t]} |a(s)-1|\eta, \ |c(xr(t)+t)/c(x)-1|\leq \eta,
$$

\noindent and

$$
B(\sup_{|y|\leq |x|} \left\vert \frac{r(yr(t)+t)}{r(t)} -1\right\vert \eta.
$$

\noindent We finally get for $t_1\leq t <uep(F)$, from Formula \ref{evtp.quot40}, that

\begin{eqnarray} 
\log \Gamma(F,x,r(t),t)&=& o(1) - \int_{t}^{xr(t)+t}\frac{a(s)}{r(s)} ds\\
&=& \int_{t}^{xr(t)+t}\frac{1}{r(s)} ds + S(1,t),
\end{eqnarray} 

\noindent where

$$
|S(1,t)|=|\int_{t}^{xr(t)+t}\frac{a(s)-1}{r(s)} ds| \leq B(t) \int_{-|x|r(t)+t}^{|x|r(t)+t}\frac{1}{r(s)} ds. 
$$

\noindent Next, by change of variables $u=(s-t)/r(t)$ in Formula \ref{evtp.quot40}, 

\begin{eqnarray}
\int_{t}^{xr(t)+t}\frac{1}{r(s)} ds&=&\int_{0}^{x}\frac{r(ur(t)+1)}{r(u)} du\\
&=:& \int_{0}^{x} du + S(2,t)\\
&=& x + S(2,t)
\end{eqnarray}

\noindent where
$$
|S(2,t)|=|\int_{0}^{x} \left(\frac{r(ur(t)+1)}{r(u)}-1\right) du|\leq A(t)|x| as \ t \rightarrow uep(F)^{-}.
$$

\noindent In the same manner

$$
|\int_{-|x|r(t)+t}^{|x|r(t)+t}\frac{1}{r(s)}| ds\leq 2 A(t)|x| as \ t \rightarrow uep(F)^{-}.
$$

\noindent In conclusion, the five last equations together yield that

$$
\Gamma(F,x,r(t),t) \rightarrow e^{-x} \ as \ t \rightarrow uep(F)^{-},
$$

\noindent for all $x\in \mathbb{R}$, which was the target. $\blacksquare$\\

\newpage 
\section{Appendix} \label{evtp.sec.append} $ $

\noindent \textbf{Proof of Part B of Lemma \ref{evtp.prop.03B}}.\\

\noindent Let us suppose that Formulas \ref{evtp.11na} and \ref{evtp.11nb} hold with $A>0$ and $B\in \mathbb{B}$. Let us break the proof into two steps.\\

\noindent \textbf{Step 1}. Let $\ell_1$ be an adherent point of $\{\alpha(t)/a(t), \ t<uep(F)\}$ as $t\rightarrow uep(F)^{-}$. Then there exits a sequence $(t_n)_{n\geq 0}$ such that $t_n \rightarrow uep(F)^{-}$ and $\alpha(t_n)/a(t_n) \rightarrow \ell_1$ as $n\rightarrow +\infty$.\\

\noindent Now let $\ell_2$ be an adherent point of $\{(\beta(t_n)-b(t_n))/a(t_n), \ n\geq 0 \}$ as $n\rightarrow +\infty$. Then there exists a sub-sequence 
$(t_{n_k})_{k\geq 0}$ of $(t_n)_{n\geq 0}$ such that
$$
(\beta(t_{n_k})-b(t_{n_k}))/a(t_{n_k}) \rightarrow \ell_2 \ as \ k \rightarrow +\infty.
$$

\noindent We also have that $\alpha(t_{n_k})/a(t_{n_k})$ converges to $\ell_1$ as a sub-sequence of $(\alpha(t_n)/a(t_n))_{n\geq 1}$.\\

\noindent We are going to see that we will necessarily have that $\ell_1=A$ and $\ell_2=B$. Let us prove but by excluding all the other possibilities. Let us give them into cases, we will show to be impossible. In the four first cases, we suppose that $\ell_1$ is infinite (equal to $+\infty$) or $\ell_2$ is. If this is impossible, we have that $\ell_1$ and $\ell_2$ are finite. Then we split the hypothesis $\ell_1\neq A$ and $\ell_2 \neq B$ into fours cases which also are shown ti be impossible. The conclusion will be that $\ell_1=A$ and $\ell_2=B$.\\

\noindent Case 1. $\ell_1=+\infty$ and $\ell_2=+\infty$. Fix $x_0>0$. It follows that for any $C>0$, there exists $k_0$ such that for any $k\geq k_0$

$$
t_{n_k})/a(t_{n_k} \geq C \ \ and \ \ \beta(t_{n_k})-b(t_{n_k}))/a(t_{n_k} \geq C,
$$

\bigskip \noindent which implies that, for $k\geq k_0$, we have

$$
\alpha(t_{n_k})x_0+\beta(t_{n_k}) \geq C a(t_{n_k}) (x_0+1)+b(t_{n_k}), 
$$

\bigskip \noindent which in turn implies for $k\geq k_0$,

$$
\Gamma(F, x_0, \alpha(t_{n_k}), \beta(t_{n_k})) \leq \Gamma(F, C(x_{0}+1), a(t_{n_k}), b(t_{n_k})).
$$

\noindent Now, by Formulas (\ref{evtp.11na}) and (\ref{evtp.11nb}), we get, as $k\rightarrow +\infty$,
$$
\exp(Ax_{0}+B) \leq \exp(- C(x_{0}+1)),
$$

\noindent for all $C>0$. The conclusion is absurd since the left-hand member of the latter inequality is fixed while the right-hand one tends to zero as 
$C\rightarrow +\infty$.\\

\noindent \textbf{Remark}. In the next cases, we will use similar methods. So, we will skip some intermediate steps and go directly to comparison to the $\Gamma$ quantities and conclude.\\

\bigskip \noindent Case 2. $\ell_1=+\infty$ and $\ell_2=-\infty$. Fix $x_{0}$ such that $-1<x_{0}<0$. It follows that for any $C>0$, there exists $k_0$ such that for any $k\geq k_0$

$$
t_{n_k})/a(t_{n_k} \geq C \ \ and \ \ \beta(t_{n_k})-b(t_{n_k}))/a(t_{n_k} \leq -C,
$$

\noindent which implies that, for $k\geq k_0$, we have

$$
\Gamma(F, x_0, \alpha(t_{n_k}), \beta(t_{n_k})) \geq \Gamma(F, C(x_{0}-1), a(t_{n_k}), b(t_{n_k})), 
$$

\noindent which implies, as $k\rightarrow +\infty$,
$$
\exp(Ax_{0}+B) \leq \exp(- C(x_{0}-1)),
$$

\noindent which implies, by letting $C\uparrow +\infty>0$, that $\exp(Ax_{0}+B) \geq +\infty$. This is absurd.\\

\bigskip \noindent Case 3. $\ell_1\in \mathbb{R}_{+}$ and $\ell_2=+\infty$. Fix $x_0>0$. It follows that for any $C>0$ and for any $\varepsilon>0$, there exists $k_0$ such that for any $k\geq k_0$

$$
t_{n_k})/a(t_{n_k} \geq (\ell_{1}-\varepsilon) \ \ and \ \ \beta(t_{n_k})-b(t_{n_k}))/a(t_{n_k} \geq C,
$$

\noindent which implies that, for $k\geq k_0$, we have

$$
\Gamma(F, x_0, \alpha(t_{n_k}), \beta(t_{n_k})) \leq \Gamma(F, ((\ell_{1}-\varepsilon)x_{0}+C), a(t_{n_k}), b(t_{n_k})) 
$$

\noindent which in turn implies for $k\geq k_0$,

$$
\exp(Ax_{0}+B) \leq \exp(-(\ell_{1}-\varepsilon)x_{0}-C),
$$

\noindent which is impossible since left-hand member of the latter inequality is fixed while the right-hand one tends to zero as 
$C\rightarrow +\infty$.\\

\bigskip \noindent Case 4. $\ell_1\in \mathbb{R}_{+}$ and $\ell_2=-\infty$. Fix $x_{0}<0$. It follows that for any $C>0$ and for any $\varepsilon>0$, there exists $k_0$ such that for any $k\geq k_0$

$$
t_{n_k})/a(t_{n_k} \leq (\ell_{1}+\varepsilon) \ \ and \ \ \beta(t_{n_k})-b(t_{n_k}))/a(t_{n_k} \leq - C,
$$

\noindent which implies that, for $k\geq k_0$, we have

$$
\Gamma(F, x_0, \alpha(t_{n_k}), \beta(t_{n_k})) \geq \Gamma(F, ((\ell_{1}+\varepsilon)x_{0}-C), a(t_{n_k}), b(t_{n_k})) 
$$

\noindent which in turn implies for $k\geq k_0$,
$$
\exp(Ax_{0}+B) \geq \exp(-(\ell_{1}+\varepsilon)x_{0}+C),
$$

\noindent which is impossible since left-hand member of the latter inequality is fixed while the right-hand one tends to $+\infty$ as 
$C\rightarrow +\infty$.\\

\bigskip \noindent Since the cases above are excluded, we necessarily have that $\ell_1$ and $\ell_2B$ are finite. Now let us consider the following cases.\\

\bigskip \noindent Case 5. $\ell_1 \neq A$. Then, we have the sub-cases :\\

\noindent (a) $\ell_1 <A$ and $\ell_2\leq A$. Fix $x_{0}=1$. Consider $\varepsilon>0$ such that $\ell_1<A-\varepsilon$. By the same, method, we get a value $k_0 \geq 0$ such that for $k\geq k_0$,

$$
t_{n_k})/a(t_{n_k} \geq (A-\varepsilon) \ \ and \ \ \beta(t_{n_k})-b(t_{n_k}))/a(t_{n_k} \geq B+\varepsilon,
$$

\noindent which leads, $k\geq k_0$, to
$$
\Gamma(F, x_0, \alpha(t_{n_k}), \beta(t_{n_k})) \geq \Gamma(F, (A x_{0}+B-\varepsilon (x_{0}+1)), a(t_{n_k}), b(t_{n_k})), 
$$

\noindent which implies $1 \geq \exp(2\varepsilon)$, which is impossible.\\

\noindent (b) $\ell_1 <A$ and $\ell_2\geq B$. Fix $x_{0}=-1$. Consider $\varepsilon>0$ such that $\ell_1<A-\varepsilon$. By the same, method, we get a value $k_0 \geq 0$ such that for $k\geq k_0$,

$$
\Gamma(F, x_0, \alpha(t_{n_k}), \beta(t_{n_k})) \leq \Gamma(F, (A x_{0}+B-\varepsilon (x_{0}-1)), a(t_{n_k}), b(t_{n_k})) 
$$

\noindent which implies $1 \leq \exp(-2\varepsilon)$, which is impossible.\\

\noindent (c) $\ell_1 >A$ and $\ell_2\geq  B$. Fix $x_{0}=1$. Consider $\varepsilon>0$ such that $\ell_1<A+\varepsilon$. By the same, method, we get a value $k_0 \geq 0$ such that for $k\geq k_0$,
$$
\Gamma(F, x_0, \alpha(t_{n_k}), \beta(t_{n_k})) \leq \Gamma(F, (A x_{0}+B+\varepsilon (x_{0}+1)), a(t_{n_k}), b(t_{n_k})) 
$$

\noindent which implies $1 \leq \exp(-2\varepsilon)$, which is impossible.\\

\noindent (d) $\ell_1 >A$ and $\ell_2\leq B$. Fix $x_{0}=-2$. Consider $\varepsilon>0$ such that $\ell_1<A+\varepsilon$. By the same, method, we get a value $k_0 \geq 0$ such that for $k\geq k_0$,
$$
\Gamma(F, x_0, \alpha(t_{n_k}), \beta(t_{n_k})) \geq \Gamma(F, (A x_{0}+B+\varepsilon (x_{0}+1)), a(t_{n_k}), b(t_{n_k})) 
$$

\noindent which implies $1 \geq \exp(\varepsilon)$, which is impossible.\\

\noindent At this step, we conclude that $\ell_1=A$. Let us suppose that $ell_2 \neq B$ in two cases.\\
 
\bigskip \noindent Case 6. $\ell_2 \neq B$. Then, we have the sub-cases :\\

\noindent (a) $\ell_2 < B$. Fix $x_{0}=1/2$. Consider $\varepsilon>0$ such that $\ell_2<B-\varepsilon$. By the same, method, we get a value $k_0 \geq 0$ such that for $k\geq k_0$,
$$
\Gamma(F, x_0, \alpha(t_{n_k}), \beta(t_{n_k})) \geq \Gamma(F, (A x_{0}+B+\varepsilon (x_{0}-1)), a(t_{n_k}), b(t_{n_k})) 
$$

\noindent which implies $1 \geq \exp(\varepsilon/2)$, which is impossible.\\

\noindent (b) $\ell_2 >B$. Fix $x_{0}=1$. Consider $\varepsilon>0$ such that $\ell_2 > B+\varepsilon$. By the same, method, we get a value $k_0 \geq 0$ such that for $k\geq k_0$,
$$
\Gamma(F, x_0, \alpha(t_{n_k}), \beta(t_{n_k})) \leq \Gamma(F, (A x_{0}+B+\varepsilon (x_{0}+1)), a(t_{n_k}), b(t_{n_k})) 
$$

\noindent which implies $1 \leq \exp(-2\varepsilon)$, which is impossible.\\

\noindent The conclusion of this first step is $\ell_1=A$. Since all adherent points of $\{\alpha(t)/a(t), \ t<uep(F)\}$ as $t\rightarrow uep(F)^{-}$ are equal to $A$. Then $\alpha(t)/a(t) \rightarrow A$, as $t\rightarrow uep(F)^{-}$.\\

\bigskip \noindent \textbf{Step 2.} Let $\ell_2$ be an adherent point of $\{(\beta(t)-b(t))/a(t), \ t<uep(F) \}$ as $t\rightarrow uep(F)^{-}$. Then there exists a sequence $(t_{n})_{n\geq 0}$ such that $t_n \rightarrow uep(F)^{-}$ and
$$
(\beta(t_{n})-b(t_{n}))/a(t_{n}) \rightarrow \ell_2 \ as \ n \rightarrow +\infty.
$$

\noindent But, by our partial conclusion, we also $\alpha(t_n)/a(t_n) \rightarrow \ell_1=A$. By the first step, we have that $\ell_2=B$.\\

\noindent This concludes the proof.$\blacksquare$

\include{04_evt_math_spr_en} 

\part{Appendix}

\chapter{Elements of Theory of Functions and Real Analysis} \label{funct}

\section{Review on limits in $\overline{\mathbb{R}}$. What should not be ignored on limits.} \label{funct.sec.1}

\noindent \textbf{Definition} $\ell \in \overline{\mathbb{R}}$ is an accumulation point of a sequence 
 $(x_{n})_{n\geq 0}$ of real numbers finite or infinitie, in $\overline{\mathbb{R}}$, if and only if there exists a subsequence $(x_{n(k)})_{k\geq 0}$ of
 $(x_{n})_{n\geq 0}$ such that $%
x_{n(k)}$ converges to $\ell $, as $k\rightarrow +\infty $.\newline

\noindent \textbf{Exercise 1 : } Set $y_{n}=\inf_{p\geq n}x_{p}$ and $%
z_{n}=\sup_{p\geq n}x_{p} $ for all $n\geq 0$. Show that :\newline

\noindent \textbf{(1)} $\forall n\geq 0,y_{n}\leq x_{n}\leq z_{n}$\newline

\noindent \textbf{(2)} Justify the existence of the limit of $y_{n}$ called limit inferior of the sequence $(x_{n})_{n\geq 0}$, denoted by $%
\liminf x_{n}$ or $\underline{\lim }$ $x_{n},$ and that it is equal to the following%
\begin{equation*}
\underline{\lim }\text{ }x_{n}=\lim \inf x_{n}=\sup_{n\geq 0}\inf_{p\geq
n}x_{p}
\end{equation*}

\noindent \textbf{(3)} Justify the existence of the limit of $z_{n}$ called limit superior of the sequence $(x_{n})_{n\geq 0}$ denoted by $%
\lim \sup x_{n}$ or $\overline{\lim }$ $x_{n},$ and that it is equal%
\begin{equation*}
\overline{\lim }\text{ }x_{n}=\lim \sup x_{n}=\inf_{n\geq 0}\sup_{p\geq
n}x_{p}x_{p}
\end{equation*}

\bigskip

\noindent \textbf{(4)} Establish that 
\begin{equation*}
-\liminf x_{n}=\limsup (-x_{n})\noindent \text{ \ \ and \ }-\limsup
x_{n}=\liminf (-x_{n}).
\end{equation*}

\newpage \noindent \textbf{(5)} Show that the limit superior is sub-additive and the limit inferior is super-additive, i.e. :  for two sequences
$(s_{n})_{n\geq 0}$ and $(t_{n})_{n\geq 0}$ 
\begin{equation*}
\limsup (s_{n}+t_{n})\leq \limsup s_{n}+\limsup t_{n}
\end{equation*}%
and%
\begin{equation*}
\lim \inf (s_{n}+t_{n})\leq \lim \inf s_{n}+\lim \inf t_{n}
\end{equation*}

\noindent \textbf{(6)} Deduce from (1) that if%
\begin{equation*}
\lim \inf x_{n}=\lim \sup x_{n},
\end{equation*}%
then $(x_{n})_{n\geq 0}$ has a limit and 
\begin{equation*}
\lim x_{n}=\lim \inf x_{n}=\lim \sup x_{n}
\end{equation*}

\bigskip

\noindent \textbf{Exercise 2.} Accumulation points of $%
(x_{n})_{n\geq 0}$.\newline

\noindent \textbf{(a)} Show that if $\ell _{1}$=$\lim \inf x_{n}$ and $\ell
_{2}=\lim \sup x_{n}$ are accumulation points of $(x_{n})_{n\geq 0}.
$ Show one case and deduce the second using point (3) of exercise 1.\newline

\noindent \textbf{(b)} Show that $\ell _{1}$ is the smallest accumulation point of $(x_{n})_{n\geq 0}$ and $\ell _{2}$ is the biggest.
(Similarly, show one case and deduce the second using point (3) of exercise 1).\newline

\noindent \textbf{(c)} Deduce from (a) that if $(x_{n})_{n\geq 0}$ has
a limit $\ell ,$ then it is equal to the unique accumulation point and so,%
\begin{equation*}
\ell =\overline{\lim }\text{ }x_{n}=\lim \sup x_{n}=\inf_{n\geq
0}\sup_{p\geq n}x_{p}.
\end{equation*}

\noindent \textbf{(d)} Combine this this result with point \textbf{(6)} of Exercise 1 to show that a sequence $(x_{n})_{n\geq 0}$ of $\overline{\mathbb{R}}
$ has a limit $\ell $ in $\overline{\mathbb{R}}$ if and only if\ $\lim \inf
x_{n}=\lim \sup x_{n}$ and then%
\begin{equation*}
\ell =\lim x_{n}=\lim \inf x_{n}=\lim \sup x_{n}
\end{equation*}

\newpage

\noindent \textbf{Exercise 3. } Let $(x_{n})_{n\geq 0}$ be a non-decreasing sequence
of $\overline{\mathbb{R}}$. Study its limit superior and its limit inferior and deduce that%
\begin{equation*}
\lim x_{n}=\sup_{n\geq 0}x_{n}.
\end{equation*}

\noindent Dedude that for a non-increasing sequence $(x_{n})_{n\geq 0}$
of $\overline{\mathbb{R}},$%
\begin{equation*}
\lim x_{n}=\inf_{n\geq 0}x_{n}.
\end{equation*}

\bigskip

\noindent \textbf{Point 4.} (Convergence criteria)\newline

\noindent \textbf{Prohorov Criterion} Let $(x_{n})_{n\geq 0}$ be a sequence of $\overline{%
\mathbb{R}}$ and a real number $\ell \in \overline{\mathbb{R}}$ such that: Every subsequence of $(x_{n})_{n\geq 0}$ 
also has a subsequence ( that is a subssubsequence of $(x_{n})_{n\geq 0}$ ) that converges to $\ell .$
Then, the limit of $(x_{n})_{n\geq 0}$ exists and is equal $\ell .$\newline

\noindent \textbf{Upcrossing or Dowcrossing Criterion. } Upcrossings and downcrossings. \newline

\noindent Let $(x_{n})_{n\geq 0}$ be a sequence in $\overline{\mathbb{R}}$ and two real numbers $a$ and $b$ such that $a<b.$
We define%
\begin{equation*}
\nu _{1}=\left\{ 
\begin{array}{cc}
\inf  & \{n\geq 0,x_{n}<a\} \\ 
+\infty  & \text{if (}\forall n\geq 0,x_{n}\geq a\text{)}%
\end{array}%
\right. .
\end{equation*}%
If $\nu _{1}$ is finite, let%
\begin{equation*}
\nu _{2}=\left\{ 
\begin{array}{cc}
\inf  & \{n>\nu _{1},x_{n}>b\} \\ 
+\infty  & \text{if (}n>\nu _{1},x_{n}\leq b\text{)}%
\end{array}%
\right. .
\end{equation*}%
.

\noindent As long as the $\nu _{j}'s$ are finite, we can define for $\nu
_{2k-2}(k\geq 2)$

\begin{equation*}
\nu _{2k-1}=\left\{ 
\begin{array}{cc}
\inf  & \{n>\nu _{2k-2},x_{n}<a\} \\ 
+\infty  & \text{if (}\forall n>\nu _{2k-2},x_{n}\geq a\text{)}%
\end{array}%
\right. .
\end{equation*}%
and for $\nu _{2k-1}$ finite, 
\begin{equation*}
\nu _{2k}=\left\{ 
\begin{array}{cc}
\inf  & \{n>\nu _{2k-1},x_{n}>b\} \\ 
+\infty  & \text{if (}n>\nu _{2k-1},x_{n}\leq b\text{)}%
\end{array}%
\right. .
\end{equation*}

\noindent We stop once one $\nu _{j}$ is $+\infty$. If $\nu
_{2j}$ is finite, then 
\begin{equation*}
x_{\nu _{2j}}-x_{\nu _{2j-1}}>b-a. 
\end{equation*}

\noindent We then say : by that moving from $x_{\nu _{2j-1}}$ to $x_{\nu
_{2j}},$ we have accomplished a crossing (toward the up) of the segment $[a,b]$
called \textit{up-crossings}. Similarly, if one $\nu _{2j+1}$
is finite, then the segment $[x_{\nu _{2j}},x_{\nu _{2j+1}}]$ is a crossing downward (downcrossing) of the segment $[a,b].$ Let%
\begin{equation*}
D(a,b)=\text{ number of upcrossings of the sequence of the segment }[a,b]\text{.}
\end{equation*}

\bigskip

\noindent \textbf{(a)} What is the value of $D(a,b)$ if \ $\nu _{2k}$ is finite and $\nu
_{2k+1}$ infinite.\newline

\noindent \textbf{(b)} What is the value of $D(a,b)$ if \ $\nu _{2k+1}$ is finite and $\nu
_{2k+2}$ infinite.\newline

\noindent \textbf{(c)} What is the value of $D(a,b)$ if \ all the $\nu _{j}'s$ are finite.%
\newline

\noindent \textbf{(d)} Show that $(x_{n})_{n\geq 0}$ has a limit iff
for all $a<b,$ $D(a,b)<\infty.$\newline

\noindent \textbf{(e)} Show that $(x_{n})_{n\geq 0}$ has a limit iff
for all $a<b,$ $(a,b)\in \mathbb{Q}^{2},D(a,b)<\infty .$\newline

\bigskip

\noindent \textbf{Exercise 5. } (Cauchy Criterion). Let $%
(x_{n})_{n\geq 0}$ $\mathbb{R}$ be a sequence of (\textbf{real numbers}).\newline

\noindent \textbf{(a)} Show that if $(x_{n})_{n\geq 0}$ is Cauchy,
then it has a unique accumulation point $\ell \in 
\mathbb{R}$ which is its limit.\newline

\noindent \textbf{(b)} Show that if a sequence $(x_{n})_{n\geq 0}\subset 
\mathbb{R}$ \ converges to $\ell \in \mathbb{R},$ then, it is Cauchy.%
\newline

\noindent \textbf{(c)} Deduce the Cauchy criterion for sequences of real numbers.

\newpage

\begin{center}
\textbf{SOLUTIONS}
\end{center}

\noindent \textbf{Exercise 1}.\newline

\noindent \textbf{Question (1) :}. It is obvious that :%
\begin{equation*}
\underset{p\geq n}{\inf }x_{p}\leq x_{n}\leq \underset{p\geq n}{\sup }x_{p},
\end{equation*}

\noindent since $x_{n}$ is an element of $\left\{
x_{n},x_{n+1},...\right\} $ on which we take the supremum or the infinimum.%
\newline

\noindent \textbf{Question (2) :}. Let $y_{n}=\underset{p\geq 0}{\inf }%
x_{p}=\underset{p\geq n}{\inf }A_{n},$ where $A_{n}=\left\{
x_{n},x_{n+1},...\right\} $ is a non-increasing sequence of sets : $\forall n\geq 0$,
\begin{equation*}
A_{n+1}\subset A_{n}.
\end{equation*}

\noindent So the infinimum on $A_{n}$ increases. If $y_{n}$ increases in $%
\overline{\mathbb{R}},$ its limit is its upper bound, finite or infinite. So%
\begin{equation*}
y_{n}\nearrow \underline{\lim }\text{ }x_{n},
\end{equation*}%
is a finite or infinite number.\newline

\noindent \textbf{Question (3) :}. We also show that $z_{n}=\sup A_{n}$ decreases and $z_{n}\downarrow \overline{\lim }$ $x_{n}$.\newline

\noindent \textbf{Question (4) \label{qst4}:}. We recall that 
\begin{equation*}
-\sup \left\{ x,x\in A\right\} =\inf \left\{ -x,x\in A\right\}. 
\end{equation*}

\noindent Which we write 
\begin{equation*}
-\sup A=\inf -A.
\end{equation*}

\noindent Thus,

\begin{equation*}
-z_{n}=-\sup A_{n}=\inf -A_{n} = \inf \left\{-x_{p},p\geq n\right\}..
\end{equation*}

\noindent The right hand term tends to $-\overline{\lim}\ x_{n}$ and the left hand to $\underline{\lim} \ -x_{n}$ and so 

\begin{equation*}
-\overline{\lim}\ x_{n}=\underline{\lim }\ (-x_{n}).
\end{equation*}

\bigskip \noindent Similarly, we show:
\begin{equation*}
-\underline{\lim } \ (x_{n})=\overline{\lim} \ (-x_{n}).
\end{equation*}

\noindent 

\noindent \textbf{Question (5)}. These properties come from the formulas, where $A\subseteq \mathbb{R},B\subseteq \mathbb{R}$ :%
\begin{equation*}
\sup \left\{ x+y,A\subseteq \mathbb{R},B\subseteq \mathbb{R}\right\} \leq
\sup A+\sup B.
\end{equation*}

\noindent In fact : 
\begin{equation*}
\forall x\in \mathbb{R},x\leq \sup A
\end{equation*}

\noindent and
\begin{equation*}
\forall y\in \mathbb{R},y\leq \sup B.
\end{equation*}

\noindent Thus 
\begin{equation*}
x+y\leq \sup A+\sup B,
\end{equation*}

\noindent where 
\begin{equation*}
\underset{x\in A,y\in B}{\sup }x+y\leq \sup A+\sup B.
\end{equation*}%
Similarly,%
\begin{equation*}
\inf (A+B\geq \inf A+\inf B.
\end{equation*}

\noindent In fact :

\begin{equation*}
\forall (x,y)\in A\times B,x\geq \inf A\text{ and }y\geq \inf B.
\end{equation*}

\noindent Thus 
\begin{equation*}
x+y\geq \inf A+\inf B.
\end{equation*}

\noindent Thus 
\begin{equation*}
\underset{x\in A,y\in B}{\inf }(x+y)\geq \inf A+\inf B
\end{equation*}

\noindent \textbf{Application}.\newline

\begin{equation*}
\underset{p\geq n}{\sup } \ (x_{p}+y_{p})\leq \underset{p\geq n}{\sup } \ x_{p}+\underset{p\geq n}{\sup } \ y_{p}.
\end{equation*}

\noindent All these sequences are non-increasing. Taking infimum, we obtain the limits superior :

\begin{equation*}
\overline{\lim }\text{ }(x_{n}+y_{n})\leq \overline{\lim }\text{ }x_{n}+%
\overline{\lim }\text{ }x_{n}.
\end{equation*}

\bigskip

\noindent \textbf{Question (6) :} Set

\begin{equation*}
\underline{\lim } \ x_{n}=\overline{\lim } \ x_{n},
\end{equation*}

\noindent Since : 
\begin{equation*}
\forall x\geq 1,\text{ }y_{n}\leq x_{n}\leq z_{n},
\end{equation*}%

\begin{equation*}
y_{n}\rightarrow \underline{\lim} \ x_{n}
\end{equation*}%

\noindent and 

\begin{equation*}
z_{n}\rightarrow \overline{\lim } \ x_{n},
\end{equation*}

\noindent we apply Sandwich Theorem to conclude that the limit of $x_{n}$ exists and :

\begin{equation*}
\lim \text{ }x_{n}=\underline{\lim }\text{ }x_{n}=\overline{\lim }\text{ }%
x_{n}.
\end{equation*}

\bigskip 
\noindent \textbf{Exercice 2}.\newline

\noindent \textbf{Question (a).}\\

\noindent Thanks to question (4) of exercise 1, it suffices to show this property for one of the limits. Consider the limit superior and the three cases:\\

\noindent \textbf{The case of a finite limit superior} :

\begin{equation*}
\underline{\lim \text{ }}x_{n}=\ell \text{ finite.}
\end{equation*}

\noindent By definition, 
\begin{equation*}
z_{n}=\underset{p\geq n}{\sup }x_{p}\downarrow \ell .
\end{equation*}

\noindent So: 
\begin{equation*}
\forall \varepsilon >0,\exists (N(\varepsilon )\geq 1),\forall p\geq
N(\varepsilon ),\ell -\varepsilon <x_{p}\leq \ell +\varepsilon .
\end{equation*}

\noindent Take less than that:

\begin{equation*}
\forall \varepsilon >0,\exists n_{\varepsilon }\geq 1:\ell -\varepsilon
<x_{n_{\varepsilon }}\leq \ell +\varepsilon.
\end{equation*}

\noindent We shall construct a subsequence converging to $\ell$.\\

\noindent Let $\varepsilon =1:$%
\begin{equation*}
\exists N_{1}:\ell -1<x_{N_{1}}=\underset{p\geq n}{\sup }x_{p}\leq \ell +1.
\end{equation*}

\noindent But if 
\begin{equation}
z_{N_{1}}=\underset{p\geq n}{\sup }x_{p}>\ell -1, \label{cc}
\end{equation}

\noindent there surely exists an $n_{1}\geq N_{1}$ such that%
\begin{equation*}
x_{n_{1}}>\ell -1.
\end{equation*}

\noindent if not we would have 
\begin{equation*}
( \forall p\geq N_{1},x_{p}\leq \ell -1\ ) \Longrightarrow \sup \left\{
x_{p},p\geq N_{1}\right\} =z_{N_{1}}\geq \ell -1,
\end{equation*}
which is contradictory with (\ref{cc}). So, there exists $n_{1}\geq N_{1}$ such that
\begin{equation*}
\ell -1<x_{n_{1}}\leq \underset{p\geq N_{1}}{\sup }x_{p}\leq \ell -1.
\end{equation*}

\noindent i.e.

\begin{equation*}
\ell -1<x_{n_{1}}\leq \ell +1.
\end{equation*}

\noindent We move to step $\varepsilon =\frac{1}{2}$ and we consider the sequence%
 $(z_{n})_{n\geq n_{1}}$ whose limit remains $\ell$. So, there exists $N_{2}>n_{1}:$%
\begin{equation*}
\ell -\frac{1}{2}<z_{N_{2}}\leq \ell -\frac{1}{2}.
\end{equation*}

\noindent We deduce like previously that $n_{2}\geq N_{2}$ such that%
\begin{equation*}
\ell -\frac{1}{2}<x_{n_{2}}\leq \ell +\frac{1}{2}
\end{equation*}

\noindent with $n_{2}\geq N_{1}>n_{1}$.\\

\noindent Next, we set $\varepsilon =1/3,$ there will exist $N_{3}>n_{2}$ such that%
\begin{equation*}
\ell -\frac{1}{3}<z_{N_{3}}\leq \ell -\frac{1}{3}
\end{equation*}

\noindent and we could find an $n_{3}\geq N_{3}$ such that%

\begin{equation*}
\ell -\frac{1}{3}<x_{n_{3}}\leq \ell -\frac{1}{3}.
\end{equation*}

\noindent Step by step, we deduce the existence of $%
x_{n_{1}},x_{n_{2}},x_{n_{3}},...,x_{n_{k}},...$ with $n_{1}<n_{2}<n_{3}%
\,<...<n_{k}<n_{k+1}<...$ such that

$$
\forall k\geq 1, \ell -\frac{1}{k}<x_{n_{k}}\leq \ell -\frac{1}{k},
$$

\noindent i.e.

\begin{equation*}
\left\vert \ell -x_{n_{k}}\right\vert \leq \frac{1}{k}.
\end{equation*}

\noindent Which will imply: 
\begin{equation*}
x_{n_{k}}\rightarrow \ell 
\end{equation*}

\noindent Conclusion : $(x_{n_{k}})_{k\geq 1}$ is very well a subsequence since $n_{k}<n_{k+1}$ for all $k \geq 1$ 
and it converges to $\ell$, which is then an accumulation point.\\

\noindent \textbf{Case of the limit superior equal $+\infty$} : 
$$
\overline{\lim} \text{ } x_{n}=+\infty.
$$
\noindent Since $z_{n}\uparrow +\infty ,$ we have : $\forall k\geq 1,\exists
N_{k}\geq 1,$ 
\begin{equation*}
z_{N_{k}}\geq k+1.
\end{equation*}

\noindent For $k=1$, let $z_{N_{1}}=\underset{p\geq N_{1}}{\inf }%
x_{p}\geq 1+1=2.$ So there exists 
\begin{equation*}
n_{1}\geq N_{1}
\end{equation*}%
such that :%
\begin{equation*}
x_{n_{1}}\geq 1.
\end{equation*}

\noindent For $k=2:$ consider the sequence $(z_{n})_{n\geq n_{1}+1}.$
We find in the same manner 
\begin{equation*}
n_2 \geq n_{1}+1
\end{equation*}%
\noindent and 
\begin{equation*}
x_{n_{2}}\geq 2.
\end{equation*}

\noindent Step by step, we find for all $k\geq 3$, an $n_{k}\geq n_{k-1}+1$ such that
\begin{equation*}
x_{n_{k}}\geq k.
\end{equation*}

\noindent Which leads to $x_{n_{k}}\rightarrow +\infty $ as $k\rightarrow +\infty $.\\

\noindent \textbf{Case of the limit superior equal $-\infty$} : 

$$
\overline{\lim }x_{n}=-\infty.
$$

\noindent This implies : $\forall k\geq 1,\exists N_{k}\geq 1,$ such that%
\begin{equation*}
z_{n_{k}}\leq -k.
\end{equation*}

\noindent For $k=1,\exists n_{1}$ such that%
\begin{equation*}
z_{n_{1}}\leq -1.
\end{equation*}
But 
\begin{equation*}
x_{n_{1}}\leq z_{n_{1}}\leq -1
\end{equation*}

\noindent Let $k=2$. Consider $\left( z_{n}\right) _{n\geq
n_{1}+1}\downarrow -\infty .$ There will exist $n_{2}\geq n_{1}+1:$%
\begin{equation*}
x_{n_{2}}\leq z_{n_{2}}\leq -2
\end{equation*}

\noindent Step by step, we find $n_{k1}<n_{k+1}$ in such a way that $x_{n_{k}}<-k$ for all $k$ bigger that $1$. So%
\begin{equation*}
x_{n_{k}}\rightarrow +\infty 
\end{equation*}

\bigskip

\noindent \textbf{Question (b).}\\

\noindent Let $\ell$ be an accumulation point of $(x_n)_{n \geq 1}$, the limit of one of its subsequences $(x_{n_{k}})_{k \geq 1}$. We have

$$
y_{n_{k}}=\inf_{p\geq n_k} \ x_p \leq x_{n_{k}} \leq  \sup_{p\geq n_k} \ x_p=z_{n_{k}}
$$

\noindent The left hand side term is a subsequence of $(y_n)$ tending to the limit inferior and the right hand side is a 
subsequence of $(z_n)$ tending to the limit superior. So we will have:

$$
\underline{\lim} \ x_{n} \leq \ell \leq \overline{\lim } \ x_{n},
$$

\noindent which shows that $\underline{\lim} \ x_{n}$ is the smallest accumulation point and $\overline{\lim } \ x_{n}$ is the largest.\\

\noindent \textbf{Question (c).} If the sequence $(x_n)_{n \geq 1}$ has a limit $\ell$, it is the limit of all its subsequences,
so subsequences tending to the limits superior and inferior. Which answers question (b).\\

\noindent \textbf{Question (d).} We answer this question by combining point (d) of this exercise and point (\textbf{6}) of the exercise \textbf{1}.\\

\noindent \textbf{Exercise 3}. Let $(x_{n})_{n\geq 0}$ be a non-decreasing sequence, we have:%
\begin{equation*}
z_{n}=\underset{p\geq n}{\sup} \ x_{p}=\underset{p\geq 0}{\sup} \ x_{p},\forall
n\geq 0.
\end{equation*}

\noindent Why? Because by increasingness,%
\begin{equation*}
\left\{ x_{p},p\geq 0\right\} =\left\{ x_{p},0\leq p\leq n-1\right\} \cup
\left\{ x_{p},p\geq n\right\}
\end{equation*}

\bigskip

\noindent Since all the elements of $\left\{ x_{p},0\leq p\leq
n-1\right\} $ are smaller than that of $\left\{ x_{p},p\geq n\right\} ,$
the supremum is achieved on $\left\{ x_{p},p\geq n\right\} $ and so 
\begin{equation*}
\ell =\underset{p\geq 0}{\sup } \ x_{p}=\underset{p\geq n}{\sup }x_{p}=z_{n}
\end{equation*}%
Thus%
\begin{equation*}
z_{n}=\ell \rightarrow \ell .
\end{equation*}

\noindent We also have $y_n=\inf \left\{ x_{p},0\leq p\leq n\right\}=x_n$ which is a non-decreasing sequence and so converges to
$\ell =\underset{p\geq 0}{\sup } \ x_{p}$. \\

\bigskip

\noindent \textbf{Exercise 4}.\\

\noindent Let $\ell \in \overline{\mathbb{R}}$ having the indicated property. Let $\ell ^{\prime }$ be a given accumulation point.%
\begin{equation*}
 \left( x_{n_{k}}\right)_{k\geq 1} \subseteq \left( x_{n}\right) _{n\geq 0}%
\text{ such that }x_{n_{K}}\rightarrow \ell ^{\prime}.
\end{equation*}

\noindent By hypothesis this subsequence $\left( x_{n_{K}}\right) $
has in turn a subsubsequence $\left( x_{n_{\left( k(p)\right) }}\right)_{p\geq 1} $ such that $x_{n_{\left( k(p)\right) }}\rightarrow
\ell $ as $p\rightarrow +\infty $.\newline

\noindent But as a subsequence of $\left( x_{n_{\left( k\right)
}}\right) ,$ 
\begin{equation*}
x_{n_{\left( k(\ell )\right) }}\rightarrow \ell ^{\prime }.
\end{equation*}%
Thus
\begin{equation*}
\ell =\ell ^{\prime}.
\end{equation*}

\noindent Applying that to the limit superior and limit inferior, we have:%
\begin{equation*}
\overline{\lim} \ x_{n}=\underline{\lim}\ x_{n}=\ell.
\end{equation*}

\noindent And so $\lim x_{n}$ exists and equals $\ell$.\\

\noindent \textbf{Exercise 5}.\\

\noindent \textbf{Question (a)}. If $\nu _{2k}$ finite and $\nu _{2k+1}$ infinite, it then has exactly $k$ up-crossings : 
$[x_{\nu_{2j-1}},x_{\nu _{2j}}]$, $j=1,...,k$ : $D(a,b)=k$.\\

\noindent \textbf{Question (b)}. If $\nu _{2k+1}$ finite and $\nu _{2k+2}$ infinite, it then has exactly $k$ up-crossings:
$[x_{\nu_{2j-1}},x_{\nu_{2j}}]$, $j=1,...,k$ : $D(a,b)=k$.\\

\noindent \textbf{Question (c)}. If all the $\nu_{j}'s$ are finite, then, there are an infinite number of up-crossings : 
$[x_{\nu_{2j-1}},x_{\nu_{2j}}]$, $j\geq 1k$ : $D(a,b)=+\infty$.\\

\noindent \textbf{Question (d)}. Suppose that there exist $a < b$ rationals such that $D(a,b)=+\infty$. 
Then all the $\nu _{j}'s$ are finite. The subsequence $x_{\nu_{2j-1}}$ is strictly below $a$. 
So its limit inferior is below $a$. This limit inferior is an accumulation point of the sequence $(x_n)_{n\geq 1}$, 
so is more than $\underline{\lim}\ x_{n}$, which is below $a$.\\

\noindent Similarly, the subsequence $x_{\nu_{2j}}$ is strictly below $b$. So the limit superior is above $a$. 
This limit superior is an accumulation point of the sequence $(x_n)_{n\geq 1}$, so it is below $\overline{\lim}\ x_{n}$, 
which is directly above $b$. Which leads to:

$$
\underline{\lim}\ x_{n} \leq a < b \leq \overline{\lim}\ x_{n}. 
$$

\noindent That implies that the limit of $(x_n)$ does not exist. In contrary, we just proved that the limit of $(x_n)$ exists, 
meanwhile for all the real numbers $a$ and $b$ such that $a<b$, $D(a,b)$ is finite.\\

\noindent Now, suppose that the limit of $(x_n)$ does not exist. Then,

$$
\underline{\lim}\ x_{n} < \overline{\lim}\ x_{n}. 
$$

\noindent We can then find two rationals $a$ and $b$ such that $a<b$ and a number $\epsilon$ such that $0<\epsilon$, all the

$$
\underline{\lim}\ x_{n} < a-\epsilon < a < b < b+\epsilon <  \overline{\lim}\ x_{n}. 
$$

\noindent If $\underline{\lim}\ x_{n} < a-\epsilon$, we can return to question \textbf{(a)} of exercise \textbf{2} and construct a subsequence of $(x_n)$
which tends to $\underline{\lim}\ x_{n}$ while remaining below $a-\epsilon$. Similarly, if $b+\epsilon < \overline{\lim}\ x_{n}$, 
we can create a subsequence of $(x_n)$ which tends to $\overline{\lim}\ x_{n}$ while staying above $b+\epsilon$. 
It is evident with these two sequences that we could define with these two sequences all $\nu_{j}$ finite and so $D(a,b)=+\infty$.\\

\noindent We have just shown by contradiction that if all the $D(a,b)$ are finite for all rationals $a$ and $b$ such that $a<b$, 
then, the limit of $(x)n)$ exists.\\

\noindent \textbf{Exercise 5}. Cauchy criterion in $\mathbb{R}$.\\

\noindent Suppose that the sequence is Cauchy, $i.e.$,
$$
\lim_{(p,q)\rightarrow (+\infty,+\infty)} \ (x_p-x_q)=0.
$$

\noindent Then let $x_{n_{k,1}}$ and $x_{n_{k,2}}$ be two subsequences converging respectively to $\ell_1=\underline{\lim}\ x_{n}$ and $\ell_2=\overline{\lim}\ x_{n}$. So

$$
\lim_{(p,q)\rightarrow (+\infty,+\infty)} \ (x_{n_{p,1}}-x_{n_{q,2}})=0.
$$

\noindent, By first letting $p\rightarrow +\infty$, we have

$$
\lim_{q\rightarrow +\infty} \ \ell_1-x_{n_{q,2}}=0.
$$

\noindent Which shows that $\ell_1$ is finite, else $\ell_1-x_{n_{q,2}}$ would remain infinite and would not tend to $0$. 
By interchanging the roles of $p$ and $q$, we also have that $\ell_2$ is finite.\\

\noindent Finally, by letting $q\rightarrow +\infty$, in the last equation, we obtain
$$
\ell_1=\underline{\lim}\ x_{n}=\overline{\lim}\ x_{n}=\ell_2.
$$

\noindent which proves the existence of the finite limit of the sequence $(x_n)$.\\

\noindent Now suppose that the finite limit $\ell$ of $(x_n)$ exists. Then

$$
\lim_{(p,q)\rightarrow (+\infty,+\infty)} \ (x_p-x_q)=\ell-\ell=0.
$$0
\noindent Which shows that the sequence is Cauchy.\\

\section{Topology and measure theory complements}

\noindent In this section, we recall some relations on measurability of real-valued applications defined on intervals of 
$\mathbb{R}$ and their continuity. We begin by the seminal Egoroff result.

\begin{theorem} \label{funct.topo.theo1} (Egoroff's Theorem) Let $E$ be a Borel subset of $\mathbb{R}$ such that
$0<\lambda (E)<+\infty$, where $\lambda$ is the Lebesgues measure on $\mathbb{R}$. Let $f_{n}:E\longmapsto \mathbb{R}$ be a sequence of finite et
measurable functions with respect to the Borel $\sigma$-algebras $\mathbb{B}(E)$ and 
$\mathbb{B}(\mathbb{R})$.\\

\noindent Suppose that this sequence converges $a.s$ to a finite function $f$. Let $\varepsilon >0$ an arbitrary positive real
number. Then, there exists a measurable subset $A$ of $E$\ such that $%
\lambda (E\setminus A)<\varepsilon $ and $f_{n}$ converge uniformly to $f$ on $A.$
\end{theorem}

\bigskip \noindent \textbf{Proof of Theorem \ref{funct.topo.theo1}}. By definition, $f_{n}$ converges to $f$ $a.s$ if and only if 
\begin{equation*}
\lambda (\{x\in E,f_{n}(x)\text{ does not converges to }f(x)\})=0.
\end{equation*}

\bigskip \noindent Since the $f_{n}$'s and $f$ are finite, we have 
\begin{equation*}
\{x\in E,f_{n}(x) \nrightarrow f(x)\}=\bigcup\limits_{\eta
>0}\bigcap\limits_{n=1}^{+\infty }\bigcup\limits_{m=n}^{+\infty }\{x\in
E,\left\vert f_{m}(x)-f(x)\right\vert \geq \eta \},
\end{equation*}

\noindent $\nrightarrow$ stands for : \textit{does not converge to}. This can be achieved on real numbers $\eta =1/k,k\geq 1$ so that 
\begin{equation*}
\{x\in E,f_{n}(x) \nrightarrow f(x)\}=\bigcup\limits_{k>0}%
\bigcap\limits_{n=1}^{+\infty }\bigcup\limits_{m=n}^{+\infty }\{x\in
E,\left\vert f_{m}(x)-f(x)\right\vert \geq 1/k\}.
\end{equation*}

\noindent So $f_{n}$ converges to $f$ $a.s$ if and only if for any $k\geq 1,$%
\begin{equation*}
\lambda (\bigcap\limits_{n=1}^{+\infty }\bigcup\limits_{m=n}^{+\infty
}\{x\in E,\left\vert f_{m}(x)-f(x)\right\vert \geq 1/k\})=0.
\end{equation*}

\noindent Put for any $k\geq 1$a and $n\geq 1$,
\begin{equation*}
A_{k,n}=\bigcup\limits_{m=n}^{+\infty }\{x\in E,\left\vert
f_{m}(x)-f-x)\right\vert \geq 1/k\}.
\end{equation*}

\noindent For a fixed $k\geq 1,$ ($A_{k,n})_{n\geq 1}$ is a nondecreasing sequence (in $n$) of measurable sets with finite measures (their measures are bounded above by $\lambda (E)$) with limit

\begin{equation*}
A_{k}=\bigcap\limits_{n=1}^{+\infty }\bigcup\limits_{m=n}^{+\infty }\{x\in
E,\left\vert f_{m}(x)-f(x)\right\vert \geq 1/k\}.
\end{equation*}

\noindent Then, we have
\begin{equation}
\lambda (A_{k,n})\downarrow \lambda (A_{k})=0.  \label{funct.topo.001}
\end{equation}

\noindent \noindent Now let $\varepsilon >0$ be fixed. For each fixed $k\geq 1$, we can use \ref{funct.topo.001}
and find an indice $n(k)$ such that
\begin{equation*}
\lambda (A_{k,n(k)})<2^{-k}\varepsilon /2.
\end{equation*}

\noindent Set
\begin{equation*}
A=\bigcap\limits_{k\geq 1}A_{k,n(k)}^{c}.
\end{equation*}

\noindent wher the compléments are taken in $E$. Then we have
\begin{equation*}
E\setminus A=\bigcup\limits_{k\geq 1}A_{k,n(k)}.
\end{equation*}

\noindent We surely have $A\subset E$ and by the $\sigma$-subadditivity of the measure $\lambda$, we see that
\begin{equation*}
\lambda (E\setminus A)=\lambda (\bigcup\limits_{k\geq 1}A_{k,n(k)})\leq
\sum_{k=1}^{+\infty }\lambda (A_{k,n(k)})\leq \varepsilon /2<\varepsilon .
\end{equation*}

\noindent To finish, let us show that $f_{n}$ converges to $f$ uniformly on $A$. Let $%
\eta >0$ and consider a value $k\geq 1$ such that $1/k_{0}<\eta $ and set $%
N=n(k_{0}).$ Now any $x\in A$ belongs to $any$ of the $A_{k,n(k)}$. Then for any $x\in A$, $x$ belongs to $A_{k_{0},n(k_{0})}$, and by this, we have 
\begin{equation*}
m\geq N=n(k_{0})\Longrightarrow \left\vert f_{m}(x)-f(x)\right\vert
<1/k_{0}<\eta .
\end{equation*}

\noindent We just proved this :

\begin{equation*}
\forall \eta >0,\exists N\geq 0,\sup_{x\in A}\sup_{n\geq N}\left\vert
f_{m}(x)-f(x)\right\vert \leq <1/k_{0}<\eta .
\end{equation*}

\noindent We conclude that $f_{n}$ converges to $f$ uniformly on $A.$\\

\bigskip We are moving now to the Lusin Theorem that will bring us back to outer measures.\\

\begin{theorem} \label{funct.topo.theo2}
(Lusin's Theorem) Let $f:\mathbb{R}\longmapsto \mathbb{R}$ be a finite and measurable
function with respect to the Borel $\sigma $-algebras $B(\mathbb{R})$ and $%
\lambda $ be the Lebesgues measure on $\mathbb{R}$. Let $\varepsilon >0$ be an
arbitrary positive real number. Then, there exists a measurable set E\ such
that $\lambda (A^{c})<\varepsilon $ and $f$ is continuous on $A.$
\end{theorem}

\bigskip \noindent To prove the Lusin'sTheorem, we need the following. The first lemma will oblige
us to go back the proof of Theorem of Caratheodory.

\bigskip

\begin{lemma} \label{funct.topo.lem1} Let be a Borel subset $A$ of $\mathbb{R}$ such that its Lebesgues
measure is finite, that is $\lambda (E)<+\infty$. Then, for any $\varepsilon
>0,$ there exists a finite union $K$ of bounded intervals such that $\lambda
(E\Delta K)<\varepsilon .$
\end{lemma}

\bigskip \noindent \textbf{Proof of Lemma \ref{funct.topo.lem1}}. We have to go back to the proof of the Theorem of Caratheodory which
justified the existence which the Lebesgues measure.\\

\noindent  Recall that the Lebesgues measure is uniquely defined on the class of intervals
\begin{equation*}
\mathcal{I}=\{]a,b],-\infty \leq a\leq b<+\infty \}
\end{equation*}

\noindent by
\begin{equation*}
\lambda (]a,b])=b-a.
\end{equation*}

\noindent This class $\mathcal{I}$ is a semi-algebra and the algebra $\mathcal{C}$\
its generates is the class of finite sums of disjoints intervals. It is
easily proved that $\lambda $ is additive on $\mathcal{I}$ and is readily
extensible to a an additive application on $\mathcal{C}$, that we always denote by $\lambda$. The broad extension of $\lambda $ to a measure 
on a $\sigma$-algebra $\mathcal{A}^{0}$ including $\mathcal{C}$ may be done by the method of the outer measure, defined as follows, for any subset
 $A$ of $\mathbb{R}$ :
\begin{equation}
\lambda ^{0}(A)=\inf \{\sum_{n=0}^{\infty }\lambda (A_{n}),\text{ }A\subset
\bigcup\limits_{n=0}^{\infty }A_{n},A_{n}\in \mathcal{C}\}.  \label{funct.topo.002}
\end{equation}

\noindent A subset of $\mathbb{R}$ is $\lambda ^{0}$-mesurable if and only for any subset$D$ of $\mathbb{R}$, we have
\begin{equation*}
\lambda ^{0}(A)=\lambda ^{0}(AD)+\lambda ^{0}(AD^{c}).
\end{equation*}

\noindent By denoting $\mathcal{A}^{0}$ \ the set of \ $\lambda ^{0}$-measurable
subsets of $\mathbb{R}$, is proved that $(\Omega ,\mathcal{A}^{0},\lambda ^{0})$ is
a measurable space and $\mathcal{C} \subset \mathcal{A}^{0}.$ This measurable space,
surely, includes the measurable space $(\Omega ,\mathcal{B}(\mathbb{R}),\lambda ^{0})$
since $\mathcal{B}(R)=\sigma (\mathcal{C})\subset \mathcal{A}^{0}.$ The measure $%
\lambda ^{0}$ is the unique extension of $\lambda $ to a measure on $\mathcal{B}(\mathbb{R})$, still
denoted by $\lambda$.\\

\noindent By using \ref{funct.topo.002} and the characterization of the infimum based on the fact that $\lambda (E)$ is finite, we conclude that for any $\varepsilon>0,$ there exists a union  $\bigcup_{n\geq 1} A_{n}$, formed by elements $A_n$ in $\mathcal{C}$ and covering $E$, such that we have 
\begin{equation}
\sum_{n=0}^{\infty }\lambda (A_{n})<\lambda (E)+\varepsilon /2.  \label{funct.topo.003}
\end{equation}

\noindent But for each $n\geq 1$,  $A_{n}$ is a finite sum of intervals - elements of $\mathcal{I}$ -  of the form 
$$
A_{n}=\sum_{1\leq j\leq p(n)}I_{n,j},
$$ 

\noindent with 
$$
\lambda (A_{n})=\sum_{1\leq j\leq p(n)}\lambda (I_{n,j}),
$$

\noindent since $\lambda $ is additive on $\mathcal{C}$. We may rephrase this by saying that 
for any $\varepsilon >0,$ there exists a union of intervals $I_{k}$, covering $E$, such that 
\begin{equation*}
\sum_{k=1}^{\infty }\lambda (I_{k})<\lambda (E)+\varepsilon /2.
\end{equation*}

\noindent Here, we see that none of these intervals is unbounded. Otherwise, for one
them, we would have  $\lambda (I_{k})=+\infty $ and Formula (\ref{funct.topo.003}) would be
impossible. Since $\sum_{k=1}^{\infty }\lambda (I_{k})$ is finite, we may
find $k_{0}$ such that
\begin{equation}
0\leq \ \sum_{k=1}^{\infty }\lambda (I_{k})-\ \sum_{k=1}^{k_{0}}\lambda
(I_{k})=\sum_{k\geq k_{0}+1}\lambda (I_{k})<\varepsilon /2.  \label{8.004}
\end{equation}

\noindent Put $K=\cup _{1\leq j\leq k_{0}} I_j$. We finally have
\begin{equation*}
K\backslash E=KE^{c}\subset \left( \bigcup\limits_{k=0}^{\infty
}I_{k}\right) E^{c}=\left( \bigcup\limits_{k=0}^{\infty }I_{k}\right)
\backslash E
\end{equation*}

\noindent with
\begin{equation*}
E\subset \left( \bigcup\limits_{k=0}^{\infty }I_{k}\right) =\left(
\bigcup\limits_{n=0}^{\infty }A_{n}\right).
\end{equation*}

\noindent Then, by \ref{funct.topo.002}, we have
\begin{equation*}
\lambda (K\backslash E)\leq \lambda \left( \bigcup\limits_{k=0}^{\infty
}I_{k}\right) -\lambda (E)\leq \left( \sum_{k=0}^{\infty }\lambda
(A_{k})\right) -\lambda (E)<\varepsilon /2.
\end{equation*}

\noindent Next, by denoting $J=\bigcup\limits_{k=0}^{\infty }I_{k},$ we have $E\subset J$ and 
\begin{equation*}
K^{c}=\left( \bigcup\limits_{k=0}^{\infty }I_{k}\backslash
\bigcup\limits_{k=0}^{k_{0}}I_{k}\right) +J^{c}
\end{equation*}

\noindent and then
\begin{equation*}
E\backslash K=EK^{c}=E\left( \bigcup\limits_{k=0}^{\infty }I_{k}\backslash
\bigcup\limits_{k=0}^{k_{0}}I_{k}\right) +EJ^{c}=E\left(
\bigcup\limits_{k=0}^{\infty }I_{k}\backslash
\bigcup\limits_{k=0}^{k_{0}}I_{k}\right) \subset
\bigcup\limits_{k_{0}+1}^{\infty }I_{k}
\end{equation*}

\noindent and
\begin{equation*}
\lambda (E\backslash K)\leq \lambda \left( \bigcup\limits_{k_{0}+1}^{\infty
}I_{k}\right) \leq \sum_{k\geq k_{0}+1}\lambda (I_{k})<\varepsilon /2.
\end{equation*}

\noindent We conclude that
\begin{equation*}
\lambda (E\Delta K)<\varepsilon .
\end{equation*}%


\bigskip \noindent  We are now able to prove Lusin's Theorem.\\

\noindent \textbf{Proof of Theorem \ref{funct.topo.theo2}}.\\

\noindent Let $\varepsilon >0.$ Fix an arbitrary $N\in \mathbb{Z}$ and put $A_{N}=(N,N+1),
$ $.$ Since $\lambda (A_{N})=1<\infty ,$ Egoroff's theorem ensures we can
find $C_{N}\subset A_{N}$ such that such that $\lambda (A_{N}\Delta C_{N})<2^{N+1}\varepsilon $ and 
\begin{equation*}
f_{n}1_{C_{N}}\rightarrow 1_{C_{N}}
\end{equation*}

\noindent uniformly. It is important that anything happens in the sets $A_{N}$ in the sequel.
In particular, the complements are meant within the sets $A_{N}.$ Since $f_{n}1_{C_{N}}$
is measurable, there exists a sequence of simple functions 
\begin{equation}
f_{n}1_{C_{N}}=\sum_{j=1}^{m(n)}\alpha _{j,n}1_{A_{n,j}} \label{funct.topo.simple}
\end{equation}

\noindent such that $f_{n}1_{C_{N}}\rightarrow f_{n}1_{C_{N}}$. We consider in (\ref{funct.topo.simple}) canonical simple
functions, i.e., we have $A_{n,1}+...+A_{n,m(n)}=C_{N}$, and the values $\alpha _{j,n}$ are finite and distinct between them.\\

\noindent By using Lemma \ref{funct.topo.lem1}, we may find for any $n\geq 1$ and for each $1\leq j\leq 1$, a finite union
of closed bounded intervals $K_{n,j}\subset A_{N}$ such that 
\begin{equation*}
\lambda (A_{n,j}\Delta K_{n,j})<2^{N+n+1}\varepsilon /n.
\end{equation*}

\noindent Put
\begin{equation*}
\widetilde{f}_{n}=\sum_{j=1}^{m(n)}\alpha _{j,n}1_{K_{n,j}}.
\end{equation*}

\noindent Remark that
\begin{eqnarray*}
\{x\in C_{N},\widetilde{f}_{n}(x)&=&f_{n}(x)\}=\sum_{j=1}^{m(n)}\{x\in C_{N},\widetilde{f}_{n}(x)=f_{n}(x)\}\cap A_{n,j}\\
&=&\sum_{j=1}^{m(n)}\{x\in C_{N},x\in A_{n,j}\cap K_{n,j}\}\cap A_{n,j}
\end{eqnarray*}

\noindent So, we have
\begin{equation*}
\{x\in C_{N},\widetilde{f}_{n}(x)\neq f_{n}(x)\}\cap A_{n,j}=A_{n,j}\cap
(A_{n,j}\cap K_{n,j})^{c}=A_{n,j}\backslash (A_{n,j}\cap K_{n,j}).
\end{equation*}

\noindent By summing over $j$, we get
\begin{equation*}
\{x\in C_{N},\widetilde{f}_{n}(x)\neq
f_{n}(x)\}=\bigcup\limits_{j=1}A_{n,j}\backslash (A_{n,j}\cap K_{n,j}).
\end{equation*}

\noindent So, we have on $\widetilde{f}_{n}=f_{n}$ : 
\begin{equation*}
B_{N}=\bigcap\limits_{n=1}^{\infty }\bigcap\limits_{j=1}^{m(n)}\left\{
(A_{n,j}\backslash (A_{n,j}\cap K_{n,j}))^{c}\right\} \cap C_{N}^{c}.
\end{equation*}

\noindent with 
\begin{eqnarray*}
\lambda (B_{N}^{c})&\leq& \sum\limits_{n=1}^{\infty }\sum_{j=1}^{m(n)}\lambda
(A_{n,j}\backslash (A_{n,j}\cap K_{n,j}))\cup C_{N}\\
&\leq & \sum\limits_{n=1}^{\infty }\sum_{j=1}^{m(n)}\lambda ((A_{n,j}\cap
K_{n,j})+\lambda (C_{N})\\
&\leq & 2^{N+1}\varepsilon +2^{N+1}\varepsilon\\
&=&2^{N}\varepsilon.
\end{eqnarray*}

\noindent \noindent We extend $\widetilde{f}_{n}$ to 
\begin{equation*}
\left( \sum_{j=1}^{m(n)}K_{n,j}\right) ^{c}
\end{equation*}

\noindent by assigning the value zero to points outside of 
$$
K_{n}=\sum_{j=1}^{m(n)}K_{n,j}.
$$

\noindent Next, each $K_{n,j}$ can be put into a sum of bounded intervals. Indeed, we have
\begin{eqnarray*}
K_{n,j}&=&\cup _{1\leq \ell \ell p(n,j)}K_{n,j,\ell}\\
&=&K_{n,j,1}+K_{n,j,1}^{c}K_{n,j,2}+...+K_{n,j,1}^{c}K_{n,j,2}^{c}...K_{n,j,p-1(n,j)}^{c}K_{n,j,p(n,j)}.
\end{eqnarray*}

\noindent The claim is true since  each $K_{n,j,\ell }^{c}$ is also a finite sum of
intervals and the class of finite sums of interavls is stable by intersection. Finally, $\widetilde{f}_{n}$ may be expressed in the form
\begin{equation*}
\widetilde{f}_{n}=\sum_{j=1}^{q(n)}\beta _{j,n}1_{B_{n,j}}
\end{equation*}

\noindent where the $B_{n,j}$ are bounded intervals.  The discontinuity points of the
functions $\widetilde{f}_{n}$ consists of some of the possibly closing ends
of the bounded intervals $B_{n,j}.$ The set $D$ of all discontinuity points
of the $\widetilde{f}_{n}$'s, $n\geq 1$, is atmost countable. And, we have on $B_{N}\cap
D^{c}$, that the functions are continuous and converge uniformly to $f$
whereas

\begin{equation*}
\lambda ((B_{N}\cap D^{c})^{c})\leq 2^{-N}\varepsilon .
\end{equation*}%

\noindent We conlude that $f$ is continuous on $B_{N}^{\ast }=B_{N}\cap D^{c}$. Now 
\begin{equation*}
B=\left( \bigcup\limits_{N=1}^{\infty }B_{N}^{\ast }\right) \cap \mathbb{Z}.
\end{equation*}

\noindent Then $f$ is continuous on $B$ with
\begin{equation*}
\lambda (B)\leq \varepsilon .
\end{equation*}

\newpage
\section{Hamel Equations} \label{funct.sec.2}
\begin{lemma} \label{funct.hamel.lem.1} Let $h:$ $\mathbb{R}\longmapsto \mathbb{R}$ be a measurable function such
that 
\begin{equation*}
\forall (x,y)\in \mathbb{R}^{2},h(x+y)=h(x)+h(y).
\end{equation*}

\noindent Then, we have 
\begin{equation*}
\forall x\in \mathbb{R},h(x)=xh(1).
\end{equation*}
\end{lemma}

\bigskip

\noindent \textbf{Proof of Lemma \ref{funct.hamel.lem.1}}. We begin to prove that $h(0)=0$. Wehave
\begin{equation*}
h(0)=h(0+0)=h(0)+h(0)=2h(0),
\end{equation*}

\noindent which implies that $h(0)=0$. Next for any $x\in \mathbb{R}$%
\begin{equation*}
0=h(0)=h(x+(-x))=h(x)+h(-x)
\end{equation*}

\noindent so that for $x\in \mathbb{R}$%
\begin{equation*}
h(-x)=-h(x).
\end{equation*}

\noindent Based on this, we may and do concentrate on positive values of $x$. We may easily see, for 
$0<p\in \mathbb{Z}$ and $0<q\in \mathbb{Z}$, that 

\begin{equation*}
h(p)=h(\underset{p\text{ times}}{\underbrace{1+...+1}})=ph(1)
\end{equation*}

\noindent and

\begin{equation*}
h(1)=h(\underset{q\text{ times}}{\underbrace{(1/q)+...+(1/q)}})=qh(1/q)
\end{equation*}

\noindent so that
\begin{equation*}
h(1/q)=1/qh(1).
\end{equation*}

\noindent Then, for any nonnegative rational $x=p/q$, we have 
\begin{equation*}
h(p/q)=h(\underset{p\text{ times}}{\underbrace{(1/q)+...+(1/q)}}%
)=ph(1/q)=(p/q)h(1).
\end{equation*}

\noindent
\noindent Now, let $x$ be any nonnegative real number. There exists a sequence of
nonnegative rational number  $x_{n}>0$ such that $x_{n}$ decreases to $x$ and a sequence of nonnegative
rational number $y_{n}>0$ such that $x_{n}$ increases.\\

\noindent If $h$ is right-continuous, we have

\begin{equation*}
h(x_{n})=x_{n}h(1)\rightarrow xh(1)\text{ and }h(x_{n})\rightarrow h(x).
\end{equation*}

\bigskip \noindent Now, if $h$ is left-continuous, we have

\begin{equation*}
h(x_{n})=y_{n}h(1)\rightarrow xh(1)\text{ and }h(y_{n})\rightarrow h(x).
\end{equation*}

\bigskip \noindent In both cases, since $h(0)=0$, we have for any $x\geq 0,$ 
\begin{equation*}
h(x)=xh(1).
\end{equation*}

\noindent For $x<0$,
\begin{equation*}
h(x)=h(-(-x))=-h(-x)=-x(-xh(1))=xh(1).
\end{equation*}

\bigskip\noindent We are going to give more general solutions.\\

\noindent We have the first general result.\\

\begin{lemma} \label{funct.hamelG.lem.1} Let $h:$ $\mathbb{R}\longmapsto \mathbb{R}$ be a function such
that 
\begin{equation}
\forall (x,y)\in \mathbb{R}^{2},h(x+y)=h(x)+h(y). \label{funct.Hamel.Eq01}
\end{equation}

\noindent If $h$ is monotone or $h$ is bounded in a neighborhood of zero, or if $h$ is right-continuous or left-contunous at one point, then we have 
\begin{equation}
\forall x\in \mathbb{R},h(x)=xh(1). \label{funct.HamelG.sol}
\end{equation}

\bigskip \noindent This conclusion is obtained if $h$ is only bounded in some right or left neighborhood of zero. 
\end{lemma}

\bigskip \noindent \textbf{Proof of lemma \ref{funct.hamelG.lem.1}}. In the proof the previous lemma, we proved that if \ref{funct.Hamel.Eq01} holds, then we have any rational number $r$,
\begin{equation*}
h(r)=rh(1).
\end{equation*}

\bigskip \noindent (a) Suppose that $h$ is nondecreasing. For any x $\in $ $\mathbb{R}$, for any $%
\delta >0$ there exist two rational numbers $r_{1}=r_{1}(\delta )$ and $%
r_{2}=r_{2}(\delta )$ such that
\begin{equation*}
r_{1}<x<r_{2}\text{ and }r_{2}-r_{1}<\delta .
\end{equation*}

\noindent We have
\begin{equation*}
h(r_{1})\leq h(x)\leq h(r_{2})\text{ and }h(r_{i})=r_{i}h(1).
\end{equation*}

\noindent Then
\begin{equation*}
0\leq h(r_{2})-h(x)\leq h(r_{2})-h(r_{1}),
\end{equation*}

\noindent that is
\begin{equation*}
0\leq r_{2}h(1)-h(x)\leq \delta .
\end{equation*}

\noindent As $\delta \downarrow 0$, $r_{2}=r_{2}(\delta \downarrow x$ and then 
\begin{equation*}
h(x)=xh(1).
\end{equation*}

\noindent Formula (\ref{funct.HamelG.sol}) is proved for $h$ non-decreasing. If $h$ is non-increasing, then its opposite $-h$ is non-decreasing
and satisfies (\ref{funct.Hamel.Eq01}). Then (\ref{funct.HamelG.sol}) holds for $-h$, and this formula is the same for $h$ and for $-h$. The proof of (\ref{funct.HamelG.sol}) is finished for a monotone function $h$ satisfying (\ref{funct.Hamel.Eq01}) .\\

\bigskip \noindent (b) Suppose that $h$ is bounded in some neighborhood of zero, say in $D=\{x,\left\vert
x\right\vert \leq \delta _{0}\},$ $\delta _{0}>0$. Set, for $0<\delta
<\delta _{0},$
\begin{equation*}
A(\delta )=\sup_{\left\vert x\right\vert \leq \delta }\left\vert
h(x)\right\vert .
\end{equation*}

\noindent Fix $\delta $ such that $0<\delta <\delta _{0}.$ Then for any $n\geq 1,$ for
any $x$ such that $\left\vert x\right\vert \leq \delta /n,$ we have

\begin{equation*}
h(nx)=nh(x)
\end{equation*}

\noindent and then

\begin{equation}
\left\vert h(x)\right\vert =\left\vert \frac{h(nx)}{n}\right\vert \leq \frac{%
A(\delta )}{n}.  \label{funct.voisin01}
\end{equation}

\noindent Now let $x$ be an arbitrary real number. For each $n\geq 1,$ let $r_{n}$ a
rational number such that 
\begin{equation}
\left\vert x-r_{n}\right\vert \leq \delta /n.  \label{funct.voisin02}
\end{equation}

\noindent We have

\begin{eqnarray*}
\left\vert h(x)-h(1)x\right\vert  &=&\left\vert h(x)-h(1)\left\{
r+(x-r)\right\} \right\vert  \\
&=&\left\vert h(x)-h(1)r+h(1)(x-r)\right\vert  \\
&=&\left\vert h(x)-h(r)+(x-r)h(1)\right\vert \text{ \ (since \ }h(r)=h(1)r%
\text{)} \\
&=&\left\vert h(x-r)+(x-r)h(1)\right\vert \text{ \ \ \ \ \ \ (By Assumption \ref{funct.Hamel.Eq01}
\ )} \\
&\leq &\frac{A(\delta )}{n}+|h(1)|\frac{\delta }{n}\text{ .\ \ \ \ \ \ \ \ \ \
\ \ \ \ \ \ \ \ \ (By (\ref{funct.voisin01}) and (\ref{funct.voisin01}) }
\end{eqnarray*}

\noindent By Letting $n\rightarrow +\infty ,$ we get%
\begin{equation*}
h(x)=h(1)x.
\end{equation*}

\bigskip \noindent (c) The proofs of (\ref{funct.hamelG.lem.1}) for a function $h$ satisfying (\ref{funct.Hamel.Eq01}) which is left-continuous or righ-continuous at some point are very similar. We only give the proof for one case. Let $h$ be right-continuous at $x_0$. Let $x$ an arbitrary real number. For any sequence of real numbers $0\leq r_n \downarrow 0$, as $n\uparrow +\infty$, we have

$$
h(x_0-x+r_n)=h(x_0+r_n)+h(-x) \rightarrow h(x_0)+h(-x)=h(x_0-x).
$$

\noindent Then $h$ is right-continuous at all points $x_0-x$, then at any point of $\mathbb{R}$. We concude by applying Lemma \label{funct.hamel.lem.1} through it right-continuity part.\\

\noindent The proof lemma (\ref{funct.hamelG.lem.1}) is complete. $\blacksquare$\\

\bigskip \noindent Here are other versions of Hamel Equations.\\

\begin{corollary} \label{funct.hamel.cor.1} \ Let $h:$ $\mathbb{R}_{+}\longmapsto \mathbb{R}_{+}$ be a measurable
function such that 
\begin{equation*}
\forall (s,t)\in \mathbb{R}_{+}^{2},k(st)=k(s)k(t).
\end{equation*}
Then we have 
\begin{equation*}
\forall t\in \mathbb{R},k(t)=t^{\log k(e)}.
\end{equation*}
\end{corollary}

\noindent \textbf{Proof of Corrolary \ref{funct.hamel.cor.1}}. Given the assumption holds, set $h(x)=\log k(e^{x}),x\in \mathbb{R}$. We have%
\begin{equation*}
h(x+y)=\log k(e^{x}e^{y})=\log k(e^{x})+\log k(e^{y})=h(x)+h(y).
\end{equation*}

\noindent Then by Lemma \ref{funct.hamel.lem.1}, we get for any $x\in \mathbb{R}$%
\begin{equation*}
h(x)=\log k(e^{x})=x\log k(e).
\end{equation*}

\noindent By the transform $t=e^{x}$, we get for any $t\geq 0,$ 
\begin{equation*}
k(t)=t^{\log k(e)}
\end{equation*}

\bigskip 

\begin{corollary} \label{funct.hamel.cor.2}  Let $\ell :$ $\mathbb{R}\longmapsto \mathbb{R}$ be a measurable function
such that 

\begin{equation*}
\forall (s,t)\in \mathbb{R}^{2},\ell (st)=\ell (s)+\ell (t).
\end{equation*}

\noindent Then we have 
\begin{equation*}
\forall t\in \mathbb{R}, \ \ell(t)=t^{\ell(e)}.
\end{equation*}
\end{corollary}

\bigskip \noindent \textbf{Proof of Collorary \ref{funct.hamel.cor.2}}. Given the assumption holds, set $h(y)=\ell (e^{y}),y>0.$ We have%
\begin{equation*}
h(s+t)=\ell (e^{s}e^{t})=\ell (e^{s})+\ell (e^{t})=h(s)+h(t).
\end{equation*}

\noindent Then 
\begin{equation*}
h(y)=yh(1)=y\ell (1),y\in \mathbb{R}
\end{equation*}

\noindent which implies

\begin{equation*}
\ell (e^{y})=y\ell (1),y\in \mathbb{R}
\end{equation*}%

\noindent Put $x=e^{y}$ to get
\begin{equation*}
\ell (x)=\ell (e)\log x,x>0
\end{equation*}

\section{Miscelleanuous facts} \label{funct.facts}

\bigskip \noindent \textbf{FACT 1}. For any $a\in \mathbb{R},$ 
\begin{equation*}
\left\vert e^{ia}-1\right\vert =\sqrt{2(1-\cos a)}\leq 2\left\vert \sin
(a/2)\right\vert \leq 2\left\vert a/2\right\vert ^{\delta }.
\end{equation*}

\bigskip \noindent This is easy for $\left\vert a/2\right\vert >1.$ Indeed
for $\delta >0,\left\vert a/2\right\vert ^{\delta }>0$ and%
\begin{equation*}
2\left\vert \sin (a/2)\right\vert \leq 2\leq 2\left\vert a/2\right\vert
^{\delta }
\end{equation*}

\bigskip \noindent Now for $\left\vert a/2\right\vert >1,$ we have the
expansion

\begin{eqnarray*}
2(1-\cos a) &=&a^{2}-\sum\limits_{k=2}^{\infty }(-1)^{2}\frac{a^{2k}}{(2k)!}%
=x^{2}-2\sum\limits_{k\geq 2,k\text{ }even}^{\infty }\frac{a^{2k}}{(2k)!}-%
\frac{a^{2(k+1)}}{(2(k+1))!} \\
&=&a^{2}-2x^{2(k+1)}\sum\limits_{k\geq 2,k\text{ }even}^{\infty }\frac{1}{%
(2k)!}\left\{ \frac{1}{a^{2}}-\frac{1}{(2k+1)((2k+2)...(2k+k)}\right\} .
\end{eqnarray*}

\bigskip \noindent For each $k\geq 2,$ for $\left\vert a/2\right\vert <1,$%
\begin{equation*}
\left\{ \frac{1}{a^{2}}-\frac{1}{(2k+1)((2k+2)...(2k+k)}\right\} \geq
\left\{ \frac{1}{4}-\frac{1}{(2k+1)((2k+2)...(2k+k)}\right\} \geq 0.
\end{equation*}

\bigskip \noindent Hence 
\begin{equation*}
2(1-\cos a)\leq a^{2}.
\end{equation*}

\bigskip \noindent But for $\left\vert a/2\right\vert ,$ the function $%
\delta \hookrightarrow \left\vert a/2\right\vert ^{\delta }$ is
non-increasing $\delta ,0\leq \delta \leq 1$. Then%
\begin{equation*}
\sqrt{2(1-\cos a)}\leq \left\vert a\right\vert =2\left\vert a/2\right\vert
^{1}\leq 2\left\vert a/2\right\vert ^{\delta }.
\end{equation*}


\end{document}